\documentclass{amsart}
\usepackage{color}

\usepackage{amsmath} 
\usepackage{amssymb}
\usepackage{mathrsfs}

\newtheorem{theorem}{Theorem}[section] 
\newtheorem{claim}[theorem]{Claim}

\newtheorem{tuc}[theorem]{The Uniqueness Claim}
\newtheorem{cc}[theorem]{Crucial Claim}

\newtheorem{conclusion}[theorem]{Conclusion}
\newtheorem{observation}[theorem]{Observation}

\theoremstyle{definition}
\newtheorem{definition}[theorem]{Definition}

\newtheorem{explanation}[theorem]{Explanation}

\newtheorem{discussion}[theorem]{Discussion}

\newtheorem{hypothesis}[theorem]{Hypothesis}

\theoremstyle{remark}
\newtheorem{remark}[theorem]{Remark}

\newtheorem{notation}[theorem]{Notation}
\newtheorem{context}[theorem]{Context}

\newcommand{\ic}{{\rm ic}}

\newcommand{\cmp}{{\rm cmp}}

\newcommand{\cor}{{\rm cor}}
\newcommand{\otp}{{\rm otp}}

\newcommand{\true}{{\rm true}}
\newcommand{\fsupp}{{\rm fsupp}}

\newcommand{\wsupp}{{\rm wsupp}}

\newcommand{\ec}{{\rm ec}}
\newcommand{\uec}{{\rm uec}}
\newcommand{\wec}{{\rm  wec}}

\newcommand{\md}{{\rm md}}
\newcommand{\suc}{{\rm suc}}

\newcommand{\pigyon}{{\dagger}}

\newcommand{\meagre}{{\rm meagre}}

\newcommand{\Hechler}{{\rm Hechler}}
\newcommand{\dep}{{\rm dp}}
\newcommand{\deq}{{\rm dq}}

\newcommand{\Ord}{{\rm Ord}}

\newcommand{\tr}{{\rm tr}}
\newcommand{\bd}{{\rm bd}}
\newcommand{\eq}{{\rm eq}}

\newcommand{\id}{{\rm id}}

\newcommand{\cov}{{\rm cov}}
\newcommand{\dom}{{\rm dom}}

\newcommand{\cf}{{\rm cf}}

\newcommand{\rang}{{\rm rang}}
\newcommand{\stt}{{\rm stt}}
\newcommand{\lqq}{{\text{`}\text{`}}}
\newcommand{\page}{{\rm pg.}}

\newcommand{\rest}{{\restriction}}

\newcommand{\wilog}{{\rm without loss of generality}}
\newcommand{\Wilog}{{\rm Without loss of generality}}

\newcommand{\then}{{\underline{then}}}
\newcommand{\when}{{\underline{when}}}

\newcommand{\oor}{{\underline{or}}}
\newcommand{\Then}{{\underline{Then}}}

\newcommand{\Iff}{{\underline{iff}}}
\newcommand{\mn}{{\medskip\noindent}}

\newcommand{\bfm}{{\mathbf m}}

\newcommand{\bfr}{{\mathbf r}}

\newcommand{\bfG}{{\mathbf G}}
\newcommand{\bfb}{{\mathbf b}}

\newcommand{\bfV}{{\mathbf V}}

\newcommand{\bfB}{{\mathbf B}}

\newcommand{\bfn}{{\mathbf n}}

\newcommand{\bfM}{{\mathbf M}}
\newcommand{\bfq}{{\mathbf q}}

\newcommand{\bfU}{{\mathbf U}}
\newcommand{\bfR}{{\mathbf R}}

\newcommand{\cB}{{\mathscr B}}

\newcommand{\gb}{{\mathfrak b}}

\newcommand{\cE}{{\mathscr E}}
\newcommand{\varp}{{\varepsilon}}
\newcommand{\gd}{{\mathfrak d\/}} 

\newcommand{\cH}{{\mathscr H}}

\newcommand{\cF}{{\mathscr F}}
\newcommand{\cG}{{\mathscr G}}
\newcommand{\cI}{{\mathscr I}}

\newcommand{\bbL}{{\mathbb L}}

\newcommand{\bbP}{{\mathbb P}}
\newcommand{\bbR}{{\mathbb R}}
\newcommand{\bbQ}{{\mathbb Q}}
\newcommand{\cP}{{\mathscr P}}

\newcommand{\cS}{{\mathscr S}}
\newcommand{\cT}{{\mathscr T}}
 
\newcommand{\cU}{{\mathscr U}}
\newcommand{\cX}{{\mathscr X}}
\newcommand{\cY}{{\mathscr Y}}

\newcommand{\fat}{{ \rm fat}}
\newcommand{\cL}{{\mathscr{L}}}

\newcount\skewfactor
\def\mathunderaccent#1#2 {\let\theaccent#1\skewfactor#2
\mathpalette\putaccentunder}
\def\putaccentunder#1#2{\oalign{$#1#2$\cr\hidewidth
\vbox to.2ex{\hbox{$#1\skew\skewfactor\theaccent{}$}\vss}\hidewidth}}
\def\name{\mathunderaccent\tilde-3 }

\newbox\noforkbox \newdimen\forklinewidth
\setbox0\hbox{$\textstyle\bigcup$}
\setbox1\hbox to \wd0{\hfil\vrule width \forklinewidth depth \dp0
                        height \ht0 \hfil}
\setbox\noforkbox\hbox{\box1\box0\relax}
\def\unionstick{\mathop{\copy\noforkbox}\limits}
\def\nonfork#1#2_#3{#1\unionstick_{\textstyle #3}#2}
\def\nonforkin#1#2_#3^#4{#1\unionstick_{\textstyle #3}^{\textstyle
    #4}#2}
\setbox0\hbox{$\textstyle\bigcup$}
\setbox1\hbox to \wd0{\hfil{\sl /\/}\hfil}
\setbox2\hbox to \wd0{\hfil\vrule height \ht0 depth \dp0 width
                                \forklinewidth\hfil}
\newbox\doesforkbox
\setbox\doesforkbox\hbox{\box1\box0\relax}
\def\nunionstick{\mathop{\copy\doesforkbox}\limits}

\def\fork#1#2_#3{#1\nunionstick_{\textstyle #3}#2}
\def\forkin#1#2_#3^#4{#1\nunionstick_{\textstyle #3}^{\textstyle
    #4}#2}

\newcommand{\stickT}{%
\setbox255=\hbox{\raise1ex\hbox{$\hspace{0.2pt}\,\bullet\,$}}
\mathord{\rlap{\hbox to\wd255{\hss\hbox{$|$}\hss}}
\box255}
}
\newcommand{\stickS}{%
\setbox255=\hbox{\raise0.6ex\hbox{$\scriptstyle\bullet$}}
\mathord{\rlap{\hbox to\wd255{\hss\hbox{$\scriptstyle|$}\hss}}
\box255}
}

\newenvironment{PROOF}[2][\proofname.]
   {\begin{proof}[#1]\renewcommand{\qedsymbol}{\textsquare$_{\rm #2}$}}
   {\end{proof}}

\usepackage[hidelinks]{hyperref}

\begin{document}
\makeatletter\def\shfiuwefootnote{\gdef\@thefnmark{}\@footnotetext}\makeatother\shfiuwefootnote{Version 2022-08-13\_2. See \url{https://shelah.logic.at/papers/1126/} for possible updates.}

\title{Corrected Iteration \\ Sh1126}

\author{Saharon Shelah}

\address{Einstein Institute of Mathematics\\
Edmond J. Safra Campus, Givat Ram\\
The Hebrew University of Jerusalem\\
Jerusalem, 9190401, Israel\\
 and \\
 Department of Mathematics\\
 Hill Center - Busch Campus \\ 
 Rutgers, The State University of New Jersey \\
 110 Frelinghuysen Road \\
 Piscataway, NJ 08854-8019 USA}
\email{shelah@math.huji.ac.il}

\urladdr{http://shelah.logic.at}

\thanks{Partially supported by European Research Council Grant No.338821,  and the Israel Science Foundation grant 1838/19. Paper 1126 on the Author's list. The author thanks Alice Leonhardt for the beautiful typing of earlier versions (up to 2019) and in later versions the  author would like to thank the typist for his work and is also grateful for the generous funding of typing services donated by a person who wishes to remain anonymous. References like \cite[2.7=La32]{Sh:945} means the label of Th.2.7 is a32.  The reader should note that the version in my website is usually more updated than the one in the mathematical archive.  First typed October 18, 2017}



\date{August 11, 2022}  

\begin{abstract}
    For $\lambda$ inaccessible, we may consider $(<\lambda)$-support iteration of some definable in fact specific $(<\lambda)$-complete $\lambda^+$-c.c. forcing notions.  But do we have ``preservation by restricting to a sub-sequence of the iterated forcing"? To regain it we ``correct" the iteration.  We prove this for a characteristic case for iterations which holds by ``nice'' for $\lambda = \aleph_{0}.$ This is done generally in a work H. Horowitz. 
    
    \noindent 
    For \cite{Sh:945} we use so called strongly $(< \lambda^{+})$-directed $\bfm.$ We could here restrict ourselves to reasonable $\bfm$ (see \ref{b36}(3)).  
\end{abstract}

\maketitle
\numberwithin{equation}{section}
\setcounter{section}{-1}

\newpage

\centerline{Annotated Content}

\S0 \quad Introduction (label z23, z26), pg.\pageref{0},

\S1 \quad Iteration Parameters, (label c0-c41), pg.\pageref{1},

\S1A \quad The frame, \page \pageref{1A},

\S1B \quad Special sufficient conditions [Used only in \S3D],  page \pageref{1B},
  
\S1C \quad On existentially closed $ \mathbf{m}$, \page\pageref{1C},

\S2 \quad The Corrected $\bbP_{\bfm}$ (label b), pg.\pageref{2},

\S3 \quad The Main Conclusion (label e), pg.\pageref{3},
 
\S3A \quad Wider $ \mathbf{m} $-s, \page \pageref{3A},

\S3B \quad Ordinal equivalence, \page \pageref{3B},

\S3C \quad  Representing a condition from $ \mathbb{P}_ \mathbf{m} [M_ \mathbf{m}]$, \page \pageref{3C},

\S3D \quad  The main result, \page \pageref{3D},

\S4 \quad General $\bfm$'s \page, \pageref{3E},

\S4A \quad Alternative proof \page, \pageref{4A},

\S4B \quad General $\bfm$'s \page, \pageref{3E}

\S4C \quad Nicely existentially closed  \page, \pageref{4C}. \\ 

\underline{What numbers are in each sub-section}

\S1A \quad numbers  \ref{c0} till  \ref{c32n},

\S1B \quad   numbers  \ref{c33n} till    \ref{c33s},

\S1C \quad numbers \ref{c34}  till \ref{c41},

\S2 \quad numbers \ref{b2} till \ref{b47},

\S3A \quad number   \ref{e4} till \ref{e16},

\S3B \quad number   \ref{e19} till \ref{e28},

\S3C \quad number   \ref{e29} till \ref{e32},

\S3D \quad number   \ref{e33} till \ref{e50},

\S4A \quad number     \ref{h2} till \ref{h22}.

\S4B  \quad number   \ref{e51} till \ref{e63},

\S4C \quad number   \ref{g0} till \ref{e70}.

\newpage  

\section{Introduction}\label{0}

This work is dedicated to proving a theorem on $(<\lambda)$-support
iterations of $(< \lambda)$-complete ``nicely" definable $\lambda^+$-c.c. forcing notions for $\lambda$ inaccessible. A nicely definable forcing notion can be, for example, random reals forcing (when $ \lambda = {\aleph_0} $).  Pedantically, at each stage it is a different forcing notion,  but it has the same definition at every step of the iteration. 
Assume $\bbQ$ is such a definition, $\langle
\bbP_\alpha,\name{\bbQ}_\beta:\alpha \le \alpha_*,\beta < \alpha_*\rangle$ is such an iteration, $\bbQ_\beta = \bbQ^{\mathbf V[\bbP_\beta]}$ has generic $\name\eta_\beta$. A question is: assuming   $\langle \name\eta_\beta:\beta < \alpha_*\rangle$ is generic for $\bbP_{\alpha_*}$, and
letting $\beta_*$ be maximal such that $2
\beta_* \le \alpha_*,$  does it follows that also the sequence   $\langle \name\eta_{2 \beta}:\beta$ satisfies $2 \beta < \alpha_*\rangle$  is
generic for the iteration $\langle \bbP_\alpha,\name{\bbQ}_\beta:\alpha \le \beta_*,\beta < \beta_*\rangle$?

The point is that in the parallel case for $\lambda = \aleph_0$ so  for FS-iterated forcing such a claim is true.  In fact,  by Judah-Shelah \cite{Sh:292}, if $\langle \bbP_\alpha, \name{\bbQ}_\beta:\alpha \le \alpha(*),\beta <  \alpha(*)\rangle$ is FS-iteration of Suslin-c.c.c. forcing notions, $\name{\bbQ}_\beta$ with the generic $\name\eta_\beta \in {}^\omega
\omega$ and for notational transparency, its definition is with no parameter and the function $\zeta:\beta(*) \rightarrow \alpha(*)$ is increasing
and $\bbP = \langle \bbP'_\alpha,\name{\bbQ}'_\beta:\alpha \le
\beta(*),\beta < \beta(*)\rangle$ is FS iteration,
$\name{\bbQ}'_\beta$ defined exactly as $\name{\bbQ}_{\zeta(\beta)}$  but now in $\mathbf V^{\bbP'_\beta}$ rather than $\mathbf V^{\bbP_{\zeta(\beta)}}$ then $\Vdash_{\bbP_{\alpha(*)}} ``\langle \name\eta_{\zeta(\beta)}:\beta < \beta(*)\rangle$ is
generic for $\bbP'_{\beta(*)}$ over $\mathbf V"$.   For CS iteration of Suslin proper forcing a weaker result holds, see \cite[\S2]{Sh:292} and  \cite{Sh:630}.

Now this is not clear to us for $(< \lambda)$-support iteration of $(< \lambda)$-strategically complete forcing notions.   The solution is
essentially to change the iteration to what we call ``corrected iteration".  We use a ``quite generic" $(<
\lambda)$-support iteration which ``includes" the one we like and use the complete sub-forcing it generates. Here we deal with a characteristic case (used in \cite{Sh:945}).  The proof applies also to  partial memory iteration.  On wide generalization (including the case  $\lambda = \aleph_0$) and application (for $\lambda = \aleph_0$)  this is continued in a work  of H. Horowitz and the author \cite{Sh:1204};   more fully \cite{Sh:1204} generalizes \S1,\S2, \S3A,\S3B, \S3D of the  present work whereas \S3C, \S3E, \S3F were added later, and \S3C  is inverse engineering of \cite[4.2,4.4]{Sh:1204}. Our main result is \ref{b35}, proving that there is ``corrected iteration", i.e.  one satisfying the promised property or see \ref{b32} and more in  \ref{b44}, \ref{b47}. 

The problem arises as follows.  In \cite{Sh:945} it is proved that for $\lambda$ inaccessible, consistently $\cov_\lambda(\meagre)$, the covering number of the meagre ideal on $\lambda$ is strictly smaller
than $\gd_\lambda$, the dominating number.  The result here is used there but the editor prefers to separate it. 
In \S3F we have  an alternative proof of the main theorem, for this we noted   in some earlier places   
what rely on what. 

We have two extreme versions of our frameworks, one we call fat, that is, in Definition \ref{c6}, $\bbP_{\bfm, t} = [u_{\bfm, t}]^{\leq \lambda}$ (used in \cite{Sh:945}). The other is the lean one when the $\cP_{\bfm, t}$ are restricted to the leaves (i.e. $t / E_{\bfm}'$). This was the original version and is the one continued in Horowitz-Shelah \cite{Sh:1204}.
 
The interest in having ``$\bfm$ is strongly $\lambda^{+}$-directed'' is that it implies $\Vdash_{\bbP_{\bfm}}$``$\{ \name{\eta}_{s}: s \in M \}$ cofinal in $\left( \Pi_{\varp < \lambda}\theta_{\varp}, <_{J_{\lambda}^{\rm{bd}}} \right)$'', by \ref{c37}. Now using $\bfm \in \bfM_{\rm{ec}}$ (being full and wide) as constructed in \S1C, does not give this, e.g. because there may be $t \in L_{\bfm}$ above all members of $M_{\bfm}.$ This is circumvented in \ref{b17}  by having, on the one hand for cofinaly many $c \in M_{\bfm}, \bfm(< c) \in \bfM_{\rm{ec}}$ and on the other hand  having ``$\bfm$ is strongly $(< \lambda^{+})$-directed'' (see \ref{b36}(2)). An alternative approach is to restrict ourselves  to the fat context.  

This work is continued in \cite{Sh:1204} and lately in \cite{Sh:F2170}, which in particular sort out when corrected iteration is necessary; we have lecture on this in the Set Theoretic Conference, in Jerusalem, July 2022.

We thank Shimoni Garti and Haim Horowitz for helpful
comments. We thank   Johannes Sch\"urz  and  Martin Goldstern for pointing out several times problem with the application to \cite{Sh:945}, in particular  in 2019 that an earlier version  of the proof of \cite[2.7=La32]{Sh:945} the statement  $ \circledast '_ 4$  was insufficient; and  later pointing out a problem in earlier version of \S3E.  We thanks Mark Po\'or for pointing out many points which need correction. 

For a reader of \cite{Sh:945} we try to give
exact references to the places here we rely 
on there (pages refer to the 2022-08  version; 
there we assume that $\bfm$ is ordinary, that is, $ L_ \mathbf{m} $  has set of elements  an ordinal $ \alpha (\mathbf{m} ) $  and $ \beta < \gamma < \alpha ( \mathbf{m} ) $  implies $ \beta < _{L_ \mathbf{m} } \gamma $).

\begin{enumerate} 
    \item[(a)] on  \cite[1.8=Lz32, page 6]{ Sh:945}, the definition  of $ \mathbf{Q}= \mathbf{Q}_ {\lambda,  \bar{ \theta} ,\alpha (*)}$,
    see here Def \ref{c6},  Claim \ref{c6d}, page  \pageref{c6}, so $ \mathbf{q} $  
    there is (essentially) $ \mathbf{q} _ \mathbf{m} $ here, and so $\mathbb{P} _{\mathbf{m} \upharpoonright L} $ here is dense in  $ \mathbb{P} ^ \mathbf{q} _{ 0, \alpha } $ 
    there when $ L = L_ \mathbf{m} \upharpoonright \alpha, $ 

    \item[(b)] on  \cite[1.9=Lz33, pag.7]{Sh:945} where $\mathbb{P} ^ \mathbf{q} _{1, \alpha }$ defined  there, is $ \mathbb{P} _ \mathbf{m} [L_ \mathbf{m} \upharpoonright \alpha] $ here; see \ref{b11}(3), page \pageref{b11},
    
    \item[(c)]  on  \cite[1.10=Lz35, pag.7]{ Sh:945}, 
    claim on the existence of  generic; 
    include  changing the  generic in $  < \lambda $ 
    places see here  \ref{c8}, \ref{c11}, page \pageref{c8}, \pageref{c11} respectively, 
    
    \item[(d)] on  \cite[1.11=Lz38, pag.8]{Sh:945} 
    see  \ref{b35} page \pageref{b35} or \ref{b38}, page \pageref{b38}, 
    
    \item[(e)] in $(\ast)_{1}$(A) in the proof of \cite[2.7 = La32, page 15]{Sh:945}, see (a)-(e) above,
    
    \item[(f)] in $(\ast)_{4}$ in the proof of \cite[2.7 = La32, page 16]{Sh:945},
    
    See \ref{b38}.
    
    \item[(g)] after $(\ast)_{7}$ in the proof of \cite[2.7 = La32, page 17]{Sh:945}
    
    See \ref{z35}(4).
    
    \item[(h)]  on $ \boxplus _1 $ inside the proof of
    Lemma  \cite[2.7=La32, pag. 17]{Sh:945},
    more details are in \ref{b35}, that is:  $ \boxplus (a)( \alpha ) $ by \ref{b35}(A)(c); $ \boxplus (a)(\beta ) $ by \ref{b35}(a)(h);  $ \boxplus (b) $ by \ref{b35}(C); $ \boxplus (c) $ by  \ref{b35}(A)(b); $ \boxplus (d)  $  by \ref{b35}(B); $ \boxplus (e) $ by \ref{b35}(A)(e),
    
    \item[(i)] on  $ \circledast '_4 $ inside the proof of  Lemma   \cite[2.7=La32, pag.18-19]{Sh:945},
    see \cite[4.12-4.27 = Le53-Le70]{Sh:1126}, 
    
    \item[(j)] In \cite[2.8 = La35, pg. 21]{Sh:945} we use \ref{e70}, page \pageref{e70}.
\end{enumerate} 

Note that even if $s \in M_{\bfm} \Rightarrow u_{s} \cap M_{\bfm} = \emptyset$ still: if $\bfm \in \bfM_{\rm{ec}}$ then $M_{\bfm} \models s < t \Rightarrow \name{\eta_{s}} < \name{\eta}_{t} \mod J_{\lambda}^{\rm{bd}},$ see \ref{c37}.  

\begin{notation}\label{z6}
    We try to use standard notation.  We use
    $\theta,\kappa,\lambda,\mu,\chi,$ for cardinals and $\alpha,\beta,\gamma,\delta,\varepsilon,\zeta, \xi$ for ordinals.  We use
    also $i$ and $j$ as ordinals.  We adopt the Cohen convention that $p \le q$ means that $q$ gives more information, in forcing notions.  The
    symbol $\triangleleft$ is preserved for ``being  a proper  initial segment".
    Also recall ${}^B A = \{f:f$ a function from $B$ to $A\}$ and let ${}^{\alpha >}A = \cup\{{}^\beta A:\beta < \alpha\}$, some prefer ${}^{<
    \alpha} A$, but ${}^{\alpha >} A$ is used systematically in the
    author's papers.  Lastly, $J^{\text{\rm bd}}_\lambda$ denotes the
    ideal of the bounded subsets of $\lambda$. 
\end{notation}

Recall from \cite{Sh:945}:

\begin{definition}\label{z23}   
    Let $\lambda$ be inaccessible, $\bar\theta = 
    \langle \theta_\varepsilon:\varepsilon
    < \lambda\rangle$ be a sequence of regular cardinals $< \lambda$ satisfying $\theta_\varepsilon > \varepsilon$.
    
    1) We define the forcing notion $\bbQ = \bbQ_{\bar\theta}$ by:
    
    \begin{enumerate}
        \item[$(\alpha)$]  $p \in \bbQ$ \underline{iff}:          
        \begin{enumerate}
            \item[(a)]  $p = (\eta,f) = (\eta^p,f^p)$, 
            
            \item[(b)]  $\eta \in \prod\limits_{\zeta < \varepsilon} \theta_\zeta$ for some $\varepsilon < \lambda$, ($\eta$ is called the trunk of $p$), 
            
            \item[(c)] $f \in \prod\limits_{\zeta < \lambda} \theta_\zeta$, 
            
            \item[(d)]   $\eta \triangleleft f$. 
        \end{enumerate}
            
        \item[$(\beta)$]  $p \le_{\mathbb Q} q$ iff: 
        
        \begin{enumerate}
            \item[(a)]  $\eta^p \trianglelefteq \eta^q$,

            \item[(b)]  $f^p \le f^q$, i.e. $(\forall \varepsilon < \lambda)
            f^p(\varepsilon) \le f^q(\varepsilon)$, 
            
            \item[(c)]  if $\ell g(\eta^p) \le \varepsilon < \ell g(\eta^q)$ then $\eta^q(\varepsilon) \in [f^p(\varepsilon),\lambda)$, actually follows.
        \end{enumerate}
    \end{enumerate}

    2) The generic is $\name\eta = \cup\{\eta^p:
    p \in \name{\mathbf G}_{\bbQ_{\bar\theta}}\}$.
\end{definition} 

The new forcing defined above is not $\lambda$-complete anymore.  By fixing a trunk $\eta$ one can define a short increasing sequence of conditions which goes up to some $\theta_\zeta$ at the $\zeta$-th coordinate and hence has no upper bound in $\prod\limits_{\zeta < \varp} \theta_\zeta$. However, this forcing is $(< \lambda)$-strategically complete since the  COM (= completeness) player can increase the trunk at each move.

\begin{remark}\label{z26}
    0) The forcing parallel to the creature forcing from \cite{Sh:326}, \cite{Sh:961} but they are ${}^{\omega}\omega$-bounding. 

    1) The forcing is \underline{parallel} to the creature forcing from
    \cite[\S1,\S2]{Sh:326}, \cite{Sh:961} 
    though they are ${}^\omega \omega$-bounding and 
    \underline{not to Hechler} forcing, whose parallel for $\lambda$ is $\bbQ^{\dom}_\lambda = \bbQ^{\Hechler}_\lambda = \{(\nu,f):f \in {}^\lambda \lambda,\nu \triangleleft f\}$, 
    ordered naturally.  We can change the definition of order, saying $ p < q $ iff  $ p=q $ or $ p \le q  \wedge \tr(p) \not=  \tr(q)$ and then  all (strictly) increasing  sequence of length $ < \lambda $ have upper bound, but the gain is doubtful as we shall use only strategic completeness for some derived forcing notions.

    2) Closer to \cite{Sh:326} we can use $\bar{\theta} = \langle \theta_{1, \varp}, \theta_{0, \varp}: \varp < \lambda \rangle$ such that $\theta_{1, \varp} \geq \theta_{0, \varp} = \cf(\theta_{\theta, \varp}) > \varp$ and $\lambda > \theta_{1, \varp},$ and let $\bbQ$ be such that:
    
    \begin{enumerate}
        \item[(a)] $p = (n, f) = (n_{p}, f_{p}) \in \bbQ_{\bar{Q}}$ \underline{iff}:
        
        \begin{itemize}
            \item $\eta \in \Pi_{\varp < \zeta} \theta_{1, \varp}, \zeta < \lambda,$
            
            \item $f \in \Pi_{\varp \in [\zeta, \lambda)}[\theta_{1, \varp}]^{< \theta_{0, \varp}}.$
        \end{itemize}
        
        \item[(b)] $\bbQ_{\bar{\theta}} \models p \leq q$ \underline{iff}: 
        
        \begin{itemize}
            \item $p, q \in \theta_{\theta},$
        
            \item $\eta_{p} \unlhd \eta_{q},$
            
            \item $\varp \in [\lg(\eta_{q}), \lambda) \Rightarrow f_{p}(\varp) \subseteq f_{q}(\varp),$
            
            \item $\varp \in [\lg(\eta_{p}), \lg(\eta_{q})) \Rightarrow \eta_{q}(\varp) \in f_{p}(\varp).$
        \end{itemize}
    \end{enumerate}
    
    Does not matter.
\end{remark}

\begin{notation}\label{z28}
    1) $L,M,N$ are linear orders and $r,s,t$ are members.
    
    2) If $ \eta \in \Pi _{ \varepsilon < \zeta  } \theta _ \varepsilon $ where $ \zeta < \lambda $ 
    then $ (\Pi _{ \varepsilon < \lambda}
    \theta _ \varepsilon )^{ [ \eta ]}$ will 
    mean $ \{ \nu \in\Pi _{ \varepsilon <  \lambda } \theta _ \varepsilon :  
    \nu $ satisfies  $ \eta \trianglelefteq 
    \nu \}  $.  

    3) For a cardinal $ \lambda $ by induction on ordinal  $ \alpha $ we define $ \beth _ \alpha ( \lambda ) $ as $ \lambda +   \Sigma _{\beta < \alpha } 2^{ \beth _ \beta ( \lambda ) } $ and $ \beth _ \alpha = \beth  ( \alpha ) = \beth _ \alpha  ( {\aleph_0} ) $. 
\end{notation}

\begin{discussion}\label{z32}
    1) Fat $\lambda^{+}$-directed $\bfm$ are helpful when we like to have $\Vdash_{\bbP_{\bfm}}$``$\{ \name{\eta}_{s}: s \in M_{\bfm} \}$ is cofinal in $(\Pi_{\varp < \lambda} \theta_{\varp}, <_{J_{\lambda}^{\rm{bd}}})$'' as in \cite{Sh:945}, see Definition \ref{c4}.
\end{discussion}

Recall,

\begin{definition}\label{z35}\ 
    1) We say that a forcing notion
    $\mathbb P$ is $\alpha$-strategically complete \underline{when}  for each $p \in \mathbb P$ in the  following game $\Game_\alpha(p,\mathbb P)$ between the players COM and INC, the player COM has a winning strategy.
    
    A play lasts $\alpha$ moves; in the $\beta$-th move, first the player
    COM chooses $p_\beta \in \mathbb P$ such that $p \le_{\mathbb P} p_\beta$
    and $\gamma < \beta \Rightarrow q_\gamma \le_{\mathbb P} p_\beta$ and
    second the player INC chooses $q_\beta \in \mathbb P$ such that $p_\beta
    \le_{\mathbb P} q_\beta$.
    
    The player COM wins a play if he has a legal move for every $\beta < \alpha$.
    
    2) We say that a forcing notion $\mathbb P$ is $(< \lambda)$-strategically
    complete \underline{when} it is $\alpha$-strategically complete for every
    $\alpha < \lambda$.
    
    Basic properties of $\bbQ_{\bar\theta}$ are summarized and proved in
    \cite[\S2]{Sh:949}. 
\end{definition}

\newpage

\section{Iteration Parameters}\label{1}

\subsection{The frame}\label{1A}

\begin{hypothesis}\label{c0}
    1) $\lambda = \lambda^{< \lambda}$ is strongly inaccessible.

    2) $\bar\theta = \langle \theta_\varepsilon:\varepsilon < \lambda\rangle$.
    
    3) $\theta_\varepsilon$ is an infinite regular cardinal $> \varepsilon$ and $< \lambda$.
    
    4) Assume $\lambda_2 \ge \lambda_1 \ge \lambda_{0} = \cf(\lambda_{0}) > \lambda$ are such that\footnote{usually $\lambda_2  = (\lambda _2 ) ^\lambda  \ge \lambda_1$ suffices but see \ref{e27}, \ref{e37}, however in \S4A we add $\lambda_{2} \geq \beth_{\lambda_{1}^{+}}.$} $(\lambda_1)^{\lambda_{0}} = \lambda_1,$   so notations should have the parameter $\bar\lambda = (\lambda_2,\lambda_1, \lambda_{0}, \lambda)$ and even\footnote{we mainly can use $\lambda_{0} = \lambda^{+},$ but when we restrict ourselves to lean $\bfm$-s, $\lambda_{0} = \lambda$ seem to suffice, see mainly \ref{c8}(f)($\gamma$), \S2, \S3C but does not seem worthwhile  to pursue.} $\bar\lambda =  (\lambda_2,\lambda_1, \lambda_{0}, \lambda,\bar\theta)$.
\end{hypothesis}

\begin{notation}
    1) $L,M$ denote partial orders, well founded if not said otherwise.

    2) Below $ \mathbf{m} , \mathbf{n} $ will be members of $ \mathbf{M} $; we may write e.g. 
    $ L, \mathbf{q} $ instead $ L_ \mathbf{m} , \mathbf{q } _\mathbf{m} $
    when $ \mathbf{m} $ is clear from the context,
    see Def \ref{c4}, \ref{c6}.

    3) We may not pedantically distinguish the subset $ L_1 $ of  $ L $ and the sub-partial order $ L_1 $ of $ L $.
\end{notation}

\begin{remark}\label{c1}
    Here there is  no harm in adding:
    
    \begin{enumerate}
        \item[(a)]  $\theta_\varepsilon > \prod\limits_{\zeta <
        \varepsilon} 2^{\theta_\zeta} + 2^{\aleph_0}$ for $\varepsilon <
        \lambda$,  and/or,
        
        \item[(b)]  for $ \mathbf{m} \in \mathbf{M} $ demanding $M_ \mathbf{m}$ is 
        a linear order, well founded 
        (it suffices to assume even $M \cong (\kappa,<),\kappa$ regular from $[\lambda_{0}, \lambda_{1})$).
    \end{enumerate}
\end{remark}

\begin{definition}\label{c2}
    1) For a partial order $L$ (not necessarily well founded) let:
    
    \begin{enumerate}
        \item[$(\alpha)$]   $\rm{dp}(L) = \cup\{\rm{dp}_L(t)+1:t \in L\}$, see below,
        
        \item[$(\beta)$]  $\rm{dp}_L(t) = \dep(t,L) \in \Ord \cup \{\infty\}$ be defined by
        $\dep_L(t) = \cup\{\dep_L(s)+1$: $s <_L t\}$.
        
        \item[$(\gamma)$]  $L_{<t} = L \rest \{s \in L:s <_L t\}$,
        
        \item[$(\delta)$]  $L_{\le t} = L \rest \{s \in L:s \le_L t\}$.
    \end{enumerate}
    
    2) Let $L^+ = L(+)$ be $L \cup \{\infty\}$ with the natural order (but we may write $t <_L \infty$ instead of $t <_{L(+)} \infty$).
    
    3) We say the set $L$ is an initial segment of the partial order $L_*,$ when: 
    
    \begin{itemize}
        \item  $L \subseteq L_*$, i.e. $s \in L \Rightarrow s \in L_*,$
        
        \item $s <_{L_*} t \wedge t \in L \Rightarrow s \in L$.
    \end{itemize}
    
    The class $\mathbf{M} $ is central in this work, see explanation \ref{c0x}, in particular, $M_{\bfm}$ is our aim, the rest ($L_{\bfm}$ first of all) are the scaffoldings.
\end{definition}

\begin{definition}\label{c4}
    1) Let $\mathbf M$ be the class of objects $\mathbf m$,  called iteration parameters, of the following form (so really $\mathbf M
    = \mathbf M[\bar\lambda]$ and if we omit sub-clauses
    $(\theta), (\iota)$ of clause (e) we may write $\mathbf M[*]$). 
    
    \begin{enumerate}
        \item[(a)]  $L$, a partial order,
        
        \item[(b)]  $M \subseteq L$, as partial orders, (in the main case $M$ is  linearly ordered), 
        
        \item[(c)]  $(\alpha) \quad  u = \langle u_t:t \in L \rangle,$ $\bar{\cP} = \langle \cP_t:t \in L\rangle,$ each $\cP_{t}$ is closed under subsets and $\cP_{t} \subseteq [u_{t}]^{\leq \lambda},$
          
        \item[${{}}$]  $(\beta) \quad u_t \subseteq \{s \in L:s <_L t\},$
         
        \item[(d)]   $\rm{dp}(L) < \infty$, that is $L$ is well founded, 
        
        \item[(e)]  
        
        \begin{enumerate}
            \item[$(\alpha)$] $E'$ is a two-place relation (on $L$),
            
            \item[$(\beta)$] $E'' := E' \rest (L \setminus M)$ is an equivalence relation on $L \setminus M,$
            
            \item[$(\gamma)$] the order $\leq_{L}$ is the transitive closure of $\bigcup \{ \leq_{L} \rest (s / E'): s \in L \setminus M \} \cup \{ \leq_{L} \rest  M\},$ equivalently (using $(\delta)$-$(\eta)$ below): 
            
            \begin{itemize}
                \item if $s, t \in L \setminus M$ are not $E''$-equivalent, then $s <_{L}$ t \underline{iff} for some $r_{1} \leq_{\bfm} r_{2},$ we have $s \leq_{\bfm} r_{2}$ from $\in s / E_{\bfm}', r_{2} \leq_{\bfm} t, r_{2} \in t / E_{\bfm}',$
                
                \item if $s \in L \setminus M$ and $t \in M,$ then $s \leq_{L} t$ \underline{iff} for some $r \in (s / E') \cap M$ we have $s < r \leq t,$
                
                \item if $s \in M$ and $t \in L \setminus M,$ then $s <_{L} t$ \underline{iff} for some $r \in (t / E') \cap M$ we have $s \leq r < t.$ 
            \end{itemize}
            
            \item[$(\delta)$] if $s E' t$ then $s \notin M \vee t \notin \in M,$
            
            \item[$(\varp)$] is $t \in L \setminus M$ then $\{ s \in L: s E' t \} = \{ s \in L: t E' s \};$ we call it $t / E';$ so $E'$ is a symmetric relation,
            
            \item[$(\zeta)$] if $s, t \in L \setminus M$ are $E''$-equivalent then $s / E' = t / E',$
            
            \item[$(\eta)$] if $t \in  L \setminus M$ then $u_{t} \subseteq t / E',$
            
            \item[$(\theta)$] if $t \in L \setminus M$ then $t / E'$ has cardinality $\leq \lambda_{2},$
            
            \item[$(\iota)$] $\Vert M \Vert \leq \lambda_{1},$
        \end{enumerate}
        
        \item[(f)] disjoint subsets $M_{\bfm}^{\rm{fat}}, M_{\bfm}^{\rm{lean}}$ of $M_{\bfm}$ such that: 
        
        \begin{itemize}
            \item if $s \in M_{\bfm}^{\rm{fat}}$ then $\cP_{\bfm, s} = [u_{t}]^{\leq \lambda},$
            
            \item if $s \in M_{\bfm}^{\rm{lean}}$ then $u \in \cP_{\bfm, s} \Rightarrow (\exists t)(u \subseteq t / E_{\bfm})$
            
            \item we let $M_{\bfm}^{\rm{non}} = M_{\bfm} \setminus (M_{\bfm}^{\rm{fat}} \cup M_{\bfm}^{\rm{lean}}).$
        \end{itemize}
    \end{enumerate}
    
    2) Saying $\bfm \in  M$ is lean means that $M_{\bfm} = M_{\bfm}^{\rm{lean}}.$ The lean context means that we restrict ourselves to lean $\bfm$: similarly for fat and neat, see below.
    
    3) We say $\bfm \in M$ is fat \underline{when} $M_{\bfm} = M_{\bfm}^{\rm{fat}}$ and moreover $t \in L_{\bfm} \Rightarrow \cP_{t} = [u_{t}]^{\leq \lambda}.$ 
    
    4) $M_{\bfm}$ is neat \underline{when} $M_{\bfm} = M_{\bfm}^{\rm{lean}} \cup M_{\bfm}^{\rm{fat}}.$  
\end{definition}

\begin{remark}
    1) We may demand $\bfm$ is strongly $(< \lambda)$-directed, see Definition \ref{b36}(2) or even reasonable, see Definition \ref{b36}(3); is harmless here and help \cite{Sh:945}.
    
    2) It may seem reasonable to demand:
    
    \begin{itemize}
        \item[$\boxplus$] is $s \in L_{\bfm} \setminus M_{\bfm}$ and $s \in A \in \cP_{t},$ then $(s / E') \cap u_{t} \in \cP_{t}.$   
    \end{itemize}
    
    However in the crucial claim \ref{e44}, \ref{e47} this cause problems for  $t \in M_{\bfm} \setminus M_{\bfm}^{\fat}.$
\end{remark}

\begin{definition}\label{c5}
    For $\mathbf m \in \mathbf M$.

    0) In \ref{c4} we let $\mathbf m = (L_{\mathbf m},M_{\mathbf m}, \bar u_{\mathbf m}, \bar{\cP}_{\bfm} ,E'_{\mathbf m}, M_{\bfm}^{\rm{lean}}, M_{\bfm}^{\rm{fat}}),$ 
    $\bar u_{\mathbf m} = \langle u_{\mathbf m,t}:t \in  L_{\mathbf m}\rangle, \bar{\cP}_{\bfm} = \langle \cP_{\bfm, t}: t \in L_{\bfm} \rangle, $ for $t \in L_{\bfm} \setminus M_{\mathbf m}$ let
    $t/E_{\mathbf m} = (t/E'_{\mathbf m})  \cup M_{\mathbf m}$ and for $t \in
    M_{\mathbf m}$ let $t/E_{\mathbf m} = M_{\mathbf m}$; \ so there is no relation $E_{\mathbf m}$ but  $t/E_{\mathbf m }$ for $t \in L_{\mathbf m,}$ is well defined.  
    
    1) In \ref{c4}, let $\dep_{\mathbf m}(t) =
    \dep_{L_{\mathbf m}}(t),\dep_{\mathbf m} =  \dep(L_{\mathbf m})$ and $\le_{\mathbf m}  \, = \,  \le_{L_{\mathbf m}}$.
    
    2) For $L \subseteq L_{\mathbf m}$:
    
    \begin{enumerate}
        \item[(a)]  let $\mathbf n = \mathbf m \rest L$ mean $\mathbf n \in \mathbf M,L_{\mathbf n} = L, \leq_{\bfn} = \leq_{\bfm} \rest L_{\bfn}, E'_{\bfn} = E'_{\bfm} \rest L, $ 
        $ u_{\mathbf n,t}  = u_{\mathbf m,t} \cap L $ and $ \cP_{\mathbf n,t} = \cP_{\mathbf m,t} \cap [L]^{\le \lambda}$ for $t \in L$ and $M_{\bfn} = M_{\bfm} \cap L, M_{\bfn}^{\rm{lean}} = M_{\bfm}^{\rm{lean}} \cap L, \, M_{\bfn}^{\rm{fat}} = M_{\bfm}^{\rm{fat}} \cap L,$ 
        
        \item[(b)]   let $\dep_{\mathbf m}(L) = 
        \dep(L_{\mathbf m} \rest L)$ and we may write $\dep(L)$ when $\bfm$ is clear from the context. 
    \end{enumerate}

    3) For $t \in L^+_{\mathbf m}$, let $\mathbf m_{<t} =  \mathbf m(<t) = \mathbf m \rest L_{<t}$ where $L_{<t} = L_{\mathbf m(<t)} = L_{\mathbf m,< t} =  \{s:s <_{\mathbf m} t\}$ so $u_{\mathbf m(<t),s} = u_{\mathbf m,s}$ for $s \in L_{<t}$, etc.
    
    3A)  Also $\mathbf m_{\le t} = \mathbf m(\le t) = \mathbf m \rest L_{\le t}$ where $L_{\le t} = L_{\mathbf m(\le t)} = L_{<t} \cup \{t\}$; let $L_{< \infty} = L,L_{\le \infty} = L^+$, etc.
    
    4) $\mathbf M_{< \mu}$ is the class of $\mathbf m \in \mathbf M$ such that $L_{\mathbf m}$ has cardinality $< \mu$.  Similarly $\mathbf M_{\le
     \mu},\mathbf M_{= \mu},\mathbf M_{> \mu},\mathbf M_{\ge \mu}$; let $\mathbf M_\mu = \mathbf M_{=\mu}$.
    
    5) For $\mathbf m,\mathbf n \in \mathbf M$ let $\mathbf m \approx \mathbf n$, and we may say $\mathbf m,\mathbf n$ are equivalent \underline{meaning} that  $L_{\mathbf m} =
    L_{\mathbf n}$ (as partial orders) and $t \in L_{\mathbf n} \Rightarrow u_{\mathbf m,t} = u_{\mathbf n,t} \wedge \cP_{\mathbf m,t} = \cP_{\mathbf n,t}$; note that there are no demands on $M$ and $E' $. 
    
    6) We say $f$ is an isomorphism from $\mathbf m_1 \in \mathbf M$ onto $\mathbf m_2 \in \mathbf M$ \when:
    
    \begin{enumerate}
        \item[(a)]  $f$ is an isomorphism from the partial order $L_{\mathbf m_1}$ onto the partial order $L_{\mathbf m_2}$,
          
        \item[(b)]  for $s,t \in L_{\mathbf m_1}$ we have $s \in u_{\mathbf m_1,t} \Leftrightarrow f(s) \in u_{\mathbf m_2,f(t)}$ and $\cP_{\bfm_{2}, f(t)} = \{ \{ f(s): s \in u \}: u \in \cP_{\bfm_{1}, t} \},$
        
        \item[(c)]  for $s,t \in L_{\mathbf m_1}$ we have $s E'_{\mathbf m_1,} t \Leftrightarrow f(s) E'_{\mathbf m_2} f(t)$, 
                
        \item[(d)]  $M_{\mathbf m_2} = \{f(s):s \in M_{\mathbf m_1}\}$.
    \end{enumerate}
    
    7) We define weak isomorphisms as in part (6) omitting clauses (c),(d).
    
    8) We say that $ \mathbf{m} $ is ordinary \when  \, the   set of elements of $ L_ \mathbf{m} $ 
    is an ordinal $ \alpha _ \mathbf{m} = \alpha ( 
    \mathbf{m})$  satisfying $ \beta < _{L_ \mathbf{m}} \gamma \Rightarrow \beta < \gamma $.

    9) For a forcing  notion $ \mathbb{P} $ we say that $ q \in \mathbb{P} $  is essentially above $ p \in \mathbb{P} $  (inside $ \mathbb{P} $) when  $ q \Vdash p \in \name{ \mathbf{G} } $.   
            
    10) We say $\bfm \in \bfM_{\rm{\bd}}$ or $\bfm$ is bounded, \underline{when}: if $s \in L \setminus M$ \underline{then} for some $t \in M$ we have $s / E' \subseteq L_{\leq t},$ or just\footnote{In the main case $M_{\bfm}$ is $\lambda^{+}$-directed, so  this does not make a difference. Also no real case when we restrict ourselves to bounded $\bfm$'s.} there is $X \in [M]^{< \lambda_{0}}$ such that $s / E' \subseteq \bigcup_{t \in X} L_{\leq t}.$   
    
    11) We say $\bfm \in \bfM_{\rm{wbd}}$ or $\bfm$ is \emph{weakly bounded} when $L_{\bfm} = \bigcup \{ L_{\bfm(\leq t)}: t \in M_{\bfm} \}.$ 
\end{definition}

\begin{discussion}\label{c5d} 
    Concerning the aim of the choice to use $u_t$ (and $\cP_t$) in \ref{c4}, note the following.

    1) By the partial order we already can get partial memory, so why not simply use only $u_t  
    := \{s:s <_L   t\}$?  After all, the index set is only partially ordered, not necessarily linearly, so these sets can be independent of each other.  The reason is that a partial order is transitive, so this simple definition would imply $s \in u_t  \Rightarrow u_s \subseteq u_ t $   which means (by definition)  the memory is transitive,  but  we do not want that to hold in general, (this is central in \cite{Sh:592}).  Here  $\bar u$ is not necessarily transitive, that is,  $s \in u_t \nRightarrow u_s \subseteq u_t$.  By a  partial order we cannot get it.

    2) In \cite{Sh:700}, \cite{Sh:700a}  we use $\cP_t$'s which are ideals, but here not
    necessarily: this helps, but has  a price; we are relying on  ``$\bbQ_{\bar\theta}$ is close to
    being $\lambda$-centered", i.e. any subset of $\{p \in \bbQ_{\bar\theta}:\tr(p) = \eta\}$ of cardinality $< \theta_{\ell g(\eta)}$ has a lub in this forcing. But for the fat context we get more  than $(< \lambda)$-complete ideal.  

    3) What is the point of ``$\bfm$ being neat''? It tells us that in that case it is easy to be an automorphism of $\bfm,$ see \ref{c11}(2), we may forget to say we use it. 
\end{discussion}

\begin{explanation}\label{c0x}   
    For $\mathbf m \in \mathbf M$:

    \begin{enumerate}
        \item[(a)]  We shall use $L_{\mathbf m}$ as the index set for the iteration; always a well founded partial order.
        
        \item[(b)]  $M_{\mathbf m}$ is the part of the index set we are really interested in, it may be $(\kappa,<)$ as in \cite{Sh:945}.
        
        \item[(c)]  The    other part in the interesting case  is ``generic enough $\mathbf m$", more accurately  existentially closed enough  so that the iteration restricted to $M$ will be  ``stabilized"    under further extensions. That is, for every $\bfm \in \bfM$ we define an iteration resulting in the  forcing $\bbP_{\bfm}$, adding a generic $\name\eta_s$ for $s \in L_{\bfm}$,  we are interested in the extension $\bfV[\langle \name\eta_s:s \in M_{\bfm}\rangle]$, it is the generic extension for the forcing we call $\bbP_{\bfm}[M_{\bfm}]$.  But, in general, even  if $\bfn \in \bfM$ extends $\mathbf{m} $ (see Definition \ref{c26} below of $\le_ \mathbf{M}$) maybe $\bbP_{\bfn}[M_{\bfm}] \ne \bbP_{\bfm}[M_{\bfm}]$.  Our aim is to define $\bfM,\le_{\bfM}$ so that  for a dense set of $\bfm$'s this holds; (done in the crucial claim \ref{c41}). So our aim is having $\bbP_{\bfm}[M_{\bfm}]$,  hence  the $s \in L_{\bfm} \setminus M_ \mathbf{m} $ serves as scaffolding, (but see \ref{b47}). 
        
        Existentially closed structures are used in model theory, but this approach gives non-well founded structures, which  is ``bad" for us.  So an essential  point  here  is to prove (under suitable definitions)  that ``generic, existentially closed enough $\bfm$" is well defined  in spite of $L_{\mathbf m}$  being required to be well founded.
        
        \item[(d)]  of course, the aim of    $ \mathbf{m}  \in \bfM$ is to be used to define the forcing, as in \ref{c6} below. 
    \end{enumerate}
\end{explanation}

\begin{definition}\label{c6}
    1) In the fat context, for $\mathbf m \in \mathbf M$ let $L = L_{\mathbf m}$ and we define the iteration $\mathbf q_{\mathbf m}$ to consist of:

    \begin{enumerate}
        \item[(a)]  a forcing notion $\bbP_t = \bbP_{\mathbf m,t}$  for $t \in L^+$; we let $\bbP_{\mathbf m} = \bbP_\infty$,
        
        \item[(b)]  $\name{\bbQ}_t$ a $\bbP_t$-name of a sub-forcing of $\bbQ_{\bar\theta}$ in the universe $\mathbf V^{\bbP_t}$, even
        $\mathbb{Q} _t\le_{\ic} \mathbb{Q} _{\bar{ \theta}}$   (i.e. $\bbQ_{t} \subseteq \bbQ_{\bar{\theta}}$ as quasi orders and incompatibility and compatibility are  preserved\footnote{But maximal anti-chains - not necessarily.  Recall that $\bbQ_{\bar\theta}$ is from  \ref{z23}, \ref{z26}.  What is   $ \name{\mathbb{Q}  }_t $?  It is implicitly defined in clause (c) and explicitly in \ref{c13}).}),

        \item[(c)]  $p \in \bbP_t$ \Iff:
        
        \begin{enumerate}
            \item[$(\alpha)$]  $p$ is a function, 
            
            \item[$(\beta)$]  $\dom(p) \subseteq L_{<t}$ has cardinality $< \lambda$, 
            
            \item[$(\gamma)$]  if $s \in \dom(p)$ then $p(s)$ consists of $\tr(p(s)) \in \prod\limits_{\varepsilon <\zeta(s)} \theta_\varepsilon$ for some  $\zeta_s = \zeta(s) < \lambda$   and  $\xi = \xi_{p(s)} = \xi(p(s)) \le \lambda$ and $\mathbf B_{p(s)}$ and $\bar r = \bar r_{p(s)} = \langle r(\zeta):\zeta < \xi_{p(s)} \rangle = \langle r_{p(s)}(\zeta):\zeta < \xi_{p(s)}\rangle \in {}^\xi(u_s)$ lists the coordinates used in computing $p(s)$  and are such that: 
            
        \begin{enumerate}
                \item[$\bullet_1$]  $\mathbf B_{p(s)}$ is a $\lambda$-Borel function\footnote{that is, a definition of one.}, $\mathbf B = \mathbf B_{p(s)}: {}^\xi(\prod\limits_{\varepsilon < \lambda} \theta_\varepsilon) \rightarrow \prod\limits_{\varp < \lambda} \theta_\varp$ moreover into $(\prod\limits_{\varp < \lambda} \theta_\varp)^{[\tr(p(s))]}$; and considering  $(d)(\alpha)$ below less pedantically $p(s) = (\tr(p(s)),\name f_{p(s)})$, where \newline $\name f_{p(s)} = $ $\mathbf B_{p(s)}(\ldots, \name\eta_{r_{p(s)}(\zeta)},\ldots)_{\zeta < \xi_{p(s)}}$  which means: absolutely, i.e. in every  forcing extension  $\mathbf V^{\bbQ}$  of $ \mathbf{V} $ where $\mathbb{Q} $ is a $(< \lambda)$-strategically complete  and is   $\lambda^+$-c.c. forcing notion, still $\mathbf B_{p(s)}$ is such a  ($\lambda$-Borel) function; we may write $\xi_{p,s}$ instead of $\xi_{p(s)}$, etc., 
        \end{enumerate}
    \end{enumerate}

    \item[(d)]
    
    \begin{enumerate}
        \item[$(\alpha)$]  $\name \eta_s$ is the $\bbP_t$-name, when $t \in
        L^+_{\mathbf m},s \in L_{< t}$ defined by 
        $\cup\{\tr(p(s)):p \in \name{\mathbf G}_{\bbP_t}\}$,
        
        \item[$(\beta)$]   For $p \in \bbP_t$ and $s \in \dom(p)$ we interpret $p(s)$ as a 
        $\bbP_s$-name \newline $(\tr(p(s)),\mathbf B_{p,s}(\ldots,\name\eta_{r_{p,s}(\zeta)},
        \ldots)_{\zeta < \xi_{p,s}})$. 
    \end{enumerate}

    \item[(e)]  $\bbP_t \models ``p \le q"$ \Iff:
        
        \begin{enumerate}
            \item[$(\alpha)$]   $p,q \in \bbP_t$, 
            
            \item[$(\beta)$]   $\dom(p) \subseteq \dom(q)$, 
            
            \item[$(\gamma)$]   if $t \in \dom(p)$ then $(q \rest L_{<t}) \Vdash_{\bbP_{\mathbf m(<t)}} ``p(t) \le_{\name{\bbQ}_{\bar\theta}} q(t)"$.
        \end{enumerate}
    \end{enumerate}
     
    2) In the general context we replace clause (c)($\gamma$) by: (so part (1) is a special case with $\iota_{p(s)} = 1, \bar{r}_{p(s), 0}  = \bar{r}_{p(s)}$).

    \begin{enumerate}
        \item[$(\gamma)$] if $s \in \dom(p)$ then $p(s)$ consists of $\tr(p(s)) \in \Pi_{\varp < \zeta(s)} \theta_{\varp}$ for some $\zeta_{s} = \zeta(s) < \lambda$ and $\varp = \varp_{p(s)} = \varp(p(s)) \leq \lambda$ and $\bfB_{p(s)}$ and $\bar{r} = \bar{r}_{p(s)} = \langle r(\zeta): \zeta < \varp_{p(s)} \rangle = \langle r_{p(s)}(\zeta): \zeta < \varp_{p(s)} \in {}^{\varp}(u_{s})$ lists the coordinates used in computing $p(s)$ and\footnote{What is the point of ``$\iota < \iota(p(s))"$? As the support is not just $u_s$ but also $\cP_s$ and $\cP_s$ is a family of suitable subsets of $u_s, p(s)$ is $(\tr(p(s)),\name f_s),f_s$ is a name of a member of  $\prod\limits_{\varp < \lambda} \theta_\varp$ such that $\tr(p(s))$ is a (proper)  initial segment.  But how is $\name f_s$ computed?  As our memory is $\cP_s  \subseteq \cP(u_s 
        )$ and not just $u_s  $ (or even a 
        $(< \lambda)$-complete ideal) $\name f_s$ is composed of  $\iota_{p(s)}$
        names each coming from $\langle \name\eta_t:t \in u_\iota\rangle,u _ \iota 
        \in \cP_s$ for $ \iota < \iota ( p (  s ) $. }  $\langle \bfB_{p(s), \iota}, \bar{r}_{p(s), \iota}: \iota < \iota(p(s) \rangle$ are such that:
    \end{enumerate}

    \begin{enumerate}
        \item[$ $] \begin{enumerate}
            \item[$\bullet_1$]  $\mathbf B_{p(s)}$ is a $\lambda$-Borel  function\footnote{that is, a definition of one}, $\mathbf B = \mathbf B_{p(s)}: {}^\xi(\prod\limits_{\varepsilon < \lambda} \theta_\varepsilon) \rightarrow \prod\limits_{\varp < \lambda} \theta_\varp$ moreover into $(\prod\limits_{\varp < \lambda} \theta_\varp)^{[\tr(p(s))]}$; and considering $(d)(\alpha)$ below less pedantically $p(s) = (\tr(p(s)),\name f_{p(s)})$, where $\name f_{p(s)} = \mathbf B_{p(s)}(\ldots, \name\eta_{r_{p(s)}(\zeta)},\ldots)_{\zeta < \xi_{p(s)}}$  which means: absolutely, i.e. in every forcing extension  $\mathbf V^{\bbQ}$  of $ \mathbf{V} $ where  $ \mathbb{Q} $ is a  $(< \lambda)$-strategically complete  and is  $\lambda^+$-c.c. forcing notion, still    $\mathbf B_{p(s)}$ is such a  ($\lambda$-Borel) function; we may write $\xi_{p,s}$ instead of $\xi_{p(s)}$, etc.,
    
            \item[$\bullet_2$]   $\iota _{p(s)} =  \iota(p(s)) < \lambda$ moreover\footnote{This and the rest of $(c)(\gamma)$ are used in the proof of \ref{e33}.  The aim is that defining $\mathbf B_{p(s)}$ from $\langle \mathbf B_{p(s),\iota}: \iota < \iota(p(s))\rangle$, the $\sup$ will not give in $\varepsilon$ the value $\theta_\varepsilon$.}  $< \theta_{\ell g(\tr(p(s))}$, 
    
            \item[$\bullet_3$] for $\iota < \iota_{p(s)}, \bar r_{p(s),\iota} =  \bar r_{p(s)} \rest w_{p(s),\iota}$  so $w_{p(s),\iota} = w(p(s),\iota) = \dom(\bar r_{p(s),\iota}) \subseteq \xi_{p(s)}$  and $\bar r_{p(s),\iota}$ is a subsequence of $\bar r_{p(s)}$,
    
            \item[$\bullet_4$]   $\mathbf B_{p(s),\iota}$ is a Borel function from ${}^{w(p(s),\iota)}(\prod\limits_{\varepsilon < \lambda} \theta_\varepsilon)$ into $(\prod\limits_{\varepsilon < \lambda} \theta_\varepsilon)^{[\tr(p(s))]}$, 
      
            \item[$\bullet_5$]  $\mathbf B_{p(s)}(\langle 
            \name\eta_{r_{p(s)}(\zeta)}:\zeta < \xi_{p(s)}\rangle) = \sup\{\mathbf B_{p(s),\iota}(\langle \eta_{r_{p(s)}(\zeta)}:\zeta \in w_{p(s),\iota}\rangle): \iota < \iota(p(s))\}$ and naturally $ \name{ f }_{p(s)} = \sup\{\name f_{p(s),\iota}:\iota \leq \iota(p(s))\},\name f_{p(s),\iota} = \mathbf B_{p(s),\iota}(\langle \name\eta_\zeta:\zeta \in w_{p(s),\iota}\rangle)$,
        
            \item[$\bullet_6$]  for each $\iota < \iota(p(s))$ for some $u \in
            \cP_{\mathbf m,s}$ we have $\{r_{p(s)}(\zeta):\zeta \in w_{p(s),\iota}\} \subseteq u$ so is a subset of $u_s$, 
    
            \item[$\bullet_7$] (follows) when $\bfm$ is lean,  if $\iota < \iota_{p(s)}$ and $\varp \in  w_{p(s), \iota}, r_{p(s)}(\varp) \in L_{\mathbf m} \backslash M_{\mathbf m}$
            \then \, $\{r_{p(s)}(\zeta):\zeta \in w_{p(s),\iota}\} \subseteq
            r_{p(s)}(\varp)/E_{\mathbf m}$, 
        \end{enumerate}
    \end{enumerate}    
    
    [Why? As Definition \ref{c4}(2) together with $\bullet_{6}$ implies $\{ r_{p(s)}()\zeta : \zeta \in W_{p(s)}, \iota\} \subseteq r_{p(s)}(\varp) / E_{\bfm}'].$
    
    \begin{enumerate} 
    
        \item[$ $]
        \begin{enumerate} 
            \item[$ \bullet _8$]  we let $ {\mathscr F} _{p(s)}$ be the set $ \{ f_{p(s), \iota }: \iota < \iota (p(s))\} $,  so we may write $ p(s)= ( \tr(p(s)), {\mathscr F} _{p(s))}).$
        \end{enumerate}
    \end{enumerate}
\end{definition}

The following matters only for \cite{Sh:945}. 

\begin{claim}\label{c6d}
    Assume $\mathbf{m} \in \mathbf{M}$  is, (see \ref{c5}(8))  ordinary\footnote{As $ L_\mathbf{m} $ is well founded, this is not a real restriction.}, that is the   set of elements of $ L_ \mathbf{m} $  is an ordinal $ \alpha _ \mathbf{m} = \alpha(\mathbf{m}) $  satisfying $ \beta < _{L_ \mathbf{m}} \gamma \Rightarrow \beta < \gamma $. 
     
    There is a unique object  $ \mathbf{q} = ( \bar{ u }, \bar{ \mathbb{P} } , \bar{ \mathbb{Q} } ,
    \bar{\eta})$ such that: 
    
    \begin{enumerate} 
        \item[(a)]   $\bar{u} = \bar{ u }_ \mathbf{m} $ so $ \alpha _ \mathbf{m} =\lg (\bar{ u } )$, 
        \item[(b)] $ \langle \mathbb{P} ^ \mathbf{q} _{0, \alpha }, \mathbb{Q} ^ \mathbf{q} _{0, \beta } : \alpha \le  \alpha _\mathbf{m}  , \beta < \alpha _ \mathbf{m} \rangle $ , the  $(  < \lambda )$-support iteration  such that: $ \name{ \mathbb{Q} }_ \alpha $ is essentially the forcing notion form  from \ref{c6}, 
        
        \item[(c)]   $ \mathbf{q} $ is  as in  \cite[1.8=Lz32, page 32]{Sh:945}.
    \end{enumerate}  
\end{claim}

\begin{PROOF}{\ref{c6d}}
    Follows from \ref{c13}  below. 
\end{PROOF} 

\begin{definition}\label{c7}
    1) For $p \in \bbP_{\mathbf m}$ let,
    
    \begin{enumerate}
        \item[(a)]   $\fsupp(p)$, the full support of $p$ be $\cup\{\{r_{p(s)}(\zeta):\zeta < \xi_{p,s}\} \cup \{ \{s\}:s  \in \dom(p)\}$

        \item[(b)]  $\wsupp(p)$, the wide support of $p$ be the set of $ s \in L_ \mathbf{m} $ such that for some $ t $  at least one of the following hold:  
        
        \begin{enumerate} 
            \item[$\bullet  _1$] $ s = t \in \fsupp(p) $, 
                
            \item[$ \bullet  _2$] $ t \in \fsupp(p) \setminus M, s \in t/E'_{\mathbf{m}}.$  
        \end{enumerate} 
    \end{enumerate}
    
    2) For $\mathbf m \in \mathbf M$ let $\bbP^{\mathbf m}_t = \bbP_{\mathbf m,t}$, etc., in Definition \ref{c6}.
    
    3) For $L \subseteq L_{\mathbf m}$ let $\bbP_{\mathbf m}(L) = \bbP_{\mathbf m} \rest \{p \in \bbP_{\mathbf m}:\fsupp(p) \subseteq L\}$, that is: 
    
    \begin{itemize}
        \item $ p \in \mathbb{P} _\mathbf{m} (L) $ iff $p \in \mathbb{P} _ \mathbf{m} $ and $ \fsupp(p) \subseteq L$,
        
        \item $ p \le _{\mathbb{P}_\bfm (L) }  q $ iff $ p \in \mathbb{P}_ \mathbf{m} (L) \wedge  q \in \mathbb{P} _\mathbf{m} (L)  \wedge  p \le _{\mathbb{P} _ \mathbf{m} } q $, 
    \end{itemize}
    
    4) For $\mathbf m \in \mathbf M$ and $t \in L_{\mathbf m}$ let\footnote{not used, could have used it in \ref{c13}}  $\name{\bbQ}_t = \name{\bbQ}_{\mathbf m,t}$ be the $\bbP_t$-name of
    $\bbQ_{\bar\theta} \rest \{(\nu,\name{f}[\bfG_{\bbP_{\bfm(< t)}}]): (\nu, \name{f})$ as in Definition \ref{c6}(c)$(\gamma)$ with $s$ there for $t$ here$\}.$ 
\end{definition}

\begin{claim}\label{c8}
    For $\mathbf m \in \mathbf M$ (so $\bbP_t = \bbP_{\mathbf m,t}$, etc.): 

    \begin{enumerate}
        \item[(a)]  the iteration $\mathbf q_{\mathbf m}$ is well defined, i.e. exists and is unique, 
        
        \item[(b)]
        
        \begin{enumerate}
            \item[$(\alpha)$]   if $t \in L^+_{\mathbf m}$ \then \, $\bbP_t$ is  indeed a forcing notion and is equal to $\bbP_{\mathbf m(< t)}$,
            
            \item[$(\beta)$] the $\bbP_t$-name $\name\eta_s$ does not depend on $t$ as long as $s <_{L_ \mathbf{m} } t \in L^+_{\mathbf m}$,
            
            \item[$(\gamma)$]   $\name\eta_t$ is a $\bbP_{\mathbf m(\le t)}$-name.  
        \end{enumerate}
        
        \item[(c)]    if $s <_L t$ are from $L^+_{\mathbf m}$ \then: 
    
        \begin{enumerate}
            \item[$(\alpha)$]   $p \in \bbP_s \Rightarrow p \in \bbP_t \wedge p \rest {L_{<s}} = p$,
    
            \item[$(\beta)$]   if $p,q \in \bbP_s$ then $\bbP_t \models ``p \le q" \Leftrightarrow \bbP_s \models ``p \le q"$,

            \item[$(\gamma)$]  if $p \in \bbP_t$ then $p \rest L_{<s} \in \bbP_s$ and $\bbP_t \models ``(p \rest L_{\bfm (<s)}) \le p"$,
            
            \item[$(\delta)$]   $\bbP_t \models ``p \le q" \Rightarrow \bbP_s \models ``p \rest L_{\bfm (<s)} \le q \rest L_{ \bfm (<s)}"$, 
    
            \item[$(\varepsilon)$]   $\bbP_s \lessdot \bbP_t$, moreover
    
            \item[$(\zeta)$]  $p \in \bbP_t \wedge (p \rest L_{\bfm(<s)}) \le q \in \bbP_s \Rightarrow q \cup (p \rest (L_{\bfm(<t)} \backslash L_{\bfm(<s)}) \in \bbP_t$ is a $\le$-lub of $p,q.$ 
        \end{enumerate}
        
        \item[(d)]  if $L$ is an initial segment of $L_{\mathbf m}$ \then \, $\bbP_{\mathbf m \rest L} = \bbP_{\mathbf m} \rest \{p \in \bbP_{\bfm} :\dom(p) \subseteq L$, equivalently $\fsupp(p) \subseteq L\}$; this holds in particular for $L_{ \bfm (\leq t)}$ and for $L_{\bfm(< t)}$. 
        
        \item[(e)]  if $L_1 \subseteq L_2$ are initial segments of $L_{\mathbf m}$, then the parallel of clause (b) holds replacing $\bbP_{\mathbf m,s},\bbP_{\mathbf m,t}$ by $\bbP_{\mathbf m \rest L_1},\bbP_{\mathbf m \rest L_2}$, respectively. Also the parallel of clause (c) holds.
         
        \item[(f)]  if $ p \in \mathbb{P} _\mathbf{m} $ \then \,: 
        
        \begin{enumerate} 
             \item[($ \alpha $)] $ \dom (p )$ has cardinality $ < \lambda$,
             
             \item[($ \beta $)]  $ \fsupp(p)$ has cardinality at most $ \lambda $,
             
             \item[($\gamma$)]
             
             \begin{itemize}
                 \item[$\bullet_{1}$] $\wsupp(p)$ is included in the union of  $\leq \lambda$ sets of the form $ t/E_{\mathbf{m}} $ or $ \{ t \}, $
                 
                 \item[$\bullet_{2}$] if $\bfm$ is lean then the union is even of $< \lambda$ such sets.
             \end{itemize}
        \end{enumerate} 
    \end{enumerate}
\end{claim}

\begin{PROOF}{\ref{c8}}
    Straightforward. For $t \in L^+_{\mathbf m}$, by induction on $\dep_{\mathbf m}(t)$, define $\bbP_t$ and prove the relevant parts of (a),(b),(c),(d),(e).
\end{PROOF}

Note the following:   

\begin{observation}\label{c9}
    If $\mathbf B$ is a $\lambda$-Borel function from ${}^\xi(\Pi \bar\theta)$ to $\cP(\lambda)$ or even  $\cH(\lambda^+)$ where $\xi \le \lambda$ \then \, there is a $\lambda$-Borel function $\mathbf B'$ from ${}^\xi(\Pi\bar\theta)$
    to $\bbQ_{\bar\theta}$ (so absolutely\footnote{That is, for every forcing notion $ \mathbb{P}$ which is $ \lambda $-strategically complete,  this property continue to hold in  $ \mathbf{V} ^\mathbb{P} $; here the property is that the range is as indicated; parallely below.  We could demand just   preserving the regularity of $ \lambda $ and the $ \theta _ \varepsilon $-s,}  to $\bbQ_{\bar\theta}$) such that for any $\bar\eta \in {}^\xi(\Pi\bar\theta)$ we have, absolutely:

    \begin{enumerate}
        \item[$\bullet$]  if $\mathbf B(\bar\eta) \in \bbQ_{\bar\theta}$ then $\mathbf B'(\bar\eta) = \mathbf B(\bar\eta)$, 
        
        \item[$\bullet$]  if $\mathbf B(\bar\eta) \notin \bbQ_{\bar\theta}$ then
        $\mathbf B'(\bar\eta) = (\emptyset,0_\lambda)$, the minimal member of $\bbQ_{\bar\theta}$.
    \end{enumerate}
\end{observation}

\begin{PROOF}{\ref{c9}}
    Just define  $\bfB'(\bar{\eta})$ as $\bfB(\bar{\eta})$ if  $\bfB(\bar{\eta}) \in \bbQ_{\bar{\theta}}$ and the trivial condition $(\langle \rangle, 0_{\lambda})$ otherwise.
\end{PROOF}

\begin{remark}\label{c10}
    1) A reader may wonder, e.g.:
    
    \begin{itemize}
        \item[$(*)$] if $\langle \bfB_{\alpha}: \alpha < \alpha_{*} \leq \lambda \rangle$ is a sequence of $\lambda$-Borel  subsets of $\Pi_{\varp < \lambda} \theta_{\varp}$ which form a partition (in $\bfV$), does they from a partition also in $\bfV^{\bbP}.$ 
    \end{itemize}
    
    In our case as $\bbP$ is $\lambda$-strategically complete (see \ref{c11}(3A)) the answer is obviously yes. 
    
    2) Note that in $(*)$ we cannot weaken the assumption too much because ``if $\bbP$ add a new subset to $\theta < \lambda$ this certainly faill''. Even $(< \lambda)$-strategically complete is not enough. Why? assume $\lambda$ is a Mahlo cardinal $S \subseteq \{ \theta < \lambda: \theta$ inaccessible$\}$ is stationary, such that (for transparency) $\diamond_{S}$ holds. We can find $\cT$ such that:
    
    \begin{enumerate}
        \item[$\boxplus$ (a)] $\cT$ a subtree of $({}^{\lambda}2, \lhd),$
        
        \item[(b)] $\cT$ with no $\lhd$-maximal nodes,
        
        \item[(c)] if $\delta \in \lambda \setminus S$ a limit ordinal, $\eta \in {}^{\delta}2$ and $\alpha < \delta \Rightarrow \eta \rest \alpha \in \cT,$ then $\eta \in \cT,$
        
        \item[(d)] $\cT$ has no $\lambda$-branch. 
    \end{enumerate}
    
    Let $\bfB_{0} = \{ \eta \in {}^{\lambda}2: \bigwedge_{\alpha < \lambda} \eta \rest \alpha \in \cT \}$ and $\bfB_{1} = {}^{\lambda}2.$
    
    In $\bfV$ those two $\lambda$-Borel sets form a partition: the first is empty and the second is all.  The forcing notion $\cT$ add a $\lambda$-branch to $\cT,$ hence $(\bfB_{0}, \bfB_{1})$ are no longer disjoint so fail to form a partition of ${}^{\lambda}2.$ Lastly, for $\alpha < \lambda$ the forcing notion $\cT$ is $\alpha$-strategically complete (just COM choose $p_{\alpha} \in \cT$ of length $> \alpha$). 
    
    3) Alternatively, if it suffice to us to have ``for $\alpha < \kappa$, COM do not lose in the game of length $\alpha$''  let $\lambda$ be inaccessible and $S$ as above or just such that  $\lambda \setminus S$ is fat i.e. for every club $E$ of $\lambda$ and $\alpha < \lambda$ there is an increasing continuous $h: \alpha \to E$ such that $S \cap \rang(h) = \emptyset.$ Let $\bbQ = \{\eta: \eta \in {}^{\lambda >}\lambda$ be increasing continuous with range disjoint to $S$ and $\sup(\rang(\eta_{i}))$ is not in $S\}.$ Let the sequence $\langle \eta_{i}: i < \lambda \rangle$ of pairwise $\unlhd$-incomparable be such that $\lg(\eta_{i}) \in S$ and $(\forall \alpha < \lg(\eta_{i}))[\eta_{i} \rest \alpha \in \bbQ]$ and it is dense in $\bbQ.$ For $i < \lambda,$ let $\bfB_{1 +i}$ be $\{ \nu \in {}^{\lambda} \lambda: \eta_{i} \lhd \nu \},$ so closed and $\bfB_{0} = \{\nu \in {}^{\lambda}\lambda: \nu$ is not increasing continuous$\},$ now $\langle \bfB_{i}: i < \lambda \rangle$ is as required.  
    
    4) Another avenue is to assume $\aleph_{0} < \theta = \cf(\theta) < \lambda, S_{0} \subseteq \{ \delta < \lambda: \cf(\lambda) < \theta \}, S \subseteq \{ \delta < \lambda: \cf(\delta) = \theta$ and $S_{0} \cap \delta$ is a stationary subset of $\delta \}.$ Now let $\bbQ = \{\eta: \eta \in {}^{\lambda >} 2$ and for no $\delta \leq \lg(\eta)$ we have $\delta \in S$ and for some club $E$ of $\delta$ do we have $\alpha \in E \cap S_{0} \Rightarrow \eta(\alpha) = 1 \}.$ Continue as in \ref{c10}(3).
    
    5) Note that if in \ref{z35}(1) we let INC to choose first, then \ref{c10}(a) does not work whereas in \ref{c10}(2), (3) this does no matters.  
    
    6) Anyhow in \ref{c9} this is not necessary; it is enough that being a member of $\bbQ_{\bar{\theta}}$ is a $\lambda$-Borel set. 
\end{remark}

\begin{claim}\label{c11}
    Let $\mathbf m \in \mathbf M$.

    1) If $L^+_{\mathbf m} \models ``s < t"$ \then:
    
    \begin{enumerate}
        \item[$(\alpha)$]  $\Vdash_{\bbP_{\mathbf m,t}} ``\name\eta_s \in
        \prod\limits_{\varepsilon < \lambda} \theta_\varepsilon"$, 
        
        \item[$(\beta)$]  if
        $\mathbf G \subseteq \bbP_t$ is generic over $\mathbf V, \eta_r = \name\eta_r[\mathbf G]$ for $r \in L_{\mathbf m,<t}, \, u \in \cP_{\bfm, t}$ and $\nu \in \Pi \bar\theta$ is from   $\mathbf V[\langle \eta_r:r \in u\rangle] \subseteq \mathbf
        V[\mathbf G],$ \then \, $\nu <_{J^{\bd}_\lambda} \eta_s$.
    \end{enumerate}

    2) $\bbP_{\mathbf m}$ satisfies the $\lambda^+$-c.c., and even 
    the $ \lambda ^+$-Knaster (and more).
    
    3) $\bbP_{\mathbf m}$ is $(< \lambda)$-strategically complete (even
    $\lambda$-strategically complete but not used\footnote{Recall that being $ \lambda$-strategically complete  means that a play of the game  lasts $ \lambda $ moves,  and the COM player to win needs to have  a legal choice  in each move. So COM needs just to have a common upper bound to suitable increasing sequences of length $ < \lambda $.}).
    
    3A) If $\bar p = \langle p_i:i < \delta\rangle$ is $\le_{\bbP_{\mathbf m}}$-increasing, $\delta < \lambda$ and $i < j < \delta \wedge t 
    \in \dom(p_i) \Rightarrow  \tr(p_i(t)) \triangleleft \tr(p_j(t))$ \then\footnote{But $\tr(p_i(t)) \trianglelefteq \tr(p_j(t))$ does not suffice, but if e.g. $ \cf (\delta ) < \theta _0$  it suffice.} \, $\bar p$ has a $\le_{\bbP_{\mathbf m}}$-upper bound $p$.  Moreover, $\dom(p) = \cup\{\dom(p_i):i < \delta\}$ and $s \in \dom(p_i) \Rightarrow \tr(p(s)) = \cup\{\tr(p_j(s)):j \in [i,\delta)\}$; in fact also $\fsupp(p) = \cup\{\fsupp(p_i):i < \delta\}$ and $p$ is a lub of $\bar p$. Also, we can weaken the demand above to $i < \delta \wedge s \in \dom(p_i) \Rightarrow \delta < \theta_{\varepsilon(s)}$ where we
    let $\varepsilon(s) = \sup\{\ell g(\tr(p_j(s))):j \in [i,\delta)\}$.
    
    3B) If $\zeta < \lambda$ and $L^+_{\mathbf m} \models ``s < t"$, 
    \then \, the following is a dense open subset of $\bbP_t$: $\cI_{s,t,\zeta} = \{p \in \bbP_t:s \in \dom(p)$ and
    $\tr(p(s))$ has length $\ge \zeta\}$.
    
    3C) If $p \in \bbP_{\mathbf m}$ and $\zeta < \lambda$  \then \, for some $q \in \bbP_{\mathbf m}$ we have $p \le q$ and $t \in \dom(p) \Rightarrow \tr(p(t)) \triangleleft \tr(q(t))$ and $t \in \dom(q) \Rightarrow \ell
    g(\tr(q(t))) > \zeta$.
    
    4) If $\name x$ is a ${\bbP}_{\mathbf m}$-name of a member of $\cH(\lambda^+)$, e.g. of 
    ${\bbQ}_{\bar \theta}$ (in $\mathbf V[\bbP_{\mathbf m}]$) \then \, for some $\xi \le \lambda$ and $\lambda$-Borel function $\mathbf
    B:{}^\xi(\Pi \bar\theta) \rightarrow \cH(\lambda^+)$ and a sequence $\langle r_\zeta:\zeta < \xi\rangle$ of members of $L_{\mathbf m}$ we have $\Vdash_{\bbP_{\mathbf m}} ``\name x = \mathbf B(\ldots,\name\eta_{r_\zeta},\ldots)_{\zeta < \xi}"$.
    
    4A) If $t \in L^+_{\mathbf m}$ and $u \subseteq L_{\bfm(<t)}$ and $\Vdash_{\bbP_t} ``\name y$ is a member of $\bbQ_{\bar\theta}$ from 
    $\mathbf V[\langle \name\eta_s:s \in u\rangle]"$, \then \, for some $\xi \le \lambda$ and $\lambda$-Borel functions as in \ref{c6}(c)
    $(\gamma)$, $\mathbf B_i:{}^\xi(\Pi \bar\theta)
    \rightarrow \bbQ_{\bar\theta}$ for $i < \xi$ 
    and sequence $\langle r_\zeta:\zeta < \xi\rangle$ of members of $u$ we have $\Vdash_{\bbP_t}$ ``for some $i < \xi$ we have $\name y = \mathbf
    B_i(\ldots,\name\eta_{r_\zeta},\ldots)_{\zeta < \xi}"$.
    
    5) If $\mathbf m,\mathbf n$ are equivalent \then \, $\bbP_{\mathbf m} = \bbP_{\mathbf n}$ and $\bbP_{\mathbf m,t} = \bbP_{\mathbf n,t}$ for $t \in L^+_{\mathbf m} = L^+_{\mathbf n}$.
    
    6) Assume that $p,q \in \bbP_{\mathbf m}$ are incompatible \then \, there are $q_1$ and $s$ such that:
    
    \begin{enumerate}
        \item[(a)]  $q_1 \in \bbP_{\mathbf m,s}$,
        
        \item[(b)]  $s \in \dom(p) \cap \dom(q)$,
    
        \item[(c)]  $(q \rest L_{\mathbf m,<s}) \le_{\bbP_{\mathbf m}} q_1$,
        
        \item[(d)]  $(p \rest L_{\mathbf m,<s}) \le_{\bbP_{\mathbf m}} q_1$,
                
        \item[(e)]  $q_1 \Vdash_{\bbP_{\mathbf m,<s}} ``p(s)$ and $q(s)$ are incompatible in $\bbQ_{\bar\theta}$ which means
        $\tr(p(s)) \perp \tr(q(s))$, i.e. they are $\trianglelefteq$-incomparable \oor \, $(\alpha) + (\beta) + (\gamma)$ where:
        
        \begin{enumerate}
            \item[$(\alpha)$]  $\ell g(\tr(q(s))) \ne \ell g(\tr(p(s)))$,
            
            \item[$(\beta)$]  if $\ell g(\tr(q(s))) < \ell g(\tr(p(s)))$ \then \, for some
            ordinal $\varepsilon,\ell g(\tr(q(s))) \le \varepsilon < \ell g(\tr(p(s)))$ and  $ q_1 \rest L_{\mathbf m(<s)} \Vdash_{\bbP_{\mathbf m(<s)}}
            \tr(p(s))(\varepsilon) < \name f_{q(s)}(\varepsilon)"$, 
            
            \item[$(\gamma)$]  if $\ell g(\tr(q(s))) > \ell g(\tr(p(s)))$ \then \, for some
            ordinal $\varepsilon,\ell g(\tr(q(s))) > \varepsilon \ge  \ell g(\tr(p( s  )))$ 
            and  $ q_1\rest L_{\mathbf m(<s)} \Vdash_{\bbP_{\mathbf m(<s)}} ``\tr(q(s))(\varepsilon) < \name f_{p(s)}(\varepsilon)"$.
        \end{enumerate}
    \end{enumerate}
    
    7) $\Vdash_{\bbP_{\mathbf m}} ``\mathbf V[\langle \name\eta_s:s \in L_{\mathbf m}\rangle] = \mathbf V[\name{\mathbf G}]"$.
    
    8) For $ t \in L^+_\mathbf{m} $ the sequence $ \langle \name{ \eta } _s : s \in L_{\mathbf{m}, < t }  \rangle $  is generic for $ \mathbb{P}_{\mathbf{m} , t }$; that is:   
    
     \begin{itemize}
         \item[$(*)$] if $ \mathbf{G} \subseteq \mathbb{P} _{\mathbf{m} ,t}$  is generic over $ \mathbf{V} $ and $ \eta _ s = \name{ \eta } _s [\mathbf{G} ]$  for $ s \in L_{\mathbf{m}, < t}$ then $ \mathbf{V} [\mathbf{G} ]=   \mathbf{V} [\langle \eta _s : s \in L_{\mathbf{m}, < t}  \rangle ]$.  
     \end{itemize}
     
    9) For $\bfm \in \bfM, \pi$ is an automorphism of $\bfm$ \underline{when}:
    
    \begin{enumerate}
        \item[(a)] $\pi$ is a permutation of $L_{\bfm},$
        
        \item[(b)] $\pi \rest M_{\bfm}$ is the identity,
        
        \item[(c)] if for every $s \in L_{\bfm}  \setminus M_{\bfm},$ for some $t \in L_{\bfm} \setminus M_{\bfm}$ we have $\pi \rest (s / E_{\bfm})$ is an isomorphism from $\bfm \rest (s / E_{\bfm})$ onto $\bfm \rest (t / E_{\bfm}).$
    \end{enumerate}
     
     10) In part (8), moreover, in $ \mathbf{V} [\mathbf{G} ]$, if $ \bar{ \eta }'= \langle \eta '_s: s \in L_{\mathbf{m},t}  \rangle $ and $ \eta '_s \in \Pi _{\varepsilon < \lambda } \theta _ \varepsilon $  and the set $ \{( s, \varepsilon ) : s \in L_{\mathbf{m}, < t},  \varepsilon < \lambda $  and $ \eta '_s (\varepsilon )\not= \eta _s (\varepsilon ) \}  $ has cardinality $ < \lambda $ \then \, also  $ \bar{ \eta }' $  is generic  (for $ \mathbb{P} (L_{\mathbf{m}, < t})  )$ and $\mathbf{V} [\bar{ \eta } '] = \mathbf{V} [\mathbf{G} ]$.   
\end{claim}

\begin{remark}\label{c12}
    What is the use of e.g. (6), (6A)?  See \ref{b35}(A)(b) and \ref{c13}.
\end{remark}

\begin{PROOF}{\ref{c11}}
    We prove all parts simultaneously  by induction on $\dep_{\mathbf m}$.
    
    1) For clause $(\alpha)$ for each $\mathbf m$, using the induction hypothesis and \ref{c8}(e),
    the problem is only when $\dep_{\mathbf m}(t) = \dep_{\mathbf m}-1$ and  use part (5A) proved below (and \ref{c8}(c)$(\zeta )$).   For
    clause $(\beta)$ use also part (6A) for $\bbP_{\mathbf m(<t)}$ proved below  in \ref{c8}(c)$(\zeta )$.    In both cases the proof of the parts quoted does not rely on
    part (1), (but may depend on the induction hypothesis).  

    2) Recall that $\lambda$ is strongly inaccessible.  If $p_\varepsilon \in \bbP_{\mathbf m}$ for $\varepsilon < \lambda^+$ then we can find by the $\Delta$-system lemma 
    a set $u$ and unbounded $S \subseteq \lambda^+$ such that $\varepsilon \ne \zeta \in S \Rightarrow \dom(p_\varepsilon) \cap \dom(p_\zeta) = u$ and $\langle \tr(p_\varepsilon(\beta)):\beta \in u \rangle$ is the same for all  $\varepsilon \in S$.  Now $p_\varepsilon,p_\zeta$ has a common upper bound for every $\varepsilon,\zeta \in u$, i.e. we
    define $r$ by:
    
    \begin{enumerate}
        \item[$\bullet$]  $\dom(r) = \dom(p_\varepsilon) \cup \dom(p_\zeta)$, 
        
        \item[$\bullet$]  $r(s) = p_\varepsilon(s)$ is $s \in
        \dom(p_\varepsilon) \backslash \dom(p_\zeta)$,
        
        \item[$\bullet$]  $r(s) = p_\zeta(s)$ if $s \in \dom(p_\zeta) \backslash \dom(p_\varepsilon)$,
        
        \item[$\bullet$]  if $s \in \dom(p_\varepsilon) \cap
        \dom(p_\zeta)$ \then \, $r(s) = (\tr(p_\varepsilon(s)),
        \max\{\name f_{p_\varepsilon(s)},\name f_{p_\zeta(s)}\})$.
    \end{enumerate}
    
    3) By (4), the second sentence + (5B) below which use only the induction hypothesis.
    
    3A) We define $p$ by:
    
    \begin{enumerate}
        \item[$\bullet$]  $\dom(p) =  \cup\{\dom(p_i):i < \delta\}$
        
        \item[$\bullet$]  $\tr(p(s)) = \cup\{\tr(p_i(s)):i < \delta$
        satisfies $s \in \dom(p_i)\}$
        
        \item[$\bullet$]  $\name f_{p(s)} = \sup\{\name f_{p_i(s)}:i < \delta$ satisfies $s \in \dom(p_i)\}$.
    \end{enumerate}
    
    Note that here having to really start with $\langle \name f_{p_i(s),\iota}:\iota < \iota(p_i(s))\rangle$ and get $\langle \name
    f_{p(s),\iota}:\iota < \iota(p(s))\rangle$, see \ref{c6}(c)$(\gamma)$ causes no problem, similarly in the proof of part (2) - just take the union.
    
    3B) Obvious by the definition of $\bbP_{\mathbf m}$ and \ref{c8}(c),  recalling  that  $ \mathbb{P} _{\mathbf{m} (< s)}$ is  $ (< \lambda )$-strategically complete, that is part (4) and (5B).  
    
    3C) The proof is by induction on $\rm{dp}_{\bfm}$ and is splitted in cases:
    
    \underline{Case 1}:  $\dep_{\mathbf m}$ is zero:
    
    So $L_ \mathbf{m}$ is empty.

    \underline{Case 2}:  $\dep_{\mathbf m} = \alpha +1$:
    
    Hence $L_2 = \{s \in L:\dep_{\mathbf m}(s)=\alpha\}$ is non-empty and
    letting $L_1 = L_{\mathbf m} \backslash L_2$; clearly $s \in L_1  \Rightarrow \dep_{\mathbf m}(s) < \alpha$, so $\dep_{\mathbf
    m \rest L_1} \le \alpha$.   Let $\zeta_* = \sup(\{\ell g(\tr(p(s))+1:s \in \dom(p)\} \cup \{\zeta +1\})$. Hence applying parts (3) 
    and (5B) to $\mathbf m \rest L_1$, i.e. the induction hypothesis we can find $q_1$ such that $\bbP_{\mathbf m \rest L_1} 
    \models ``p \rest L_1 \le q_1"$ and $[s \in \dom(q_1) \Rightarrow \ell g(\tr(q_1(s)) > \zeta_*]$ and $q_1$ forces a value to 
    $\name f_{p(s)} \rest \zeta_*$, call it $\rho_{s}$ for $s \in \dom(p) \cap L_2.$
    
    Define $q \in \bbP_{\mathbf m}$ by $\dom(q) =
    \dom(q_1) \cup (L_2 \cap \dom(p)),
    q \rest L_1 = q_1$ and if $s \in L_2 \cap \dom(p)$ then $q(s) = (\rho_s,\rho_s \char 94 (\name f_{p(s)} \rest [\zeta_*,\lambda))$, fully $\iota(q(s)) = \iota(p(s)), \bar{s}_{q(s), \iota} = \bar{s}_{p(s), \iota}$ and $\bfB_{q(s), \iota}$ is like $\bfB_{p(s), \iota}$ only restricting  the range to $(\Pi_{\varp < \lambda} \theta_{\varp})^{{\rm{tr}(q(s))}}$
    
    Easily $q$ is as required.
    
    \underline{Case 3}:  $\delta = \dep_{\mathbf m}$ is a limit ordinal of cofinality $\ge \lambda$:
    
    So $\alpha = \sup\{\dep_{\mathbf m}(s)+1:s \in \dom(p)\}$ is an ordinal $< \delta$ and let $L = \{s \in L_{\mathbf m}:\dep_{\mathbf m}(s) <
    \alpha\}$, so $L$ is an initial segment of $L_{\mathbf m}$ and applying the induction hypothesis to $\mathbf m \rest L,p$ we get $q$ as required in $\bbP_{\mathbf m \rest L}$ hence in $\bbP_{\mathbf m}$.
    
    \underline{Case 4}:  $\delta = \dep_{\mathbf m}$ is a limit ordinal of
    cofinality $< \lambda$:
    
    Let $\langle \alpha_i:i < \cf(\delta)\rangle$ be increasing continuous with limit $\delta$, let $\alpha_{\cf(\delta)} = \delta$ and for $i \le
    \cf(\delta)$ let $L_i := \{s \in L_{\mathbf m}:\dep_{\mathbf m}(s) < 1 + \alpha_i\}$.
    
    Now we choose $(p_i,\zeta_i)$ by induction on $i < \cf(\delta)$ such that:
    
    \begin{enumerate}
        \item[(a)]  $p_i \in \bbP_{\mathbf m \rest L_i}$, 
        
        \item[(b)]  $\bbP_{\mathbf m \rest L_i} \models ``(p \rest L_i) \le p_i$
        and $p_j \le p_i"$ when $j<i$, 
        
        \item[(c)]  if $i$ is a limit ordinal then $p_i$ is gotten from
        $\langle p_j:j<i\rangle$ as in part (4), 
        
        \item[(d)]  if $s \in \dom(p_i)$ then $\ell g(\tr(p_i(s))) \ge \zeta_i$, 
        
        \item[(e)]   $\langle \zeta_j:j < i\rangle$ is an increasing continuous sequence of ordinals $< \lambda$ and if $i$ is non-limit then $\zeta_i$ is $> \zeta$ and
        $\ge  |\dom(p)|$ and $> \sup(\{\ell
        g(\tr(p_j(s))):j<i$ and $s \in p_j\} \cup \{\ell g(\tr(p(s))):s \in \dom(p)\})$.
    \end{enumerate}
    
    Using \ref{c8} and the induction hypothesis this is easy.
    
    4) For transparency assume $\Vdash ``\name y \in
    \prod\limits_{\varepsilon < \lambda} \theta_\varepsilon"$ or just $\in {}^\lambda \mathbf V$. By parts (4) + (5B), i.e. part (3), 
    for each $\zeta < \lambda$ the following subset of $\bbP_{\mathbf m,t}$ is open and dense: $\cI_\zeta = \{p \in \bbP_{\mathbf m,t}$: for some $\nu \in  \prod\limits_{\varepsilon < \zeta} \theta_\varepsilon$ or $\in
    {}^\zeta \mathbf V$ (from $\mathbf V!$) we have 
    $p \Vdash_{\bbP_{\mathbf m,t}} ``\name y \rest \zeta =  \nu"\}$.  Clearly there is a maximal antichain $\langle p_{\zeta,\varepsilon}:\varepsilon < \xi_\zeta\rangle$ of $\bbP_{\mathbf m,t}$ included in $\cI_\zeta$ and by part (2) \wilog \,
    $\xi_\zeta \le \lambda$, the rest should be clear.  In the general case we can code $\name y$ as a subset of $\lambda$, etc.
    
    4A) This too should be clear as $\bbP_t$ satisfies the $\lambda^+$-c.c.
    
    5) Look at the definitions.  
    
    6) Using parts (4) and (5B) and the definition this is easy.
    
    7) Suppose toward contradiction that $\mathbf G_1 \ne \mathbf G_2$ are generic subsets of $\bbP_{\mathbf m}$ but $s \in L_{\mathbf m} \Rightarrow \name\eta_s[\mathbf G_1] = \eta_s = \name\eta_s[\mathbf G_2]$.
    
    Let $p_1 \in \mathbf G_1 \backslash \mathbf G_2$ hence there is $p_2 \in \mathbf G_2$ such that $p_2 \Vdash_{\bbP_{\mathbf m}} ``p_1 
    \notin \name{\mathbf G}_2"$ hence $p_1,p_2$ are incompatible.  Let $L_* = \{s \in
    L_{\mathbf m}:\mathbf G_1 \cap \bbP_{\le s} = \mathbf G_2 \cap \bbP_{\le s}\}$ so
    $L_*$ is an initial segment of $L_{\mathbf m}$.  If $L_* = L_{\mathbf m}$ we can
    easily get a contradiction, so $L_* \ne L_{\mathbf m}$ and let $r \in L_{\mathbf m} \backslash L_*$ be such that $L_{ \bfm(< r)} \subseteq L_*$.   Now as in part (8) we can get a contradiction having found a common
    upper bound to $p_1,p_2$.
    
    Alternatively use part (6).
    
    8), 9), 10) Easy too.  
\end{PROOF}

\begin{conclusion}\label{c13}
    Let $\mathbf m \in \mathbf M$ and for notational transparency  is ordinary (see \ref{c5}(8), 
    which means that  for some ordinal $\beta(*),t \in L_{\mathbf m} \Leftrightarrow t \in \beta(*)$
    and $s <_{\mathbf m} t \Rightarrow s < t$.)  
    \Then \, $\mathbf q$ is essentially\footnote{In particular - $ \mathbb{P} _{\bfm, \alpha } $ is a sub-forcing of the one we get by the iteration.} a $(< \lambda)$-support iteration of length $\beta(*)$  with $\name{\bbQ}_\alpha =
    \{(\nu,f) \in \bbQ_{\bar\theta}^{\mathbf V[\langle \name\eta_\beta:\beta < \alpha\rangle]}:\nu \triangleleft f,f = \sup\{f_\iota:\iota < \iota(\alpha)\},\iota(\alpha) < \lambda,\nu \triangleleft f_\iota$ and
    $\{f_\iota:\iota < \iota(\alpha)\} \subseteq \cup\{\bbQ_{\bar\theta}^ {\mathbf V[\langle \name\eta_\alpha:\alpha \in u\rangle]}:u \in 
    \cP_{\mathbf m,\alpha}\}\}$ with the natural order, i.e. the order of $\bbQ^{\mathbf V[\bbP_\alpha]}_{\bar\theta}$ restricted to this set. 
\end{conclusion}

\begin{PROOF}{\ref{c13}}
    Should be clear by \ref{c11}.
\end{PROOF}

Till now $(E'_{\mathbf m},M_{\mathbf m})$ have played no role and we could have omitted them.

\begin{definition}\label{c26}
    1) We define the two-place relation $\le \, = \, \le_{\mathbf M}$ on $\mathbf M$
    as follows: $\mathbf m \le \mathbf n$ \Iff:
    
    \begin{enumerate}
        \item[(a)]  $L_{\mathbf m} \subseteq L_{\mathbf n}$, as partial orders of course,
        
        \item[(b)]  $M_{\mathbf m} = M_{\mathbf n}$, (yes! equal), and $M_{\bfm}^{\rm{fat}} = M_{\bfn}^{\rm{fat}}, M_{\bfm}^{\rm{lean}} = M_{\bfn}^{\rm{lean}},$ 
        
        \item[(c)]   $u_{\mathbf m,t} = u_{\mathbf n,t} \cap L_{\mathbf m}$ and\footnote{This is the parallel in clause (d) are covered by clause (f) but see part (2).} $\cP_{\bfm, t} = \{ u \cap L_{\bfm}: u \in \cP_{\bfn, t} \}$ for $t \in M_{\mathbf m}$,
        
        \item[(d)]  $u_{\mathbf m,t} = u_{\mathbf n,t}$ and $\cP_{\bfm, t} = \cP_{\bfn, t}$ for $t \in L_{\mathbf m} \backslash M_{\mathbf m},$
        
        \item[(e)]  if  $t \in L_{\mathbf m } \backslash M_{\mathbf m}$  then $t/E'_{\mathbf m,} = t/E'_{\mathbf n,} $ hence  $E'_{\mathbf m} = E'_{\mathbf n} \rest L_{\mathbf m}.$
        
        \item[(f)] Hence,
        
        \begin{enumerate}
            \item[$\bullet_{1}$] if $t \in L_{\bfm} \setminus M_{\bfm}$ then $\cP_{\bfm, t} = \cP_{\bfn, t},$
            
            \item[$\bullet_{2}$]if $t \in M_{\bfm}$ and $s \in L_{\bfm} \setminus M_{\bfm}$ then $\{ u \in \cP_{\bfm, t}: u \subseteq s / E_{\bfm} \} = \{ u \in \cP_{\bfn, t}: u \subseteq s / E_{\bfn} \},$
            
            \item[$\bullet_{3}$] if $t \in M_{\bfm}$ then $\{ u \in \cP_{\bfm, t}: u \subseteq M_{\bfm} \} = \{ u \in \cP_{\bfn, t}: u \subseteq M_{\bfm} \}$ 
        \end{enumerate}
    \end{enumerate}
    
    2) We define the two-place relation $\le_* = \le^*_{\mathbf M}$ as in part (1) omitting clauses (b),(e) and (f); natural but not used 
    here.
     
    3) We define the two-place relation $\leq_{\bfM}^{\rm{bd}}$ by $\bfm \leq_{\bfM}^{\rm{bd}} \bfn$ iff $\bfm \leq_{\bfM} \bfn$ and both are bounded, see \ref{c5}(10). 
\end{definition}

\begin{claim}\label{c28}
    1) $\le_{\mathbf M}$ is a partial order or $\bfM$ and $\leq_{\bfM}^{\rm{bd}}$ a partial order on $\bfM_{\rm{\bd}}$ in fact is $\leq_{\bfM} \rest \bfM_{\rm{bd}}$.
    
    2) If $\langle \mathbf m_\alpha:\alpha < \delta\rangle$ is $\le_{\mathbf M}$-increasing, \then \, its union $\mathbf m_\delta$ (naturally defined) is a $\le_{\mathbf M}$-lub and $|L_{\mathbf m_\delta}| \le \Sigma\{|L_{\mathbf m_\alpha}|:\alpha < \delta\}$.
    
    2A) Similarly for $\bfM_{\rm{bd}}.$
    
    2B) We can restrict ourselves to any of the context (see \ref{c4})(2) including the fat context (there for $t \in M_{\bfm_{0}}$, $\cP_{t}$ should be $\cP(u_{\bfm_{\delta}, \ast})$, which may be different then $\bigcup \{ \cP(u_{\bfm_{\alpha}, \ast}): \alpha < \delta \}$).  
    
    3) If $\mathbf m \le_{\mathbf M} \mathbf n$ and $L \subseteq L_{\mathbf m}$ \then \, $p \in \bbP_{\mathbf m}(L) \Leftrightarrow p \in \bbP_{\mathbf n}(L)$ for
    every $p$.
    
    4) If $\bfm \le_{\mathbf{M} } \bfn$ and $\bbP_{\mathbf m} \lessdot \bbP_{\mathbf n}$ and $L \subseteq L_{\mathbf m}$ \then \, $\bbP_{\mathbf m}(L) = \bbP_{\mathbf n}(L)$ 
    as quasi orders.
    
    5) if $\bfm \leq_{\bfM} \bfn$ \underline{then}: 
    
    \begin{itemize}
        \item $\bfm$ is lean \underline{iff} $\bfn$ is lean, 
        
        \item $\bfm$ is fat \underline{if}\footnote{Why not ``iff''? because maybe $\bfm$ is fat but for some $t \in L_{\bfn} \setminus L_{\bfm}, \cP_{t} \neq [u_{t}]^{\leq \lambda}.$} $\bfn$ is fat, 
        
        \item $\bfm$ is neat \underline{if} $\bfn$ is neat,
        
        \item $\bfm$ is bounded \underline{if} $\bfn$ is. 
    \end{itemize}
\end{claim}

\begin{PROOF}{\ref{c28}}
    Easy.  
    
    1) Obvious. 
    
    2) Why is $L_{\mathbf m_\delta} := \cup\{L_{\mathbf m_\alpha}: \alpha <
    \delta\}$ well founded?  Toward contradiction assume $\bar t = \langle t_n:n < \omega\rangle$ is $<_{L_{\mathbf m_\delta}}$-decreasing.  We can
    replace $\bar t$ by any infinite sub-sequence.  So \wilog: 
    
    \begin{enumerate}
        \item[$(*)$]  either $(\alpha)$ or $(\beta)$,  where:

        \begin{enumerate}
            \item[$(\alpha)$]  for every $n < m$ there is $s_{n,m} \in M_{\mathbf m_0}$ such that $t_m <_{L_\delta} s_{n,m} <_{L_\delta} t_n,$
            
            \item[$(\beta)$]  for no $n < m$ this holds.
        \end{enumerate}
    \end{enumerate}
    
    If clause $(\alpha)$ holds, then $\langle s_{n,n+1}:n < \omega\rangle$ is a $<_{M_{\bfm_0}}$-decreasing sequence contradiction. If clause $(\beta)$ holds, then for $n < \omega$, let $\alpha(n) = \min\{\alpha:t_n \in L_{\mathbf m_\alpha}\}$; \wilog \, the sequence $ \langle \alpha ( n ) : n < \omega  \rangle $  is monotonically increasing or constant; so as $M_{\mathbf m_{\alpha(n)}} = M_{\mathbf m_0}$,  by \ref{c26}(1)(e) we get $t_n/E_{\mathbf m_{\alpha(n+1)}} =  t_{n+1}/E_{\mathbf m_{\alpha(n+1)}}$  (recalling part (1)),   hence $t_{n+1} \in  L_{\mathbf m_{\alpha(n)}}$ hence $\alpha(n+1) \le \alpha(n)$.
    So $\{t_n:n < \omega\} \subseteq M_{\bfm_{\alpha(0)}}$ hence as 
    $L_{\mathbf m_{\alpha(n)}}$ is well founded we are done.
    
    The proofs of (2A) and (2B) are easy too. 
     
    Finally for (3), (4) and (5), see the proof of $\boxplus_\alpha$ in the proof of \ref{c33s}.
\end{PROOF}

\begin{claim}\label{c31}
    $(\mathbf M,\le_{\mathbf m})$ has amalgamation. That is, if $\mathbf m_0 \le_{\mathbf M} \mathbf m_1,\mathbf m_0 \le_{\mathbf M} \mathbf
    m_2$ and $L_{\mathbf m_1} \cap L_{\mathbf m_2} = L_{\mathbf m_0}$ \then \, there
    is $\mathbf m \in \mathbf M$ such that $\mathbf m_1 \le_{\mathbf M} \mathbf
    m,\mathbf m_2 \le_{\mathbf M} \mathbf m$ and $L_{\mathbf m} = L_{\mathbf m_1}
    \cup L_{\mathbf m_2}$.  In fact, $ \mathbf{m} $ is unique, so we call it  $ \mathbf{m}_1 \oplus _{\mathbf{m}_0} \mathbf{m}_2$ 
\end{claim}

\begin{PROOF}{\ref{c31}}
    Note that by clause $(e)(\gamma)$ of Definition \ref{c4} and clause (e) of Definition \ref{c26}(1):
    
    \begin{enumerate}
        \item[$(*)_1$]  assume $(s_1 \in L_{\mathbf m_1}  \backslash L_{\mathbf m_0})
        \wedge  (s_3 \in L_{\mathbf m_2} \backslash L_{\mathbf m_0})$ and $s_2 \in
        L_{\mathbf m_0}$;
        
        \begin{enumerate}   
            \item[$\bullet$]  if $(s_1 <_{\mathbf m_1} s_2) \wedge (s_2 <_{\mathbf m_2} s_3),$ \then \,  for some $s'_1,s'_2 \in M_{\mathbf m_0}$ we have $s_{1}' \in (s / E_{\bfm}') \cap M_{\bfm}, s_{3}' \in t / E_{\bfm}' \cap M_{\bfm}, s_1 < _{\mathbf m_1} s'_1 \le_{\mathbf m_0} s_2 $  and $ s_2 \le _{\mathbf m_1} s'_2  < _{\mathbf m_2} s_3$,
            
            \item[$\bullet$]  if $s_3 <_{\mathbf m_2} s_2 \wedge s_2 <_{\mathbf m_1} s_1,$ \then \, for some $s'_1,s'_2 \in M_{\mathbf m_0}$ we have $s_{1}' \in (s / E_{\bfm}') \cap M_{\bfm}, s_{3}' \in (t / E_{\bfm}') \cap M_{\bfm},$  $s_3 < _{\mathbf m_2} s'_2 \le  _{\mathbf m_1} s_2$ and  $s_2 \le _{\mathbf m_2} s'_1   <_{\mathbf m_1} s_1$. 
        \end{enumerate}
    \end{enumerate}
    
    We now define $\mathbf m$ by:
    
    \begin{enumerate}
        \item[$(*)_2$]    
        
        \begin{enumerate}
                \item[(a)]
                    \begin{enumerate}
                        \item[$(\alpha)$]  $t \in L_{\mathbf m}$  \underline{iff} $t \in L_{\mathbf m_1} \vee t \in L_{\mathbf m_2}$, 
                        
                        \item[$(\beta)$]  $M_{\mathbf m} = M_{\mathbf m_0}$ and $M_{\bfm}^{\rm{fat}} = M_{\bfm_{0}}^{\rm{fat}}, M_{\bfm}^{\rm{lean}} = M_{\bfm_{0}}^{\rm{lean}}.$
                    \end{enumerate}
                    
                \item[(b)]  $s <_{\mathbf m} t$ \Iff \, one of the following occurs:
                
                \begin{enumerate}
                    \item[$(\alpha)$]  $s <_{\mathbf m_1} t$,
                
                    \item[$(\beta)$]  $s <_{\mathbf m_2} t$, 
            
                    \item[$(\gamma)$]  $s \in L_{\mathbf m_1} \backslash 
                    L_{\mathbf m_0}$ and $t \in L_{\mathbf m_2} \backslash L_{\mathbf m_0}$ and 
                    for some    $ r \in M_{\mathbf m_0}$,   $ s \le_{\mathbf m_1} r \wedge r \le_{\mathbf m_2} t$,
                    
                    \item[$(\delta)$]  $s \in L_{\mathbf m_2} \backslash L_{\mathbf
                    m_0}$ and $t \in L_{\mathbf m_1} \backslash L_{\mathbf m_0}$ 
                    and for some $r \in M_{\mathbf m_0},s \le_{\mathbf m_2} r \wedge r \le_{\mathbf m_1} t$.
                \end{enumerate}
                
            \item[(c)]  $u_{\mathbf m,t}$ is:
            
            \begin{enumerate}
                \item[$(\alpha)$]  $u_{\mathbf m_1,t} \cup u_{\mathbf m_2,t}$ if\footnote{but recall that for $ {\ell} \in \{ 1,2\} $ we have: 
                $t \in L_{\mathbf m_0} \backslash M_{\mathbf m_0} \Rightarrow u_{\mathbf m_\ell,t} = u_{\mathbf m_0,t} \wedge \cP_{\mathbf m_\ell,t} = \cP_{\mathbf m_0,t}.$} $t \in L_{\mathbf m_0}$, 
                 
                \item[$(\beta)$]  $u_{\mathbf m_1,t}$ if $t \in L_{\mathbf m_1} \backslash L_{\mathbf m_0}$.
                  
                \item[$(\gamma)$]  $u_{\mathbf m_2,t}$ if $t \in L_{\mathbf m_2} 
                \backslash L_{\mathbf m_0}$.
            \end{enumerate}
            
            \item[(d)] $E'_{\mathbf m} = E'_{\mathbf m_1} \cup E'_{\mathbf m_2}$.
            
            \item[(e)] $\cP_{\bfm, t}$ is: 
            
            \begin{enumerate}
                \item[$(\alpha)$] $\cP_{\bfm_{1}, t}$ if $t \in L_{\bfm_{1}} \setminus L_{\bfm_{0}},$
                
                \item[$(\beta)$] $\cP_{\bfm_{2}, t},$ if $t \in L_{\bfm_{2}} \setminus L_{\bfm_{0}},$
                
                \item[$(\gamma)$] $\cP_{\bfm_{1}, t} \cup \cP_{\bfm_{1}, t},$ if $t \in M_{\bfm_{0}}^{\rm{lean}},$ 
                
                \item[$(\delta)$] $\{ u_{1} \cup u_{2}: u_{1} \in \cP_{\bfm_{1}, t}, u_{2} \in \cP_{\bfm_{2}, t} \}$ if $t \in M_{\bfm_{0}}^{\rm{fat}},$
                
                \item[$(\varp)$] $\cP_{\bfm_{1}, t} \cup \cP_{\bfm_{1}, t}$ if $t \in \bfM_{\bfm}^{\rm{non}}.$
            \end{enumerate}
        \end{enumerate}
    \end{enumerate}

    Clearly,
    
    \begin{enumerate}
        \item[$\odot$]  $\mathbf m \in \mathbf M$ and $\mathbf m_1 \le_{\mathbf M} \mathbf m$ and $\mathbf m_2 \le_{\mathbf M} \mathbf m$.
    \end{enumerate}
     
    So we are done  proving the existence of $ \mathbf{m} $, the  uniqueness is obvious.
\end{PROOF}
  
\begin{observation}\label{c32n}
    1) For $ p,q \in \mathbb{P} _ \mathbf{m} $ we have:  $ \mathbb{P} _ \mathbf{m} \models \lqq p \le q" $  \underline{iff} $ \dom (p) \subseteq   \dom(q)$  and  $ q $ is essentially above $ p $ inside  $ \mathbb{P} _ \mathbf{m} $, (see \ref{c5}(9) or below).  
        
    2) For $p,q \in \bbP_{\mathbf m}$ the following conditions are equivalent:
    
    \begin{enumerate}
        \item[(a)]  $q \Vdash ``p \in \name{\mathbf G}_{\bbP_{\mathbf m}}"$,  that is $ q $ is essentially above $ p$, see \ref{c5}(9), 
        
        \item[(b)]  if $s \in \dom(p)$ \then \, either $s \in \dom(q)$ and $(q \rest L_{\mathbf m,< s}) \Vdash_{\bbP_{\mathbf m,<s}} ``p(s) \le q(s)"$ \underline{or} $s \notin \dom(q),\tr(p(s)) = \emptyset$ and $q
        \rest L_{\mathbf m,<s} \Vdash_{\bbP_{\mathbf m,<s}} ``p(s)$ is trivial,  i.e. $\name f_{p(s)}$ is constantly zero",
        
        \item[(c)]  $\bbP_{\mathbf m} \models ``p \le q^+"$ where $\dom(q^+) = \dom(q) \cup \dom(p)$ and $q^+(s)$ is: 
        
        \begin{enumerate}
            \item[$(\alpha)$]  $q(s)$ if $s \in \dom(q)$, 
            
            \item[$(\beta)$]  the trivial condition if $s \in \dom(p) \backslash \dom(q)$; note that $\fsupp(q^+) = \fsupp(q) \cup
            \fsupp(p)$.
        \end{enumerate}
    \end{enumerate}
\end{observation}

\begin{remark}
    We shall use this freely.
\end{remark}

\begin{PROOF}{\ref{c32n}}
    1) Easy but we shall elaborate.
     
    Let $ p, q \in \mathbb{P} _ \mathbf{m} $. 
    If $ p \le q $ then clearly $\dom(p) \subseteq \dom(q)$  and $ q \Vdash _{\mathbb{P} _ \mathbf{m}} \lqq p \in \name{ \mathbf{G} }"$, that is  $ q $  is essentially above $p$. 
    
    For the other direction assume 
    $\dom(p) \subseteq \dom(q)$ but $ \mathbb{P} _ \mathbf{m} \models \neg ( p \le q )$ and we shall prove that $ q $ is not essentially  above $p,$ this suffices.  By the present assumption there is $ s \in \dom(p) $ (hence $ s \in \dom(q)$) but $ q \rest L_{ \bfm (< s)} \not\Vdash  \lqq p(s) \le q(s)"$. 
    
    Hence there is $q_{1} \in \bbP_{\bfm (< s)}$ above $q \rest L_{\bfm(<s)}$ such that $q_{1} \Vdash_{\bbP_{\bfm}(< s)}$``$\neg (p(s) \leq q(s))$''. By the properties of $\bbQ_{\bar{\theta}}$ (and $\bbQ_{s}'$, \ref{c11}(6)) there are $q_{2}, q'$ such that: 
    
    \begin{itemize}
        \item[$(*)_{1}$]
        
        \begin{enumerate}
            \item[(a)] $q' \in \bbP_{\bfm}, \dom(q') = \{ s\},$
            
            \item[(b)] $q_{1} \leq q_{2}$ in $\bbP_{\bfm (< s)},$
            
            \item[(c)] $q_{2} \Vdash_{\bbP{\bfm(s)}}$``$q(s) \leq q'(s)$ but $q'(s), p(s)$ are incompatible''
        \end{enumerate}
    \end{itemize}
    
    Lastly, choose the function $q_{3}$ by: 
    
    \begin{itemize}
        \item[$(*)_{2}$] 
        
        \begin{enumerate}
            \item[(a)] $\dom(q_{3}) = \dom(q_{2}) \cup \dom(q),$
            
            \item[(b)] $q_{3} \rest \dom(q_{2}) = q_{2},$
            
            \item[(c)] $q_{3}(s) = q'(s),$
            
            \item[(d)] $q_{2}(t) = q(t)$ if $t \in \dom(p) \setminus (\dom(q_{2}) \cup \{ s \}). $
        \end{enumerate}
    \end{itemize}
    
    Clearly $q_{3} \in \bbP_{\bfm}, q \leq q_{3}$ and $q_{3} \Vdash_{\bbP_{\bfm}}$``$p \notin \name{\bbQ}_{\bbP_{\bfm}}$'' so we are done. 
    
    2) \underline{(a) implies (c)}:
    
    By the choice of $ q^+ $ we have $ q \le q^+$, so clause (a) implies  that $ q $ is essentially above $ p $ hence by part (1) in $ \mathbb{P} _\mathbf{m} $ we have $ p \le q^+ $ so clearly
    clause (c) holds.
     
    \underline{(c) implies (a)}:
     
    Easy.
    
    \underline{(c) iff (b)}:
    
    Obvious recalling the properties of $\bbQ_{\bar\theta}$.
\end{PROOF}

\subsection{Special sufficient conditions} \label{1B}

\begin{claim}\label{c33n}
    For $\mathbf m \in \mathbf M$, recalling \ref{c7}(3), we have $\bbP_{\mathbf m}(L_1) \lessdot \bbP_{\mathbf m}(L_3)$ \when \,:
    
    \begin{enumerate}
        \item[$(*)$]  
        
        \begin{enumerate}
            \item[(a)]  $L_2 \subseteq L_3$ are initial segments of $L_{\mathbf m},$
            
            \item[(b)]  $L_1 \subseteq L_3$ and $L_0 = L_1 \cap L_2,$
            
            \item[(c)]  $L_0$ is an initial segment of $L_1$, (follows),  
            
            \item[(d)]  $\bbP_{\mathbf m}(L_0) \lessdot \bbP_{\mathbf m}(L_2),$
            
            \item[(e)]  $L_1 \backslash L_0$ is disjoint to $M_{\mathbf m},$
            
            \item[(f)] if $t \in L_1 \backslash L_0$ then  $(t/E_{\mathbf m}) \cap L_{\mathbf m,<t} \subseteq L_1$.
        \end{enumerate}
    \end{enumerate}
\end{claim}

\begin{remark} \label{c35p}
    1) We may phrase it differently.  Recall that assuming $\bbP' \lessdot \bbP$, we say $p' \in \bbP'$ is a reduction of $p \in \bbP$ where every condition $r \in \bbP'$ stronger     than $p'$ (in $ \mathbb{P} ' $)  is still compatible (in $\bbP$) with $p$.   Let $ \mathbb{P} _ {\ell} = \mathbb{P}_\mathbf{m}   (L_ {\ell} )$.  
    Now the statement is: to find a reduction of $p_3$ from  $ \mathbb{P} _3 $ to $ \mathbb{P} _ 1 $  first consider $p_2 =$ the reduction of $ p_3 $ to $\mathbb{P} _ 2 $,    then let $p_0$ be a reduction of $p_2$ from  $ \mathbb{P} _2 $ to $ \mathbb{P}_{0}$ and finally extend $p_0$ to a condition $p_1$ by appending the information from $p_3$ on ($L_{1}$ minus $L_{0}$).
     
    2) Claim \ref{c33n} is used only in the
    proof of \ref{c33s} which is used only in the
    proof of \ref{e35} and \ref{e37}. 
\end{remark}

\begin{PROOF}{\ref{c33n}}
    As $\dep_{\mathbf m}(L_1) < \infty$ it suffices to  prove by induction on the ordinal $\gamma$ that:
    
    \begin{enumerate}
        \item[$\boxplus_\gamma$]  if $\langle L_\ell:\ell \le 3\rangle$ satisfies $(*)$ of the claim and $\dep_{\mathbf m}(L_1) \le \gamma$ \then:
        
        \begin{enumerate}
            \item[$\bullet_1$]  we have $p_1 \in \bbP_{\mathbf m}(L_1)$ and $p_1 \le q_1 \in \bbP_{\mathbf m}(L_1) \Rightarrow p_3,q_1$ are compatible in
            $\bbP_{\mathbf m}(L_3)$ \when:

            \begin{enumerate}
                \item[$(a) $]  $p_3 \in \bbP_{\mathbf m}(L_3),$
                
                \item[$(b) $]  $p_0 \in \bbP_{\mathbf m}(L_0),$
                
                \item[$(c) $]   if $p_0 \le q_0 \in \bbP_{\mathbf m}(L_0)$ 
                \then \, $p_2 := p_3 \rest L_2$ and $q_0$ are compatible in
                $P_{\mathbf m}(L_2),$
                
                \item[$(d) $] $p_1 = p_0 \cup (p_3 \rest (L_1 \backslash L_0))$.
            \end{enumerate}
            
            \item[$\bullet_2$]  $\bbP_{\mathbf m}(L_1) \lessdot \bbP_{\mathbf m}(L_3).$
        \end{enumerate}
    \end{enumerate}

    Why this holds?  Assume we have arrived to $\gamma$.
    
    \underline{Clause $\bullet_1$}:  (notice that here we do not  use the induction hypothesis): 
    Recalling clause (f) of the assumption,
    indeed, $p_1 = p_0 \cup (p_3 \rest (L_1 \backslash L_0))  \in \bbP_{\mathbf m}(L_1)$ by the definitions (clauses $\bullet_1 (a)), (b), (d) $  of $\boxplus_\gamma$), e.g. why
    $\fsupp(p_1) \subseteq L_1$?  Note that if $s \in \dom(p_3 \rest (L_1 \backslash L_0)$ then $s \in L_1 \backslash L_0 \subseteq L_1$ and
    $\{r_{p_3(s)}(\zeta):\zeta < \xi_{p(s)}\}$ is included in $L_3$ because $p \in \bbP_{\mathbf m}(L_3)$ and in $L_{<s}$ by Definition
    \ref{c6}.  As $s \in L_1 \backslash L_0$ by $(*)(e)$ we have $s \notin M_{\mathbf m}$ hence by Definition \ref{c6} we have $\{r_{p_3(s)}(\zeta):\zeta < \xi_{p(s)}\} \subseteq u_s \subseteq s/E_{\mathbf m}$.  By $(*)(f)$ we have $(s/E_{\mathbf m})
    \cap L_{\mathbf m <t} \subseteq L_1$ hence together $\{r_{p_3(s)}(\zeta):\zeta < \xi_{p(s)}\} \subseteq L_1$, and we are done proving $\fsupp(p_1) \subseteq L_1$.
    
    So the first statement in $\boxplus_\gamma \bullet_{1}$ holds; what about the second?  Toward contradiction assume $q_1$ contradicts the
    desired conclusion. Then by \ref{c11}(6) there are $s$ and $p^+_3$ such that:
    
    \begin{enumerate}
        \item[$\oplus$]  $(a) \quad s \in \dom(q_1) \cap \dom(p_3)$,
        
        \item[${{}}$]  $(b) \quad p^+_3 \in \bbP_{\mathbf m}(L_{\mathbf m,<s})$, 
        
        \item[${{}}$]  $(c) \quad p^+_3$ is above $p_3 \rest L_{\mathbf m,<s}$ and above $q_1 \rest L_{\mathbf m,<s}$,  
          
        \item[${{}}$]  $(d) \quad p^+_3 \Vdash_{\bbP_{\mathbf m,<s}} ``p_3(s),q_1(s) \in \bbQ_{\bar\theta}$ are incompatible (in
        $\bbQ_{\bar\theta}$)".
    \end{enumerate}
    
    So $s \in \dom(q_1) \subseteq L_1$ and as $L_2$ is an initial segment of $L_{\mathbf m}$ and
    clause  (c) of  $ \bullet _2 $ 
    (of $\boxplus_\gamma$),  clearly $s \in L_0$ is impossible, so $s \in \dom(q_1) \backslash L_0
    \subseteq L_1 \backslash L_0$. As $\bbP_{\mathbf m} \models ``p_1 \le q_1"$, necessarily $q_1 \rest L_{\mathbf m,<s} \Vdash_{\bbP_{\mathbf m,<s}} ``p_1(s) \le q_1(s)"$, so as
    $q_1 \rest L_{\mathbf m,<s} \le p^+_3 \rest 
    L_{\mathbf m,<s}$ (by $\oplus(c))$, 
    also $p^+_3 \rest L_{\mathbf m,<s} \Vdash_{\bbP_{\mathbf m,<s}} 
    ``p_1(s) \le q_1(s)"$.  As $s \notin L_0$ clearly $p_1(s) = p_3(s)$ by clauses $\boxplus_\gamma  \bullet _2 (b), (d) $, so  $p^+_3 \rest L_{\mathbf m,<s} \Vdash_{\bbP_{\mathbf m,<s}} 
    ``p_3(s) \le q_1(s)"$ and again easy contradiction to $\oplus(d)$.
    
    \underline{Clause $\bullet_2$}:  
    
    Clearly $\bbP_{\mathbf m}(L_1) \subseteq \bbP_{\mathbf m}(L_3)$ as quasi orders.  Next we shall prove $\bbP_{\mathbf m}(L_1) \le_{\ic}
    \bbP_{\mathbf m}(L_3)$, so assume $q_1,q_2 \in \bbP_{\mathbf m}(L_1)$ has a common upper bound $p_3$ in $\bbP_{\mathbf m}(L_3)$, and we should find one in $\bbP_{\mathbf m}(L_1)$.  Hence (see \ref{c6}(e)$(\beta)$) we have  $\dom(q_1) \cup \dom(q_2) \subseteq \dom(p_3)$.
    
    As $p_3 \rest L_2 \in \bbP_{\mathbf m}(L_2)$ by $(*)(a)$ and  we are assuming $\bbP_{\mathbf m}(L_0) \lessdot \bbP_{\mathbf m}(L_2)$,
    see $(*)(d)$ there is $p_0 \in \bbP_{\mathbf m}(L_0)$ such that $p_0 \le q \in
    \bbP_{\mathbf m}(L_0) \Rightarrow q,p_3 \rest L_2$ are compatible in $\bbP_{\mathbf m}(L_2)$ and let $p_1 = p_0 \cup (p_3 \rest (L_1 \backslash L_0))$.  By  $\boxplus_\gamma(b)$, which we have proved noting that clauses (a)-(d) of $\boxplus_\gamma   \bullet _2 $  holds, we know that  $p_1 \in \bbP_{\mathbf m}(L_1)$ and $p_1 \le p'_1 \in \bbP_{\mathbf m}(L_1) \Rightarrow p_3,p'_1$ are compatible in
    $\bbP_{\mathbf m}(L_3)$. It suffices to prove that $p_1$ is a common upper
    bound of $q_1,q_2$.
    
    We could have replaced $p_0$ by $p'_0$ whenever $p_0 \le p'_0 \in  \bbP_{\mathbf m}(L_0)$.  So \wilog \, for $\ell=1,2$ we have  $\dom(q_\ell) \cap L_0 \subseteq \dom(p_0)$ hence $\subseteq \dom(p_1)$, also recall $\dom(q_\ell) \backslash L_0 \subseteq \dom(p_3) \cap L_1 \backslash L_0$ and by the choice of $p_1$ we have $\dom(p_3) \cap L_1 \backslash L_0 \subseteq \dom(p_1) \backslash L_0$.
    
    So recalling $\dom(q_\ell) \subseteq L_1$
    together $\dom(q_\ell) \subseteq \dom(p_1)$.
    
    As we are assuming $\bbP_{\mathbf m}(L_0) \lessdot \bbP_{\mathbf m}(L_2)$ 
    \wilog \,  $p_0$ is above\footnote{Why?  It suffices to prove that there is $p'_0 \in \bbP_{\mathbf m}(L_0)$ above $p_0$ and above $q_\ell \rest L_0$. So toward contradiction assume this fails hence there  is $p^+_0 \in \bbP_{\mathbf m}(L_0)$  above $p_0$ incompatible with $q_\ell \rest L_0$.  By the choice 
    of $p_0$ we know that $p^+_0,(p_3 \rest L_2)$ are compatible, so  let $p^+_3 \in \bbP_{\mathbf m}(L_2)$ be a common upper bound.   Now $L_2$ is an initial segment of $L_{\mathbf m}$ by $(*)(a)$ and $p_3$ is above $q_\ell$ hence $p_3 \rest L_2$ is above $q_\ell \rest L_2$
    and as $q_\ell \in \bbP_{\mathbf m}(L_1),L_0 = L_1 \cap L_2$ we have $q_\ell \rest L_2 = q_\ell \rest L_0,p_3 \rest L_2$ is above $q_\ell
    \rest L_0$ but $p^+_3$ is above $p_3 \rest L_2$ hence  $p^+_3$ is above $q_\ell \rest L_2$. 
    Also $p^+_3$ is above $p^+_0$ which forces $q_\ell \rest L_0 \notin \name{\mathbf G}_{\bbP_{\mathbf m}(L_0)}$, equivalently
    $q_\ell \rest L_0 \notin \name{\mathbf G}_{\bbP_{\mathbf m}(L_2)}$, contradiction.} 
    $q_\ell \rest L_0$.  If toward contradiction we assume that $\ell \in \{1,2\}$ and $q_\ell \nleq p_1$ then for some $s \in \dom(q_\ell)$ we have $(q_\ell \rest L_{\mathbf m,<s}) \le (p_1 \rest L_{\mathbf m,<s})$ but $p_1 \rest L_{\mathbf m,<s} \nVdash_{\bbP_{\mathbf m}(L_{\mathbf m,<s})} ``q_\ell(s) \le  p_1(s)"$.  Clearly, $s \in L_0$ is impossible so $s \in L_1 \backslash L_0$ hence $s \notin M_{\mathbf m}$ by clause $(*)(e)$.
     
    Let $L'_0 = L_0,L'_1 = L_0 \cup (L_1 \cap L_{\mathbf m,<s}),L'_2 = L_2,L'_3 = L_3$ so $(L'_0,L'_1,L'_2,L'_3)$ satisfies the assumptions
    of the present claim and $\dep_{\mathbf m}(L'_1) < \gamma$,  hence by the induction hypothesis, $\bbP_{\mathbf m}(L'_1) \lessdot \bbP_{\mathbf m}(L'_3)$.
    
    Recall $s \in L_1 \backslash L_0$ hence $(s/E_{\mathbf m}) \cap L_{\mathbf m,<s}
    \subseteq L_1$ by clause (f) of the assumption of the claim, so $\fsupp(p_1 \rest \{s\}) \backslash \{s\},\fsupp(q_\ell \rest \{s\}) \backslash \{s\}$ are $\subseteq L'_1$ hence $p_1(s),q_\ell(s)$ are $\bbP_{\mathbf m}(L'_1)$-names.  So recalling $p_1 \rest L_{\mathbf m,<s} \nVdash_{\bbP_{\mathbf m}(L_{\mathbf m,<s})} ``q_\ell(s) \le p_1(s)"$ and $\bbP_{\mathbf m}(L'_1) \lessdot
    \bbP_{\mathbf m}(L'_3)$ and $L_{\mathbf m,<s} \subseteq L_3 = L'_3$ we  have $p_1 \rest L'_1 \nVdash_{\bbP_{\mathbf m}(L'_1)} ``q_\ell(s) \le
    p_1(s)"$. Hence there is $p^+_1$ such that $p_1 \rest L'_1 \le p^+_1 \in \bbP_{\mathbf m}(L'_1)$ such that $p^+_1  \Vdash_{\bbP_{\mathbf m}(L'_1)} ``q_\ell(s) \nleq p_1(s)"$ so recalling
    $\bbP_{\mathbf m}(L'_1) \lessdot \bbP(L'_3)$ we have $p^+_1 \Vdash_{\bbP_{\mathbf m}(L'_3)} ``q_\ell(s) \nleq p_1(s)"$.
    
    But by $\boxplus_{\gamma_1}  \bullet _2$ 
    for $\gamma_1 = \dep_{\mathbf m}(L'_1)$, we
    know that $p^+_1$ and $p_3 \rest L_{\mathbf m,<s}$  are compatible (in $\bbP_{\mathbf m}$,
    equivalently $\bbP_{\mathbf m}(L_{\mathbf m,<s}))$ so let $p^+_3 \in \bbP_\mathbf{m}  (L_{\mathbf m,<s})$ be a common upper bound of $p^+_1,p_3 \rest L_{\mathbf m,<s}$.  Now $p^+_3 \Vdash_{\bbP_{\mathbf m}(L'_3)} ``q_\ell(s)
    \le p_1(s)"$ because: $q_\ell \le p_3$ by the choice of $p_3$; $p_1(s) = p_3(s)$ by the choice of $p_1$ and $p_3 \le p^+_3$, see above.
    However, $p^+_3 \Vdash_{\bbP_{\mathbf m}(L'_3)}  ``q_\ell(s) \nleq p_1(s)"$ as $p^+_1
    \le p^+_3$, see above.
    
    So we have proved $\bbP_{\mathbf m}(L_1) \le_{\ic} \bbP_{\mathbf m}(L_3)$.
    
    To finish proving clause $\boxplus_\gamma  \bullet_{1},$  that is, $\bbP_{\mathbf m}(L_1) \lessdot \bbP_{\mathbf m}(L_3)$ note that clause $\boxplus_\gamma  \bullet_{1}$  does this as for every $p_3 \in  \bbP_{\mathbf m}(L_3)$ there is $p_0$ as in $\boxplus_\gamma  \bullet _1 (b), (c)$  by clause (d) of the claim's assumption and let $p_1$ be as defined in $\boxplus_\gamma   \bullet_{1} (d)$. 
\end{PROOF}

\begin{claim}\label{c33s}
    We have $\bbP_{\mathbf m_1}(L_1) = \bbP_{\mathbf m_2}(L_1)$ (i.e. as quasi orders) and $\bbP_{\mathbf m_\ell}(L_1) \lessdot \bbP_{\mathbf m_\ell}$ for $\ell=1,2$ \when:
    
    \begin{enumerate}
        \item[$\boxdot$]
        
        \begin{enumerate}
            \item[(a)]  $\mathbf m_1 \le_{\mathbf M} \mathbf m_2$, 
                
            \item[(b)]  $L_0 \subseteq L_1 \subseteq L_{\mathbf m_1}$, 
                
            \item[(c)]  $L_0$ is an initial segment of $L_1$, 
                
            \item[(d)]  $\bbP_{\mathbf m_1}(L_0) = \bbP_{\mathbf m_2}(L_0)$, 
                
            \item[(e)]  $\bbP_{\mathbf m_\ell}(L_0) \lessdot \bbP_{\mathbf m_\ell}$ for $\ell=1,2$, 
                  
            \item[(f)]   if $t \in L_1 \backslash L_0$ \then \, $t \notin M_{\mathbf m_2}$ and  $L_{\mathbf m_1,<t} \cap (t/E_{\mathbf m_1}) = L_{\mathbf m_2,<t} \cap (t/E_{\mathbf m_2}) \subseteq L_1$.
        \end{enumerate}
    \end{enumerate}
\end{claim}

\begin{remark} \label{c33t}
    Used only in the proof of $\boxplus_{4.4}$ inside the proof of \ref{e35}, so we could have used $M_\beta,\cE$ from there.
\end{remark}

\begin{PROOF}{\ref{c33s}}
    For $\ell \in \{1,2\}$ let $\bar L_\ell = \langle L_{\ell,i}:i < 4 \rangle$ be defined by:
    
    \begin{enumerate}
        \item[$\oplus_1$]  $(a) \quad L_{\ell,0} = L_0$,
        
        \item[${{}}$]  $(b) \quad L_{\ell,1} = L_1$, 
        
        \item[${{}}$]  $(c) \quad L_{\ell,2} = \{s \in L_{\mathbf m_\ell}:s \le_{\mathbf m_\ell} t$ for some $t \in L_0\}$,
    
        \item[${{}}$]  $(d) \quad L_{\ell,3} = L_{\mathbf m_\ell}$.
    \end{enumerate}
    
    Clearly,
    
    \begin{enumerate}
        \item[$\oplus_2$]  $(a) \quad (\mathbf m_\ell,\bar L_\ell)$ satisfies the assumptions of \ref{c33n} hence,
        
        \item[${{}}$]  $(b) \quad \bbP_{\mathbf m_\ell}(L_{\ell,1}) \lessdot \bbP_{\mathbf m_\ell}(L_{\ell,3})$ which means $\bbP_{\mathbf m_\ell}(L_1) \lessdot \bbP_{\mathbf m_\ell}$ for $\ell=1,2$.
    \end{enumerate}

    Why $\oplus_2$? Clearly it suffices to prove clause (a), so we just have to check clauses  $(*)(a)-(f)$ of \ref{c33n}.
    
    \underline{Clause $(*)(a)$}:
    
    By $\oplus_1(d),L_{\ell,3} = L_{\mathbf m_\ell}$ hence is an initial segment of $L_{\mathbf m_\ell}$ and by $\oplus_1(c),L_{\ell,2}$  is an initial segment of $L_{\mathbf m_\ell}$ which is $L_{\ell,3}$ so $L_{\ell,2} \subseteq L_{\ell,3}$.
    
    \underline{Clause $(*)(b)$}:
    
    For the first statement, $L_{\ell,1} \subseteq L_{\ell,3}$ is trivial by $\oplus_1(d) + \oplus_1(b) + \boxdot(a),(b)$.   The second statement says $L_{\ell,0} = L_{\ell,1}
    \cap L_{\ell,2}$.  Now $L_{\ell,0} \subseteq L_{\ell,1}$ by $\boxdot(a),(b)$ of the claim and $\oplus_1(a),(b)$.  Also $L_{\ell,0}
    \subseteq L_{\ell,2}$ holds by $\oplus_1(c)$ (and $\oplus_1(a)$). Together $L_{\ell,0} \subseteq L_{\ell,1} \cap L_{\ell,2}$; to prove
    the inverse inclusion assume $s \in L_{\ell,2} \cap L_{\ell,1}$, so as $s \in L_{\ell,2}$ by $\oplus_1(c)$ there is $t \in L_0$ such that $s
    \le_{\mathbf m_\ell} t$.  But $s \in L_{\ell,1} = L_1$ so by $\boxdot(c)$ of the claim we have $s \in L_0 = L_{\ell,0}$ as promised.
    
    \underline{Clause $(*)(c)$}:
    
    Holds by  condition $\boxdot(c)$ of the claim.
    
    \underline{Clause $(*)(d)$}:
    
    By clause $\boxdot(f)$ of the claim and $\oplus_1(c),L_{\ell,2}$ is an initial segment of $L_{\mathbf m_\ell}$, hence by \ref{c8}(e) we have $\bbP_{\mathbf m_\ell}(L_{\ell,2}) \lessdot \bbP_{\mathbf m_\ell} = \bbP_{\mathbf m_\ell}(L_{\ell,3})$.  By $\boxdot(e) \, 
    \bbP_{\mathbf m_0}(L_{\ell,0}) \lessdot \bbP_{\mathbf m_\ell}$; so together as
    $L_{\ell,0} \subseteq L_{\ell,2}$, we have $\bbP_{\mathbf m_\ell}(L_0) \lessdot \bbP_{\mathbf m_\ell}(L_{\ell,2})$. 
    
    \underline{Clauses $(*)(e),(f)$}:
    
    Hold by condition $\boxdot(f)$ of the claim.  
         
    So $\oplus_2$ holds indeed.  So now we deal with the other half.
     
    \underline{Proof of}: $\bbP_{\mathbf m_1}(L_1) = \bbP_{\mathbf m_2}(L_1)$.
     
    Let $\langle s_\alpha:\alpha < \alpha(*)\rangle$ list $L_1 \backslash L_0$ such that $s_\alpha \le_{L_{\mathbf m}} s_\beta \Rightarrow \alpha
    \le \beta$.  This is possible as $L_{\mathbf m_2}$ is well founded.
     
    Now, 
     
    \begin{enumerate}
        \item[$\oplus_3$]  for $\ell=1,2$ and $\alpha \le \alpha(*)$ let $\bar L^*_{\ell,\alpha} = \langle L^*_{\ell,\alpha,i}:i < 4\rangle$ be (but we can omit $\ell$) where:
        
        \begin{enumerate}
            \item[(a)]  $L^*_{\ell,\alpha,0} = L_0,$
            
            \item[(b)]  $L^*_{\ell,\alpha,1} = L_0 \cup \{s_\beta:\beta < \alpha\},$
            
            \item[(c)]  $L^*_{\ell,\alpha,2} = \{s \in L_{\mathbf m_\ell}:s \le_{\mathbf m_\ell} t$ for some $t \in L_0\},$
            
            \item[(d)]  $L^*_{\ell,\alpha,3} = L_{\mathbf m_\ell},$
        \end{enumerate}
        
        \item[$\oplus_4$]
            \begin{enumerate}
                \item[(a)]  $(\bar{\mathbf m}_\ell,\bar L^*_{\ell,\alpha})$ satisfies the assumption of \ref{c33n},
                
                \item[(b)]  $\bbP_{\mathbf m_\ell}(L^*_{\ell,\alpha,1}) \lessdot \bbP_{\mathbf m_\ell}(L^*_{\ell,\alpha,3})$.
            \end{enumerate}
    \end{enumerate}
    
    [Why?  Note the $\mathbf m_\ell,\langle L^*_{\ell,\alpha,i}:i < 4\rangle$ satisfies the assumptions of \ref{c33n}, hence $\oplus_2$
    holds for $\mathbf m_\ell,\bar L_{\ell,\alpha}$ for $\alpha  \le \alpha(*)$.]
    
    Now by induction on $\alpha \le \alpha(*)$ we prove that:
    
    \begin{enumerate}
        \item[$\boxplus_\alpha$]   $\bbP_{\mathbf m_1}(L^*_{\alpha,1}) =  \bbP_{\mathbf m_2}(L^*_{\alpha,1})$.
    \end{enumerate}
    
    \underline{Case 1}:  $\alpha = 0$:
    
    As $L^*_{1,\alpha,1} = L_0 = L^*_{2,\alpha,1}$, clause $\boxdot(d)$  of the assumption gives $\boxplus_\alpha$ as promised.
    
    \underline{Case 2}:  $\alpha$ a limit ordinal:
    
    Easy by the definition of the iteration.  That is, first,  if $ \dom(p) \in [L_{\mathbf{m} _2}]^{\le \lambda } $  then  we know $p
    \in \bbP_{\mathbf m_1}(L^*_{\alpha,1}) \Leftrightarrow \bigwedge\limits_{\beta < \alpha}
    [p \rest L^*_{\beta,1} \in \bbP_{\mathbf m_1}(L^*_{\beta,1})] \Leftrightarrow \bigwedge\limits_{\beta < \alpha}
    [p \rest L^*_{\beta,1} \in \bbP_{\mathbf m_2}(L^*_{\beta,1})] \Leftrightarrow p \in \bbP_{\mathbf m_2}(L^*_{\alpha,1})$; 
    second, for $p,q \in \bbP_{\mathbf m_1}(L^*_{\alpha,1})$  by the definition of the order and the induction hypothesis, $\bbP_{\mathbf m_1}(L^*_{\alpha,1}) \models ``p \le q"$ \Iff \, $\bigwedge\limits_{\beta < \alpha} [\bbP_{\mathbf m_1}(L^*_{\beta,1}) \models ``p \rest L^*_{\beta,1} \le q
    \rest L^*_{\beta,1}"]$ \Iff \, $\bigwedge\limits_{\beta < \alpha}
    [\bbP_{\mathbf m_2}(L^*_{\beta,1}) \models ``p \rest L^*_{\beta,1} \le q
    \rest L^*_{\beta,1}"]$ \Iff \,
    $\bbP_{\mathbf m_2}(L^*_{\alpha,1}) \models ``p \le q"$.
    
    So $\boxplus_\alpha$ holds.
    
    \underline{Case 3}:  $\alpha = \beta +1$:
    
    Clearly,

    \begin{enumerate}
        \item[$(*)_1$]  $p \in \bbP_{\mathbf m_1}(L^*_{\alpha,1}) \Leftrightarrow 
        p \in \bbP_{\mathbf m_2}(L^*_{\alpha,1})$.
    \end{enumerate}
    
    Next,
    
    \begin{enumerate}
        \item[$(*)_2$]  assume $p,q \in \bbP_{\mathbf m_1}(L^*_{\alpha,1})$ and we shall prove that $\bbP_{\mathbf m_1}(L^*_{\alpha,1}) \models ``p \le q"$ implies $\bbP_{\mathbf m_2}(L^*_{\alpha,1}) \models ``p \le q"$.
    \end{enumerate}
    
    [Why? If $s_\beta \notin \dom(p)$ this is obvious by the induction hypothesis.  Hence we can assume $s_\beta \in \dom(p)$, so as we are assuming $\bbP_{\mathbf m_1}(L^*_{\alpha,1}) \models ``p \le q"$, clearly $s_\beta \in \dom(q)$ hence $s_\beta \in \dom(p) \cap \dom(q)$.  First, similarly $\bbP_{\mathbf m_1}(L^*_{\beta,1}) \models ``(p \rest L^*_{\beta,1}) \le (q \rest L^*_{\beta,1})"$ and $(q \rest L^*_{\beta,1}) {\Vdash_{\bbP_{\bfm_{1} (<s_{\beta})}}} ``p(s_\beta) \le_{\bbQ_{\bar\theta}} q(s_\beta)"$ by the definition of $\bbP_{\mathbf m_1}(L^*_{\beta,1})$.   Second, as $q \rest L^*_{\beta,1} \in \bbP_{\mathbf m_1}(L^*_{\beta,1}) =  \bbP_{\mathbf m_2}(L^*_{\beta,1})$ and $\bbP_{\mathbf m_1}(L^*_{\beta,1}) \lessdot
    \bbP_{\mathbf m_1  }$  by $\oplus_4$ and $\bbP_{\mathbf m_2}(L^*_{\beta,1})
    \lessdot \bbP_{\mathbf m_2}$ by $\oplus_4 $ 
    and $p(s_\beta), q(s_\beta)$ are  $\bbP_{\mathbf m_1}(L^*_{\beta,1})$-names  (as $\fsupp(p(s_\beta),\fsupp(q(s_\beta))
    \subseteq L^*_{\beta,1})$ necessarily we have $q \rest L^*_{\beta,1}  \Vdash_{\bbP_{\mathbf m_2}} ``p(s_\beta) \le_{\bbQ_{\bar\theta}} 
    q(s_\beta)"$.  Third, as $\bbP_{\mathbf m_1}(L^*_{\beta,1})  \models ``p \rest L^*_{\beta,1} \le q \rest L^*_{\beta,1}"$, 
    by the induction hypothesis $\bbP_{\mathbf m_2}
    (L^*_{\beta,1}) \models ``p \rest L^*_{\beta,1} \le q \rest L^*_{\beta,1}$.  Fourth, by the last two sentence and the definition of
    the order in $\bbP_{\mathbf m_2}$ we have $\bbP_{\mathbf m_2} \models ``p \le
    q"$ so the conclusion of $(*)_2$ holds also in this case.
    
    Note that if  $s_\beta \in \dom(p) \backslash \dom(q)$ then $p \nleq q$, so we are done
    proving $(*)_2$.]
     
    \begin{enumerate}
        \item[$(*)_3$]  if $p,q \in \bbP_{\mathbf m_1}   (L^*_{\alpha,1})$ and $\bbP_{\mathbf m_2}(L^*_{\alpha,1}) \models ``p \le q"$ then $\bbP_{\mathbf m_1}(L^*_{\alpha,1}) \models ``p \le q"$.
    \end{enumerate}
    
    [Why?  Similar to the proof of $(*)_2$.]
    
    By $(*)_1,(*)_2,(*)_3$ clearly $\boxplus_\alpha$ holds.  So we carried the induction so $\boxplus_\alpha$ holds for every $\alpha \le \alpha(*)$ and for $\alpha = \alpha(*)$ we get $\bbP_{\mathbf m_1}(L_1) = \bbP_{\mathbf m_2}(L_2)$.  Together with $\oplus_2(b)$ in the beginning of the proof we are done.
\end{PROOF}

\subsection{On existentially closed $ \mathbf{m}'$s}\label{1C}
 
\begin{definition}\label{c34}
    0) For $L_{\mathbf{m}} \in \mathbf M$ let: 
    
    \begin{enumerate}
        \item[(a)]  $\dep^*_{\mathbf m}(L) = \cup\{\dep_{M_{\mathbf m}}(t)+1:t \in L     \cap   M_{\mathbf m}\},$ for $L \subseteq L_{\mathbf m},$
        
        \item[(b)]  $L^{\dep}_{\mathbf m,\gamma} = \{t \in L_{\mathbf m}:  t \in M_{\bfm} \Rightarrow \dep_{M_{\mathbf m}}(t) <
        \gamma$ and   $ t \in L_ \mathbf{m} \setminus M_ \mathbf{m} \Rightarrow  \sup\{\dep_{M_{\mathbf m}}(s):s \in M_{\mathbf m}$ and $s <_{L_{\mathbf m}} t\} < \gamma\}.$ So,  
        
        \begin{itemize}
            \item $L_{\bfm, \gamma}^{\rm{dp}}$ is an initial segment of $L_{\bfm},$
           
            \item $L_{\bfm, \gamma}^{\rm{dp}}$ is $\subseteq$-increasing continuous with $\gamma$ and is equal to $L_{\bfm}$ for $\gamma = \rm{dp}_{\bfm}^{*}(M_{\bfm}),$ or for $\gamma = \rm{dp}_{\bfm}^{\ast}(M_{\bfm}) +1$ (if $(\exists t \in L \setminus M)(\forall s \in M_{\bfm})(t > s)$).
        \end{itemize}
          
        \item[(c)] $L^{\deq}_{\mathbf m,\gamma} = \{t \in L_{\mathbf m}:  t \in M_ \mathbf{m} , \dep _{M_ \mathbf{m} }(t) < \gamma $  \underline{or}   $ t \in L_ \mathbf{m} \setminus M_ \mathbf{m} $  and $ \min \{\dep_{M_ \mathbf{m} }(s):  s \in M_\mathbf{m} \cup \{ \infty\},  t < s \}\le \gamma \} $, note that  (we mean):
        
            \begin{itemize}
                \item for $ \gamma = 0 $ this is 
                $ \{t \in L_ \mathbf{m} :$  if $(\exists s \in M_{\bfm})(t \leq s)$ then for some  $ s \in M_ \mathbf{m} $ we have $ t < s  $ and  $ \dep_{M_ \mathbf{m} }(s)=0\} $,
                
                \item each $L_{\bfm, \gamma}^{\rm{dq}}$ is an initial segment of $L_{\bfm},$
                
                \item  the set $L^{\deq}_{\mathbf m,\gamma} $ is $ \subseteq $-increasing with $\gamma,$ but not necessarily continuous,
                
                \item (meaningful only if we do not assume $\bfm$ is bounded, see \ref{c5}(10)) if $ t \in L_ \mathbf{m}$  then we have: for  
                no $ s \in M_ \mathbf{m}$ do we have
                $ t \leq  s $ iff \,  $ t  \in  L^\deq_{\mathbf{m}, \gamma} \setminus \cup \{L^\deq _{M_\mathbf{m} , \beta }: \beta < \gamma  \} $ 
                for $ \gamma = \rm{dp}_{\bfm}^{*}(M_{\bfm}) = \cup \{\dep_{M_ \mathbf{m} } (s)+ 1 : s \in M_ \mathbf{m} \} $. 
            \end{itemize}
    \end{enumerate}
    
    \begin{enumerate}
        \item[1) $ $ $ $ (a)] For an ordinal $\gamma$ let $\mathbf M^{\rm{bec}}_\gamma$ (here $\rm{bec}$ stands for bounded existentially closed) be the class of $\mathbf m \in \mathbf M_{\bd}$ such that, recalling Definition \ref{c7}(3):
    \end{enumerate}

    \begin{enumerate}
        \item[$ $]
        
        \begin{enumerate}
            \item[$(*)$] if $\mathbf m \le_{\mathbf M} \mathbf m_1 \le_{\mathbf M}  \mathbf m_2$ and $\bfm_{1}, \bfm_{2}$ are bounded, \then \, $\bbP_{\mathbf m_1}(L^{\dep}_{\mathbf m_1,\gamma})
            \lessdot \bbP_{\mathbf m_2}(L^{\dep}_{\mathbf m_2,\gamma})$ hence $L \subseteq L^{\dep}_{\mathbf m_1,\gamma}$ implies  $\bbP_{\mathbf m_1}(L) = \bbP_{\mathbf m_2}(L)$ (by \ref{c28}(4)).
        \end{enumerate}
    \end{enumerate}
    
    \begin{enumerate}
        \item[(b)] Let $\bfM_{\gamma}^{\rm{uec}}$ (where ueb stand for unbounded existentially closed) is defined similarly omitting ''bounded''.
        
        \item[(c)] Let $\bfM_{\gamma}^{\wec}$ (where wec stand for weakly bounded existential closed) is defined similarly replacing ``bounded'' by ``weakly bounded''. 
         
        \item[(d)] We may write $\bfM_{\gamma}^{\ec}$ for $\bfM_{\gamma}^{\rm{uec}}.$
    \end{enumerate}
    
    2) Let $\mathbf M_{\ec} = \bfM^{\ec}_\infty$ be the class of $\mathbf m$ which $\in \mathbf M^{\ec}_\gamma$ for every ordinal $\gamma;$ similarly $\bfM_{\rm{bec}} = \bfM_{\infty}^{\rm{bec}}$. 
    
    3) Let $\mathbf M^{\ec}_{\chi,\gamma} = \{\mathbf m \in \mathbf M^{\ec}_\gamma:|L_{\mathbf m}| \le \chi\}$, similarly $\mathbf M^{\ec}_{\chi,\infty}$ and for bec. 
\end{definition}

\begin{observation}\label{c37}
    1) Of course, $\mathbf M^{\ec}_{\gamma_2} \subseteq \mathbf M^{\ec}_{\gamma_1}$ and $L^{\dep}_{\mathbf m,\gamma_1} \subseteq L^{\dep}_{\mathbf m,\gamma_2}$ are initial segments of $L_{\mathbf m}$ when $\gamma_1 \le \gamma_2$.
    
    2) In \ref{c34}(1), the following are equivalent:
    
    \begin{enumerate}
        \item[(a)]  $\bbP_{\mathbf m_1}(L^{\dep}_{\mathbf m_1,\gamma}) \lessdot
        \bbP_{\mathbf m_2}(L^{\dep}_{\mathbf m_2,\gamma})$ for every $\gamma,$

        \item[(b)]  $\bbP_{\mathbf m_1} \lessdot \bbP_{\mathbf m_2}$.
    \end{enumerate}
     
    3) If $\bfm \in \bfM_{\rm{ec}}$ and $M_{\bfm} \models$``$s < t$'' (in particular, $s, t \in M_{\bfm}$) then $\Vdash_{\bbP_{\bfm}}$``$\name{\eta}_{s} < \name{\eta}_{t} \mod J_{\lambda}^{\rm{bd}}$''. Moreover, if $M_{\bfm} \models s_{i} < t$ for $i < i_{*} \leq \lambda $ and $ \bfB$ is an $i_{*}$-place $\lambda$-Borel function from $\Pi_{\varp} \theta_{\varp}$ into $\Pi_{\varp < \lambda} \theta_{\varp},$ then $\Vdash_{\bbP_{\bfm}}$``$\bfB(\cdots, \name{\eta}_{s_{i}}, \cdots)_{i < i_{*}} < \eta_{i} \mod J_{\lambda}^{\rm{bd}},$

    4) If for every $L \in [L_{\bfm}]^{\leq \lambda}$ for some $t \in M_{\bfm}$ we have $L \in \cP_{\bfm, t}$ \underline{then} (see \ref{b36}(3)) $\Vdash_{\bbP_{\bfm}}$``$\{ \name{\eta}_{t}: t \in M_{\bfm} \}$ is cofinal in $\left( \Pi_{\varp < \lambda} \theta_{\varp} \right)$''.
\end{observation}

\begin{remark}
    Recall if $\bfm$ is fat, then $L \in \cP_{\bfm, t}$ means $L \subseteq u_{\bfm, t}.$
\end{remark}

\begin{PROOF}{\ref{c37}}
    1) Easy.
     
    2) First, concerning  $(a) \Rightarrow (b)$, note that for $\gamma$ large enough we
    have $L^{\dep}_{\bfm_\ell,\gamma} = L_{\bfm_\ell}$  hence $\bbP_{\mathbf m_\ell}(L^{\dep}_{\mathbf m_\ell,\gamma}) = \bbP_{\mathbf m_\ell}$, so clear.  Second, assume (b), note that $L^{\dep}_{\mathbf m_\ell,\gamma}$ is an initial segment of $L_{\mathbf m_\ell}$ hence $\bbP_{\mathbf m_\ell}(L^{\dep}_{\mathbf m_\ell,\gamma})  \lessdot \bbP_{\mathbf m_\ell}$ for $\ell=1,2$ by \ref{c8}(c), hence we have $\bbP_{\mathbf m_1}(L^{\dep}_{\mathbf m_1,\gamma}) \lessdot \bbP_{\mathbf m_1} \lessdot \bbP_{\mathbf m_2}$, but $\lessdot$ is transitive, hence $\bbP_{\mathbf m_1}(L^{\dep}_{\mathbf m_1,\gamma}) \lessdot \bbP_{\mathbf m_2}$.  Also  $\bbP_{\mathbf m_2}(L^{\dep}_{\mathbf m_2,\gamma}) \lessdot \bbP_{\mathbf m_2}$ and $L^{\dep}_{\mathbf m_1,\gamma} \subseteq L^{\dep}_{\mathbf m_2,\gamma}$ by the definition.  Hence by the definition  $p \in \bbP_{\mathbf m_1}(L^{\dep}_{\mathbf m_1,\gamma}) \Rightarrow p \in \bbP_{\mathbf m_2}(L^{\dep}_{\mathbf m_2,\gamma})$; but lastly $(\bbQ_1 \lessdot \bbP \wedge \bbQ_2 \lessdot \bbP \wedge (\forall p)(p \in \bbQ_1 \Rightarrow p \in \bbQ_2)  \Rightarrow \bbQ_1 \lessdot \bbQ_2$ so we are done.

    3) Easy, as $\bfm \in \bfM_{\rm{ec}}$ its suffice to find $\bfn$ such that $\bfm \leq_{\bfM} \bfn$ and $\bfn$ satisfies the conclusion. So given $i_{*}, t, s_{i}$ such that $s_{i} <_{\bfm} t$ (for $i < i_{*}$) we define $\bfn \in \bfM$ as follows: 
    
    \begin{enumerate}
        \item[(a)] the set of elements of $L_{\bfn}$ are those of $L_{\bfm}$ and $r_{*},$ a new element,
        
        \item[(b)] the order $<_{\bfm}$ is defined by: $r_{1} <_{\bfm} r_{2}$ \underline{iff} $r_{1} <_{\bfm} r_{2}$ or $r_{1} \leq_{\bfm} s_{i} \wedge r_{2} = r_{*}$ for some $i < i_{*}$ or $r_{1} = r_{\ast} \wedge t \leq_{\bfm} r_{2},$
        
        \item[(c)] $M_{\bfn} = M_{\bfm},$  
        
        \item[(d)] $E_{\bfn}' = \{ (r_{1}, r_{2}): (r_{1}, r_{2}) \in E_{\bfm}' \ \text{\underline{or}} \ r_{1} = r_{*} \wedge r_{2} \in \{ s_{i}: i < i_{*} \} \cup \{ t \} \ \text{\underline{or}} \ r_{2} = r_{*} \wedge r_{1} \in \{ s_{i}: i < i_{*} \} \cup \{ t \} \},$
        
        \item[(e)] $u_{\bfn, r}$ is: 
        
        \begin{itemize}
            \item $u_{\bfm, r}$ \underline{if} $r\in L_{\bfm} \setminus \{ t \},$
            
            \item $u_{\bfm, r} \cup \{ r_{\ast} \}$ \underline{if} $r = t,$
            
            \item $\{ s_{i}: i < i_{*} \}$ \underline{if} $r = r_{*}.$
        \end{itemize}
        
        \item[(f)] $\cP_{\bfn, r}$ is: 
        
        \begin{itemize}
            \item $\cP_{\bfm, r}$ \underline{if}  $r \in L_{\bfm} \setminus \{  t \}$,
            
            \item $\cP_{\bfm, r} \cup \{ \{ r_{*} \}  \}$ \underline{if} $r = t,$ except when $t \in M_{\bfm}^{\rm{fat}},$ in which case it us $\cP(u_{\bfn, t}),$ 
            
            \item $\cP(\{ s_{i}: i < i_{*} \})$ \underline{if} $r = r_{*}.$
        \end{itemize}
    \end{enumerate}
    
    4) Easy by \ref{c11}(1)($\beta$). 
\end{PROOF}

\begin{definition}\label{c39}
    Let $\bfm \in \bfM.$
        
    1) We say $\mathbf m$ is $\mu$-wide\footnote{No real harm if we demand $ \mu \ge \lambda_{0}$ and use $ \lambda {}^{ + }.$ in part (1A).} \when \, 
    $ \mu \ge \lambda_{0}$ and
    for every $t \in L_{\mathbf m}\backslash M_{\mathbf m}$ there are $t_\alpha \in L_{\mathbf m} \backslash M_{\mathbf m}$ for $\alpha < \mu$ such that:
    
    \begin{enumerate}
        \item[(a)]  $\mathbf m \rest (t_\alpha/E_{\mathbf m})$ is isomorphic 
        to $\mathbf m \rest (t/E_{\mathbf m})$ over $M_{\mathbf m},$

        \item[(b)]  $\beta < \gamma < \mu \Rightarrow t_\beta/E''_{\mathbf m} \ne t_\gamma/E''_{\mathbf m}$.
    \end{enumerate}
    
    1A) We say $\bfm$ is wide \underline{when} it is $\lambda_{0}$-wide, see \ref{c0}. We say $\bfm$ is very wide when it is $\vert L_{\bfm} \vert$-wide.
    
    2) We say $\mathbf m$ is \emph{full} \when \,: if $\mathbf m \rest M_{\mathbf m} \le_{\mathbf M} \mathbf n$ and $E''_{\mathbf n}$ has exactly one equivalence class \then \, for some $t \in L_{\mathbf m} \backslash M_{\mathbf m}$, we have $\mathbf n$ is isomorphic to $\mathbf m \rest (t/E_{\mathbf m})$ over $M_{\mathbf m}$. Similarly for $\bfM_{\rm{wbd}}.$
    
    3) We say $\bfm$ is \emph{$\mu$-wide} or \emph{full inside} $\bfM_{\bd}$ when we restrict ourselves to $\bfM_{\bd}.$ 
\end{definition}
   
\begin{cc}\label{c41}
    1) If $\chi  = \chi  ^ \lambda  \ge 2^{\lambda_2}$ (see \ref{c0}) and $\mathbf m \in \mathbf M_{\le \chi}$ \then \, for some $\mathbf n$ we have $\mathbf m \le_{\mathbf M} \mathbf n \in \mathbf M_\chi$ and $\mathbf n \in \mathbf M_{\uec}$.
    
    2) If in addition $\bfm$ is bounded, \underline{then} for some $\bfn$ we have $\bfm \leq_{\bfM} \bfn \in \bfM_{\chi}$ and $\bfn \in \bfM_{\rm{bec}}.$
\end{cc}
   
\begin{PROOF}{\ref{c41}}
    Let $x = u$ for part (1) and $x = b$ for part (2). Let $\mathscr{X} = \mathscr{X}_{\bfm} = \{\bfn: \bfn$ is bounded if $x = b;$ and $(\bfm \rest M_{\bfm}) \leq_{\bfM} \bfn$ and $L_{\bfn} \setminus M_{\bfm} = t / E_{\bfn}''$ for some $t,$ hence $\Vert L_{\bfn} \Vert \leq \lambda_{2}\}.$

    We define a two-place relation $\cE$ on $\cX$:

    \begin{enumerate}
        \item[$(*)_0$]  $\mathbf n_1 \cE \mathbf n_2$ \Iff \, $(\mathbf n_1,\mathbf n_2 \in \cX$ and) there is an isomorphism $h$ from $\mathbf n_1$ onto $\mathbf n_2$ over $\mathbf m \rest M_{\mathbf m}$, 
        that is: an isomorphism from $L_{\mathbf n_1}$  onto $L_{\mathbf n_2}$ over $M_{\mathbf m}$ (as partial orders) such that: 
        
        \begin{enumerate} 
            \item[(a)] $t \in L_{\mathbf n_1} \Rightarrow u_{\mathbf n_2,h(t)} = \{h(s):s  \in u_{\mathbf n_1,t}\}$, 
            
            \item[(b)] $t \in L_{\mathbf n_1} \Rightarrow \cP_{\mathbf n_2,h(t)} = \{\{h(s):s \in u\}:u \in \cP_{\mathbf n_1,t}\}$, 
            
            \item[(c)] $s,t \in L_{\mathbf n_1} \Rightarrow (s E'_{\mathbf n_1 } t \Leftrightarrow h(s) E'_{\mathbf n_2 } h(t))$.   
        \end{enumerate}
    \end{enumerate} 

    Clearly $\cE$ is an equivalence relation.
    
    By our assumptions $\chi \ge 2^{\lambda_2}$ and $\mathbf n \in \cX \Rightarrow |L_{\mathbf n}| \le \lambda_2 \wedge (\forall t \in L_{\bfn})(\cP_{\bfn, t} \subseteq [L_{\bfn, < t}]^{\leq \lambda}),$  hence recalling $\lambda_2 = (\lambda_2)^\lambda$  clearly $\cE$ has $\le 2^{\lambda_2}$ equivalence classes and let $\langle \mathbf n_\alpha:\alpha <
    2^{\lambda_2} \rangle$ be a set of representatives (not necessary, but no harm
    in allowing repetitions).  
    
    By \ref{c28}(2) and \ref{c31} we can find $\mathbf n$ such that:
    
    \begin{enumerate}
        \item[$(*)_1$]  $(a) \quad \mathbf m \le_{\mathbf M} \mathbf n \in \mathbf M_\chi,$
        
        \item[${{}}$]  $(b) \quad$ for every $\alpha < 2^{\lambda_2}$ we can find $\langle t_{\alpha,i}:i < \chi\rangle$ such that:
        
        \begin{enumerate}
            \item[${{}}$]  $(\alpha) \quad t_{\alpha,i} \in L_{\mathbf n} \backslash L_{\mathbf m},$
            
            \item[${{}}$]  $(\beta) \quad (\alpha \ne \beta) \vee (i \ne j)  \Rightarrow t_{\alpha,i}/E_{\mathbf n} \ne t_{\beta,j}/E_{\mathbf n},$
            
            \item[${{}}$]  $(\gamma) \quad \mathbf n \rest  (t_{\alpha,i}/E_{\mathbf n})$ is $\cE$-equivalent to $\mathbf n_\alpha$,
            see \ref{c5}(0) on $t_{\alpha,i}/E_{\mathbf n}$.
        \end{enumerate}
    \end{enumerate}

    We shall now  prove that $\mathbf n$ is as required. Let $\mathbf n \le_{\mathbf M} \mathbf n_1 \le_{\mathbf M} \mathbf n_2, $ and $ \bfn_{1}, \bfn_{2}$ are bounded when $x = b$   and define $\cF$ as the set of functions $f$ such that some $L_1, L_2$ satisfy:
    
    \begin{enumerate}
        \item[$(*)_2$]  $(a) \quad L_\ell \subseteq L_{\mathbf n_\ell}$,
        
        \item[${{}}$]    $(b) \quad M_ \mathbf{m} = M_ \mathbf{n} \subseteq L_1  \cap L_2$,
        
        \item[${{}}$]  $(c) \quad L_\ell$ has cardinality $\le \lambda_2$,
        
        \item[${{}}$]  $(d) \quad L_\ell$ is $E_{\mathbf n_{\ell} }$-closed, i.e.  $M_ \mathbf{m} \subseteq L_{\ell}$ and $ t
        \in L_\ell \backslash M_{\mathbf m} \Rightarrow t/E_{\mathbf n_2} \subseteq L_\ell$, 
          
        \item[${{}}$]  $(e) \quad f$ is an isomorphism from  $\mathbf n_2 \rest L_1$ onto $\mathbf n_2 \rest L_2$  over $M_{\mathbf m}$, i.e.:  
        
        \begin{enumerate}
        \item[${{\bullet_1}}$]  $ \quad f$ is a one-to-one mapping from $L_1$ onto $L_2$, 
        
        \item[${{\bullet _2}}$]  $\quad f \rest M_{\mathbf m}$ is the identity,
        
        \item[${{\bullet _3}}$]  $ \quad f$ maps $\le_{\mathbf n_2} \rest L_1$ onto $\le_{\mathbf n_2} \rest L_2$, 
        
        \item[${{\bullet _4}}$] $\quad $  $ s E'_{\mathbf{n} _1}  t  \Leftrightarrow f(s) E'_{\mathbf n_2} f(t)$,

        \item[${{\bullet _5}}$] $\quad$ for $s,t \in L_1$ we have $s \in u_{\mathbf n_2,t} \Leftrightarrow f(s) \in u_{\mathbf n_2,f(t)}.$
          
        \item[$\bullet_{6}$] for $t \in L_{1}$ we have $\cP_{\bfn_{2}, f(t)} = \{ \{ f(s): s \in u \}: u \in \cP_{\bfn_{1}, t}, u \subseteq L_{1} \}.$ 
        \end{enumerate}
    \end{enumerate}
    
    Clearly,
    
    \begin{enumerate}
        \item[$(*)_3$]  if $f \in \cF$ and $L' \subseteq L_{\mathbf n_1}, L'' \subseteq L_{\mathbf n_2}$ and $|L'| + |L''| \le \lambda_2$ \then \, for some $g \in \cF$ extending $f$ we have: 
        
        \begin{enumerate} 
            \item[(a)] $L' \subseteq \dom(g) $,
            
            \item[(b)] $ L'' \subseteq \rang(g)$,  
            
            \item[(c)] $\rang(g) \backslash (L'' \cup  \rang(f)) \subseteq L_{\mathbf n_{2}}$,
            
            \item[(d)] $\dom(g) \backslash (L'   
            \cup  \dom(f)) \subseteq L_{\mathbf n_{1}}.$
        \end{enumerate} 
    \end{enumerate}
    
    We can finish as in the parallel of the Tarski-Vaught criterion for $\bbL_{\infty,\lambda^+_2}$  but we shall elaborate.  That is, first we can prove by
    induction on the ordinal $\gamma < |L_{\mathbf n_2}|^+$ (and in fact just $\gamma < \|M_{\mathbf n_2}\|^+$) that $ (*)_4-(*)_6 $ below holds: 
    
    \begin{enumerate}
        \item[$(*)_4$]  letting $L_\gamma = L^{\dep}_{\mathbf n_2,\gamma}$, if 
        $g \in \cF$ \then: 
        
        \begin{enumerate}
            \item[(a)] $g$ maps $\dom(g) \cap L_\gamma$ onto $\rang(g) \cap L_\gamma$, 
            
            \item[(b)]  $g$ induces an isomorphism $\hat g$ from  $\bbP_{\mathbf n_2}(\dom(g) \cap L_\gamma)$ onto $\bbP_{\mathbf n_2}(\rang(g) \cap L_\gamma)$, that is: $\hat
            g(p)=q$ \underline{iff}:
            
            \begin{enumerate}
                \item[$(\alpha)$]  $p \in \bbP_{\mathbf n_2}(\dom(g) \cap L_\gamma)$, 
                
                \item[$(\beta)$]  $q \in \bbP_{\mathbf n_2}(\rang(g) \cap L_\gamma)$, 
                
                \item[$(\gamma)$]  $g$ maps $\dom(p)$ onto $\dom(q)$ and $s \in \dom(p) \Rightarrow \tr(p(s)) = \tr(q(g(s)))$, 
                
                \item[$(\delta)$]   if $s \in \dom(g),g(s) = t \in \rang(g)$ and $\name f_{p(s)} = \mathbf B_{p(s)}(\ldots, \name\eta_{r_{p(s)}(\zeta)},\ldots)_{\zeta < \xi_{p(s)}}$ and $f_{q(t)} = \mathbf B_{q(t)}(\ldots,\name\eta_{r_{q(t)}(\zeta)},\ldots)_{\zeta < \xi_{q(t)}}$ \then \, $\xi_{q(t)} = \xi_{p(s)},\mathbf B_{q(t)} = \mathbf
                B_{p(s)}$ and $\zeta < \xi_{p(s)} \Rightarrow r_{q(t)}(\zeta) = g(r_{p(s)}(\zeta))$, 
                
                \item[$(\varepsilon)$]   moreover in $(\delta)$ we have $\iota(s,p) = \iota(t,q)$ and if $\iota < \iota(s,p)$ then $w_{p,s,\iota} = w_{q,t,\iota},\mathbf B_{p(s),\iota} = \mathbf B_{q(t),\iota}.$
            \end{enumerate}
        \end{enumerate}
    \end{enumerate}
    
    [Why? We use freely \ref{c11}(9). Let $ \chi_* $ be such that $\gamma,  g, \mathbf{n}, \mathbf{n}_1,\mathbf{n} _2  \in {\mathscr H} (\chi _*)$.  Let $ \mathfrak{A} \prec ({\mathscr H} (\chi_*), \in ) $ be such that  $\gamma, g, \mathbf{n}, \mathbf{n}  _1,\mathbf{n} _2  \in \mathfrak{A} , \| \mathfrak{A} \| = \chi , \chi + 1  \subseteq  \mathfrak{A} $  and $ [\mathfrak{A} ]^{\le \lambda }\subseteq \mathfrak{A} $, (hence $ \mathfrak{A} \prec _{\mathbb{L} _{\lambda ^+, \lambda ^+}} ({\mathscr H} (\chi_*), \in $)). 
    
    For $ {\ell} = 1,2 $ let $ L^*_ {\ell} = L_{\mathbf{n} _ {\ell} }\cap \mathfrak{A}$ and $ \mathbf{n} ^*_{\ell} = \mathbf{n} _ {\ell} \upharpoonright L^*_{\ell} $, so by absoluteness $ \mathbb{P} _{\mathbf{n} ^*_{\ell}}  (L_{\mathbf{n} ^*_{\ell}}) = \mathbb{P} _{\mathbf{n} _{\ell} }( L_{\mathbf{n}^* _{\ell} })$ hence  $ \mathbb{P} _{\mathbf{n} ^*_{\ell} }(L_{\mathbf{n} ^*_{\ell} })  \lessdot \mathbb{P} _{\mathbf{n} _{\ell} }(L_{\mathbf{n} _{\ell} })$.
    By the choice of $ \mathbf{n} $ as very wide and full (see Definition \ref{c39}), also $\mathbf{n} \upharpoonright (\mathfrak{A} \cap L_ \mathbf{n} )$ is very wide and full of cardinality $ \chi  $. But we have  $ \mathbf{n} \upharpoonright (\mathfrak{A} \cap L_ \mathbf{n} ) 
    \le \mathbf{n} _2 \upharpoonright (\mathfrak{A} \cap L_ {\mathbf{n}_2} ) $   both of cardinality $ \chi $ hence  also $ \mathbf{n} ^*_2 $   is very wide  and full (see Definition \ref{c4}) of cardinality $ \chi$.  Now  easily $ g $ can be extended to an automorphism of  $\mathbf{n} ^*_2 $. The promised statement now follows.]
    
    Second, 
     
    \begin{enumerate}
        \item[$(*)_5$]  $\bbP_{\mathbf n_2}(L_\gamma \cap L_{\mathbf n_1})  \lessdot \bbP_{\mathbf n_2}(L_\gamma)$.
    \end{enumerate}
    
    [Why? By\footnote{Can repeat the proof of $(*)_4 $  but for variety we give another proof.} the definitions and the induction hypothesis $\bbP_{\mathbf n_2} (L_\gamma \cap L_{\mathbf n_1}) \subseteq  \bbP_{\mathbf n_2}(L_\gamma)$ as quasi orders.  
    
    Also if $p_1,p_2 \in \bbP_{\mathbf n_2}(L_\gamma \cap L_{\mathbf n_1})$ are compatible in
    $\bbP_{\mathbf n_2}(L_\gamma)$ let $q \in \bbP_{\mathbf n_2}(L_\gamma)$  be a common upper bound there.  We can find an $E_{\mathbf n_2}$-closed $L' \subseteq  L_{\mathbf n_1}$ of cardinality $\le \lambda_2$ (recalling $\mathbf n \in \cX \Rightarrow |L_{\mathbf n}| \le
    \lambda_2$) such that $p_1,p_2 \in \bbP_{\mathbf n_1}(L')$ and  $E_{\mathbf n_2}$-closed $L'' \subseteq L_{\bfn_{1}}$ of cardinality $\le \lambda_2$ such that $L' \subseteq L''$ and $q \in \bbP_{\mathbf n_2}(L'')$. Now we can find $f_1 \in \cF$ such that $\dom(f_1) = \cup\{t/E_{\mathbf n_2}:t \in L'\} \cup M_{\bfm}$ recalling that $t/E_{\mathbf m} \supseteq M_{\mathbf m}$, see \ref{c5}(0) and $f_1$ is the identity. Then by $(*)_3$ we can find $f_2 \in \cF$ extending $f_1$ with $\dom(f_2) = \cup\{t/E_{\mathbf n_2}:t \in L''\}$ and $\rang(f_2) \backslash \rang(f_1) \subseteq L_{\mathbf n_1}$.  So recalling $(*)_4(b)$ applied to $f$ we have $\bbP_{\mathbf n_2}
    \models ``(p_1 \le \hat f_2(q)) \wedge (p_2 \le
    \hat f_2(q))"$ and $\hat f_2(q) \in \bbP_{\mathbf n_2}(L_\gamma \cap L_{\mathbf n_1})$ recalling
    $(*)_4$.  So $p_1,p_2$ are compatible also in $\bbP_{\mathbf n_2}(L_\gamma \cap L_{\mathbf n_1})$.  
    Obviously, if $p_1,p_2 \in \bbP_{\mathbf n_2}(L_\gamma \cap L_{\mathbf n_1})$ are compatible in $\bbP_{\mathbf n_2}(L_\gamma \cap L_{\mathbf n_1  }),$ say, $q$ witnesses, then
    $q$ is a common upper bound of  $p_1,p_2$ in $\bbP_{\mathbf n_2} (L_\gamma)$.
    
    So every antichain of $\bbP_{\mathbf n_2}(L_\gamma \cap L_{\mathbf n_1})$ is an antichain of $\bbP_{\mathbf n_2}(L_\gamma)$.  Similarly to the above every maximal
    antichain of $\bbP_{\mathbf n_2}(L_\gamma \cap L_{\mathbf n_1})$ is a maximal antichain of $\bbP_{\mathbf n_2}(L_\gamma)$; similarly for the
    other direction.  So we are done.]
    
    \begin{enumerate}
        \item[$(*)_6$]  $\bbP_{\mathbf n_1}(L_\gamma \cap L_{\mathbf n_1}) = \bbP_{\mathbf n_2}(L_\gamma \cap L_{\mathbf n_1}) \lessdot \bbP_{\mathbf n_2}(L_\gamma)$.
    \end{enumerate}
    
    [Why?  We prove this by induction on $\gamma$, as in proving the Tarski-Vaught criterion is sufficient (we shall  elaborate  later in the proof of \ref{e35}, more specifically $\boxplus_4$ proves a similar statement in detail with weaker assumptions).]
    
    Hence (using $\gamma = |L_{\mathbf n_2}|^+$),
    
    \begin{enumerate}
        \item[$(*)_7$]  $\bbP_{\mathbf n_1} \lessdot \bbP_{\mathbf n_2}$.
    \end{enumerate}
     
    Hence for every $L \subseteq L_{\mathbf n_1}$ by \ref{c28}(4)  we have $\bbP_{\mathbf n_1}(L) = 
    \bbP_{\mathbf n_2}(L)$ as required for $\bfn \in \bfM_{\ec}$, see Definition \ref{c34}.
\end{PROOF}

\begin{definition}\label{c58}
    1) For $\bfm \in M,$ let $\bfn = \bfm^{[\rm{bd}]}$ be $\bfm \rest L_{\bfm}^{\rm{bd}},$ where $L_{\bfm}^{\bd} = \{s \in L_{\bfm}:$ for some $t \in M_{\bfm}$ we have $s / E_{\bfm}' \subseteq L_{\bfm(\leq t)}$ or just for some $\cX \in [M_{\bfm}]^{\leq \lambda}$ we have $s / E_{\bfm}' \subseteq \bigcup \{ L_{\bfm(\leq t)}: t \in X \} \}.$ 
    
    1A) For $\bfm \in \bfM,$ let $\bfn := \bfm^{[\rm{wbd}]}$ be $\bfm \rest L_{\bfm}^{\rm{wbd}},$ where $L_{\bfm}^{\rm{wbd}} := \bigcup \{ L_{\bfm(\leq t)}: t \in M_{\bfm} \}.$ 
    
    2) Assume $\bfn_{1} \leq_{\bfM} \bfm_{1}, \bfn_{1} \leq_{\bfM} \bfn_{2}$ and $L_{\bfn_{2}} \cap L_{\bfm_{1}} = L_{\bfn_{1}}.$ Then let $\bfm_{2} = \bfn_{2} \underset{ \bfn_{1} }\oplus \bfm_{1}$ be defined by: 
    
    \begin{enumerate}
        \item[(a)] the set of elements of $L_{\bfm_{2}}$ is $L_{\bfm_{1}} \cup L_{\bfm_{2}},$
        
        \item[(b)] $<_{\bfm_{2}}$ is the transitive closure of $<_{\bfn_{2}} \cup <_{\bfm_{1}},$ 
        
        \item[(c)] $E_{\bfm_{2}}' = E_{\bfn_{2}}' \cup E_{\bfm_{1}}',$
            
        \item[(d)] $u_{\bfm_{2}, t}$ is: 
            
        \begin{itemize}
            \item $u_{\bfm_{2}, t}$ if $t \in L_{\bfn_{2}} \setminus L_{\bfn_{1}},$
                
            \item $u_{\bfm_{1}, t}$ if $t \in L_{\bfm_{1}} \setminus L_{\bfn_{1}},$
                
            \item $u_{\bfn_{2}} \cup u_{\bfm_{1}, t}$ if $t \in L_{\bfm_{1}}$ (so in $u_{\bfn_{1}, t}$ if $L \in L_{\bfm_{1}} \setminus M_{\bfm_{0}}$).
        \end{itemize}
            
        \item[(e)] $\cP_{\bfm, t}$ is defined naturally, that is: 
            
        \begin{itemize}
            \item $\cP_{\bfm_{2}, t}$ if $t \in L_{\bfn_{2}} \setminus L_{\bfn_{1}},$
                
            \item $\cP_{\bfn_{2}, t}$ if $t \in L_{\bfn_{2}} \setminus L_{\bfn_{1}},$
                
            \item $\cP_{\bfn_{2}} \cup \cP_{\bfm_{1}, t}$ if $t \in L_{\bfm_{1}}$ except when $t \in M_{\bfm_{1}}^{\fat}$ (so in $\cP_{\bfn_{1}, t}$ if $L \in L_{\bfm_{1}} \setminus M_{\bfm_{0}}$),
                
            \item $[u_{\bfm_{2}, t}]^{\leq \lambda}$ if $t \in M_{\bfm_{1}}^{\fat}.$
        \end{itemize}
    \end{enumerate}
\end{definition}

\begin{claim}\label{c61}
    1) In \ref{c58}(1) indeed $\bfm^{[\bd]} \in \bfM$ and moreover it is bounded. 
    
    2) If $\bfm \in \bfM, \bfm$ is bounded \underline{iff} $\bfm = \bfm^{[\bd]}.$
    
    3) In \ref{c58}(2) indeed $\bfm_{2} = \bfn_{2} \underset{\bfn_{1}} \oplus \bfm_{1}$ belongs to $M, \bfm_{1} \leq_{\bfM} \bfm_{2}$ and $\bfn_{1} = \bfm_{1}^{[\bd]} \leq_{\bfM} \bfn_{2}^{\rm{bd}} \Rightarrow \bfm_{2}^{[\bd]} = \bfn_{2}^{[\bd]}.$ 
    
    4) In \ref{c58}(1) we can add $\bfn^{[\bd]} \in M_{\rm{bec}}.$
\end{claim}

\begin{PROOF}{\ref{c61}}
    Easy, e.g.
    
    For part (3) we are given $\bfm \in \bfM$ and let $\bfn$ be as constructed above for $x = u.$ Clearly $\bfn^{[\bd]}$ is as constructed above for $ x = b,$ so we are done. 
\end{PROOF}

\newpage

\section{The Corrected $\bbP_{\mathbf m}$}\label{2}

\begin{discussion}\label{b2}
    Here  for $L \subseteq L_{\bfm},$ we define $\bbP_{\bfm}[L]$, the complete subforcing of the completion of $\bbP_{\bfm}$ generated by $\langle \name\eta_s:s \in L\rangle$, the central case is $L=M_{\bfm}$, of course.
\end{discussion}

\begin{definition}\label{b5}
    Let $\bbP$ be a forcing notion and $Y \subseteq \bbP$ and $\chi$ a regular cardinal.
    
    1) Let $\bbL_\chi(Y)$ be the set of sentences formed from $\{p:p \in \bbP\}$ closing under the operations $\neg p$ and $\bigwedge\limits_{i < \alpha} p_i$, for $\alpha < \chi$; so (infinitary) propositional logic.
    
    2) For a directed  $\mathbf G \subseteq \bbP$ and $\psi \in \bbL_\chi(Y)$ we define
    the truth value $\psi[\mathbf G]$ naturally (by induction on $\psi$ starting with $p[\mathbf G] = \true \Leftrightarrow p \in \mathbf G$).
    
    3) Let $\bbL^+_\chi(Y,\bbP)$, the $\bbL_\chi$-closure of $Y$ for $\bbP$, (where $Y \subseteq \bbP$; if $Y = \bbP$ we may omit $Y$) be the following partial order:
    
    \begin{enumerate}
        \item[$\bullet$]   set of elements $\{\psi \in \bbL_\chi(Y,\bbP): \, \nVdash_{\bbP} ``\psi[\name{\mathbf G}] =$ false"$\},$
        
        \item[$\bullet$]  the order $\psi_1 \le \psi_2$ iff $\Vdash_{\bbP}$ ``if $\psi_2[\name{\mathbf G}] =$ true then $\psi_1[\name{\mathbf G}] =$ true".
    \end{enumerate}
    
    4) The completion of $\bbP$ is the $\bbL_\chi$-closure of $\bbP$ which is  $\bbL_\chi^{+}(\bbP)= \mathbb{L} ^+_\chi (\mathbb{P}, \mathbb{P} )$  where $\chi$ is minimal such that $\bbP$ satisfies the $\chi$-c.c.
\end{definition}

\begin{claim}\label{b8}
    For a cardinal $\chi$ and forcing notion $\bbP$ and $Y \subseteq \bbP$ we have:

    \begin{enumerate}
        \item[(a)]  $\bbL^+_\chi(Y,\bbP)$ is a forcing notion, 
        
        \item[(b)]  $\bbP \lessdot \bbL^+_\chi(\bbP)$ under the natural identification\footnote{Pedantically $\bbP \lessdot' \bbL^+_{\chi_{1}}[\bbP]$, see \ref{b11}(8), because $\bbL^+_\chi[\bbP] \models ``p \le q"$ \Iff \, $q \Vdash_{\bbP} ``p \in \name{\mathbf G}_{\bbP}"$.},
        
        \item[(c)]  $\bbL^+_\chi(Y,\bbP) \lessdot \bbL^+_\chi(\bbP)$, 
        
        \item[(d)]  $\bbL^+_{\chi_1}(Y,\bbP) \lessdot \bbL^+_{\chi_2}(Y,\bbP)$ when $\chi_1 \le \chi_2$ are regular, 
        
        \item[(e)]  if $\bbP$ satisfies the $\chi_1$-c.c. and $\chi_1 < \chi_2$ are regular, \then \, $\bbL^+_{\chi_1}(Y,\bbP)$ is essentially equal to $\bbL^+_{\chi_2}(Y,\bbP)$, i.e. up to the natural equivalence of elements in a quasi order, 
        
        \item[(f)]  if $Y = \bbP$ then $\bbP$ is a dense subset of $\bbL^+_\chi(\bbP)$.
    \end{enumerate}
\end{claim}

\begin{PROOF}{\ref{b8}}
Easy.
\end{PROOF} 

\begin{definition}\label{b11}
    Let $\mathbf m \in \mathbf M$.
    
    1) For $t \in L_{\mathbf m},\varepsilon < \lambda$ and $\eta \in \prod\limits_{i < \varepsilon} \theta_i$ let $p=p^*_{t,\eta} \in
    \bbP_{\mathbf m}$ be the function with domain $\{t\}$ such that $p(t) = (\eta,\eta \char 94 0_\lambda)$, i.e. $f_{p(t)} \in \prod\limits_{i <\lambda} \theta_i$ is defined by $f_{p(t)}(\varepsilon)$ is $\eta(\varepsilon)$ if $\varepsilon < \ell g(\eta)$ and is zero otherwise.
    
    2) For $L \subseteq L_{\mathbf m}$ let $Y_L = Y_{\mathbf m,L} = \{p^*_{t,\eta}:t \in L$ and $\eta \in \prod\limits_{\varepsilon < \zeta} \theta_\varepsilon$ for some $\zeta < \lambda\}$.
    
    3) For $L \subseteq L_{\mathbf m}$ let $\bbP_{\mathbf m}[L]$ be $\bbL^+_{\lambda_{0}} (Y_L,\bbP_{\mathbf m})$, see Definition \ref{b5}(3) and Hypothesis \ref{c2}(4) on $\lambda_{0}.$ 
    
    4) For $L \subseteq L_{\mathbf m}$ let $\bbP_{\mathbf m}(L) =  \bbP_{\mathbf m} \rest \{p \in \bbP_{\mathbf m}:\fsupp(p) \subseteq L\}$, see Definition \ref{c7}(1), recalling \ref{c7}(2),(3).
    
    5) $\bbP'_{\mathbf m}$ is the partial order with the same set of elements as $\bbP_{\mathbf m}$ and $\le_{\bbP'_{\mathbf m}} = \{(p,q):p,q \in \bbP_{\mathbf m}$ and no $r$ above $q$ is incompatible with $p\}$ and $\bbP'_{\mathbf m}(L) = \bbP'_{\mathbf m}  \rest \{p \in \bbP_{\mathbf m}:\fsupp(p) \subseteq L\}$, we may ``forget" the distinction\footnote{Really the only difference is the possibility that $\dom(p) \nsubseteq \dom(q)$, see \ref{c32n}.}.
    
    6) For quasi orders $\bbQ_1,\bbQ_2$ let $\bbQ_1 
    \subseteq' \bbQ_2$ mean that:
    
    \begin{enumerate}
        \item[(a)] $s \in \bbQ_1 \Rightarrow s \in \bbQ_2$
        
        \item[(b)] $s \le_{\bbQ_1} t \Rightarrow s \le_{\bbQ_2} t$.
    \end{enumerate}
    
    7) For quasi orders $\bbQ_1,\bbQ_2$ let $\bbQ_1 \subseteq'_{\ic} \bbQ_2$ means that $\bbQ_1 \subseteq' \bbQ_2$ and
    
    \begin{enumerate}
        \item[(c)]   if $s,t \in \bbQ_1$ are incompatible in $\bbQ_1$ then they
        are incompatible in $\bbQ_2$.
    \end{enumerate}
    
    8) We define $\lessdot'$ similarly, that is $ \mathbb{Q}_1  \subseteq  '_{\ic} \mathbb{Q} _2 $ and every  maximal antichain of $\mathbb{Q} _1$ is a maximal antichain of $ \mathbb{Q} _2$. 
    
    9) Let $\bbQ_1 \subseteq'_{\eq} \bbQ_2$ means that $\bbQ_1 \lessdot' \bbQ_2$ and for every $p \in \bbQ_2$ there is $q \in \bbQ_1$ equivalent to it which means $\Vdash_{\bbQ_2} ``p \in \name{\bfG}_{\bbQ_2}$ \Iff \, $q \in \name{\bfG}_{\bbQ_1}$. 
\end{definition}

\begin{claim}\label{b14}
    Let $\mathbf m \in \mathbf M$ and $L \subseteq L_{\mathbf m}$.
     
    1) $\bbP_{\mathbf m}[L_{\mathbf m}]$ is equivalent to $\bbP_{\mathbf m}$ as forcing notions, in fact, $\bbP_{\mathbf m} = \bbP_{\mathbf m}(L_\bfm) \lessdot \bbP_{\mathbf m}[L_{\mathbf m}]$ and is a dense subset of it under the natural identification (see \ref{b5}(1)), but we should pedantically use $\bbP'_{\mathbf m}(L_{\mathbf m})$ or use $\lessdot'$.
    
    2) $\bbP_{\mathbf m}[L_{\mathbf m}]$ is $(< \lambda)$-strategically complete and is $\lambda^+$-c.c.
    
    3) $\bbP_{\mathbf m}(L) \subseteq \bbP_{\mathbf m}[L]$ as sets and  $\bbP_{\mathbf m}[L] \lessdot \bbP_{\mathbf m}[L_{\mathbf m}]$ and $\bbP_{\mathbf m}(L) \subseteq' \bbP_{\mathbf m}[L]$.
    
    4) If $\mathbf G  \subseteq \bbP_{\mathbf m}$ is generic over $\mathbf V$ and $\eta_t = \name\eta_t[\mathbf G]$ for $t \in L_ \mathbf{m} $  and $\mathbf G^+  =   \{\psi \in \bbL_{\lambda^+}(Y_{L_{\bfm}}, \mathbb{P}_\mathbf{m} ): \psi[\mathbf G] = \true\}$, see \ref{b5}(2)(3), \then \, $\mathbf V[\mathbf G] = \mathbf V[\mathbf G^+] = \mathbf
    V[\langle \eta_t:t \in L_{\mathbf m}\rangle]$.
    
    5) In part (4), moreover  $\mathbf G^+$ is a subset of $\bbP_{\mathbf m}[L_{\bfm}]$ generic over $\mathbf V$.
    
    6) $\bbP_{\mathbf m}(L_1) \subseteq \bbP_{\mathbf m}(L_2)$ and $\bbP_{\mathbf m}[L_1] \lessdot \bbP_{\mathbf m}[L_2]$ \when \, $L_1
    \subseteq L_2 \subseteq L_{\mathbf m}$.
    
    7) If $\mathbf m,\mathbf n \in \mathbf M$ are equivalent  \then \, $\bbP_{\mathbf m}[L] = \bbP_{\mathbf n}[L]$ and $\bbP_{\mathbf m}(L) = \bbP_{\mathbf n}(L)$ for $L \subseteq L_{\bfm}$.
     
    8) [($ > \lambda $)-continuity]  Assume $I_*$ to be a $\lambda^+$-directed partial order and $\bar L = \langle L_r:r \in I_*\rangle$ be such that $r \in I_* \Rightarrow L_r \subseteq L_{\mathbf m}$ and $r <_{I_*} s \Rightarrow L_r \subseteq L_s$ and $L = \cup\{L_r:r \in I_*\}$.  \Then \,, as sets and moreover as partial orders $\bbP_{\mathbf m}[L] = \cup\{\bbP_{\mathbf m}[L_r]:r \in I_*\}$ and $\bbP_{\mathbf m}(L) = \cup\{\bbP_{\mathbf m}(L_r):r \in I_*\}$.
    
    9) If $\mathbf m \in \mathbf M_{\ec}$ and $\mathbf m \le_{\mathbf M} \mathbf m_1
    \le_{\mathbf M} \mathbf m_2$ \underline{then} $\bbP_{\mathbf m_1}[L_{\mathbf m}] = 
    \bbP_{\mathbf m_2}[L_{\mathbf m}]$.
     
    10) The sequence $ \name{\bar{\eta } }_L = \langle \name{ \eta } _ s : s \in L \rangle $ 
    is a generic for $ \mathbb{P} _ \mathbf{m} [L] $, that is: if $ \mathbf{G} \subseteq \mathbb{P}_ \mathbf{m} [L] $ is generic over $ \mathbf{V} $ and  $ \nu _s = \name{ \eta } _s [ \mathbf{G} ] $ for $ s \in L $ then: 
     
    \begin{enumerate} 
        \item[(a)] $ \mathbf{V} [ \mathbf{G} ] = \mathbf{V} [ \langle  \nu _s : s \in L \rangle ]$,   
         
        \item[(b)] $\bar{ \nu}  = \langle \nu _s : s \in L \rangle$  determines $ \mathbf{G} $ uniquely.
    \end{enumerate} 
\end{claim}

\begin{remark}\label{b17}
    What about $\bbP_{\mathbf m}(L) \subseteq'_{\ic} \bbP_{\mathbf m}[L]$  and $\bbP_{\mathbf m}(L) \lessdot' \bbP_{\mathbf m}[L]$?
    
    Concerning the second, there may be a maximal antichain $\langle p_{i}: i < i_{*} \rangle $ of $\bbP(L),$ but some $q \in \bbP_{\bfm}$ is incompatible with $p_{i}$ for $i < i_{*}.$ This witness $\neg (\bbP_{\bfm}(L) \lessdot \bbP_{\bfm})$ hence $\neg(\bbP_{\bfm}(L) \lessdot \bbP_{\bfm}[L]).$ Concerning the first ($\bbP_{\bfm}(L) \subseteq_{\rm{ic}}' \bbP_{\bfm}[L]$) easily it holds. Note that $(\bbP_{\bfm}(L) \subseteq \bbP_{\bfm}[L])$ may fail as explained earlier as maybe $q \Vdash_{\bbP_{\bfm}}$``$p \in \name{\bfG}$'' but $\nleq_{\bbP_{\bfm}} q,$ see \ref{c5}(9) and \ref{c32n}.
\end{remark}

\begin{PROOF}{\ref{b14}}
    1) Easy.
     
    2) Follows by part (1) and \ref{c11}.
     
    3) The first statement by their definitions, the second statement by part (1).
     
    For the third clause, $ \lqq \mathbb{P} _ \mathbf{m}[L] \subseteq '\mathbb{P} _\mathbf{m}(L) "$, note that:  
        
    \begin{enumerate} 
        \item[$(*)_{1}$]  if $ p, q \in \mathbb{P}_{\mathbf{m}} (L)$, then
        $ \mathbb{P}_{\bfm}(L) \models \lqq p \le q "$ iff $\bbP_{\bfm} \models$``$ p \leq q$'' which implies $ \mathbb{P}_{\mathbf{m}}[L] \models \lqq p \le q "$  by the definition of $ \mathbb{P} _ \mathbf{m} [L].$  
    \end{enumerate} 
    
    \begin{enumerate}
        \item[$(\ast)_{2}$] if $p, q \in \bbP_{\bfm}(L)$ and $\dom(p) \subseteq \dom(q),$ then $\bbP_{\bfm}(L) \models$``$p \leq q$'' iff $\bbP_{\bfm} \models$``$p \leq q$'' iff $\bbP_{\bfm}[L] \models$``$p \leq q$''. 
    \end{enumerate}
    
    [The first ``iff'' by the definition of $\bbP_{\bfm}(L),$ the second ``iff'' by \ref{c32n}.]
     
    4), 5), 6) Should be clear recalling \ref{c11}(7).
    
    7) Easy, recalling \ref{c11}(5).
    
    8), 9)  Easy.
    
    10) By the definition of $ \mathbb{P} _ \mathbf{m} [L] $.
\end{PROOF}

\begin{tuc}\label{b20}
    There is an isomorphism from $\bbP_{\mathbf m_1}[M_1]$ onto $\bbP_{\mathbf m_2}[M_2]$ which (recalling Definition \ref{b11}(1)) maps  
    $p^*_{t,\eta}$ to $p^*_{h(t),\eta}$ for $t \in M_1,\eta \in  \cup\{\prod\limits_{\varepsilon < \zeta} \theta_\varepsilon:\zeta < \lambda\}$ \when:
    
    \begin{enumerate}
        \item[$\boxplus$]  $(a) \quad \mathbf m_\ell \in \mathbf M^{\ec}_\infty$ for $\ell=1,2$, 
        
        \item[${{}}$]  $(b) \quad M_\ell = M_{\mathbf m_\ell}$
        for $\ell=1,2$,
        
        \item[${{}}$]  $(c) \quad h$ is an isomorphism from $\mathbf m_1
        \rest M_1$ onto $\mathbf m_2 \rest M_2$.
    \end{enumerate}
\end{tuc}

\begin{PROOF}{\ref{b20}}
    By renaming \wilog \, $M_1 = M_2$ call it $M$ and $h$ is the identity and $L_{\mathbf m_1} \cap L_{\mathbf m_2} = M$.  Let $\mathbf m_0 = \mathbf m_1 \rest M = \mathbf m_2 \rest M$ so $\mathbf m_0 \le_{\mathbf M} \mathbf m_\ell$ 
    for $\ell=1,2$ and $L_{\mathbf m_0} = L_{\mathbf m_1} \cap L_{\mathbf m_2}$. 
     
    By \ref{c31}, there is $\mathbf m$ such that $\mathbf m_1 \le_{\mathbf M} \mathbf m$ and $\mathbf m_2 \le_{\mathbf M} \mathbf m$. 
    As $\mathbf m_1,\mathbf m_2 \in \mathbf M^{\ec}_\infty$ by \ref{b14}(9) 
    we have $\bbP_{\mathbf m_1}[M] = 
    \bbP_{\mathbf m}[M]$ and $\bbP_{\mathbf m_2}[M] = \bbP_{\mathbf m}[M]$ so together we get the desired conclusion.
\end{PROOF}

\begin{definition}\label{b23}
    1) We call $\mathbf m \in \mathbf M$ reduced when $L_{\mathbf m} = M_{\mathbf m}$.  We call $ \mathbf{m} $ unary when the equivalence  relation
    $ E''_ \mathbf{m} $ has exactly one equivalence class.
    
    2) For $\mathbf m \in \mathbf M$ let $\bbP^{\cor}_{\mathbf m}$ be $\bbP_{\mathbf n}[L_{\mathbf m}]$ and $\bbP^{\cor}_{\mathbf m}[L]$ be $\bbP_{\mathbf n}[L]$ for $L \subseteq L_{\mathbf m}$  \when \, $\mathbf m \le_{\mathbf M} \mathbf n \in \mathbf M_{\ec}$.
\end{definition}

\begin{remark}\label{b26}
    1) Why is $\bbP^{\cor}_{\mathbf m}[L]$ well defined? see below.
    
    2) Here ``cor" stands for corrected.
\end{remark}

The interest in the definition is because:

\begin{claim}\label{b29}
    1) If $\mathbf m \in \mathbf M$ and $L \subseteq L_{\mathbf m}$ \then \, $\bbP^{\cor}_{\mathbf m}[L]$ is well defined. 
    
    2) $\bbP^{\cor}_{\mathbf m}[M_{\mathbf m}]$ is well defined and depends only on $\mathbf m \rest M_{\mathbf m}$.
    
    3) If $\mathbf m \le_{\mathbf M} \mathbf n$ and $L_1 \subseteq L_2  \subseteq L_{\mathbf m}$ \then \, $\bbP^{\cor}_{\mathbf m}[L_1] = 
    \bbP^{\cor}_{\mathbf n}[L_1] \lessdot \bbP^{\cor}_{\mathbf n}[L_2]
    \lessdot \bbP^{\cor}_{\mathbf n}$.
    
    4) Assume $\bfm$ is bounded and $\bfm \leq_{M} \bfn \in \bfM_{\rm{bec}}.$ If $L \subseteq L_{\bfm}$ then $\bbP_{\bfm}^{\rm{cor}}[L] = \bbP_{\bfn}[L].$
    
    5) Assume $\bfm$ is weakly bounded and $\bfm \leq_{\bfM} \bfn \in \bfM_{\rm{wec}}.$ If $L \subseteq L_{\bfm}$ then $\bbP_{\bfm}^{\rm{cor}}[L] = \bbP_{\bfn}[L].$
    
    6) If $\bfn \in \bfM_{\rm{wec}}$ then $\bfm \in \bfM_{\ec}.$
\end{claim}

\begin{PROOF}{\ref{b29}}
    1) By \ref{c41}, $\bbP^{\cor}_{\mathbf m}[L]$ has at least one definition so it suffices to prove uniqueness.  So  assume $\mathbf m \le_{\mathbf M} \mathbf m_\ell \in \mathbf M_{\ec}$ for $\ell=1,2$ and we should prove that $\bbP_{\mathbf m_1}[L] = \bbP_{\mathbf m_2}[L]$. \Wilog \, $L_{\mathbf m_1} \cap L_{\mathbf m_2} = L_{\mathbf m}$.  Now by \ref{c31}  we can find $\mathbf n \in \mathbf M$  such that $\mathbf m_1 \le_{\mathbf M} \mathbf n$ and $\mathbf m_2 \le_{\mathbf M} \mathbf n$; as $\mathbf m_\ell \in \mathbf M_{\ec}$ see Definition \ref{c34} we have $\bbP_{\mathbf m_\ell} \lessdot \bbP_{\mathbf n}$ for $\ell=1,2$.  As in the end of the proof of \ref{b20} we are done.
    
    2) By \ref{b20}.

    3) Follows from Definition \ref{c34}(2) and \ref{b23}(2).
    
    4) On the one hand, we can find $\bfm_{1} \in \bfM_{\rm{bec}}$ such that $\bfm \leq_{\bfM} \bfm_{1}$ by \ref{c41}(2). On the other hand, can find $\bfm_{3} \in \bfM_{\ec}$ such that $\bfm_{1} \leq_{\bfM} \bfm_{3}$ by \ref{c4}(1). Let $\bfm_{2} = \bfm_{3}^{[\bd]}$ and let $\bfm_{0} = \bfm$ so $\bfm_{0} \leq_{\bfM} \bfm_{1} \leq_{\bfM} \bfm_{2} \leq_{\bfM} \bfm_{3}.$ By the choice of $\bfm_{1}$ we have 
    
    \begin{itemize}
        \item $\bbP_{\bfm_{1}}[L] = \bbP_{\bfm_{2}}[L] \lessdot \bbP_{\bfm_{2}}.$
    \end{itemize}
    
    As $L_{\bfm_{2}}$ is an initial segment of $L_{\bfm_{3}},$ clearly, 
    
    \begin{itemize}
        \item $\bbP_{\bfm_{2}} \lessdot \bbP_{\bfm_{3}}$ so $\bbP_{\bfm_{2}}[L] = \bbP_{\bfm_{3}}[L].$
    \end{itemize}
    
    Lastly as $\bfm_{3} \in \bfM_{\ec}, \bbP_{\bfm_{3}}[L] = \bbP_{\bfm}^{\cor}[L].$ Together we are done. 
    
    5) Similarly to part (4). 
    
    6) Easy. 
\end{PROOF}

\begin{discussion}\label{b32}  
    1) But we like to prove for reduced $\mathbf m \in \mathbf M$ and $M \subseteq M_{\mathbf m}$ that $\bbP^{\cor}_{\mathbf m \rest M} \lessdot \bbP^{\cor}_{\mathbf m}$, this is the whole point of the  corrected iteration.  This is delayed to \ref{e50}.   We now prove that this  suffices.

    2) Conclusion \ref{b35}  below is the desired conclusion but it relies on  \S 3, specifically  on 
    \ref{e50} (or \S4A).
    
    3) The reader may understand \ref{b35} without reading the rest of \S2, \S3 by ignoring clause (A)(d), \oor \, reading \ref{b5}, \ref{b8}.
    
    4) By \ref{b29}(4) we may restrict ourselves to $\bfM_{\bd}.$ We use it freely. 
\end{discussion}

\begin{conclusion}\label{b35}
    For every ordinal $\delta_*$ there is $\mathbf q = \langle \bbP_\alpha,\name\eta_\beta:\alpha \le \delta_*, \beta < \delta_{*} \rangle$ such that:

    \begin{enumerate}
        \item[(A)]
            \begin{enumerate}
                \item[(a)]  $\langle \bbP_\alpha:\alpha \le \delta_* \rangle$ is $\lessdot$-increasing sequence of forcing notions, 
                
                \item[(b)]  $\name\eta_\alpha$ is a $\bbP_{\alpha +1}$-name of a member of $\prod\limits_{\varepsilon < \lambda} \theta_\varepsilon$ which dominates $(\prod\limits_{\varepsilon < \lambda} \theta_\varepsilon)^{\mathbf V[\bbP_\alpha]}$,
                
                \item[(c)]  $\name\eta_\alpha$ is a generic for $\bbP_{\alpha +1}/\bbP_\alpha$, moreover $\langle
                \name\eta_\beta:\beta < \alpha\rangle$ is a generic 
                for $\bbP_\alpha$, 
                
                \item[(d)]  $\bbP_\alpha \lessdot' \bbL_{\lambda_{0}}^{+}(Y_\alpha,\bbP_\alpha)$ in fact $\bbP_{\alpha}$ is dense in $\bbL_{\lambda_{0}}^{+}(Y_{\alpha}, \bbP_{\alpha})$  where $Y_\alpha$ is defined as in \ref{b11}(2) with $\alpha$ here standing for $L$ there and see \ref{b5},
                
                \item[(e)]  $\bbP_\alpha$ is $(< \lambda)$-strategically
                complete and $\lambda^+$-c.c., 
                
                \item[(f)]  if $\delta \le \delta_*$ has cofinality $>
                \lambda$  \then \, $\bbP_\delta = \cup\{\bbP_\alpha:\alpha < \delta\}$,
                if $ \cf (\delta) = \lambda $ then the union  is just a dense subset of $ \mathbb{P}_\delta $, 
                
                \item[(g)]  $\bbP_{\delta_*}$ has cardinality  $|\delta_*|^\lambda$. 
            \end{enumerate}
            
        \item[(B)]  if $\cU \subseteq \delta_*$ \then \, the complete sub-forcing generated by $\langle \name\eta_\alpha:\alpha \in
        \cU\rangle$ is isomorphic to $\bbP_{\otp(\cU)}$,
        
        \item[(C)]  if $\mathbf G \subseteq \bbP_{\delta_*}$ is generic over
        $\mathbf V$ and $\eta_\alpha = \name\eta_\alpha[\mathbf G]$ for $\alpha <
        \delta_*$ and $\eta'_\alpha \in \prod\limits_{\varepsilon <
        \lambda} \theta_\varepsilon$ for $\alpha < \delta_*$ and $\{(\alpha,\varepsilon):\alpha < \delta_*,\varepsilon < \lambda$ and
        $\eta'_\alpha(\varepsilon) \ne \eta_\alpha(\varepsilon)\}$ has
        cardinality $< \lambda$ \then \, also $\langle \eta'_\alpha:\alpha <
        \delta_*\rangle$ is a generic for $\bbP_{\delta_*}$, determining a possibly
        different $\mathbf G'$ but $\mathbf V[\mathbf G'] = \mathbf V[\mathbf G]$,
        
        \item[(D)]  in clause (B), 
        moreover if $\cU \subseteq \delta_*$ and
        $\langle \alpha_i:i < \otp(\cU)\rangle$ list $\cU$ in increasing order \then \, for some unique $\mathbf G'' \subseteq \bbP_{\otp(\cU)}$ generic over $\mathbf V,i < \otp(\cU) \Rightarrow \eta'_{\alpha_i} =
        \name\eta_i[\mathbf G'']$.
    \end{enumerate}
\end{conclusion}

\begin{PROOF}{\ref{b35}}
    \Wilog \, $\lambda_1 \ge |\delta_*|;$ we can use only $\bfm \in \bfM_{\bd}$ (by \ref{b29}(4)). We define $\mathbf m \in \mathbf M$ by:
    
    \begin{enumerate}
        \item[$(*)$]
        
        \begin{enumerate}
            \item[(a)]  $L_{\mathbf m} = \delta_*,$
            
            \item[(b)]  $M_{\mathbf m} = \delta_*$ and\footnote{Other reasonable choice is  $M_{\bfm}^{\fat} = \emptyset, M_{\bfm}^{\rm{lean}} = \delta_{\ast}$ and $M_{\bfm}^{\rm{fat}} = \emptyset = M_{\bfm}^{\rm{lean}}$.} $M_{\bfm}^{\fat} = \delta_{\ast},$ 
        
            \item[(c)]  $u_{\mathbf m,\alpha} = \alpha$ and $\cP_{\mathbf m,\alpha} = [\alpha]^{\le \lambda}$ for $\alpha < \delta_*,$
            
            \item[(d)]  $ E'_{\mathbf m} = \emptyset$.
        \end{enumerate}
    \end{enumerate}
    
    It is easy to check that indeed $\mathbf m \in \mathbf M$ and let $\mathbf n \in \mathbf M_{\ec}$ be such that $\mathbf m \le_{\mathbf M} \mathbf n$, exists by the Crucial Claim \ref{c41} and let $\bbP_\alpha =  \bbP_{\mathbf n}[\{i:i < \alpha\}]$ for $\alpha \le \delta_*$. 
    
    Now clearly clause (A) holds and $\bbP_\delta = \bbP^{\cor}_{\mathbf m}$ by \ref{b23}(2), \ref{b29}(1) and  e.g.  clause (A)(b) holds by \ref{c11}(4A).
     
    As for clause (B), first note that for every $L \subseteq \delta_*$,  the sequence $\bar\eta_L = \langle \name\eta_\alpha:\alpha \in L\rangle$ is generic for $\bbP_{\mathbf m}[L]$ by Definition \ref{b11}.
    
    Second, for $M \subseteq \delta_*$ let $\alpha = \otp(M)$ and $h:M \rightarrow \alpha$ be $h(i) = \otp(i \cap M)$ so $h$ is an isomorphism from
    $\mathbf m \rest M$ onto $\mathbf m \rest \alpha$ hence by \ref{e50}(2) below, with $\mathbf m,\mathbf m \rest \alpha,M,\alpha$ here standing
    for $\mathbf m_1,\mathbf m_2,M_1,M_2$ there we have $h$ induces an isomorphism from $\bbP^{\cor}_{\mathbf m}[M]$ onto $\bbP^{\cor}_{\mathbf m \rest \alpha}[L_{\mathbf m \rest \alpha}]$.  In particular, $\id_\alpha$ induces an isomorphism from $\bbP^{\cor}_{\mathbf m \rest \alpha}$ onto $\bbP^{\cor}_{\mathbf m}[\alpha]$.
    
    Together we get clause (B). Also Clause (C)  holds by \ref{c11}(8) and  clause  (D) follows so we are done. 
\end{PROOF}

\begin{definition}\label{b36}
    1) We say $\bfm$ is essentially $(< \mu)-$directed (if $\mu = \aleph_{0}$ we may omit it) \underline{when}: if $L \subseteq M, \vert L \vert < \mu$ then for some $t \in M_{\bfm},$ we have: 
    
    \begin{itemize}
        \item[$\bullet$]  $s \in L \Rightarrow s <_{\bfm} t \wedge s \in u_{t}$ so $M_{\bfm}$ is directed\footnote{Why not add $\{ s \} \in \cP_{\bfm, t}$? See \ref{c37}(13).}). 
    \end{itemize}
    
    [Note that it follows because $\bfm$ is bounded.] 

    2) We say $\bfm$ is strongly $\mu$-directed (or $(< \mu)$-directed; if $\mu = \aleph_{0}$ we may omit it) \underline{when}: for every $L \subseteq L_{\bfm}$ of cardinality $< \mu$ there is $t \in M_{\bfm}$ such that $L \in \cP_{\bfm, t}$ (the condition implies ``$\bfm$ is weakly bounded'' and ``$\bfm$ is not lean, $t \notin M_{\bfm}^{\rm{lean}}$ when $E_{\bfm}''$ has at least two equivalence classes''). 
    
    3) We say $\bfm$ is reasonable \underline{when}: 
    
    \begin{enumerate}
        \item[$(\alpha)$] $\bfm$ is strongly $\lambda^{+}$-directed and $M_{\bfm}^{\rm{fat}}$ is cofinal in $M_{\bfm},$
        
        \item[$(\beta)$] $\bfm(\leq t) \in  \bfM_{\rm{ec}}$ for every $t \in M_{\bfm},$
        
        \item[$(\gamma)$] $\bfm$ is wide and bounded (see Definition \ref{c5}(10) and Definition \ref{c39}(1A)). 
    \end{enumerate}
\end{definition}

Similarly we can deal with such iterations with partial memory and spell out how $\bbP^{\cor}_{\mathbf m}[L]$ is defined from a $(<
\lambda)$-support iteration with partial memory. This is used in \cite{Sh:945}, but we need more: see \S3.

\begin{conclusion}\label{b38}
    Assume $M$ is a well founded partial order and $\bar u' = \langle u'_t:t \in M\rangle,u'_t \subseteq M_{<t}$ and $\bar{\cP}' = \langle \cP'_t: t \in M\rangle$ with $\cP'_t \subseteq [u'_t]^{\le \lambda}$ is closed under subsets.  \Then \, we can find $\beta(*),h,\bbP_\beta = \bbP_{0,\beta}, \name{ \mathbb{Q} }_\beta = \name{ \mathbb{Q} }_{0, \beta }, \bbP_{1,\beta}, \name{\bbQ}_{1,\alpha },
    \name\eta_\alpha,\name\eta'_s$  and  $\bbP_{1, v}, \bbP'_u$ (for $\beta \le \beta(*), \alpha < \beta(*),s \in M$  and $v \subseteq \beta(\ast), u \subseteq M$) and $h,   \bar u,\bar{\cP}$ such that:
    
    \begin{enumerate}
        \item[(A)]  
        
            \begin{enumerate}
            \item[(a)]  $\langle \bbP_\beta,\name{\bbQ}_\alpha:\beta \le
            \beta(*),\alpha < \beta(*)\rangle$ is $(< \lambda)$-support\footnote{This will be $ \mathbf{q} _ \mathbf{m} $, well up to equivalence, see \S1.} iteration,
            
            \item[(b)] 
            
            \begin{enumerate}
                \item[$(\alpha)$]  $\bar u = \langle u_\beta:\beta < \beta(*)\rangle$ such that $u_\beta \subseteq \beta$, 
                
                \item[$(\beta)$]  $\bar{\cP} = \langle \cP_\beta:\beta < \beta(*)\rangle$ such that $\cP_\beta \subseteq [u_\beta]^{\le \lambda}$ is closed under subsets, 
            \end{enumerate}
            
            \item[(c)]  $\name\eta_\alpha$ is a $\bbP_{\alpha +1}$-name of a member of $\prod\limits_{\varepsilon < \lambda}
            \theta_\varepsilon$, 
            
            \item[(d)]  $\langle \name\eta_\alpha:\alpha < \beta\rangle$ is generic for $\bbP_\beta$,  
            \item[(e)]  $\name{\bbQ}_\alpha$ is defined as in  Definition \ref{c7}(4),
            
            \item[(f)]  $\Vdash_{\bbP_{\beta(*)}} ``\name\eta_\beta \in \prod\limits_{\varepsilon < \lambda} \theta_\varepsilon$ dominates
            every $\nu \in \prod\limits_{\varepsilon < \lambda}  \theta_\varepsilon$ from $\mathbf V[\langle \name\eta_\alpha:\alpha \in u\rangle]$ when $u \in \cP_\beta$.  
        \end{enumerate}
        
        \item[(B)]
        
        \begin{enumerate}
            \item[(a)]  $h$ is a one-to-one function from $M$ into\footnote{In general not onto!} $\beta(*)$; stipulate
            $h(\infty) = \beta(*)$, 
            
            \item[(b)]  $s <_M t \Rightarrow h(s) < h(t)$, 
            
            \item[(c)]  $u_{h(t)} \cap \rang(h) = \{h(s):s \in u'_t\}$, 
            
            \item[(d)]  $\cP_{h(t)} \cap [\rang(h)]^{\le \lambda} = 
            \{\{h(s):s \in u\}:u \in \cP'_t\}$, 
        \end{enumerate}

        \item[(C)]
        
        \begin{enumerate}
            \item[(a)]  $\bbP_{1,\beta} = \bbL^+_{\lambda^+}(Y_\beta,\bbP_\beta)$ where we let $Y_\beta = \{p^*_{\alpha,\nu}:\alpha < \beta,\nu \in \prod\limits_{\varepsilon < \zeta} \theta_\varepsilon$ for some $\zeta < \lambda\}$, see \ref{b5}, \ref{b11}(1), 
            
            \item[(b)]  $\bbP_{1,u} = \bbL^+_{\lambda_{0}}(Y_u,\bbP_\beta)$, where $Y_u$ is defined similarly when $u \subseteq \beta(*)$, 
            
            \item[(c)]  $\bbP'_{u}$ is a forcing notion for $u \subseteq M$ and $\name\eta'_s$ is a $\bbP'_{{\{s\}}}$-name for $s \in M,$
        
            \item[(d)]  $h$ induces an isomorphism from  $\bbP'_{u}$ onto $\bbP_{1, \{h(s):s \in u\}}$ for $u \subseteq M$ and $\name\eta'_s$ to $\name\eta_{h(s)}$ for $s \in M$, 
            
            \item[(e)]  $\langle \name\eta_{h(s)}:s \in u\rangle$ is generic for $\bbP'_{u}$ for $u \subseteq M$, 
        \end{enumerate}

        \item[(D)]
        
        \begin{enumerate}
            \item[(a)]  $\bbP'_{u} \lessdot \bbP'_{v}$ when $u \subseteq v \subseteq M$, 
            
            \item[(b)]  $\bbP_\beta,\bbP_{1,u},\bbP'_{1, u}$ are $(< \lambda)$-strategically complete and $\lambda^+$-c.c., 
            
            \item[(c)]   if $M_1,M_2 \subseteq M$ and $f$ is an isomorphism from $M_1$ onto $M_2$ as partial orders  such that $t \in M_1 \Rightarrow u'_{f(t)} \cap M_2 = \{f(s):s \in u'_t \cap M_1\}$ and $t \in M_{1}  \Rightarrow \cP_{f(t)}' \cap [M_{2}]^{\leq \lambda} = \{ f(s): s \in u \cap M_{1} \}: u \in \cP_{t}' \}$ \underline{then}  the mapping $h(s) \mapsto h(f(s))$ induces 
            an isomorphism from the forcing notion $\bbP'_{1, {M_1}}$ onto $\bbP'_{1, {M_2}}$.
        \end{enumerate}
        
        \item[(E)]  if $M$ is $(< \lambda^{+})$-directed and  the set $ Y \subseteq M $  is cofinal in $M,$ \then \, the set  $\{ \name{ \eta }_{h(s)}:  s \in Y \} $ is cofinal in $\{ \name{\eta}_{\beta}: \beta < \beta(\ast) \}$and even in $\Pi_{\varp < \lambda} \theta_{\varp}$ in $\bfV^{\bbP_{\beta(\ast)}} $ (see \ref{c37}(3)).
    \end{enumerate}
\end{conclusion}

\begin{PROOF}{\ref{b38}}
    Easy. We can assume $\lambda_{1} \geq \vert M \vert. $ Similarly  to the proof of \ref{b35}, the proof of clause (E) is easy by \ref{e37}
\end{PROOF}

\begin{claim}\label{b41}
    If $ \mathbf{m}_1 \le_{\mathbf{M} } 
    \mathbf{m}_2 \le_{\mathbf{M} } \mathbf{n}$ 
    and $ \mathbb{P} _{\mathbf{m} _ {\ell} }    \lessdot \mathbb{P} _ \mathbf{n} $ 
    for $ {\ell} = 1,2$  then $ \mathbb{P} _{\mathbf{m}_1} \lessdot \mathbb{P} _{\mathbf{m}_2}$
\end{claim}

\begin{PROOF}{\ref{b41}}
    Easy.
\end{PROOF} 

The following will be used in \ref{b47}.

\begin{claim}{\label{b44}}
    If (A) then (B), where:
    
    \begin{enumerate} 
        \item[(A)]  
      
         \begin{enumerate} 
            \item[(a)] $ \mathbf{m}_{0}, \mathbf{m}_1, \mathbf{m} _2 \in \mathbf{M} $,
             
            \item[(b)] $ L_* $ is an initial segment of $ L_{\mathbf{m} _ 1} $,
            
            \item[(c)] $ L_* = L_{\mathbf{m} _1}\cap L_{\mathbf{m} _2}$,
            
            \item[(d)] $ \mathbf{m} _0 =  \mathbf{m} _1 \rest L_* \le _ \mathbf{M} = \mathbf{m} _2 $, 
         \end{enumerate}

        \item[(B)]  there is $ \mathbf{m} \in \mathbf{M} $ such that:
        
         \begin{enumerate}  
            \item[(a)] $ \mathbf{m}_{1} \le _ \mathbf{M} \mathbf{m},$
            
            \item[(b)] $ \mathbf{m} _2 = \mathbf{m} \rest L_{\mathbf{m} _2}.$
         \end{enumerate} 
     \end{enumerate}

     2) If $L_{1} \subseteq L_{2}$ are initial segments of $L_{\bfm}$ and $\bfm \rest L_{2} \in \bfM_{\ec}$ then $\bfm \rest L_{1} \in \bfM_{\rm{ec}}.$
     
     3) In part (1) we may add (e) to clause (A) and (c), (d) to clause (B), where: 
     
     \begin{enumerate}
         \item[(A)(e)] $L_{\ast} \subseteq L_{\bfm_{1}(< t_{\ast})},$ where $t \in M_{\bfm_{1}},$
         
         \item[(B)(c)] if $s \in L_{\bfm} \setminus L_{\bfm_{1}}$ then $s <_{\bfm_{1}} t_{\ast},$
         
         \item[$ $ (d)] if $s \in M_{\bfm_{1}} \setminus M_{\bfm_{0}}$ and $t_{\ast} \leq_{\bfm_{1}} s$ then $u_{\bfm, s} = u_{\bfm_{1}, s} \cup ((L_{\bfm} \setminus L_{\bfm_{1}}) \cap L_{\bfm(< s)}).$ 
     \end{enumerate}
\end{claim} 

\begin{PROOF}{\ref{b44}}
    1) Easy but we elaborate. 
    We define $ \mathbf{m} $ as follows:
    
    \begin{enumerate} 
    \item[$(*)_1$]   
    
        \begin{enumerate} 
            \item[(a)]  $ L_ \mathbf{m} $ as a set is  $ L_{\mathbf{m} _1} \cup L   _{\mathbf{m}_2}$,
            
            \item[(b)] $\le _ \mathbf{m} $ is the transitive closure of $ \{ (s,t):  L_{\mathbf{m} _1 } \models s < t $ or $  L_{\mathbf{m} _2 } \models s < t   \} $,
            
            \item[(c)]  $ M_ \mathbf{m} = M_  {\mathbf{m} _1} $, $M_{\bfm}^{\rm{lean}} = M_{\bfm_{1}}^{\rm{lean}}, M_{\bfm}^{\rm{fat}} = M_{\bfm_{1}}^{\rm{fat}},$
            
            \item[(d)] $ u_{\mathbf{m}, t }$ is:
               \begin{enumerate} 
                    \item[$ (\alpha ) $] $ u_{\mathbf{m} _1, t }$   when $ t \in L_{\mathbf{m} _1} \setminus L_*,$  and,
                    
                    \item[$ ( \beta ) $]  $ u_{\mathbf{m} _2, t}$  when  $ t \in L_{\mathbf{m}_2} \setminus M_{\bfm_{0}},$
                    
                    \item[$(\gamma)$] $u_{\bfm_{1}, t} \cup u_{\bfm_{2}, t}$ if $t \in M_{\bfm_{0}}.$
                \end{enumerate}
           
            \item[(e)] $\cP_{\bfm, t}$ is:
            
            \begin{enumerate}
                \item[$(\alpha)$] $\cP_{\bfm_{1}, t}$ when $t \in L_{\bfm_{1}} \setminus L_{*}$ and,
                
                \item[($\beta$)] $\cP_{\bfm_{2}, t}$ when $t \in L_{\bfm_{2}} \setminus M_{\bfm_{0}},$
                
                \item[$(\gamma)$] $[u_{\bfm, t}]^{\leq \lambda}$ if $t \in M_{\bfm_{0}}^{\fat},$
                
                \item[$(\delta)$] $\cP_{\bfm_{1}, t} \cup \cP_{\bfm_{2}, t}$ if $t \in M_{\bfm_{0}} \setminus M_{\bfm_{\gamma}}^{\fat}.$
            \end{enumerate}
            
            \item[(f)]   We define $ E'_\mathbf{m} $ by: for $ s,t \in L_ \mathbf{m} $,  we have $ s E'_{\bfm} t $  \underline{iff} $ s E'_{\mathbf{m} _1} t$  or $ s E'_{\mathbf{m} _2 } t.$
        \end{enumerate} 
    \end{enumerate}
     
     As  $L _* $  is an initial segment of $ L_{\mathbf{m} _1 }$  we have:
                
    \begin{enumerate} 
        \item[$(*)_2 $] $ L_ \mathbf{m} \models \lqq s \le t  "$  \underline{iff}  $L_{\mathbf{m} _2} \models \lqq s \le t  "$   \underline{or}  $ s \in L_{\mathbf{m} _2 },   t \in L_{\mathbf{m} _ 1} \setminus L_* $  and for some $ r \in L_*$ we have $ L_ {\mathbf{m} _2} \models \lqq s \le r"$  and  $ L_{\mathbf{m} _1 }\models \lqq r \leq t " $.  
                 
        \item[$(*)_{3}$] $L_{\bfm_{2}}$ is an initial segment of $L_{\bfm}.$        
    \end{enumerate}
    
    Now check that $ \mathbf{m} $ is as required.
        
    2) Follows.
    
    3) Easy (changing $(\ast)_{1}$ above naturally). 
\end{PROOF} 

Sometime we would like to have in addition to being in $ \mathbf{M} _{\ec}$  that $ \{ \name{ \eta }  _ s: s \in M \} $  be cofinal  in $ (\Pi _{\varepsilon < \lambda } \theta _ \varepsilon,   \le _{J^{\bd}_\lambda })$ in $\bfV^{\bbP_{\bfm}}$. Toward this we use the following claim:

\begin{claim}\label{b47}
    Assume $\mathbf{m} \in \mathbf{M}.$ 

    1) A sufficient condition for $ \mathbf{m} \in  \bfM_{\rm{bec}}$ is:
    
    \begin{enumerate} 
         \item[$(*)_\mathbf{m} $] For some $\delta, 
         \bar{L}, \bar{ c }$  we have:
         
         \begin{enumerate} 
             \item[(a)]  $ \bar{c} = \langle c_ \alpha: \alpha < \delta  \rangle 
             \in {}^{ \delta } (M_{\mathbf{m}}),$ 
             
             \item[(b)]  $ \bar{ L } = \langle L_ \alpha : \alpha < \delta  \rangle $, 
             
             \item[(c)] $ \mathbf{m} \rest L_{\alpha} $ belongs to   $ \mathbf{M} _{\ec}$  for every $ \alpha < \delta $,
                    
             \item[(d)]  $ L_ \alpha \subseteq L_{\mathbf{m}, < c_ \alpha },  L_ \alpha \subseteq  u_{\mathbf{m} , c_ \alpha }$, $M_{\bfm (< c_{\alpha})} \subseteq L_{\alpha}$ and if $t \in L_{\alpha} \setminus M_{\bfm}$ then $L_{\alpha} \cap (t / E_{\bfm})$ is an initial segment of $t / E_{\bfm},$             
            
             \item[(e)] $ \delta $ has cofinality $> \lambda,$
             
             \item[(f)] $ \bar{ c } $ is  increasing  and cofinal in $ L_ \mathbf{m}$,
             
             \item[(g)] $ \bar{L} $  is $ \subseteq $-increasing with union $ L_ \mathbf{m} $. 
         \end{enumerate} 
    \end{enumerate} 
    
    2)  A sufficient condition for $ \mathbf{m} \in \mathbf{M}_{\rm{bec}}$ is:
    
    \begin{enumerate} 
         \item[$(*)'_\mathbf{m} $] For some $ \bar{ c}, \bar{L}$  we have:   
             \begin{enumerate} 
                \item[(a)-(e)]  as above,
                
                \item[(f)]  if $ L \subseteq L_ \mathbf{m} $ has cardinality  $ \le \lambda $ then for some $ \alpha < \delta $ we have $ L \subseteq L_{\alpha }$,   
             \end{enumerate} 
     \end{enumerate} 

    3)  For  $ L_* \subseteq L_ \mathbf{m} $  we have $(A)_{L_*} \Rightarrow (B)_{L_*},$ where:
        
    \begin{enumerate} 
        \item[$(A)_{L_*}$]   if $ L \subseteq L_* $ has cardinality  $\le \lambda $   and $ \mathbf{m} \le _\mathbf{M} \mathbf{n} $  then  $ \mathbb{P} _{\mathbf{n} }[L] = \mathbb{P} _ \mathbf{m} [L],$
        
        \item[$(B)_{L_*}$]  if $ \mathbf{m} \le _\mathbf{M} \mathbf{n} $  then  $ \mathbb{P} _{\mathbf{n} }[L_*] = \mathbb{P} _ \mathbf{m} [L_*]$,   
    \end{enumerate} 

    4) If $ c \in L_ \mathbf{m}, $  $ L_*  \subseteq u_{\mathbf{m}, c }$,   $ \mathbf{m} \rest L_{*} \in \mathbf{M} _ {\ec},$ $M_{\bfm(< c)} \subseteq L_{*}$ and  $t \in L_{\bfm} \setminus M_{\bfm} $ implies $ \ L_{*} \cap  (t / E_{\bfm})$ is an initial segment of $t / E_{\bfm}$ \underline{then} clause $ (B)_{L_*}$ above holds, 
    
    5) We have $ (a) \Rightarrow (b),$  when:
    
    \begin{enumerate} 
        \item[(a)]  we have: 
        
            \begin{enumerate} 
                \item[$ ( \alpha ) $] $\bfm$ is strongly $(<\lambda^{+})$-directed,
                
                \item[$ ( \beta ) $] for every $t \in M_{\bfm}$ (or just for cofinally many $t \in M_{\bfm}$) we have $\bfm(\leq t) \in \bfM_{\rm{ec}}.$
            \end{enumerate}
        
        \item[(b)] $\bfm \in \bfM_{\rm{ec}}.$ 
    \end{enumerate} 
    
    5A) Similarly for $\bfM_{\rm{bec}}.$ 
    
    6) If $ M $ is a $  < \lambda^{+}$-directed  well founded partial order  of cardinality $ \le  \lambda _1 $, for example,  $ M = (\kappa, < ), \kappa = \cf(\kappa) \in  (\lambda, \lambda _1]$, our main case,  \underline{then} there is a strongly $\lambda^{+}$-directed $ \mathbf{m}  \in \mathbf{M}$ such that $ M_ \mathbf{m} = M$  and $ (*)'_\mathbf{m} $  from part (2)  holds, (hence $ \mathbf{m} \in \mathbf{M} _{\rm{bec}}$   and $  \{ \name{ \eta }_s : s \in M_ \mathbf{m} \} $ is cofinal in $\left(  \Pi_{\varp < \lambda} \theta_{\varp}, <_{J_{\lambda}^{\rm{bd}}} \right)$ in the universe $\bfV^{\bbP_{\bfm}}.$
\end{claim}  

\begin{PROOF}{\ref{b47}}
    Straightforward (recalling \ref{b29}(4)), i.e.
    
    1) By (2). 
    
    2) By (3) and (4).
    
    3) Obvious, see \ref{b14}(8).
    
    4) Clear.
    
    5) Easy.
     
    6) Choose $ \bar{c} $ such that:
    
    \begin{enumerate}
        \item[$(*)_1$]  
        
        \begin{enumerate}
            \item[(a)]  $ \bar{ c } \in {}^{ \delta } (M_ \mathbf{m})$   for some ordinal $ \delta $,
            
            \item[(b)] if $ \alpha < \beta < \delta $ then $c_ \beta \not\le_{M }  c_ \alpha $, 
            
            \item[(c)] $\bar{ c }$ lists $M_ \mathbf{m}$,
            
            \item[(d)] (follows),  if $ L \subseteq L_ \mathbf{m} $ has cardinality $ \le \lambda $ then for some $ \alpha < \delta $  the element $ c_ \alpha $ is an upper bound of $ L $, moreover $ L \in {\mathscr P} _{\mathbf{m}, < c_ \alpha }$.
        \end{enumerate}
    \end{enumerate}

    Now we choose fat bounded $\mathbf{m}_\alpha$ by induction on $\alpha \le \delta $ such that:
    
    \begin{enumerate}
        \item[$(*)_2$] 
        
        \begin{enumerate}
            \item[(a)]  $ \langle \mathbf{m}_ \beta: \beta \le \alpha \rangle $  is $ \le _ \mathbf{M} $-increasing continuous,

            \item[(b)] $ L_{\mathbf{m}_0 } = M$ and  $ u_{\mathbf{m}_0, s}  = M_{< s }, $ (hence $ {\mathscr P}_{\mathbf{m}_0, s }=  [u_{\mathbf{m}_0, s }]^{\le \lambda }$ recalling $\bfm_{0}$ being fat) for $ s \in M$,

            \item[(c)]  for every $ s \in L_{ \mathbf{m}_ \alpha } \setminus M$  for some $\beta < \alpha $ we have $ L_{\mathbf{m}_ \alpha } \models s \le c_ \beta $,

            \item[(d)] if $ \gamma \in [ \alpha, \delta ) $ then  $ u_{\mathbf{m} _ \alpha,  c_\gamma}  =  L_{\mathbf{m} _ \alpha , < c_ \gamma }$,

            \item[(e)]  if $ \alpha = \beta + 1 $ then $\mathbf{m}_{\alpha}(< c_ \beta) \in \mathbf{M} _{\ec},$
            
            \item[(f)] $L_ {\mathbf{m} _ \alpha }$ has cardinality at most $ 2^{ \lambda _2}$  or even $ \lambda _2 $,  but this does not matter, 
            
            \item[(g)] if $ t \in L_{\mathbf{m} _ \alpha} $ then for some  $\beta < \alpha $ we have  $ t/ E''_{\mathbf{m} _ \alpha } \subseteq L_{\mathbf{m} _{\beta + 1}}  \setminus L_{\bfm_{\beta}}.$  
        \end{enumerate} 
    \end{enumerate} 
    
    There is no problem to carry the definition; as:
    
    For $ \alpha = 0 $ we have defined $ \mathbf{m} _0$ in clause (b) of $(*)_2$ above.
    
    For $\alpha $ a limit ordinal use \ref{c28}(1), so in particular $L_{\mathbf{m}} =  \cup \{ L_{\mathbf{m}_ \beta }: \beta < \alpha\} $.  
    
    For $ \alpha = \beta + 1 $  by \ref{c41}  there is $ \mathbf{n} _ \beta \in \mathbf{M} _{\ec}$  such that $ \mathbf{m} _ \beta( < c_\beta)   \le_ \mathbf{M} \mathbf{n} _ \beta $,   \wilog \, we have $ L_{\mathbf{m}_ \beta }  \cap L_{\mathbf{n}_ \beta }= L_{\mathbf{m}( < c_\beta)}$.

    By \ref{e37} below \wilog \,  the cardinality of $ L_ {\mathbf{n} _ \beta} $  is at most $ \lambda _2$. Now apply \ref{b44}(3)  with $ \mathbf{m} _ \beta , L_{\mathbf{m} _\beta , < c_ \beta}, \mathbf{n}_\beta$   here standing for  $ \mathbf{m} _2, L_*, \mathbf{m}_1$ there.  

    So we have carried the induction. Now clearly $ \mathbf{m} _ \delta $ is as promised,  That is, $ (*)_{\mathbf{m} _\delta }$ from part (2) of the claim holds,  hence $ \mathbf{m} \in \mathbf{M} _{\ec}$ by part (2) being cofinal holds by \ref{c37};  so we are done. 
\end{PROOF}  

\newpage

\section {The main conclusion}\label{3}

\subsection{Wider $\bfm$'s} \label{3A}\

Recall that in this section our main interest is in restricting ourselves to lean $\bfm,$ but in \S3C we do not assume this and in \S3A, \S3B, \S3D we rely on \S1, \S2, in particular \S1B 

In \S3B, \S3D we restrict ourselves to lean $\bfm$, but not in \S3A, however the projection defined in \ref{e4}(1) are helpful only in the lean case. 

Note that here we fulfil the promises from \S2,
Now in \S4A we rely on \S3A, \S3C, but we do not rely on \S3B,  \S3D. Lastly, \S4A gives alternative proof of the promises from \S2 proved in \S3D, it  relies on \S3A, \S3C but not on  \S3B, \S3D (except Def \ref{e44}). In \S4B and in \ref{b47} we fulfil additional promises from \cite{Sh:945}. 

We have a debt from \S2, i.e. see discussion \ref{b32}. Toward this we explicate what appear in the proof of \ref{c41}.  We use mainly the  notions of wide, full and ``being in $ \bfM_{\ec}$".

Note that \ref{e4}(2), (4) and \ref{e5n}(3), (4) are of interest exceptionally only for the neat context. 

\begin{definition}\label{e4}
    Let $\mathbf m \in \mathbf M$.
    
    1) For $L \subseteq L_{\bfm}$ we say $p \in \bbP_{\bfm}(L)$ is the projection (to $L$) of $q \in \bbP_{\bfm}(L_{\bfm})$ and write $p = q \upharpoonleft L$ \underline{when}:
    
    \begin{enumerate}
        \item[(a)] $\dom(p) = \dom(q) \cap L,$
        
        \item[(b)] if $s \in \dom(p)$ then:
        
        \begin{enumerate}
            \item[$(\alpha)$] $\tr(p(s)) = \tr(q(s)),$
            
            \item[$(\beta)$] $\{ \name{f}_{p(s), \iota}: \iota < \iota(p(s)) \} = \{ \name{f}_{q(s), \iota}: \iota < \iota(q(s)) \ \text{and} \ \bar{r}_{p(s), \iota} \ \text{is a sequence of}\\ \text{members of} \ L \},$ see Definition \ref{c6}(2). 
        \end{enumerate}
    \end{enumerate}

    2) Let $\cF_{{\mathbf m}, \mu}$ be the set of the functions $f$ such that for
    some $L_1,L_2$:
    
    \begin{enumerate}
        \item[(a)]   $f$ is an isomorphism from $\mathbf m \rest L_1$ onto $\mathbf m \rest L_2,$
          
        \item[(b)]   $L_\ell$ is a subset of $L_{\mathbf m}$ for $\ell=1,2,$
        
        \item[(c)]  $M_{\mathbf m} \subseteq L_\ell$ for $\ell=1,2$ and $f \rest M_{\mathbf m}$ is the identity,
        
        \item[(d)]  $L_\ell$ is $E_{\mathbf m}$-closed, i.e. $M_{\mathbf m} \subseteq L_\ell$ and
        if $t \in L_{\mathbf m} \backslash M_{\mathbf m}$ and $t \in L_\ell$ then $t/E_{\mathbf m} \subseteq L_\ell$ for $\ell=1,2,$
        
        \item[(e)]  $\{t/E''_{\mathbf m}:t \in L_\ell \backslash M_{\mathbf m}\}$
          has cardinality $< \mu $. 
    \end{enumerate}
    
    2A) Let $ {\mathscr F} _ \mathbf{m} = {\mathscr F} _{\mathbf{m}, \lambda_{0}}.$
    
    3) If $L_1,L_2 \subseteq L_{\mathbf m}$ and $f$ is an isomorphism from $\mathbf m \rest L_1$ onto $\mathbf m \rest L_2$ \then \, we let $\hat f$
    be the one-to-one mapping\footnote{We have not said ``order preserving"! still it is  a function from  $ \mathbb{P} _ \mathbf{m} (L_1)$ onto $ \mathbb{P} _ \mathbf{m} (L_1)$  by the way we have defined the $ \mathbb{P} _ \mathbf{m} (L)$-s  and because of \ref{c4}(e)$(\chi )$.}  from $\bbP_{\mathbf m}(L_1)$ onto $\bbP_{\mathbf m}(L_2)$ as in $(*)_4(b)$ of the proof of \ref{c41}.
    
    4) Let $\bbP^-_{\mathbf m}(L)$ be $\{p \in \bbP_{\mathbf m}(L):\fsupp(p) \subseteq L $ and    $ \iota (p(\alpha )) \le 1 $ for every $ \alpha \in  \dom ( p) \}$ with the order inherited from $\bbP_{\mathbf m}$.
\end{definition}

\begin{observation}\label{e5n}
    Let $\mathbf m \in \mathbf M$ and $L \subseteq L_{\mathbf m}$.

    1) The projection of $q \in \bbP_{\bfm}$ to $L$ is well defined and $\in \bbP_{\bfm}(L).$
    
    2) Moreover, it is unique. 
    
    3) If $p \in \bbP_{\bfm}(L)$ is the projection of $q \in \bbP_{\bfm}(L_{\bfm})$ to $L$ \underline{then} $p \leq q$ in $\bbP_{\bfm}.$
    
    4) Each $p \in \bbP_{\bfm}$ is equivalent to $\cS := \{ (p \rest \{ t \}) \upharpoonleft L: t \in \dom(p) \wedge L \in \cP_{\bfm, \leq t}\} \cup \{ p \upharpoonleft M_{\bfm}\};$ the equivalence means $\Vdash_{\bbP_{\bfm}}$``$p \in \name{G}_{\bbP_{\bfm}}$ \underline{iff} $\cS_{p} \subseteq \name{\bfG}_{\bbP_{\bfm}}$''. More specifically it is equivalent to $\cS_{p} = \{ (p \rest \{ t \}) \upharpoonleft L: t \in \dom(p) \wedge L \in \cL_{t} \}$ \underline{when} $\cL_{t}$ satisfies: if $\iota < \iota_{p(s)}$ then for some $L \in \cL_{t},$ (recalling \ref{c6}) we have $\rang(\bar{r}_{p(t), \iota})  \subseteq L.$
    
    5) For every $p \in \bbP_{\mathbf m},p$ is equivalent to $\cS'_p := \{p^{[t]}: t \in \dom(p)  \}$ where $p^{[t]} \in \bbP_{\mathbf m}$ has domain $\{t\}$ and $p(t) = (\tr(p_t),\mathbf B_{p(t)}(\langle \eta_{r_{p(t)}(\zeta)}:\zeta \in w_{p(t)}\rangle)$; recall Definition \ref{c6} for the meaning of $\bfB_{p(t)}$, etc. 
\end{observation}

\begin{remark}\label{e5p}
    1) Note that the choice in Definition \ref{c6}$(c)(\gamma)$ to require such $\langle \name f_{p(t),\iota}:\iota <  \iota(p_t)\rangle$ exists, is necessary for \ref{e5n}(4), which is crucial in the proof of \ref{e50}.
    
    2) In Definition \ref{c39}(1A) we choose ``wide means $\lambda$-wide" as when applying it, if $X = \fsupp(p)$ then for some $Y \subseteq L_{\mathbf m}$ of cardinality $< \lambda, X \subseteq \cup\{t/E_{\mathbf m}:t \in Y\}$.
\end{remark}

\begin{PROOF}{\ref{e5n}}
    Easy e.g.
    
    4) Now if $q \in \cS_{p}$ then $q$ has the form $(p \rest \{ t \}) \upharpoonleft L$ where $L \in \cP_{\bfm, t}$ hence $\Vdash ``p \in \name{\mathbf G}$ implies $q \in \name{\mathbf G}",$ hence $\Vdash$``$p \in \bfG$ implies $\cS_{p} \subseteq \name{\bfG}$''.
    
    For the other direction assume $q \in \bbP_{\mathbf m}$  forces $\cS_p \subseteq \name{\mathbf G} \subseteq \bbP_{\mathbf m}$ and we shall prove that $q$ is compatible with $p$, this suffices, so toward contradiction assume $q,p$ are incompatible.
    
    \Wilog \, $\dom(p) \subseteq \dom(q)$ and recalling $t \in \dom(p) \Rightarrow q \Vdash ``p
    \upharpoonleft (t/E_{\mathbf m}) \in \name{\mathbf G}"$ clearly $s \in \dom(p)  \Rightarrow q \Vdash
    ``{\rm tr}(p(s)) \subseteq \name \eta_s"$  so necessarily $s \in \dom(p) \Rightarrow \tr(p(s))  \trianglelefteq  \tr(q(s))$.  Recalling \ref{c11}(6),  as $p,q$ are incompatible
    there are $s \in \dom(p) \cap \dom(q)$ and $q_1$ such that $q \rest L_{\mathbf m,< s} \le q_1 \in \bbP_{\mathbf m}(L_{\mathbf m,<s})$ and 
    $q_1 \Vdash ``q(s),p(s)$ are incompatible in $\bbQ_{\bar\theta}"$.
    
    As $\tr(p(s)) \trianglelefteq  \tr(q(s))$ this implies $q_1 \Vdash ``\tr(q(s)),p(s)$ are incompatible, so recalling $q \Vdash$``$\tr(p(s)) \subseteq \eta_{s}$'' this implies $f_{p(s)} \rest \ell g(\tr(q(s))) \ntrianglelefteq \tr(q(s))"$.   Recalling Definition \ref{c6}$(2)(c)(\gamma)$,
    $q_1 \Vdash _{\mathbb{P} _{\mathbf{m}, s }}$ 
    ``there is $\iota < \iota(s,p)$
    such that $f_{p(s),\iota},\tr(q(s))$ are incompatible". Possibly increasing $q_1$, we can fix $\iota$.   But letting $u \in \cP_{\bfm{, s}}$ be such that $\bar r_{p(s),\iota} \subseteq u$ this implies that $q_1 \Vdash ``(p \rest \{ s \}) \upharpoonleft u \notin \name{\mathbf G}$ or $\tr(q(s)) \nsubseteq \name\eta_s"$. However, $q_1,q$ are compatible and this contradicts the choice of $q$.
\end{PROOF}

\begin{claim} \label{e7}
    1) For $ \chi \ge 2^{\lambda_2}$ the
    $\mathbf{n} \in \mathbf{M} _ \chi $ constructed in \ref{c41} satisfies: if $\mathbf n \le_{\mathbf M} \mathbf n_1$ \then \,  $\mathbf n_1$ is full and  wide,  even $ \lambda _2$-wide and if $\bfn_{1} \in \bfM_{\chi}$ even very wide.
    
    2) If $\mathbf n \in \mathbf M_{\ec}$ and 
    $\mathbf n \le_{\mathbf M} \mathbf n_1$ \then \, 
    $\mathbf n_1 \in \mathbf M_{\ec}$.
    
    3) If $ \mathbf{m} \in \mathbf{M} _ \chi $ is
    full and very wide (or just $ \lambda_2 $-wide and even $\lambda_{0} $-wide), \then \, $ \mathbf{m} \in \mathbf{M} _{\ec}.$
    
    4) If $ \mathbf{m} \in \mathbf{M},$ \then \,  there is a very wide full $ \mathbf{n} \in \mathbf{M}$ such that $ \mathbf{m} \le_{\bfM} \mathbf{n}$. 
\end{claim}

\begin{PROOF}{\ref{e7}}
    1) Holds by the proof of \ref{c41}.
    
    2) Holds by Definition \ref{c34}(1),(2).
    
    3),4)  By the proof of \ref{c41}.
\end{PROOF}

\begin{claim}\label{e10}
    Assume $\mathbf m$ is $ \mu $-wide where $ \mu \ge \lambda_{0}.$

    1) If $f \in \cF_{{\mathbf m}, \mu }$ and $X \subseteq L_{\mathbf m}$ has cardinality
    $< \mu,$  \then \, there is $g$ such that:
    
    \begin{enumerate}
        \item[(a)]  $g \in \cF_{{\mathbf m}, \infty}$ and even belongs  to $ {\mathscr F} _{\mathbf{m}, \mu }$,
        
        \item[(b)]  $f \subseteq g$, 
        
        \item[(c)]  $\dom(g) = \rang(g)$,  
        
        \item[(d)]  $X \subseteq \dom(g)$.
    \end{enumerate}
    
    2) If $g \in \cF_{\mathbf m, \mu }$  and $\dom(g) = \rang(g)$ \then \, $g^{+ \mathbf m} = g \cup \id_{L_{\mathbf m} \backslash \dom(g)}$ is an automorphism of $\mathbf m$.
    
    3) If $f$ is an automorphism of $\mathbf m$ \then \, it naturally induces an automorphism $\hat f$ of $\bbP_{\mathbf m}(L_{\mathbf m})$ similarly to $\hat f$ from $(*)_4(b)$ of the proof of \ref{c41} and it induces an automorphism of $ \mathbb{P} _ \mathbf{m} [L_ \mathbf{m} ]$ as well; abusing our notation we denote both by  $ \hat{f}.$
    
    4) If $f \in \cF_{\mathbf m , \mu }$  \then \, it induces an isomorphism $\hat f$ from  $ \mathbb{P} _ \mathbf{m} [\dom (f)]$  onto  $ \mathbb{P} _ \mathbf{m} [\rang(f)]$  hence (as above) from 
    $\bbP_{\mathbf m}(\dom(f))$ onto $\bbP_{\mathbf m}(\rang(f))$.
    
    5) If $ p \in \mathbb{P} _ \mathbf{m} $ \then \, 
    the set $ \{ t/ E_ \mathbf{m} : t \in \wsupp(p) \} $  has cardinality $ < \lambda $.   
\end{claim}

\begin{PROOF}{\ref{e10}}
    1) Easy by the definition of wide in \ref{c39}(1) and of $\cF_{\mathbf m}$  in \ref{e4}(2), in particular clause (e).

    2) Just read the definition of $\mathbf m \in \mathbf M$ and of $f \in \cF_{\mathbf m}$, in particular:
    
    \begin{enumerate}
        \item[(a)]  if $t_1,t_2 \in L_{\mathbf m} \backslash M_{\mathbf m}$ are not $E'_{\mathbf m}$-equivalent then $(t_1/E_{\mathbf m}) \cap (t_2/E_{\mathbf m}) = M_{\mathbf m}$ and $\le_{\mathbf m} \rest (t_1/E_{\mathbf m} \cup t_2/E_{\mathbf m})$ is determined by $\le_{\mathbf m} \rest (t_1/E_{\mathbf m}),\le_{\mathbf m} \rest (t_2/E_{\mathbf m}),$
        
        \item[(b)]  $g \rest M_{\mathbf m} = \id_{M_{\mathbf m}}$.
    \end{enumerate}
    
    3) Naturally by the definition.  
    
    4) Let $g \in \cF$ be as in part (1) and let $h = g^{+ \mathbf m}$ so an automorphism of $\mathbf m$ which extends $g$ as in part (2).  So $\hat h$ is an automorphism of $\bbP_{\mathbf m}(L_{\mathbf m})$ and clearly $\hat f = \hat h \rest \bbP_{\mathbf m}(\dom(f))$ is as required.
   
    5) Is  clear,  see \ref{c8}(f).
\end{PROOF}

\begin{claim}\label{e16}
    Let $\mathbf m \in \mathbf M$  and $ \mu \ge \lambda_{0}.$ 
    
    If $f_1,f_2 \in \cF_{\mathbf m, \mu }$ \then:  
    
    \begin{enumerate}
        \item[(a)]  $f_1 \subseteq f_2 \Rightarrow \hat f_1 \subseteq  \hat f_2,$
        
        \item[(b)]   $f_1 = f^{-1}_2 \Rightarrow \hat f_1 = (\hat f_2)^{-1}$.
    \end{enumerate}
\end{claim}

\begin{PROOF}{\ref{e16}}
    Just consider the definition, see \ref{e4}(3) and $(*)_4(b)$ of the proof of \ref{c41}.
\end{PROOF}

\subsection {Ordinal equivalence} \label{3B}\

\begin{context}\label{e18}
    All $\bfm$-s are lean\footnote{So maybe we can use $\lambda_{0} = \lambda.$}. 
\end{context}

\begin{observation}\label{e19}
    1) $\bbP_{\bfm}^{-}(L) \subseteq \bbP_{\bfm}(L),$ see Definition \ref{e4}(4).
    
    2) For every $p \in \bbP_{\bfm}$ there is a sequence $\langle p_{i}: i < i(*) \rangle$ of $< \lambda$ members of $\bbP_{\bfm}^{-},$ (see \ref{e4}(6)) such that $\Vdash_{\bbP_{\bfm}(L_{\bfm})}$``$p \in \name{\bfG} \iff \{ p_{i}: i <i(*) \} \subseteq \name{\bfG}$''.
\end{observation}

\begin{PROOF}{\ref{e19}}
    1) By their definitions.
    
    2) Should be clear, see Definition \ref{e4}(4) and \ref{e5n}(3).  
\end{PROOF}
    
\begin{remark}
    \label{e21}
    1) Observation \ref{e19} is not used.
    
    2) Probably we can avoid using ``wide" and prove earlier the density of $\mathbf M_{\ec}$
    with smaller cardinality but the present way seems more transparent.
\end{remark}

\begin{definition}\label{e24}
    Assume $\mathbf m \in \mathbf M$.
    
    1) Let $\cY_{\mathbf m}$ be the set of pairs $(t,\bar s)$ such that  $t \in L_{\mathbf m} \backslash M_{\mathbf m}$ and $\bar s \in 
    {}^\zeta(t/E''_{\mathbf m})$ for some $\zeta < \lambda^+$;  we may write $\bar s$ instead of $(t,\bar s)$  as usually $\bar s$ determines $t/E''_ \mathbf{m} $,  but this is the only information  about $t$  that matter.  We could have used instead pairs  $ (t/E''_\mathbf{m} , \bar{s})$.
          
    2)  By induction on the ordinal $\gamma$ we define when  $(t_1,\bar s_1),(t_2,\bar s_2)$ are $\gamma$-equivalent in $\mathbf m$ or are $(\mathbf m,\gamma)$-equivalent:
    
    \begin{enumerate}
        \item[(a)]  if $\gamma=0$, \then \, letting $L_\ell = (M_{\mathbf m} \cup \rang(\bar s_\ell))$ for $\ell=1,2$ there is $h$ such that: 
        
        \begin{enumerate}
            \item[$(\alpha)$]  $h$ is an isomorphism from $\mathbf m \rest L_1$ onto $\mathbf m \rest L_2$, 
            
            \item[$(\beta)$]  $h$ maps $\bar s_1$ to $\bar s_2$,
            
            \item[$(\gamma)$]  $h \rest M_{\mathbf m}$ is the identity,
            
            \item[$(\delta)$]  $h$ induces an isomorphism from  $\bbP_{\mathbf m}(L_1)$ onto $\bbP_{\mathbf m}(L_2)$ (as defined in \ref{c4}$(*)_4(b)$),  
            
            \item[$(\varepsilon)$]  moreover, $h$ induces an isomorphism from $\bbP_{\mathbf m}[L_1]$ onto $\bbP_{\mathbf m}[L_2]$, as defined in \ref{b20}, so $p^*_{t,\eta} \mapsto p^*_{h(t),\eta}$, see \ref{b11}(3), 
        \end{enumerate}
        
        \item[(b)]  if $\gamma = \beta +1$ \then \, for every  $\ell \in \{1,2\}$ for every 
        $\varepsilon < \lambda^+$ and  $\bar s'_\ell \in {}^\varepsilon(t_\ell/E''_{\mathbf m})$  there is  $\bar s'_{3 - \ell} \in {}^\varepsilon(t_{3 - \ell}/E''_{\mathbf m})$ such that $(t_1,\bar s_1 \char 94 \bar s'_1),(t_2,\bar s_2 \char 94 \bar s'_2)$ are $\beta$-equivalent, 
        
        \item[(c)]  if $\gamma$ is a limit ordinal \then \, $(t_1,\bar s_1), (t_2,\bar s_2)$ are $\beta$-equivalent for every $\beta < \gamma$.
    \end{enumerate}
\end{definition}

\begin{remark}\label{e26}
    1) Note above that if $\bar s_\ell$ is the empty sequence \then \, $t_\ell$ would not be determined by $\bar s_\ell$, still in those cases the
    equivalence just means $\bar s_1 = \bar s_2$.

    2) We can use $t/E_{\mathbf m}$ or $t/E'_{\mathbf m}$ instead of  $t/E''_{\mathbf m}$ as everything is over $M_{\mathbf m}$.
\end{remark}

\begin{claim}\label{e27}
    For $\mathbf m \in \mathbf M$ and ordinal $\alpha$ the number of equivalence classes of ``being $(\mathbf m,\alpha)$-equivalent" is $\le \beth_{1 +\alpha +1}(\lambda_1)$. 
\end{claim}

\begin{PROOF}{\ref{e27}}
    By induction on $\alpha$.

    \underline{Case 1}:  $\alpha = 0$:
    
    Note that the set of elements of $\bbP_{\mathbf m}(M_{\mathbf m} \cup \rang(\bar s))$ has cardinality $\le 2^{\lambda_1}$ (and even $\le
    (\lambda_1)^\lambda$) and depends just on $\mathbf m \rest (M_{\mathbf m} \cup \rang(\bar s))$ but there are $\beth_2(\lambda_1)$  possibilities for the quasi order on $\bbP_{\mathbf m}(L_1)$ and even for $\bbP_{\mathbf m}[L_1]$.  
    
    \underline{Case 2}:  $\alpha$ is a limit ordinal:
    
    By clause (c) of Definition \ref{e24}, the number of $\alpha$-equivalence classes is $\le \prod\limits_{\beta < \alpha}$ (the
    number of $\beta$-equivalence classes) 
    $\le \prod\limits_{\beta < \alpha} \beth_{1 + \beta +1} (\lambda_1) \le  (\beth_{1 + \alpha +1}(\lambda_1))^{\beth_{1 + \alpha}} = 
    \beth_{1 + \alpha +1}(\lambda_1)$.
    
    \underline{Case 3}:  $\alpha = \beta +1$:
    
    Clearly every $\alpha$-equivalence class can be coded as a set of $\beta$-equivalence classes hence the number of $\alpha$-equivalence
    classes is $\le 2^{\beth_{1 + \beta +1}(\lambda_1)}  = \beth_{1 + \beta +2}(\lambda_1) = \beth_{1 + \alpha +1}(\lambda_1)$, as promised.
\end{PROOF}

\begin{definition}\label{e28}
    For an ordinal $\beta$, let $\cG_{\mathbf m,\beta}$ be the set of  functions $f$ such that for some $t^\ell_i,\bar s^\ell_i$ for $i<i(*)$ and $\ell \in \{1,2\}$ we have:
    
    \begin{enumerate}
        \item[(a)]  $i(*) < \lambda^+$,  
        
        \item[(b)]  $\langle t^\ell_i:i < i(*)\rangle$ is a sequence of pairwise non-$E''_{\mathbf m}$-equivalent members of $L_{\mathbf m} \backslash M_{\mathbf m}$,  
          
        \item[(c)]  $\bar s^\ell_i \in {}^{\zeta(i)}(t^\ell_i/E''_{\mathbf m})$ where $\zeta(i) < \lambda^+$, 
          
        \item[(d)]  $(t^1_i,\bar s^1_i),(t^2_i,\bar s^2_i)$ are $\beta$-equivalent (members of $\cY_{\bfm}$), 
          
        \item[(e)]  $f$ is an isomorphism from $\mathbf m \rest L_1$ onto $\mathbf m \rest L_2$ when $L_\ell = \cup\{\rang(\bar s^\ell_i):i < i(*)\} \cup M_{\mathbf m}$, 
          
        \item[(f)]  $f \rest M_{\mathbf m} =$ the identity, 
        
        \item[(g)]  $f$ maps $\bar s^1_i$ to $\bar s^2_i$ for $i<i(*)$.
    \end{enumerate}

    2) For $f \in \cG_{\mathbf m,0}$ we define $\hat f$ as the mapping from   $\bbP_{\mathbf m}(\dom(f))$ onto $\bbP_{\mathbf m}(\rang(f))$ induced by $f$; see clause \ref{e24}(2)(a)$(\varepsilon)$;  (clearly well defined 1-to-1 function, but does it preserve the order? we shall return to this in \ref{e33}). 
\end{definition}

\subsection {Representing $p \in \bbP_{\bfm}[M_{\bfm}]$}\label{3C}\

Applying this subsection in \S3D we may assume all $\bfm$-s are lean and so maybe $\lambda_{0} = \lambda$ is O.K., but certainly not applying it in \S4.  

\begin{claim}\label{e29}
    Assume $\mathbf m$ is  $ \mu $-wide and $ \mu \ge \lambda_{0}.$
    
    1) The conditions $p,q \in  \bbP_{\mathbf m}(L_\mathbf{m})$  are compatible \when \, for some  $\psi$ the following condition holds:
    
    \begin{enumerate}
        \item[$(\stt)_{p,q,\psi}$]  $(a) \quad \psi \in \bbP_{\mathbf m}[M_{\mathbf m}]$, 
        
        \item[${{}}$]  $(b) \quad p, q \in \bbP_{\bfm}(L_{\bfm})$ and $\wsupp(p) \cap \wsupp(q) \subseteq M_{\mathbf m}$, see Definition \ref{c7}(1)(b), equivalently   $\big(  s \in \fsupp(p) \backslash M_{\mathbf m}\big) \wedge  \big( t \in \fsupp(q) \backslash M_{\mathbf m}\big) 
        \Rightarrow \neg(s E''_{\mathbf m} t)$,
        
        \item[${{}}$]  $(c) \quad$ if $\psi \le \varphi \in \bbP_{\mathbf m}[M_{\mathbf m}]$ then $\varphi,p$ are compatible in 
        $\bbP_{\mathbf m}[L_{\mathbf m}]$,
        
        \item[${{}}$]  $(d) \quad \psi,q$ are compatible in $\bbP_{\mathbf m}[L_{\mathbf m}]$, equivalently $q \nVdash_{\bbP_{\mathbf m}} ``\psi[\name{\mathbf G}]=$ false".
    \end{enumerate}

    2) For a dense set of $\psi \in \bbP_{\bfm}[M_{\bfm}]$ there are $\bar L,\bar p$
    such that:
     
    \begin{enumerate}
        \item[(a)]  $\bar p = \langle p_\varp:\varp < \mu  \rangle \in {}^{\mu}  (\bbP_{\bfm})$, 
        
        \item[(b)]  $\bar{ L }  = \langle L_\varp:\varp <  \mu  \rangle$ where $ \fsupp (p_\varepsilon ) \subseteq  L_ \varepsilon $, 
        
        \item[(c)]  $\bfm \rest L_\varp \le_{\bfM} \bfm$ so in particular $t \in L_\varp \backslash M_{\bfm} \Rightarrow  t/E_{\bfm} \subseteq L_\varp$, 
        
        \item[(d)]  $\langle L_\varp \backslash M_{\bfm}:\varp <  \mu  \rangle$ are pairwise disjoint,  
        
        \item[(e)]  $(L_\varp \backslash M_{\bfm})/E''_{\bfm}$  has cardinality $ <   
        \lambda_{0},$ (see \ref{c2}(4) and \ref{c8}(f)($\gamma$)),    
        
        \item[(f)]  for every permutation $\pi$ of  $ \mu $  there is an automorphism $\hat f$ of $\bfm$ mapping  $(L_\varp,p_\varp)$ to $(L_{\pi(\varp)},p_{\pi(\varp)})$
        for $ \varepsilon < \mu,$
    
        \item[(g)]  if $u \subseteq \ \mu $  has cardinality $\lambda$ then $\psi,\bigvee\limits_{\varp   \in u }  p_\varp$ are equivalent in  $\bbP_{\bfm}[L_{\bfm}]$, i.e.  $\psi \le \bigcup\limits_{\varp \in u }  p_\varp \le \psi$.
    \end{enumerate}
    
    3) Assume that  $ L $ is a $ \mu $-wide  initial segment of $ L_ \mathbf{m} $ and $ \psi _0 \in \mathbb{P} _ \mathbf{m}   [ M_ \mathbf{m} \cap L ]$.  \underline{Then} there is a pair  $ ( \psi, \bar{ p } ) $ satisfying $ \psi _0 \le \psi \in \mathbb{P} _ \mathbf{m} [ M_ \mathbf{m} \cap L ] $ and clauses (a)-(g) above  hold   and:
    
    \begin{enumerate}
        \item[(h)]  if $ \varepsilon < \mu $  then  $ p_ \varepsilon \in \mathbb{P} _ \mathbf{m}(L). $
    \end{enumerate}
    
    Also we can add:     
    
    \begin{enumerate}
        \item[(i)]  the sequence $ \langle \name{ \eta} _s: s \in L \cap M_ \mathbf{m} \rangle $ is a generic for $ \mathbb{P} _ \mathbf{m}  [ L \cap M_ \mathbf{m} ]$, that is it determines $ \name{ \mathbf{G}} _{\mathbb{P}  _\mathbf{m} } [ L \cap M_{\mathbf{m} }]$.
    \end{enumerate}    
\end{claim} 

\begin{remark}\label{e31} 
    1) In \ref{e29}(1)  instead of $ \stt _{p,q,\psi}$ we  can use the stronger statement: 
    
    \begin{enumerate} 
        \item[$(\stt)'_{p,q,\psi}$] as there but  omit
        clause (d) and add to clause (c): also  $\varphi,q$ are compatible in $\bbP_{\mathbf m}[L_{\mathbf m}]$, 
     \end{enumerate} 
     
    But the present choice is more convenient in the proof of \ref{e29}(1).
    
    2) We use $\lambda > \aleph_0$ in the proof, to eliminate it we can imitate the completeness theorem\footnote{but we give details. First as a warm up notice that (for $\lambda = \aleph_{0})$: 

    \begin{enumerate}
        \item[$(\ast)$] if $r \in \bbP_{\bfm}$ \underline{then} we can find $\cT$ and $\bar{r}, \bar{s}$ such that:  
    \end{enumerate}
    
    \begin{enumerate}
        \item[(a)]  
        
        \begin{enumerate}
            \item[$(\alpha)$] $\cT$ is a sub-tree of ${}^{\omega >} \omega$ which is well founded,
        
            \item[($\beta$)] if $\eta \in \cT,$ \underline{then} $\suc_{\cT}(\eta)$ is empty or is $\omega.$  
        \end{enumerate}
        
        \item[(b)] $\bar{r} = \langle r_{\eta}: \eta \in \cT \rangle$ and $r_{\langle \rangle} = r,$
        
        \item[(c)] $r_{\eta} \in \bbP_{\bfm}$ and $r_{\eta} \subseteq r_{\nu}$ for $\eta \unlhd \nu \in \cT,$
        
        \item[(d)] $\bar{s} = \langle s_{\eta}: \eta \in \cT \setminus \max(\cT) \rangle$ such that $\eta \lhd \nu \Rightarrow s_{\nu} \nleq s_{\eta},$
        
        \item[(e)] if $\eta = \nu^{\smallfrown} \langle k \rangle \in \cT,$ \underline{then} $s_{\nu} \in \dom(r_{\nu}) \cap M_{\bfm}$ and $r_{\nu} \rest (\dom(r_{\nu}) \setminus L_{\bfm(\leq s_{\nu})}) = r_{\eta} \rest (\dom(r_{\nu}) \setminus L_{\bfm(s_{\nu})}),$
        
        \item[(f)] if $ \eta \in \max(\cT),$ \underline{then} $\dom(r_{\eta}) \cap M_{\bfm} = \{ s_{\eta \rest \ell}: 0 < \ell \leq \lg(\eta)  \},$
        
        \item[(g)] if $\eta \in \cT \setminus \max(\cT),$ \underline{then} for some $k$ we have: 
        
        \begin{itemize}
            \item if $\ell \geq k,$ then $\tr(r_{\eta^{\smallfrown} \langle i \rangle}(s_{\eta}))$ has length $\ell,$
            
            \item if $\ell \geq k, \varrho = \tr(r_{\eta^{\smallfrown} \langle i \rangle}(s_{\eta}))$ for some $\varrho \in \Pi_{\varp < \ell} \theta_{\varp},$ \underline{then} for every $\rho \in \Pi_{\varp < \ell} \theta_{\varp}$ satisfying $\varrho \leq \ell$ and $\tr(r_{\eta}(s_{\eta})) \unlhd \rho$ for some $j < \omega$ we have $\rho = \tr(r_{\eta^{\smallfrown} \langle j \rangle}(s_{\eta})).$ 
        \end{itemize}
    \end{enumerate}
    
    This can be proved by induction on $\sup \{ \rm{rk}(M_{\bfm}(s)) + 1: s \in \dom(r) \cap M_{\bfm} \}.$
    
    Let $\langle s_{i}: i < i_{\ast} \rangle$ lists $M_{\bfm}$ such that $s_{i} <_{\bfm} s_{j} \Rightarrow i < j,$ and let $s_{i_{\ast}} = \infty.$ For $i \leq i_{\ast}$ let $L_{i} = \bigcup \{ L_{\bfm(\leq s_{j})}: j < i \},$ it is as an initial segment of $L_{\bfm}.$ We prove by induction on $i \leq i_{\ast}$ that the statement holds when $p, q \in \bbP_{\bfm}(L_{i}).$ For $i = 0$ this is trivial and limit $i$ it is. So assume $i = j+1,$ now if $s_{i} \notin \dom(p) \cup \dom(q)$ this is trivial and if $s_{j} \in \dom(p) \setminus \dom(q)$ this is obvious. Similarly if $s_{j} \in \dom(q) \setminus \dom(p).$ So assume $s_{j} = \dom(p) \cap \dom(q).$ as in the proof of \ref{e29}(1), without loss of generality $\tr(q(s_{j})) \unlhd \tr(p(s_{j})).$ As in the proof of \ref{e29}(1), for some $q_{1} \in \bbP_{\bfq},$ we have $(\ast)_{q_{1}, p, \psi}$ and $q \leq q_{1}$ and $\lg(\tr(q_{1}(s))) > \lg(\tr(p(s))),$ hence $\tr(p(s)) \lhd \tr(q_{1}(s)).$ 
    
    Clearly $(\ast)_{q_{1} \rest L_{j}, p \rest L_{j}, \psi}$ holds, therefore $q_{1} \rest L_{s}, p \rest L_{s}$ are compatible in $\bbP_{\bfm},$ hence in $\bbP_{\bfm}(L_{j}),$ and let $r \in \bbP_{\bfm}(L_{j})$ be a common upper bound. Now, $r$ forces (i.e. $\Vdash_{\bbP_{\bfm}(L_{j})}$) then $\name{f}_{q(s)} \rest ( \lg(\tr(q(s))), \lg(\tr(p(s)))) \leq \tr(q_{1}(s)),$ hence $r \Vdash_{\bbP_{\bfm}(L_{j})}$``$p(s), q(s)$ are compatible in $\name{\bbQ}_{s_{j}}$'', therefore $r, p, q$ have a common upper bound. So we are done.} for $\bbL_{\aleph_1,\aleph_0}.$
\end{remark}

\begin{PROOF}{\ref{e29}}
    1) We choose $(p_n,q_n,\psi_n)$ by induction on $n$ such that:
    
    \begin{enumerate}
        \item[$\boxplus_n$]  
        
        \begin{enumerate}
            \item[(a)] 
            
            \begin{enumerate}
                \item[($\alpha$)] $(\stt)_{p_n,q_n,\psi_n}$ holds if $n$ is even,
                
                \item[($\beta$)] $(\stt)_{q_n,p_n,\psi_n}$ holds if $n$ is odd,
            \end{enumerate}
            
            \item[(b)] $(p_0,q_0,\psi_0) = (p,q,\psi)$,
             
            \item[(c)] if $n=2m+1$ and $s \in \dom(p_{2m}) \cap M_{\mathbf m},$ \then \, $s \in \dom(q_{2m+1})$,  and $\tr(p_{2m}(s)) \triangleleft \tr(q_{2m+1}(s))$, 
            
            \item[(d)] if $n=2m+2$ and $s \in \dom(q_{2m+1}) \cap M_{\mathbf m},$ \then \, $s \in \dom(p_{2m+2})$   and $\tr(q_{2m+1}(s)) \triangleleft \tr(p_{2m+2}(s))$, 
            
            \item[(e)] if $n=m+1$ then $p_m \le p_n,q_m \le q_n$.
        \end{enumerate}
    \end{enumerate}
    
    \underline{Case 1}:  For $n=0$ use clause (b).
    
    \underline{Case 2}:   $n=2m+1$.
    
    So the triple $(p_{2m},q_{2m},\psi_{2m})$ is well defined, let $u_{2m}= \dom(p_{2m}) \cap M_{\mathbf m}$ and let $\bar\nu = \langle \nu_s:s \in u_{2m}\rangle$ be defined by $\nu_s = \tr(p_{2m}(s))$.
    
    Clearly,
    
    \begin{enumerate}
        \item[$(*)_1$]  $\psi_{2m} \Vdash p^*_{s,\nu_s}$ for $s \in u_{2m}$.
    \end{enumerate}
    
    [Why?  Clearly $p_{2m} \Vdash_{\bbP_{\bfm}} p^*_{s,\nu_s}$, i.e. $p^*_{s,\nu_s} \le p_{2m}$ in $\bbP_{\bfm}(L_{\bfm}),$ hence in $\bbP_{\mathbf m}[L_{\mathbf m}]$ and therefore, if $\psi_{2m} \nVdash p^*_{s,\nu_s},$ then $\psi' = \psi_{2m} \wedge \neg p^*_{s,\nu_s} \in \bbP_{\mathbf m}[M_{\mathbf m}]$ is $\ge \psi_{2m}$, hence compatible with $p_{2m}$, contradiction, see clause (c) in $(\stt)_{p,q,\psi}$  which holds by $ \boxplus _{2m}(a)(\alpha ) $.]

    \begin{enumerate}
        \item[$(*)_2$] there is $q'_{2m} \in \bbP_{\mathbf m}(L_{\mathbf m})$ which is above $q_{2m}$ and above $\psi_{2m}$ and naturally $u_{2m} \subseteq \dom(q'_{2m})$ hence $s \in u_{2m}$ implies $\nu_s \subseteq \tr(q'_{2m}(s)).$
    \end{enumerate}
    
    [Why?  By clause (d) of $(\stt)_{p_{2m},q_{2m},\psi_{2m}}$ which holds by $\boxplus_{2m}(a)(\alpha)$ recalling $\bbP_{\mathbf m}(L_{\mathbf m})$ is dense $\bbP_{\mathbf m}[L_{\mathbf m}]$; the ``hence" by $(*)_1$.]
    
    \begin{enumerate}
        \item[$(*)_3$]  there is $\psi'_{2m} \in \bbP_{\mathbf m}[M_m]$ such that:
        
        \begin{enumerate}
            \item[(a)]  if $\psi'_{2m} \le \varphi \in  \bbP_{\mathbf m}[M_{\mathbf m}]$ then $\varphi,q'_{2m}$ are  compatible in $\bbP_{\mathbf m}[L_{\mathbf m}]$, 
            
            \item[(b)]  if $s \in u_{2m}$ then $\psi'_{2m} \Vdash p^*_{s,\nu_s}$, 
            
            \item[(c)]  $\psi_{2m} \le \psi'_{2m}$.
        \end{enumerate}
    \end{enumerate}
    
    [Why?  Obvious using the $\lambda^+$-c.c., i.e. $\psi'_{2m} = \psi_{2m} \wedge \neg(\vee \{\varphi:\varphi \in \cI\})$  where $\cI$ is a maximal anti-chain  of members  $ \varphi 
    \in \bbP_{\mathbf m}[M_{\mathbf m}]$ satisfying
    $ \varphi  \perp q'_{2m}$ in $\bbP_{\mathbf m}[L_{\mathbf m}];$ see more in \ref{e29}.]
    
    \begin{enumerate}
        \item[$(*)_4$]  \wilog \, $\wsupp(q'_{2m}) \cap \wsupp(p_{2m}) \subseteq  M_{\mathbf m}$.
    \end{enumerate}
    
    [Why?  As $\mathbf m$ is  $ \mu $-  
    wide using an automorphism of $\mathbf m$
    which is the identity on $ \wsupp ( q_{2m}),$  i.e. by \ref{e10}. Even if $\bfm$ is fat this is fine.]
    
    Lastly, let $p_n =  p'_{2m},q_n = q'_{2m},\psi_n = \psi'_{2m}$ and check.
    
    \underline{Case 3}:  $n=2m+2.$
    
    Similar to case 2 with  the roles of the $p$'s and the $q$'s interchanged. 

    Having carried the induction we  can define $p_*$ as the upper bound of, in fact the union of 
    $\{p_n:n < \omega\}$ as in \ref{c11}(3A), in particular:
    
    \begin{enumerate}
        \item[$(*)_7$ (a)]  $(\dom(p_*) = \underset{n}\bigcup \, \dom(p_n)$; in fact, also $\fsupp(p_*) = \bigcup\limits_{n} \fsupp(p_n)$ and $ \wsupp ( p_* ) = \underset{n} \bigcup \,  \wsupp (p_n)$, 
        
        \item[(b)] if  $ s \in \dom(p_*)$  and $ n $ is minimal  such that  $s \in \dom(p_n)$ then $\tr(p_*(s)) = \bigcup\limits_{k \ge n} \tr(p_k(s))$  and  $ \{  f_{p_* , \iota }: \iota < \iota (p_*)(s)\} $  is equal to $ \{\tr(p_*(s))\cup f_{p_k , \iota }\rest [\lg(\tr(p_*(s))), \lambda ) : \iota < \iota (p_k(s))$ for some $ k \in [n, \omega ) \} $.
    \end{enumerate}
    
    Similarly let $q_*$ be the upper bound of,
    in fact the union of  $\{q_n:n < \omega\}$ as in
    \ref{c11}(3A), so again, in particular:
    
    \begin{enumerate}
        \item[$(*)_8$ (a)]  $\dom(q_*) = \bigcup\limits_{n} \dom(q_n)$,  and  also $\fsupp(q_*) = \bigcup\limits_{n} \fsupp(q_n)$ and $ \wsupp (q_*)= \cup _n \wsupp (q_n)$, 
        
        \item[(b)] if  $s \in \dom(p_*)$    and $ n $ is minimal  such that
        $s \in \dom(q_n)$ \then:
        
        \begin{enumerate} 
            \item[$ \bullet _1$] $\tr(q_*(s)) =
            \bigcup\limits_{k \ge n} \tr(q_k(s))$, 
            \item[$ \bullet _2$] $ \{  f_{q_* , \iota }: \iota < \iota (q_*)(s)\} $  is equal to $ \{\tr(p_*(s))\cup  f_{p_k , \iota }\rest [\lg(\tr(p_*(s))), \lambda ) : \iota < \iota (q_k(s))$ for some $ k \in [n, \omega ) \} $. 
        \end{enumerate} 
    \end{enumerate}
    
    Hence, 
    
    \begin{enumerate}
        \item[$(*)_9$]  $(\rm{a}) \quad p_*,q_* \in \bbP_{\mathbf m}$, 
        
        \item[${{}}$]  $(\rm{b}) \quad \dom(p_*) \cap \dom(q_*) \subseteq M_{\mathbf m}$,  moreover, $\wsupp(p_*) \cap \wsupp(q_*) \subseteq M_{\mathbf m}$, 
        
        \item[${{}}$]  $(\rm{c}) \quad \dom(p_*) \cap M_{\mathbf m} =  \dom(q_*) \cap M_{\mathbf m}$, 
        
        \item[${{}}$]  $(\rm{d}) \quad$ if $s \in \dom(p_*) \cap M_{\mathbf m}$, equivalently, $s \in \dom(p_*) \cap \dom(q_*)$ then:  $\tr(p_*(s)) = \tr(q_*(s))$.
    \end{enumerate}
    
    [Why?  Clause (a) by properties of $\bbP_{\mathbf m}$ and $p_n \le p_{n+1},q_n \le q_{n+1}$ see above, clause (b) as $\dom(p_{2m}) \cap \dom(q_{2m}) \subseteq M_{\mathbf m}$ as $(\stt)_{p_{2m},q_{2m},\psi_{2m}}$. Clause (c) by
    $\boxplus_n(c),(d)$, the first conclusion and clause (d) by $\boxplus_n(c),(d)$, the second conclusion.]
    
    It follows that $p_*,q_*$ are compatible in $\bbP_{\mathbf m}$ but $p=p_0 \le p_*,q = q_0 \le q_*$, so $p,q$ are compatible as promised.
     
    2) Let $\psi_0 \in \bbP_{\bfm} [M_ \mathbf{m} ]$ be given.  Let $p \in \bbP_{\bfm}$ be such that $p \Vdash_{\bbP_{\bfm}} ``\psi_0[\name{\bfG}] = \true"$.
     
    Let  $ {\mathscr I }_0 =  \{ \varphi :  \varphi 
    \in \mathbb{P} _ \mathbf{m} [ M_ \mathbf{m} ]$ 
    and $ \varphi , p $ are incompatible in   $ \mathbb{P} _ \mathbf{m} [L_ \mathbf{m} ] \} $
    and let $ {\mathscr I } _ 1 $ be a maximal set of pairwise  incompatible members of $ {\mathscr I }_0 $. As $ \mathbb{P} _ \mathbf{m} [L_ \mathbf{m} ] $ satisfies the $ \lambda ^ + $-c.c., clearly 
    $ {\mathscr I } _1 $ has cardinality at most $ \lambda $ and let $ \psi = \wedge \{  \neg  \varphi : \varphi \in {\mathscr I } _1\} $. Clearly we have:
    
    \begin{enumerate}
        \item[$(*) _ 1$] $\psi \in \bbP_{\bfm}[M_{\bfm}]$ and:
    
        \begin{enumerate} 
            \item[(a)]  if $ \psi \le \varphi \in  \mathbb{P} _ \mathbf{m} [M_ \mathbf{m} ]  $,  then $ p, \varphi $ are compatible in $ \mathbb{P} _ \mathbf{m} [L_ \mathbf{m} ] $,
            
            \item[(b)] $ \psi _0 \le \psi $ in  $\mathbb{P} _ \mathbf{m}[M_ \mathbf{m}]$,
        
            \item[(c)] $ \psi \le p $ in $ \mathbb{P} _ \mathbf{m} [L]$. 
        \end{enumerate} 
    \end{enumerate} 
    
    Let $L_0 = \cup\{t/E_{\bfm}:t \in \fsupp(p)\}
    \cup M_ \mathbf{m},$ so   $ (L_0 \setminus M_ \mathbf{m} ) /E''_\mathbf{m} $  has cardinality $ < \lambda_{0}$ and as $\bfm$ is  $ \mu $-wide, we can find $L_\varp,(\varp \in [1, \mu   ))$ as required, that is, choose an automorphism
    $\pi_\varp$ of $\bfm$ for $\varp <  \mu $ 
    such that $\pi_\varp \rest M_{\bfm}$  is the identity, $\langle \pi_\varp(L_0 ) \backslash M_{\bfm}:\varp < \mu  \rangle$ are pairwise disjoint  where we let  $\pi_0$  be  the identity 
    and so  $L_\varp = \pi_\varp(L)$, and let   $ p_\varp = \hat\pi(p)$ for $ \varepsilon <  \mu $.  
    Note:
    
    \begin{enumerate}
        \item[$(*)_2 $]  if $\varphi_1 \in \bbP_{\bfm}[L_{\bfm}]$,  and $\bbP_{\bfm}[L] \models ``\psi  \le \varphi_1"$ then for for all but $ < \lambda $ ordinals $\varp < \mu$,  the conditions  $p_\varp,\varphi_1$ are compatible.
    \end{enumerate}
    
    [Why?  Let $q \in \bbP_{\bfm} (L_{\bfm})$  be above $\varphi_1$ in $ \mathbb{P} _ \mathbf{m} [L_ \mathbf{m} ] $,  so the set 
    $\{t/E_{\bfm}:t \in \fsupp(q)\}$ has cardinality $< \lambda_{0}$.
    
    So for every  $\varp < \mu $ except $ < \lambda_{0}$ many,  the sets $ \wsupp(q) = \cup\{t/E_{\bfm}:t  \in \fsupp(q)\}$ and $L_\varp \backslash M_{\bfm}$ are disjoint.  Now for  every such $\varp$, the triple $(p_\varp,q,\psi)$ satisfies the assumptions of part (1), hence $p_\varp,q$ are compatible hence $p_\varp,\varphi_1$ are compatible,  so $ (*)_2 $ holds indeed].
    
    Now clearly $\langle (p_\varp,L_\varp):
    \varp <  \mu  )$ satisfies clauses (a)-(f)
    of part  (2), so we are left with clause (g), that is:
    
    \begin{itemize}
        \item if $ u \in [\mu]^{\lambda}$  then $\psi,\bigvee\limits_{\varp  \in u } p_\varp$ are equivalent in  $\bbP_{\bfm}[L_{\bfm}]$, i.e. $\psi \le \bigvee\limits_{\varp  \in u} p_\varp \le \psi$.
    \end{itemize}
    
    Why this holds?  First by the choice of $ \psi $, 
    that is by $(*)_1 $ clearly   $ p \Vdash _ { \mathbb{P} _ \mathbf{m} [L_ \mathbf{m} ]}   `` \psi \in \name{ \mathbf{G} }"$   hence for $ \varepsilon < \mu $ by the choice of $p_ \varepsilon $ also $ p_ \varepsilon \Vdash _{\mathbb{P} _ \mathbf{m} [L_ \mathbf{m} ]} ``  \psi \in \name{ \mathbf{G} }"$  hence $ \psi \le p_ \varepsilon $ in  $ \mathbb{P} _ \mathbf{m} [L_ \mathbf{m} ]  $ hence  $ \psi \le \vee_{ \varepsilon \in u } p_ \varepsilon$ $ \mathbb{P}_ \mathbf{m} [L_ \mathbf{m}]$.
     
     Second, for the other inequality, just note that:
     
     $(*)_3$ if $ q \in \mathbb{P} _ \mathbf{m} [L_ \mathbf{m} ]$ and $ \mathbb{P} _ \mathbf{m} [L_ \mathbf{m} ] \models  ``  \psi \le q " $ then $ q $ is compatible with  $ p_ \varepsilon  $ for every $  \varepsilon < \mu $ except $ < \lambda$ many.
    
    [Why does $(*)_3 $  holds? as in the proof of $ (*)_2 $.]
    
    3) We use part (2) on $\mathbf{n} =
    \mathbf{m} \restriction L$;  so find $ \psi \in \mathbb{P} _ \mathbf{n} [L_ \mathbf{n} ]$
    above $ \psi _0 $ satisfying clauses (a)-(g), but $ \mathbb{P} _ \mathbf{n} [L_ \mathbf{n} ] = \mathbb{P} _ \mathbf{m} [L_ \mathbf{n} ]= \mathbb{P}_ \mathbf{m} [L]$, and so   clause (h) is obvious  and clause (i) holds by the definition of  $\mathbb{P} _ \mathbf{m} [L_ \mathbf{m} ].$
\end{PROOF}

\begin{claim}\label{e32}
    The set $\{\psi_i:i < i(*)\} \cup \{\psi_*\}$ has a common upper bound in $\bbP_{\mathbf m}[L_{\mathbf m}]$ \when:
    
    \begin{enumerate}
        \item[$(*)$ (a)]  $\mathbf m \in \mathbf M$ is $  \mu $-wide  and $ \mu \ge \lambda_{0}$,
        
        \item[(b)]  $i(*) < \lambda $ or just $i_{\ast} < \lambda_{0},$ 
        
        \item[(c)]  $L_{i}  \subseteq L_{\bfm} $  for $i < i(*)$, 
        
        \item[(d)] $L_ i \cap L_j  =   M_ \mathbf{m} $    for $i \not= j<i(*)$,

        \item[(e)]  $ \psi_* \in \bbP_{\mathbf m}[M_{\mathbf m}]$, 
        
        \item[(f)] $t \in L _i    \Rightarrow (t  /E_ \mathbf{m}) \subseteq L_i $, 
    
        \item[(g)] $\psi_i \in \bbP_{\mathbf m}[L_i]$, 
        
        \item[(h)] if $\bbP_{\mathbf m}[M_{\mathbf m}] \models ``\psi_* \le \varphi"$ and $i<i(*)$ then $\psi_i,\varphi$ are compatible in $\bbP_{\mathbf m}[L_{\mathbf m}],$ equivalently in  $\bbP_{\mathbf m}[L_i ] $. 
    \end{enumerate}
\end{claim}

\begin{PROOF}{\ref{e32}} 
    We can  for $ i < i(*) $ replace $ L_i $ by $ L'_i$ when $M_ \mathbf{m} \subseteq L'_i \subseteq L_i $  and the parallel of clauses (f), (g)  of $ (*) $ hold. Hence \wilog:
    
    \begin{enumerate} 
        \item[$(*)_1$] the set $ \{t/E''_\mathbf{m}: t \in L_i \setminus  M_ \mathbf{m}\}$ has cardinality $ < \lambda_{0}$. 
    \end{enumerate} 
    
    As $\psi_* \in \bbP_{\mathbf m}[M_{\mathbf m}]$, there is $p \in \bbP_{\mathbf m}$ such that $p \Vdash_{\bbP_{\mathbf m}} ``\psi_*[\name{\mathbf G} _{\bbP_{\mathbf m}}]=$ true".  As $\mathbf m$ is $ \mu $-wide,  by \ref{e10} there is an automorphism  $f$ of $\mathbf m$  over $ M_ \mathbf{m} $ such that $i<i(*) \Rightarrow f''(\wsupp(p)) \cap L_i   \subseteq M_{\mathbf m}$, hence \wilog \, $i<i(*) \Rightarrow \wsupp(p) \cap L_i   \subseteq M_{\mathbf m}$.  Now we choose $p_i$ by induction  on $i \le i(*)$ such that:
    
    \mn
    \begin{enumerate}
        \item[$\boxplus$]
        \begin{enumerate}
            \item[(a)]  $p_i \in \bbP_{\mathbf m},$

            \item[(b)]  $\langle p_j:j \le i\rangle$ is increasing,

            \item[(c)]  if $s \in \dom(p_i),i<i(*)$ \then \, $\ell g(\tr(p_{i+1}(s)) > i(*),$

            \item[(d)]  $p_0 = p,$

            \item[(e)]   if $i=j+1$ then $p_i \Vdash
            ``\psi_j[\name{\mathbf G}_{\bbP_{\mathbf m}}]=$ true",

            \item[(f)]  $\wsupp(p_i)$ hence also $\fsupp(p_i)$ is  disjoint to $\cup\{L_j \backslash M_{\mathbf m}:j \in [i,i(*))\}$.
        \end{enumerate}
    \end{enumerate}
    
    This is sufficient for the claim as $p_{i(*)}$ is as required.  So let us carry the induction.  For $i=0$ use clause (d), for $i$ limit by \ref{c11}(3A) we know that $\langle p_j:j<i\rangle$ has  a $\le_{\bbP_{\mathbf m}}$-upper bound $p_i$ with domain  $\cup\{\dom(p_j):j<i\}$ satisfying  
    $\wsupp(p_i) \subseteq \cup\{\wsupp(p_j):j<i\}$ by \ref{c11}(3A), hence $p_i$ is as required, in particular as in clause (f).
    
    Recall $p_{j}$ is above $p_{0} = p$ hence above $\psi_{\ast}$ (in $\bbP_{\bfm}[L_{\bfm}]$). As in the proof of $(*)_{3}$ inside \ref{e29}(1) (or see \ref{e52}(1) below)  there is $\varphi_{j} \in \bbP_{\bfm}[M_{\bfm}]$ such that: 
    
    \begin{itemize}
        \item[$\bullet_{1}$] $\psi{*} \leq \varphi_{j},$
        
        \item[$\bullet_{2}$] if $\varphi_{j} \leq \varphi  \in \bbP_{\bfm}[M_{\bfm}]$ then $p_{j}, \varphi$ are compatible. 
    \end{itemize}
    
    Lastly, assume $i=j+1$, by $(*)$(h)  there is $ q_j \in \mathbb{P} _\mathbf{m} $ above $\varphi _j \wedge \psi_{j} $. Because $ \mathbf{m} $ is $ \mu $-wide there is an automorphism $ \pi $ of $ \mathbf{m} $ over $ M_ \mathbf{m}$  satisfying $\pi \rest L_j $ is the identity, so  $ \pi'' (\dom(q_j) \setminus M_ \mathbf{m} $   is disjoint to $ \wsupp( p_i )$ and to   $ L_ \varepsilon $ for $ \varepsilon \in i_* \setminus \{ j\}.$ So \wilog: 
    
    \begin{enumerate} 
        \item[$(*)_2$] $ q_j $ itself satisfies this.
    \end{enumerate} 
    
    Now  the statement $(\stt )_{p_j, q_j, \varphi_{j}}$ holds.
    
    [Why? because  $ \wsupp(p_j) \cap \wsupp(q_j) \subseteq M_ \mathbf{m} $  by $(*)_2$,  the choice of $\varphi_{j}$ and $q_{j}$ above.]
     
    Hence by \ref{e29} $p_j,q_j$ has a common upper bound called $p_i$.   As $\mathbf m$ is wide,  for some automorphism $ \pi $  of $ \mathbf{m}  $ over $ M_\mathbf{m} $  such that  $ \pi \rest \wsupp(p_j)$ is the identity and $\pi"\wsupp(p_j)$ is disjoint to $ \cup \{ L_ \varepsilon: \varepsilon \in  [i, i_* )\} $, hence by renaming \wilog: 
            
    \begin{enumerate} 
        \item[$(*)_3$] $ \wsupp(p_i) \setminus M_ \mathbf{m} $   is disjoint to  $ \cup \{ L_ \varepsilon: \varepsilon \in  [i, i_* )\},$
    \end{enumerate} 

    Clearly $p_i$ is as required so we have finished  proving $ (*)_3$. 
    
    So we have finished proving the last case in the 
    the induction.
    
    So we are done.
\end{PROOF}

\subsection{The main result}\label{3D}\  

Here we  continue \S3A, \S3B, and in particular prove 
the main result, it does not rely on \ref{3C}.  
Concerning \S1B, we rely on it only  in one point: quoting \ref{c33s} while   proving $ \boxplus _{4.4}$  and the beginning of Case 3   inside the proof of \ref{e35},  this can be avoided using \S4A. We have not work out if e.g. \S3D works for the fat context.

\begin{hypothesis}\label{e33z}
    We are in the lean context (for this subsection).
\end{hypothesis}

\begin{conclusion}\label{e33}
    If $\beta \ge 0$ and $\mathbf m$ is wide and 
    $f \in \cG_{\mathbf m,\beta}$ and $L_1,L_2$ its domain and range respectively  \then \, $f$ induces an isomorphism $\hat f$ from $\bbP_{\mathbf m}(L_1)$ onto $\bbP_{\mathbf m}(L_2)$.
\end{conclusion}

\begin{remark}\label{e34}
    1) See Definition \ref{e4}(3); note that this claim is not covered by Definition \ref{e4}(2).
    
    2) Here we use \ref{e5n}(4), so the choice in Definition \ref{c6}(c)($\gamma$) is justified  (see Remark \ref{e5p}(1) used below in the proof).
    
    3) We could have separated the definition of ``analyze" and its properties.
    
    4) Note that in Definition \ref{e24}, we deal only with $L_1 \subseteq t/E_{\mathbf m}$ for some $t$.
     
    5) How come even $\beta=0$ is suitable for \ref{e33}?  The point is clause $(a)(\varp)$ of Definition \ref{e24}(2).  But  there is no real harm using larger $\beta$.
\end{remark}

\begin{PROOF}{\ref{e33}} 
    By the definitions, clearly $\hat f$ is a one-to-one function  from $\bbP_{\mathbf m
    }(L_1)$ onto $\bbP_{\mathbf m}(L_2)$.  
    Next assume $p_1,q_1 \in \bbP_{\mathbf m}(L_1),\dom(p_1) \subseteq \dom(q_1)$ and let $p_2 := \hat f(p_1),q_2 := \hat f(q_1)$;
    clearly they belong to $\bbP_{\mathbf m}(L_2)$.  
    We shall prove that $\bbP_{\mathbf m}
    \models ``p_1 \le q_1"$ iff $\bbP_{\mathbf m} \models ``p_2 \le q_2"$.
    
    Let  $ i(*) <   \lambda $ and  $\bar{ t } _1 = \langle t^1_i:i<i(*)\rangle$ be such that:
    
    \begin{enumerate}
        \item[$\oplus_1$]
        
        \begin{enumerate}
            \item[(a)]  $t^1_i \in \fsupp(q_1) \backslash M_{\mathbf m} \subseteq L_1$ such that $\fsupp(q_1)$ is included in $\cup\{t^1_i/E_{\mathbf m}:i < i(*)\}$,

            \item[(b)]  $\langle t^1_i:i < i(*)\rangle$ are pairwise non-$E''_{\mathbf m}$-equivalent. 
        \end{enumerate}
    \end{enumerate}

    Next let,

    \begin{enumerate}
        \item[$\oplus_2$]
            \begin{enumerate}
                \item[(c)] let  $t^2_i = f(t^1_i)$   for $ i < i(*)$  and  let $\bar t_2 = \langle t^2_i:i <i(*)\rangle$,

                \item[(d)] $\fsupp(p_\ell) \subseteq  \cup\{t^\ell_i/E''_{\mathbf m}: i < j(*)\} \cup M_{\mathbf m}$, so $j(*) \le i(*)$,  for $ {\ell} = 1,2$. 
            \end{enumerate}
    \end{enumerate}
    
    For $i<i(*)$ let $\psi^*_{1,i} \in \bbP_{\mathbf m}[M_{\mathbf m}]$ be such that: $\vartheta \in \bbP_{\mathbf m}[M_{\mathbf m}]$ is compatible with $q_{1,i} := q_1 \upharpoonleft (t^1_i/E_{\mathbf m})$ (the projection!) \Iff
    \, $\vartheta \wedge \psi^*_{1,i} \in \bbP_{\mathbf m}[M_{\mathbf m}]$;
    clearly exists as $\bbP_{\mathbf m}$ satisfies the $\lambda^+$-c.c.  Clearly $ \mathbb{P} _ \mathbf{m} [L_ \mathbf{m} ]  \models  \lqq \psi ^*_{1, i}  \le q_{1.i} \le q_1 $  for $ i < i(*) $ and let $\psi^*_1 = \wedge\{\psi^*_{1,i}:i <i(*)\}$.
    
    Now $\psi^*_1 \in \bbP_{\mathbf m}[M_{\mathbf m}]$ as $q_1 \Vdash ``\psi^*_1[\name{\mathbf G}_{\bbP_{\mathbf m}}] =$ true".  We will say
    ``$\psi^*_1,\bar\psi^*_1 = \langle \psi^*_{1,i},q_{1,i}: i<i(*)\rangle$ analyze $q_1$ or $(q_1,\bar t_1)$" \when \, the above holds.
    
    Next choose $\varphi^*_1,\langle \varphi^*_{1,i},p_{1,i}: i <j(*)\rangle$
    which analyze $p_1,\langle t^1_i:i<j(*)\rangle$   where \wilog \, $ j(*) \le i(*)$. Why possible?  As above  recalling $ p_1 \le q_1 \Rightarrow  \fsupp(p_1 ) \subseteq \fsupp(q_1 ) $. 
    
    Lastly, let  $\psi^*_{2,i} = \check f(\psi^*_{1,i}), p_{2,i} =  \hat f(p_{1,i}),   \psi^*_2 = \check f(\psi^*_1), \varphi^*_{2,i} = \check f(\varphi^*_{1,i}), q_{2,i} = \hat f(q_{1,i}),   \varphi^*_2 = \check f(\varphi^*_1)$   where $\check f$ is the function from $\bbL_{\lambda_{0}}(Y_{L_1},\bbP_{\mathbf m})$ onto $\bbL_{\lambda_{0}}(Y_{L_2},\bbP_{\mathbf m})$ induced by $f$, i.e. where $\check f$ is the one-to-one function with domain $\bbL_{\lambda^+}[Y_{L_1}]$ defined by $p^*_{t,\eta} \mapsto p^*_{f(t),\eta}$. Now,
    
    \begin{enumerate}
        \item[$(*)$]  for $\ell=1,2$ the sequence $(p_\ell,q_\ell,\bar\psi^*_\ell, \bar\psi^*_\ell,\varphi^*_\ell,\bar\varphi^*_\ell)$ where $\bar\psi^*_\ell  = \langle \psi^*_{\ell,i},q_{\ell,i}:i < i_\ell(*)\rangle,\bar\varphi^*_\ell = \langle
        \varphi^*_{\ell,i},p_{\ell,i}:i<i(*)\rangle$  satisfy the same demands as listed above for $\ell=1,2$, that is
        
        \begin{enumerate}
            \item[(a)] $(\psi^*_\ell,\bar\psi^*_\ell)$ analyze $(q_\ell,\bar t_\ell)$ for $\ell=1,2$
            
            \item[(b)] $(\varphi^*_\ell,\bar\varphi^*_\ell)$ analyze $(p_\ell, \bar t_\ell \rest j(*))$ for $\ell=1,2$.
        \end{enumerate}
    \end{enumerate}
    
    [Why?  Think, recalling  $f \rest (t^1_i/E_{\mathbf m})$ is an isomorphism from $\mathbf m \rest ((t^1_i/E_{\mathbf m}) \cap L_1)$ onto $\mathbf m \rest ((t^2_i/E_{\mathbf m}) \cap L_2)$,  inducing an isomorphism between  $
    \mathbb{P} _ \mathbf{m} [ (t^1_i /E_ \mathbf{m} ) \cap L_1$  and  $ \mathbb{P} _\mathbf{m} [(t^2_i /E_ \mathbf{m} ) \cap L_2] $  by \ref{e24}(a)($\delta  $)  and $ \psi ^*_ 2 = \wedge \{\psi ^*_{2,i}: i < i(*) \} $ is because each  function $ f \rest (t^1_i /E_ \mathbf{m} ) $ induces the identity mapping on $ \mathbb{P} _ \mathbf{m} [M_ \mathbf{m} ]$.]
    
    Next,
    
    \begin{enumerate}
        \item[$\boxplus$]  for $\ell=1,2$ we have $(A)_\ell \Leftrightarrow (B)_\ell$ where:  
        
        \begin{enumerate}
            \item[$(A)_\ell$]  $\bbP_{\mathbf m} \models ``p_\ell \le q_\ell"$,  
            
            \item[$(B)_\ell$]  for every $i<j(*)$ we have $\bbP_{\mathbf m}[t^\ell_i/E_{\mathbf m}] \models ``(\varphi^*_\ell \wedge p_{\ell,i}) \le (\psi^*_\ell \wedge q_{\ell,i})"$.
        \end{enumerate}
    \end{enumerate}
    
    Why?  First, assume that the condition $(B)_\ell$ fails, say for $i$, hence there is $\vartheta \in \bbP_{\mathbf m}[t^\ell_i/E_{\mathbf m}]$ such that $\bbP_{\mathbf m}[t^\ell_i/E_{\mathbf m}] \models ``(\psi^*_\ell \wedge q_{\ell,i}) \le \vartheta"$,  and  $\varphi^*_\ell \wedge p_{\ell,i} \wedge \vartheta \notin \bbP_{\mathbf m}[t^\ell_i/E_ \mathbf{m} ]$.  So by claim \ref{e32} there is $q^+_\ell \in \bbP_{\mathbf m}$ such that $q^+_\ell \in \bbP_{\mathbf m}[L_{\mathbf m}]$ is above $\vartheta$,   hence above  $\psi^*_\ell$ and above $q_{\ell,j} =
    q_\ell \upharpoonleft (t^\ell_j/E_{\mathbf m})$ for $j<i(*)$. That is, first get $\psi \in \bbP_{\bfm}[M_{\bfm}]$ such that $\psi \ge
    \psi^*_\ell$ and [$\psi \le \psi' \in \bbP_ \mathbf{m}  [M_\mathbf{m} ] \Rightarrow   \psi',\vartheta$ are compatible] (using $\vartheta \ge \psi^*_\ell$).    Then apply \ref{e32} to 
    $(\{q_{\ell,j}:j <  i(*), j \not= i \} \cup \{\vartheta\}) \cup  \{\psi\}$ to get $q^+_\ell$. We have used   $ i (*) < \lambda $.  
    
    Hence by \ref{e5n}(4) the condition $q^+_\ell$ is above $q_\ell$,   but $q^+_\ell
    \Vdash ``\varphi^*_\ell \wedge p_{\ell,i}[\name{\mathbf G}] =$ false" as
    $q^+_\ell$ is above $\vartheta$.  However, $p_\ell \Vdash_{\bbP_{\mathbf m}[L_{\mathbf m}]} ``p_{\ell,i} \in \name{\mathbf G}$ and
    $\varphi^*_\ell \in \mathbf G"$.  By the last two sentences $q^+_\ell,p_\ell$ are incompatible   
    in $\bbP_{\mathbf m}[L_{\mathbf m}]$ equivalently in $\bbP_{\mathbf m}$.  So indeed $\neg(B)_\ell \Rightarrow \neg(A)_\ell$.
    
    For the other direction assume condition $(B)_\ell$ holds, but condition $(A)_\ell$ fails and we shall get a contradiction.  So there
    is $q^+_\ell \in \bbP_{\mathbf m}$ above $q_\ell$ incompatible with
    $p_\ell$.
    
    For each $i<i(*)$ as $(\psi^*_\ell,\langle \psi^*_{\ell,j},q_{\ell,j}: j<i(*)\rangle)$ analyze  $q_\ell$, clearly $\bbP_{\mathbf m}[L_{\mathbf m}] \models ``(\psi^*_\ell
    \wedge q_{\ell,i}) \le q_\ell"$ but $q_\ell \le q^+_\ell$ hence  $\bbP_{\mathbf m}[L_{\mathbf m}] \models ``(\psi^*_\ell \wedge q_{\ell,i})
    \le q^+_\ell"$, and as we are assuming clause $(B)_\ell$ we have $\bbP_{\mathbf m}[L_{\mathbf m}] \models ``(\varphi^*_\ell \wedge p_{\ell,i})
    \le q^+_\ell"$.  Hence by \ref{e5n}(4), $q^+_\ell$ is above $p_\ell$ in $ \mathbb{P}_\mathbf{m} [L_ \mathbf{m} ]$  hence they are compatible in $ \mathbb{P}_\mathbf{m}$,  contradiction.  
    So indeed $(B)_\ell \Rightarrow (A)_\ell$. Together, $\boxplus$ holds].
    
    Now clearly $(B)_1 \Leftrightarrow (B)_2$, see Definition \ref{e24}, \ref{e28}; so by $\boxplus$ we have $(A)_1 \Leftrightarrow (A)_2$ which is the desired conclusion.
\end{PROOF}

\begin{claim}\label{e35}
    We have $\bbP_{\mathbf m_1} \lessdot \bbP_{\mathbf m}$ \when:
    
    \begin{enumerate}
        \item[(a)]  $\mathbf m_1 \le_{\mathbf M} \mathbf m$, 
        
        \item[(b)]  if $t \in L_{\mathbf m} \backslash M_{\mathbf m_1}$ and $\bar s \in {}^\zeta(t/E''_{\mathbf m}),\zeta < \lambda^+$ \then \, we can find $t_i,\bar s_i$ for $i < \lambda^+$ such that:
        
        \begin{enumerate}
            \item[$(\alpha)$]   $t_i \in L_{\mathbf m_1} \backslash M_{\mathbf m_1}$, 
            
            \item[$(\beta)$]  $t_i/E''_{\mathbf m_1}   \ne t_j/E''_{\mathbf m_1}$ when $i \ne j < \lambda^+$, 
            
            \item[$(\gamma)$]  $\bar s_i \in {}^\zeta(t_i/E''_{\mathbf m_1})$, 
            
            \item[$(\delta)$]  $(t_i,\bar s_i)$ is $\xi$-equivalent to $(t,\bar s)$ in $\mathbf m$ where\footnote{no real harm in using larger $\xi$.} $\xi = 1$.
        \end{enumerate} 
        
        \item[(c)]  $\mathbf m$ is wide.
    \end{enumerate}
\end{claim}

\begin{remark}\label{e36}  
    In the proof we use conclusion \ref{e33} but not clause $(a)(\varepsilon)$ of Definition \ref{e24}(2).
\end{remark}

\begin{PROOF}{\ref{e35}}
    \begin{enumerate}
        \item[$\boxplus_1$]   for $ \beta \ge 0 $ and $f \in \cG_{\mathbf m,\beta}$, 
        
        \begin{enumerate}
        \item[(a)] $\hat f$ preserves ``$p_2$ is above
          $p_1$ in $\bbP_{\mathbf m}"$, and its negations,  
        
        \item[(b)] if $\beta > 0$ then $\hat f$ preserves also incompatibility in $\bbP_{\mathbf m}$.
        \end{enumerate}
    \end{enumerate}

    [Why? Clause (a) holds by \ref{e33}.  For clause (b) use  clause (a) and Definitions \ref{e24} and \ref{e28} or see the proof of $\boxplus_2$.]

    \begin{enumerate}
        \item[$\boxplus_2$]  if $p_i \in \bbP_{\mathbf m_1}$ for $i < i(*) < \lambda^+$ and $p \in \bbP_{\mathbf m}$ \then \, there is $p^*$ such that:
        
        \begin{enumerate}
            \item[(a)]  $p^* \in \bbP_{\mathbf m_1}$, equivalently $p^* \in \bbP_{\mathbf m}(L_{\mathbf m_1})$,
            
            \item[(b)]  $\bbP_{\mathbf m_1} \models ``p_i \le p^*"$ iff $\bbP_{\mathbf m} \models ``p_i \le p"$,  
            
            \item[(c)]  $\bbP_{\mathbf m_1} \models ``p_i,p^*$ are compatible" iff $\bbP_{\mathbf m} \models ``p_i,p$ are compatible".
        \end{enumerate}
    \end{enumerate}

    [Why?  Let $q_i \in \bbP_{\mathbf m}$ be such that: if $p_i,p$ are compatible in $\bbP_{\mathbf m}$ then $p_i \le q_i \wedge p \le q_i$. We can find $L_1 \subseteq L_2$ such that
   
    \begin{enumerate}
        \item[$\bullet$]  $M_{\mathbf m} \subseteq L_1 \subseteq L_{\mathbf m_1}, |L_1 \backslash M_{\mathbf m}| \le \lambda$, 
        
        \item[$\bullet$]  $\{p_i:i < i(*)\} \subseteq \bbP_{\mathbf m}(L_1)$, 
        
        \item[$\bullet$]  $L_1 \subseteq L_2 \subseteq L_{\mathbf m},|L_2 \backslash M_{\mathbf m}| \le \lambda$ and $p,q_i \in \bbP_{\mathbf m}(L_2)$ for $i<i(*)$.
    \end{enumerate}
    
    By the assumption of the claim there is $f \in \cG _{\mathbf m,1}$ such that:
    
    \begin{enumerate}
        \item[$\bullet$]  $\dom(f) \subseteq \cup\{(t/E''_{\mathbf m}) \cap L_2:
        t \in L_2\} \cup M_{\mathbf m}$,  
      
        \item[$\bullet$]  $t \in L_1 \Rightarrow f \rest (t/E_{\mathbf m} \cap L_2) = \id_{(t/E_{\mathbf m}) \cap L_2}$, 
    
        \item[$\bullet$]  if $q \in \{q_i:i < i(*)\} \cup \{p\} \cup \{p_i:i < i(*)\}$ and $t \in \dom(q) \backslash M_{\mathbf m}$ then $\fsupp(q(t)) \subseteq \dom(f)$, 
          
        \item[$\bullet$]  $\rang(f) \subseteq L_{\mathbf m_1}$.
    \end{enumerate}
    
    Let $p^* = \hat f(p)$: by $\boxplus_1(a)$ clearly clauses (a),(b) of $\boxplus_2$ holds; and the choice of the $q_i$'s   (and as $ p \le q_2 \Rightarrow  \hat{f}(p ) \le \hat{f}(q_i)$)   
    also the implication ``if" of clause (c).  The
    ``only if" of clause (c) holds by $\boxplus_1(b)$ so we are done.]

    \begin{enumerate}
        \item[$\boxplus_3$]  if $p \in \bbP_{\mathbf m}$ then $p \in \bbP_{\mathbf m_1}$ \Iff \, $\fsupp(p) \subseteq L_{\mathbf m_1}$.
    \end{enumerate}

    [Why?  Obvious.]
    
    Recalling Definition \ref{c34}(0)(c): 
    
    \begin{enumerate}
        \item[$\boxplus_4$]  for every ordinal $\gamma$, we have  $\bbP_{\mathbf m_1}(L^{\deq}_{\mathbf m_1,\gamma}) \lessdot 
        \bbP_{\mathbf m}(L^{\deq}_{\mathbf m,\gamma})$.
    \end{enumerate}

    [Why?  We shall prove this by induction on $\gamma$ using  $\boxplus_2 + \boxplus_3$.
    
    Note that:  
    
    \begin{enumerate}
        \item[$\boxplus_{4.1}$] 
        \begin{enumerate}
            \item[(a)]  $L^{\deq}_{\mathbf m,\gamma} \cap L_{\mathbf m_1} = L^{\deq}_{\mathbf m_1,\gamma}$, 
            
            \item[(b)]   if $f \in \cG   _{\mathbf m,\beta}, s \in \dom(f)$ and $\beta$ is an ordinal then:
            
            \begin{itemize}
                \item   $s \in L^{\deq}_{\mathbf m_1,\gamma}
                  \Leftrightarrow f(s) \in L^{\deq}_{\mathbf m,\gamma}$, 
            \end{itemize}
            
            \item[(c)]   the parallel of $\boxplus_2$ holds  replacing  the pair $ (\mathbb{P} _{\mathbf{m} _1},\mathbb{P} _ \mathbf{m}  )$ by the pair  $ (\mathbb{P} _{\mathbf{m} _1}(L^\deq _{\mathbf{m} _1, \gamma }), \mathbb{P} _ \mathbf{m}(L^\deq _{\mathbf{m} , \gamma }))$;  so  e.g.    
            $p^* \in  \bbP_{\mathbf m}(L^\deq _{\mathbf m_1,\gamma})$, 
            
            \item[(d)]  $L^{\deq}_{\mathbf m,\gamma}$ is an initial segment of $L_{\mathbf m}$, 
            
            \item[(e)]  $L^{\deq}_{\mathbf m_1,\gamma}$ is an initial segment of $L_{\mathbf m_1}$, 
            
            \item[(f)]  $\bbP_{\mathbf m_1}(L^{\deq}_{\mathbf m_1,\gamma}) \lessdot \bbP_{\mathbf m_1}(L_{\mathbf m_1})$, similarly for $\mathbf m$.
        \end{enumerate}
    \end{enumerate}

    We shall use this freely.  The inductive proof on $\gamma$ splits to
    three cases.

    \underline{Case 1}:  $\gamma =0.$

    So, 
    
    \begin{itemize}
        \item  $E = E''_{\mathbf m} \rest L^{\deq}_{\mathbf m,\gamma}$ is an
        equivalence relation on $L^{\deq}_{\mathbf m,\gamma}$, 
        
        \item   $E \rest L^{\deq}_{\mathbf m_1,\gamma} = E''_{\mathbf m_1} \rest  L^{\deq}_{\mathbf m_1,\gamma}$, 
        
        \item   if $t \in L^{\deq}_{\mathbf m_1,\gamma}$ then $t \notin M_{\mathbf m_1},t/E'_{\mathbf m_1} = t/E'_{\mathbf m},(t/E'_{\mathbf m_1}) \cap L^{\deq}_{\mathbf m_1,\gamma} = (t/E_{\mathbf m_1}) \cap
        L^{\deq}_{\mathbf m_1,\gamma} = (t/E'_{\mathbf m}) \cap L^{\deq}_{\mathbf m,\gamma}$ initial segment of $L_{\mathbf m_1}$ and of $L_{\mathbf m}$ and $\bbP_{\mathbf m}((t/E_{\mathbf m_1}) \cap L^{\deq}_{\mathbf m_1,\gamma}) = \bbP_{\mathbf m_1}((t/E_{\mathbf m_1}) \cap L^{\deq}_{\mathbf m_1,\gamma})$, 
        
        \item $\bbP_{\mathbf m}(L^{\deq}_{\mathbf m,\gamma})$ is the product with $(< \lambda)$-support of $\{\bbP_{\mathbf m}
        ((t/E_{\mathbf m_1}) \cap L^{\deq}_{\mathbf m_1,\gamma}):t \in L^{\deq}_{\mathbf m,\gamma}\}$, 
        
        \item similarly for $\mathbf m_1$.
    \end{itemize}
    
    So the result should be clear.

    \underline{Case 2}:  $\gamma = \beta +1$
    
    Let $M_\beta = \{s \in M_{\mathbf m}:\dep  _{\mathbf m}(s)=\beta\}$, clearly:
    
    \begin{enumerate}
        \item[$\boxplus_{4.2}$]
        
        \begin{enumerate}
            \item[(a)]  $M_\beta$ is a set of pairwise incomparable elements,
            
            \item[(b)] $s \in M_\beta \Rightarrow L_{\mathbf m_1,<s} \subseteq L^{\deq}_{\mathbf m_1,\beta} \wedge L_{\mathbf m,<s} \subseteq L^{\deq}_{\mathbf m_1,\beta}$, 
            
            \item[(c)]  $M_\beta$ is disjoint to $L^{\deq}_{\mathbf m_1,\beta},L^{\deq}_{\mathbf m,\beta}$, 
            
            \item[(d)]  $M_\beta \subseteq L^{\deq}_{\mathbf m_1,\gamma}$, 
            
            \item[(e)]  $L^{\deq}_{\mathbf m,\beta} \cup M_\beta$ is an initial segment of $L_{\mathbf m}$, 
            
            \item[(f)]  $L^{\deq}_{\mathbf m_1,\beta} \cup M_\beta$ is an initial segment of $L_{\mathbf m_1}$.
        \end{enumerate}
    \end{enumerate}

    As first half we prove:
    
    \begin{enumerate}
        \item[$\boxplus_{4.3}$]  $\bbP_{\mathbf m_1}(L^{\deq}_{\mathbf m_1,\beta} \cup M_\beta) \lessdot \bbP_{\mathbf m}(L^{\deq}_{\mathbf m,\beta} \cup M_\beta)$.
    \end{enumerate}
    
    Why?  Recalling $\boxplus_{4.1}(a)$, note
    
    \begin{enumerate}
        \item[$(a)^\pigyon$]   for $p,q \in \bbP_{\mathbf m_1}(L^{\deq}_{\mathbf m_1,\beta} \cup M_\beta)$ we have $\bbP_{\mathbf m_1}(L^{\deq}_{\mathbf m_1,\beta} \cup M_\beta) \models ``p \le q"$ iff 
        $\bbP_{\mathbf m}(L^{\deq}_{\mathbf m , \gamma} \cup M_\beta) \models ``p \le q"$.
    \end{enumerate}
    
    [Why?  Immediate by the definition of the order and the induction hypothesis.]

    \begin{enumerate}
        \item[$(b)^\pigyon$]  if $p_1,p_2 \in \bbP_{\mathbf m_1}(L^{\deq}_{\mathbf m_1,\beta} \cup M_\beta)$ \then \, $p_1,p_2$ are compatible in  $\bbP_{\mathbf m_1}(L^{\deq}_{\mathbf m_1,\beta} \cup M_\beta)$ iff they are compatible in $\bbP_{\mathbf m}(L^{\deq}_{\mathbf m,\beta} \cup M_\beta)$. 
    \end{enumerate}
    
    [Why?  The implication $\Rightarrow$ holds by clause $(a)^\pigyon$.  
        
    So assume $p_3 \in \bbP_{\mathbf m}(L^{\deq}_{\mathbf m,\beta} \cup M_\beta)$ is a
    common upper bound of $p_1,p_2$ in $\bbP_{\mathbf m}(L^{\deq}_{\mathbf m,\beta} \cup M_\beta)$ equivalently in $\bbP_{\mathbf m}$. 
    
    Now (by clause (b) of the claim assumption) there is $f \in \cG_{\mathbf m,1}$ (actually $ \cG _{\mathbf{m}, 0 }$ suffices here) such that:  
    
    \begin{enumerate}
        \item[$\bullet$]   $f \rest (\fsupp(p_1) \cup
        \fsupp(p_2))$ is the identity, moreover
        
        \item[$\bullet$]   $s \in \wsupp(p_1) \cup
        \wsupp(p_2) \wedge s \in \dom(f) \Rightarrow f(s) = s$,
        
        \item[$\bullet$]  $\dom(f) = \cup\{\fsupp(p_\ell):\ell = 1,2,3\}$
        
        \item[$\bullet$] $\rang(f) \subseteq L_{\mathbf m_1}$.
    \end{enumerate}
    
    Hence clearly $f \rest M_\beta = \id_{M_\beta}$ so by  $\boxplus_{4.1}(b)$ we have $\rang(f) \subseteq L^{\deq}_{\mathbf m_1,\beta} \cup M_\beta$ so $\hat f(p_3) \in  \bbP_{\mathbf m}(L^{\deq}_{\mathbf m_1,\beta} \cup M_\beta)$.
    
    By $\boxplus_1$ the condition $\hat f(p_3)$ is a common  upper bound of $p_1,p_2$ in $\bbP_{\mathbf m}$ and by the previous sentence also in $\bbP_{\mathbf m}(L^{\deq}_{\mathbf m_1,\beta} 
    \cup M_\beta)$, so by clause $(a)^\pigyon$  
    the conclusion of $(b)^\pigyon$  holds.] 
    
    \begin{enumerate} 
    \item[$(c)^\pigyon$]  
        If $\cI$ is a maximal antichain in $\bbP_{\mathbf m_1} (L^{\deq}_{\mathbf m_1,\beta} \cup M_\beta)$ \then \, $\cI$ 
        is a maximal antichain of $\bbP_{\mathbf m}
        (L^{\deq}_{\mathbf m,\beta} \cup M_\beta)$.
    \end{enumerate}
    
    [Why?  As in the proof of clause $ (b)^\pigyon $  and of $\boxplus_2$.]
    
    So we are done proving $\boxplus_{4.3}$.
    
    Now we return to proving $\boxplus_4$,  note
    
    \begin{enumerate}
        \item[$\boxplus_{4.4}$]  let $\cE = \{(s_1,s_2):s_1,s_2 \in L_*$ and
          $s_1/E_{\mathbf m} = s_2/E_{\mathbf m}\}$ where $L_* = L^{\deq}_{\mathbf m,\gamma} \backslash (L^{\deq}_{\mathbf
          m,\beta} \cup M_\beta)\}$,  \then:
        
        \begin{enumerate}
            \item[(a)]  $\cE$ is an equivalence relation on $L_*,$
            
            \item[(b)]  if $s_1,s_2 \in L_*$ and $s_1 \le_{L_{\mathbf m}} s_2$ then $s_1 \cE s_2,$
            
            \item[(c)]  if $s_1,s_2 \in L_*$ and $s_1 \cE s_2$ then  $s_1 \in L^{\deq}_{\mathbf m_1,\gamma} \Leftrightarrow s_2 \in 
            L^{\deq}_{\mathbf m_1,\gamma}$ (and both $\notin M_\beta$),
            
            \item[(d)]  if $s \in L_*$ then $L_{\mathbf m,<s} \subseteq L^{\deq}_{\mathbf m,\beta} \cup M_\beta \cup (s/\cE),$
            
            \item[(e)]  if $s \in L_* \cap L_{\mathbf m_1}$ then  $L_{\mathbf m_1,<s} \subseteq L^{\deq}_{\mathbf m_1,\beta} \cup M_\beta \cup (s/\cE)$.
        \end{enumerate}
    \end{enumerate}

    Hence let $L_0 = L^{\deq}_{\mathbf m_1,\beta} \cup M_\beta$ and $L_1 = L^{\deq}_{\mathbf m_1,\gamma} = L^{\deq}_{\mathbf m_1} \cup M_\beta$ they satisfy  all the assumptions of \ref{c33s} hence its conclusion, so we are done easily proving Case 2 of $\boxplus_4$.
    
    \underline{Case 3}:  $\gamma$ is a limit ordinal
    
    Note that in this case  the set $L^\deq_{\mathbf{m} , \gamma }  \setminus \cup \{ L^\deq_{\mathbf{m}, \beta }: \beta < \gamma \} $  consist of $ s \in  L^{\deq } _{\mathbf{m},\gamma } \setminus M_ \mathbf{m} $ which are not below any elements  from   $ L^\deq_{\mathbf{m}, < \gamma } =  \cup \{L^\deq_{\mathbf{m}, \beta }: \beta < \gamma \} $  hence as in case 2 we can treat them as in the proof of  $ \boxplus _{4.4}$, citing \ref{c33s},  so we shall ignore them below.
    
    Clearly $p \in \bbP_{\mathbf m_1}(L^{\deq}_{\mathbf m_1,<  \gamma})$ iff $p \in
    \bbP_{\mathbf m}(L^{\deq}_{\mathbf m_1,< \gamma})$; also each of them implies $p \in \bbP_{\mathbf m}(L^{\deq}_{\mathbf m,< \gamma})$. Also for $p,q \in \bbP_{\mathbf m_1}(L^{\deq}_{\mathbf m_1,< \gamma})$ we have  
    $\bbP_{\mathbf m_1}(L^{\deq}_{\mathbf m_1,< \gamma}) \models ``p \le q"$ \Iff \,
    $\bbP_{\mathbf m}(L^{\deq}_{\mathbf m,< \gamma}) \models ``p \le q"$ by the definition of the order and the induction hypothesis.  Together
    $\bbP_{\mathbf m_1}(L^{\deq}_{\mathbf m_1,< \gamma}) \subseteq  \bbP_{\mathbf m}(L^{\deq}_{\mathbf m,< \gamma})$, (as partial orders).
    
    Next assume that $q_1,q_2 \in \bbP_{\mathbf m_1}(L^{\deq}_{\mathbf m_1,< \gamma})$ and $p_3$ is a common upper bound of $q_1,q_2$ in  $\bbP_{\mathbf m}(L^{\deq}_{\mathbf m,< \gamma})$.
    
    We shall find $p_1 \in \bbP_{\mathbf m_1}(L^{\deq}_{\mathbf m_1,< \gamma})$ such that:
    
    \begin{enumerate}
        \item[$(*)_1$]  $(\rm{a}) \quad p_1$ is above $q_1,q_2$ (in  $\bbP_{\mathbf m_1}(L^{\deq}_{\mathbf m_1,< \gamma})$ or equivalently in $\bbP_{\mathbf m}(L^{\deq}_{\mathbf m_1,< \gamma}))$,
        
        \item[${{}}$]  $(\rm{b}) \quad$ if $p_1 \le p'_1 \in  \bbP_{\mathbf m}(L^{\deq}_{\mathbf m_1,< \gamma})$ then $p'_1,p_3$ are compatible in $\bbP_{\mathbf m}(L^{\deq}_{\mathbf m,< \gamma})$.
    \end{enumerate}
    
    This clearly suffices; why? e.g. if $\{r_i:i < i(*)\} \subseteq \bbP_{\mathbf m_1}(L^{\deq}_{\mathbf m_1,< \gamma})$ is a maximal antichain of $\bbP_{\mathbf m_1}(L^{\deq}_{\mathbf m_1,< \gamma})$ but not of
    $\bbP_{\mathbf m}(L^{\deq}_{\mathbf m,< \gamma})$, let $q_1 = q_2 = \emptyset$ and $p_3 \in \bbP_{\mathbf m}(L^{\deq}_{\mathbf m,< \gamma})$ be incompatible with every $r_i$; let $p_1$ be as in $(*)_1$, it gives a contradiction.
    
    If $\cf(\gamma) \ge \lambda$ then for some $\gamma_1 < \gamma$ we have $q_1,q_2 \in \bbP_{\mathbf m}(L^{\deq}_{\mathbf m_1,\gamma_1})$ and $\fsupp(p_3) \cap L^{\deq}_{\mathbf m,< \gamma} \subseteq L^{\deq}_{\mathbf m,\gamma_1}$ and use the induction hypothesis on $\gamma_1$ for
    clause (a) of $(*)_1$; for clause (b) of $(*)_1$ we also recall \ref{c11}(6); (alternatively imitate the case $\cf(\gamma) < \lambda$, choosing ``changing our minds" $\gamma_\varp < \gamma$ with the induction).  So assume $\aleph_0 \le \cf(\gamma) < \lambda$ and let $\langle \gamma_\varepsilon:\varepsilon < \cf(\gamma)\rangle$ be increasing continuous with limit $\gamma$.
    
    Now we choose $p_{1,\varepsilon}$ by induction on $\varepsilon \le \cf(\gamma)$ such that:
    
    \begin{enumerate}
        \item[$(*)_2$]
    
        \begin{enumerate}
            \item[(a)]  $p_{1,\varepsilon} \in 
            \bbP_{\mathbf m}(L^{\deq}_{\mathbf m_1,\gamma_\varepsilon})$, 
            
            \item[(b)]  $(\gamma_\varepsilon,q_1 \rest L^{\deq}_{\mathbf m,\gamma_\varepsilon},q_2 \rest L^{\deq}_{\mathbf m,\gamma_\varepsilon},p_3 \rest L^{\deq}_{\mathbf m,\gamma_\varepsilon},p_{1,\varepsilon})$ are like $(\gamma,q_1,q_2,p_3,p_1)$ in $(*)_1$, 
            
            \item[(c)]  $p_{1,\zeta} \le p_{1,\varepsilon}$ for $\zeta < \varepsilon$, 
            
            \item[(d)]  if $\varepsilon = \zeta +1$ and $s \in \dom(p_{1,\zeta})$ then $\ell g(\tr(p_\varepsilon(s)) > \cf(\gamma)$. 
        \end{enumerate}
    \end{enumerate}
    
    So we are done proving $\boxplus_4$.]

    \begin{enumerate}
        \item[$\boxplus_5$]   $\bbP_{\mathbf m_1} \lessdot \bbP_{\mathbf m}$.
    \end{enumerate}
    
    [Why?  By $\boxplus_4$ for $\gamma$ large enough.]
    
    So we are done.
\end{PROOF}

\begin{claim}\label{e37}
    If $\mathbf m \in \mathbf M$ is reduced or just $L_{\mathbf m}$ has cardinality $\le \lambda_2$  \then \, there is $\mathbf n \in \mathbf M_{\ec}$ of cardinality $\le \lambda_2$ such that $\mathbf m  \le_{\mathbf M} \mathbf n$.
\end{claim}

\begin{remark}\label{e39}
    By this we may restrict ourselves to $\mathbf M_{\le \lambda_2}$  (but then similarly in the end of \S2).
\end{remark}

\begin{PROOF}{\ref{e37}}
    We choose $\chi$ large enough and $\mathbf m_* \in \mathbf M_\chi$ which is wide, belongs to $\mathbf M_{\ec}$ and $\mathbf m \le_{\mathbf M} \mathbf m_*$; moreover is full and very wide (see \ref{e4}(1), as constructed in \ref{c41}).
    
    We can choose $\mathbf n$ such that:
    
    \begin{enumerate}
        \item[$(*)$]
        
        \begin{enumerate}
            \item[(a)]  $\mathbf n \in \mathbf M$ and $\mathbf n$ is wide and $|L_{\mathbf n}| = \lambda_2,$

            \item[(b)]  $\mathbf m \le_{\bfM} \mathbf n \le_{\mathbf M} \mathbf m_*,$

            \item[(c)]  $(\mathbf n,\mathbf m_*)$ satisfies  the criterion from \ref{e35}, with $\mathbf m_1,\mathbf m$ there 
            standing for $\mathbf n,\mathbf m_*$ here.
        \end{enumerate}
    \end{enumerate}
    
    [Why?  Let $\xi = 1$ and recalling Definition \ref{e24}(1) choose $\langle (t_\alpha,\bar s_\alpha):\alpha < \lambda_2)\rangle$ 
    such that $(t_\alpha,\bar s_\alpha) \in \cY_{\mathbf m_*},t_\alpha \in L_{\mathbf m_*} \backslash M_{\mathbf m_*},\langle t_\alpha/E_{\mathbf m}:\alpha < \lambda_2\rangle$ are pairwise distinct and for every $(t,\bar s) \in \cY_{\mathbf m_*}$ there are $\lambda^+$ ordinals $\alpha < \lambda_2$ such that $(t,\bar s),(t_\alpha,\bar s_\alpha)$ are $\xi$-equivalent, possible by \ref{e27} recalling $\lambda_2 \ge \beth_3(\lambda_1)$.   Let $L' = \cup\{t_\alpha/E_{\mathbf m_*}:\alpha <
    \lambda_2\} \cup L_{\mathbf m}$ and for each $t \in L' \backslash M_{\mathbf m_*}$ let $\langle s_{t,\alpha}:\alpha < \lambda^+\rangle$ be
    such that $s_{t,\alpha} \in L_{\mathbf m_*} \backslash M_{\mathbf m_*}$ and $\mathbf m_* \rest (s_{t,\alpha}/E_{\mathbf m_*})$ is isomorphic to
    $\mathbf m_* \rest (t/E_{\mathbf m_*})$ over $M_{\mathbf m}$.  Let $L = L' \cup \{s_{t,\alpha}:\alpha < \lambda^+,t \in  L' \backslash M_{\mathbf m_*}\}$ and $\mathbf n = \mathbf m_* \rest L$.  Now it is easy to check that $\mathbf n$ is as required.] 
    
    It suffices to prove that $\mathbf n$ belongs to $\mathbf M_{\ec}$, let
    $\mathbf n \le_{\mathbf M} \mathbf n_1 \le_{\mathbf M} \mathbf n_2$.
    
    \Wilog \, $L_{\mathbf n_2}$ has cardinality  $\le 2^{\lambda_2}$, by the LST argument;  (what is the LST argument here? let $ \chi _*$  be large enough such that $ \lambda, \mathbf{m}, \mathbf{m}_*,  \mathbf{n}_1, \mathbf{n} _2 $ all   belong to $ {\mathscr H} ( \chi _*)$ and let $ \mathfrak{A} \prec  ({\mathscr H}(\chi _*), \in ) $  be of cardinality $2^\lambda_{2}$    such that all  the above belong to it and  $ u \subseteq  \mathfrak{A} \wedge |u| \le  \lambda _2 \Rightarrow u \in \mathfrak{A} $. Now replace  $ \mathbf{n} _1, \mathbf{n} _2 $ by their restriction  to $ \mathfrak{A}  $).

    Now as $\mathbf m_*$ is very wide and full \wilog \, $\mathbf n_2 \le_{\mathbf M} \mathbf m_*$.  Now $(\mathbf n_1,\mathbf m_*)$ 
    satisfies the criterion from \ref{e35} hence $\bbP_{\mathbf n_1} \lessdot \bbP_{\mathbf m_*}$.
    
    Also the pair $(\mathbf n_2,\mathbf m_*)$ satisfies the criterion from \ref{e35} looking at the criterion.  Hence by \ref{e35} we have $\bbP_{\mathbf n_2} \lessdot \bbP_{\mathbf m_*}$.
    
    As $\mathbf n_1 \le_{\mathbf M} \mathbf n_2 \le_{\mathbf M} \mathbf m_*$ from the last two sentences it easily follows that $\bbP_{\mathbf n_1} \lessdot \bbP_{\mathbf n_2}$, so we are done.
\end{PROOF}

\begin{discussion}\label{e41}
    In what way does this proof help?  Will it not be simpler to omit in Definition \ref{c6} clause (c) the $\iota_{p(s)}, \mathbf B_{p(s),\iota}$, etc.?
    
    In this case in \ref{e4} we cannot define the projection directly hence we should look for projection as in general forcing, but then we
    run into problems of absoluteness.  More specifically, \ref{e35} seems to be problematic; anyhow this does not matter.
\end{discussion} 

\begin{definition}\label{e44}
    For $\mathbf m \in \mathbf M$ and $M$ is a subset of $M_{\bfm}$ so of cardinality $\le \lambda_1$ we define  $\mathbf n := \mathbf m \langle M\rangle \in \mathbf M$ as follows:
    
    \begin{enumerate}
        \item[(a)]  $L_{\mathbf n} = L_{\mathbf m}$ even as a partial order,

        \item[(b)]  $\bar u_{\mathbf n} = \bar u_{\mathbf m}$ and $\bar{\cP}_{\mathbf n} = \bar{\cP}_{\mathbf m}$,

        \item[(c)]  $M_{\mathbf n} = M$; yes $ M $ not $M_{\mathbf m}$!

        \item[(d)]  $E'_{\mathbf n} = \{(s,t):s,t \in L_{\mathbf m,}$ and $\{s,t\} \nsubseteq M\}.$  
    \end{enumerate}
\end{definition}

\begin{claim}\label{e47}
    Assume $\mathbf m \in \mathbf M_{\le \lambda_2}$ and $M $  is a subset of $ M_{\mathbf m}$.
    
    1) $\mathbf n := \mathbf m \langle M \rangle$ indeed belongs to $\mathbf M$ and is equivalent to $\mathbf m$ hence $\bbP_{\mathbf m}(L_{\mathbf m}) = \bbP_{\mathbf n}(L_{\mathbf m})$ i.e. $\bbP_{\bfm} = \bbP_{\bfn}.$ 
    
    2) If $\mathbf n  = \mathbf{m} \langle M \rangle  
    \le_{\mathbf M} \mathbf n_1$ \then \, for some $\mathbf m_1$ we have $\mathbf m \le_{\mathbf M} \mathbf m_1$ and $\mathbf m_1,\mathbf n_1$ are equivalent.
    
    3) If $\mathbf m \in \mathbf M_{\ec}$ and $\mathbf n = \mathbf m \langle M \rangle$ \then \, $\mathbf n \in \mathbf M_{\ec}$.
    
    4) If $\bfm \in \bfM_{\wec}$ and $\bfn = \bfm \langle M \rangle$ \underline{then} $\bfn \in \bfM_{\wec}.$ 
\end{claim}

\begin{PROOF}{\ref{e47}}
    1) Check, noting that $t \in L_{\mathbf n} \backslash M_{\mathbf n} \Rightarrow t \in L_{\mathbf m} \backslash M \Rightarrow |t/E'_{\mathbf n}| \le |L_{\mathbf n}| = |L_{\mathbf m}| \le \lambda_2$ and $|M_{\mathbf m}| = |M| \le |M_{\mathbf m}| \le \lambda_1$, (in fact, here $ M \subseteq M_ \mathbf{m} $ is not necessary,  only $ `` M $ has cardinality $ \le \lambda_1 "$).
    
    2) Given such $\mathbf n_1$ we now define $\mathbf m_1 \in \mathbf M$ by: 
    
    \begin{enumerate}
        \item[$(*)_1$]
        \begin{enumerate}
            \item[(a)]  $L_{\mathbf m_1} = L_{\mathbf n_1}$,
            
            \item[(b)]  $\bar u_{\mathbf m_1} = \bar u_{\mathbf n_1}$ and $\bar{\cP}_{\mathbf m_1} = \bar{\cP}_{\mathbf n_1}$, 
            
            \item[(c)]  $M_{\mathbf m_1} = M_{\mathbf m}$,
            
            \item[(d)]  $E'_{\mathbf m_1} = \{(s,t):s E'_{\mathbf m} t$ \underline{or}  $\{s,t\} \nsubseteq L_{\mathbf m}  $ but $ \{  s, t \} \subseteq L_{ \mathbf{n}_1}$ and $s E'_{\mathbf n_1} t\}$.
        \end{enumerate}
    \end{enumerate}
    
    Clearly:

    \begin{enumerate}
        \item[$(*)_2$] 
        \begin{enumerate}
            \item[(a)]  $\langle M_{\mathbf m}\rangle \char 94 \langle s/E''_{\mathbf m}:s \in L_{\mathbf m_1} \backslash M_{\mathbf m} \rangle \char 94 \langle t/E''_{\mathbf n_1}:t \in L_{\mathbf n_1} \backslash L_{\mathbf n}\rangle$ is a partition of $L_{\mathbf m_1} = L_{\mathbf n_1}$,  
            
            \item[(b)]  $E''_{\mathbf m_1} = E'_{\mathbf m_1} \rest \{(s,t) \in E'_{\mathbf m_1}$ and $s,t \notin M_{\mathbf m}\}$  is an equivalence relation, its equivalence classes  being the sets listed in clause (a) except $M_{\mathbf m}$,  
            
            \item[(c)]  $\mathbf m_1$ satisfies clause $(e)(\gamma)$ of Definition \ref{c4}. 
        \end{enumerate}
    
        \item[$(*)_3$]  
        
    \begin{enumerate}
        \item[(a)]  if $s \in L_{\mathbf m} \backslash M_{\mathbf m}$ then:  
            
        \begin{enumerate}
            \item[$(\alpha)$]  $s \in L_{\mathbf m_1} \backslash M_{\mathbf m_1}$, 
                
            \item[$(\beta)$]  $s/E'_{\mathbf m_1} = s/E'_{\mathbf m}$, 
                
            \item[$(\gamma)$]   $u_{\mathbf m_1,s} =
                u_{\mathbf n_1,s} = 
                u_{\mathbf n,s} = u_{\mathbf m,s}$, 
                
            \item[$(\delta)$]  $\cP_{\mathbf m_1,s} = \cP_{\mathbf m,s}= {\mathscr P} _{\mathbf{n}_1 , s }= 
            {\mathscr P} _{\mathbf{n}, s}$. 
        \end{enumerate}
        
        \item[(b)]  if $s \in L_{\mathbf m_1} \backslash L_{\mathbf m}$ then:  
        
        \begin{enumerate}
            \item[$(\alpha)$]  $s \in L_{\mathbf n_1} \backslash L_{\mathbf n}$,  
            
            \item[$(\beta)$]  $s/E'_{\mathbf m_1} = s/E'_{\mathbf n_1}$, 
            
            \item[$(\gamma)$]  $u_{\mathbf m_1,s} = u_{\mathbf n_1},s$, 
            
            \item[$(\delta)$]  $\cP_{\mathbf m_1,s} = \cP_{\mathbf n_1,s}$.  
        \end{enumerate}
        
        \item[(c)]  if $s \in M_{\mathbf m_1}$, 
        equivalently  $s \in M_{\mathbf m}$ then
        
        \begin{enumerate}
            \item[$(\alpha)$]  $u_{\mathbf m_1,s} = u_{\mathbf n_1,s}$
            
            \item[$(\beta)$]  $\cP_{\mathbf m_1,s} = \cP_{\mathbf n_1,s}  = 
            {\mathscr P} _{\mathbf{n} ,s} \cup 
            ({\mathscr P} _{\mathbf{n} _1 , s}
            \setminus {\mathscr P} _{\mathbf{n}  , s })$  
        \end{enumerate}
    \end{enumerate}
    \end{enumerate}
    
    and easily, 
    
    \begin{enumerate}
        \item[$(*)_4$]
        \begin{enumerate}
            \item[(a)] indeed $\mathbf m_1 \in \mathbf M$,
            
            \item[(b)]  $\mathbf m \le_{\mathbf M} \mathbf m_1$,
            
            \item[(c)]  $\mathbf m_1,\mathbf n_1$ are equivalent.
        \end{enumerate}
    \end{enumerate}
    
    So we are done.
    
    3) Assume $\mathbf n \le_{\mathbf M} \mathbf n_1 \le_{\mathbf M} \mathbf n_2$, as in the proof of part (2) there are $\mathbf m_1,\mathbf m_2$ such that $\mathbf m \le_{\mathbf M} \mathbf m_1 \le_{\mathbf M} \mathbf m_2$ and $\mathbf m_\ell,\mathbf n_\ell$ are equivalent for $\ell=1,2$. As $\mathbf m \in \mathbf M_{\ec}$ we have $\bbP_{\mathbf m_1} \lessdot \bbP_{\mathbf m_2}$ but this means $\bbP_{\mathbf n_1} \lessdot \bbP_{\mathbf n_2}$, as required.
    
    4) Similarly because $\bfm \in  \bfM_{\rm{wbd}} \Rightarrow \bfm \langle M \rangle \in \bfM_{\rm{wbd}}$. 
\end{PROOF}

\begin{conclusion}\label{e50}
    1) If $\mathbf m \in \mathbf M,M $ is a subset of $M_{\mathbf m}$ and $\mathbf n = \mathbf m \rest M$ \then \, $\bbP^{\cor}_{\mathbf n} \lessdot \bbP^{\cor}_{\mathbf m}$.

    2) If $\mathbf m_\ell \in \mathbf M$ and $M_\ell $ is a subset of $ M_{\mathbf m_\ell}$ for $\ell=1,2$ and $h$ is an isomorphism from $\mathbf m_1 \rest M_1$ onto $\mathbf m_2 \rest M_2$ \then \, $h$ induces  an isomorphism from $\bbP^{\cor}_{\mathbf m_1}[M_1]$ onto  $\bbP^{\cor}_{\mathbf m_2}[M_2]$.
\end{conclusion}

\begin{PROOF}{\ref{e50}}  
    1)  \Wilog \,  $\mathbf m   \in \mathbf M_{\le \lambda_2}$;  (why? because trivially $ \mathbf{n} \in \mathbf{M} _{\le \lambda_1 }$ and letting $ \mathbf{m} _1= \mathbf{m} \rest M_ \mathbf{m}$  we have  $\mathbf{m}_1    \le_\mathbf{M}  \mathbf{m}$ and  $ \mathbb{P} ^{\cor}_{\mathbf{m} _1}=  \mathbb{P} ^{\cor}_{\mathbf{m} }[M_ \mathbf{m} ]$ and  $\mathbf{n} = \mathbf{m} _1 \rest M _ \mathbf{m} $). By \ref{e37}  there is $\mathbf m_* \in \mathbf M^{\ec}_{\lambda_2}$ such that $\mathbf m \le_{\mathbf M} \mathbf m_*$ hence
    by \ref{b29}(2)  $\bbP^{\cor}_{\mathbf m} = \bbP_{\mathbf m_*}[M_{\mathbf m}]$.
    
    Let $\mathbf n_* = \mathbf m_* \langle 
    M \rangle$, see Definition \ref{e44}, so $\mathbf n_* \rest M = \mathbf n$ and by  \ref{e47}(3) 
    we have $\mathbf n_* \in \mathbf M_{\ec}$, hence $\bbP_{\mathbf n_*}[M  _{\mathbf n}] = \bbP^{\cor}_{\mathbf n}$. But $\mathbf n_*,\mathbf m_*$ are equivalent, hence $\bbP_{\mathbf n_*} = \bbP_{\mathbf m_*}$ hence $\bbP_{\mathbf n_*}[L] = \bbP_{\mathbf m_*}[L]$ for every $L \subseteq L_{\mathbf m_*}$ hence by \ref{b29}(3)
    $\bbP^{\cor}_{\mathbf n} = \bbP_{\mathbf n_*}[M 
    _{\mathbf n}] \lessdot \bbP_{\mathbf n_*}[M _{\mathbf m}] = \bbP_{\mathbf m_*}[M 
    _{\mathbf m}] = \bbP^{\cor}_{\mathbf
    m}$.  So the conclusion holds.

    2) Easy, too.
\end{PROOF}

\newpage 

\section{General $\bfm$'s}

This section depend on \S1A, \S1C, \S2, \S3A, \S3C but not on \S1B, \S3B, \S3D. 

\subsection{Alternative proof}\label{4A}

\begin{hypothesis}\label{h0}
    We are in the general context. 
\end{hypothesis}

This sub-section plays a double role. First, we give an alternative proof of the main results, they may be simpler but we lose some information and we are \underline{assuming} $\lambda_{2} \geq \beth_{\lambda_{1}^{+}}.$ Second, it give proof which works also for the fat context and even the neat and general contexts not just the lean context (as in \S3D). Specifically, 

\begin{enumerate}
    \item[$\boxplus$] in this version:  
    
    \begin{enumerate}
        \item[(a)] we ignore \S1B, that is \ref{c33n}, \ref{c33s},
        
        \item[(b)] we ignore \S3B that is \ref{e19}-\ref{e28}
        
        \item[(c)] we ignore or replace almost all \S3D, that is:
        
        \begin{enumerate} 
            \item[$(\alpha )$]   we ignore Claims 
            \ref{e33}, \ref{e34}, \ref{e35},
            
            \item[$(\beta )$]  we replace Claim \ref{e37} by   \ref{h2}(2), 
 
            \item[$(\gamma )$] Def \ref{e44}, \ref{e47}(1),(2), (3) remains,
      
            \item[$(\delta )$] Claim \ref{e50} is replaced by \ref{h26}  (whose proof just say ``repeat the proof of \ref{e50}'').
        \end{enumerate} 
    \end{enumerate}       
\end{enumerate}

\begin{definition}\label{h2}
    1) Let $\Omega_{\bfm}^{1} := \bigcup \{ \Omega_{\bfm, t}^{1}: t \in L_{\bfm} \},$ where for $t \in L_{\bfm}, \, \Omega_{\bfm, t}^{1}$ is the set of $\bfb = \langle t, \bfB, \bar{c}, d, c, \iota, g \rangle$ such that: 
    
    \begin{enumerate}
        \item[(a)] $\bar{c} = \langle c_{i}: i < i_{\bfb} \leq \lambda \rangle,$
        
        \item[(b)] $d \subseteq c = \bigcup \{ c_{i}: i < i_{\bfb} \} \subseteq \lambda,$
        
        \item[(c)] $\bfB$ is a Borel function from ${}^{c} \cH(\lambda)$ into $\cH_{< \lambda(+)} (\cB_{t}),$ so if $\rho \in ({^{c}}\cH(\lambda))^{\bfV[\bbR]},$ then $\bfB(\rho)$ belongs to $({}^{c}\cH(\lambda))^{\bfV[\bbR]}$ but not necessarily to $\bfV,$
        
        \item[(d)] $\iota < i_{\bfb},$
        
        \item[(e)] $g$ is a function from $c_{\iota}$ into $L$  such that $\varp \in c_{\iota} \Rightarrow [g(\varp) \in M_{\bfm} \equiv (\varp \in d)]$ and $\rang(g)$ is included in some $L \in \cP_{t}.$
    \end{enumerate}
    
    2) Let $\Omega_{\bfm}^{2} := \bigcup \{\Omega_{\bfm, t}^{2}: t \in L_{\bf,}  \},$ where for $t \in L_{\bfm}, \, \Omega_{\bfm, t}^{2}$ is the set of $\bar{\bfb}$ such that:
    
    \begin{enumerate}
        \item[(a)] $\bar{\bfb} = \langle \bfb_{j}: j < \lg(\bar{\bfb}) \leq \lambda \rangle,$
        
        \item[(b)] $\bfb_{j} \in \Omega_{\bfm, t},$
        
        \item[(c)] $t_{\bfb_{j}} = t, \iota_{\bfb_{j}} = j$ for $j < \lg(\bar{\bfb}),$
        
        \item[(d)] $\langle 1, t_{\bfb_{j}}, \bar{c}_{\bfb_{j}}, d_{\bfb_{j}} \rangle$ is the same for all $j < \lg(\bar{\bfb}),$
    \end{enumerate}
    
    2A) For $t \in L_{\bfm},$ we say $\bar{\bfb} \in \Omega_{\bfm, t}^{2}$ \emph{strictly represent} $p(t)$ \underline{when}: 
    
    \begin{enumerate}
        \item[(a)] $p = p \rest \{ t \} \in \bbP_{\bfm}$ and in Definition \ref{c6}(2) we have $\iota_p(s) = 1,$
        
        \item[(b)] $p(s)$ is $\bfB(\dots, g_{\bfb_{j}(\zeta)}, \dots)_{j \in \lg(\bar{b}), \zeta \in c_{\bfb[j]}, i}.$
        
        
            
            
            
            
            
            
            
            
    \end{enumerate}
    
    2B) We let $\Omega_{\bfm}^{3} := \bigcup \{ \Omega_{\bfm, t}^{3}: t \in L_{\bfm} \},$ where for $t \in L_{\bfm}$ we let $\Omega_{\bfm, t}^{3}$ be the family of $\gb$ such that $\gb$ is a subset of $\Omega_{\bfm, t}^{2}$ of cardinality $< \lambda.$ 
    
    2C) We say $\gb \in \Omega_{\bfm}^{3}$ represents $p(s)$ \underline{when}: 
    
    \begin{enumerate}
        \item[(a)] $p \in \bbP,$
        
        \item[(b)] $s \in \dom(p),$
        
        \item[(c)] $p(s) = \sup_{\varp < \varp_{\ast}} (\eta, \name{f}_{\varp}),$ where $\gb = \{ \bar{\bfb}_{\varp}: \varp < \varp_{\ast} \},$ each $\bar{\bfb}_{\varp} \in \Omega_{\bfm, \ast}^{2}$ and $(\eta, \name{f}_{\varp})$ is strictly represented by $\bar{\bfb}_{\varp.}$
    \end{enumerate}
    
    3) For $\bfm \in \bfM$ we shall define a model $\md(\bfm),$ pedantically it is \ $\md_{\bar{t}}(\bfm),$ where\footnote{So instead we can use $\langle t_{\alpha} \rest E_{\bfm}'': \alpha < \alpha_{\ast} \rangle.$} $\bar{t} = \langle t_{\alpha} = t_{\bfm, \alpha}: \alpha \leq \alpha_{\bfm} = \alpha(\bfm) \rangle, t_{\alpha(\bfm)}$ is a fix member of $M_{\bfm}$ $\bar{t} \rest \alpha_{\bfm}$  is a maximal sequence of pairwise non-$E_{\bfm}''$-equivalent members of $L_{\bfm} \setminus M_{\bfm}$ (so below $A_{\alpha}^{\ell} = A_{\bfm, \alpha}^{\ell}$ for $\alpha \leq \alpha_{\bfm}), t_{\alpha(\bfm)} \in M_{\bfm}$ and stipulate $M_{\bfm} = t_{\alpha(\bfm)} / E_{\alpha(\bfm)}''$ ignoring the case $M_{\bfm} = \emptyset$:
        
    \begin{enumerate}
        \item[(A)] The set of elements of $\md(\bfm)$ is the disjoint union of the following sets; below $\alpha < \alpha(\bfm)$: 
            
        \begin{enumerate}
            \item[(a)] $A_{\alpha}^{1} = \{ (1, t_{\alpha}, s): s \in t_{\alpha} / E_{\bfm}' \},$ see \ref{c4}(e)($\varp$), ($\zeta$),
                
            \item[(b)] $A_{\alpha}^{2} = \{ (2, t_{\alpha}, p): p \in \bbP_{\bfm}(t_{\alpha} / E_{\bfm}') \} \cup \{ (2, t_{\alpha}, p, q): \bbP_{\bfm}[t_{\alpha} / E_{\bfm}] \models p \leq q \},$  see Definition \ref{c7}, central for the lean context, 
                
            \item[(c)] $A_{\alpha}^{3} = \{ (3, t_{\alpha}, s, \bfb): s \in t_{\alpha} / E_{\bfm}'', \, \bfb \in \Omega_{\bfm}^{1} \ \text{and} \ \rang(g_{\bfb}) \subseteq t_{\alpha} / E_{\bfm}' \},$ 
            
            \item[(d)] $A_{\alpha}^{4} = \{ (4, t_{\alpha}, \psi): \psi \in \bbP_{\bfm}[t_{\alpha} / E_{\bfm}] \} \cup \{ (4, t_{\alpha}, \psi, \varphi): \bbP_{\bfm}[t_{\alpha} / E_{\bfm}]  \models \psi \leq \varphi \},$ 
                
            \item[(e)] $A_{\alpha(\bfm)}^{1} = \{ (1, \alpha_{\bfm}, s, \ell): s \in M_{\bfm}$ and $\ell = 1 \Rightarrow s \in M_{\bfn}^{\rm{lean}}, \, \ell = 2 \Rightarrow s \in M_{\bfm}^{\fat}, \, \ell = 0 \Rightarrow s \in M_{\bfm}^{\rm{non}} \},$
                
            \item[(f)] $A_{\alpha(\bfm)}^{2} = \{ (2, t_{\alpha_{\bfm}}, p): p \in \bbP_{\bfm}(L_{\bfm}) \} \cup \{ (2, \alpha_{\bfm}, p, q): \bbP_{\bfm}(L_{\bfm}) \models p \leq q \},$ 
                
            \item[(g)] $A_{\alpha(\bfm)}^{3} = \{ (3,  t_{\alpha_{\bfm}}, s, \bfb): s \in M_{\bfm}, \, \bfb \in \Omega_{\bfm}^{1}, \ \text{and} \ \rang(g_{\bfb}) \subseteq M_{\bfm} \},$ 
                
            \item[(h)] $A_{\alpha}^{4}(\bfm) = \{ (4, \alpha, \psi): \psi \in \bbP_{\bfm}[t_{\alpha} / E_{\bfm}] $ and $\alpha = \alpha_{\bfm} \} \cup \{ (4, \alpha_{\bfm}, \psi, \varphi): \bbP_{\bfm}[M_{\bfm}] \models \psi \leq \varphi \},$ 
            
            \item[(i)] notation: for $\alpha < \alpha(\bfm),$ $A_{\alpha} = A_{\alpha}^{1} \cup A_{\alpha}^{2} \cup A_{\alpha}^{3} \cup A_{\alpha}^{4} \cup A_{\alpha(\bfm)}^{1} \cup A_{\alpha(\bfm)}^{2} \cup A^{3}_{\alpha(\bfm)} \cup A_{\alpha(\bfm)}^{4}.$
        \end{enumerate}
            
        \item[(B)] The relations of $\md(\bfm)$ are the relations $R$ on $\md(\bfm)$ such that: 
            
        \begin{enumerate}
            \item[(a)] $R = \bigcup \{ R \rest A_{\alpha}: \alpha < \alpha_{\bfm} \},$
                
            \item[(b)] (an overkill) $R$ is first order definable in $(\cH(\chi_{\bfm}), \in, <^{*}_{\chi[\bfm]}, \bfm),$ where $<^{*}_{\chi[\bfm]}$ is a well ordering of $\cH(\chi_{\bfm}).$
        \end{enumerate}
            
        \item[(C)] In particular there is an individual constant for each $c \in A_{\alpha(\bfm)}^{1} \cup A_{\alpha(\bfm)}^{2} \cup A_{\alpha(\bfm)}^{3} \cup A_{\alpha(\bfm)}^{4},$ or code them by unary relations. 
    \end{enumerate}
    
    4) For $s \in L_{\bfm} \setminus M_{\bfm},$ the model $\md(\bfm) \rest (s / E_{\bfm})$ is naturally defined as (when $s E_{\bfm}' t_{\alpha}$) the restriction of the model $\md(\bfm)$ to $\bigcup \{ A_{\alpha}^{\ell}: \ell = 1, 2, 3, 4 \} \cup \{ A_{\alpha(\bfm)}^{\ell}: \ell = 1, 2, 3, 4 \}.$  
\end{definition}

\begin{definition}\label{h7}
    1) We say that $\bar{a} \in {}^{\lambda}(\md(\bfm))$ represents $p \in \bbP_{\bfm}$ \underline{when} for some $\bar{\alpha}$ we have: 
        
    \begin{enumerate}
        \item[(a)] $\bar{\alpha}$ is an sequence of ordinals $\leq \alpha(\bfm)$ of length $\zeta_{p} < \lambda,$
            
        \item[(b)] we let $\bar{\alpha} =  \langle \alpha(\varp) = \langle \alpha_{\varp}: \varp < \zeta_{p} \rangle,$
            
        \item[(c)] $\wsupp(p)$ is equal to $\bigcup \{ t_{\bfm, \alpha_{\varp}} / E_{\bfm}: \varp < \zeta_{p} \} \cup M_{\bfm},$
            
        \item[(d)] if $s \in \dom(p),$ then the following set $\gb$ represents $p(s),$ where $\gb$ is the set of $\bar{\bfb} \in \Omega_{\bfm, s}^{2}$ such that for each $i < \lg(\bar{\bfb})$ we have $t_{\bfb_{i}} = s$ and there is an $\varp < \zeta_{p}$ such that:  
        
        \begin{itemize}
            \item $a_{2 \varp} = (1, t_{\alpha_{\varp}}, s) \in A_{\alpha}^{1}$ and $\alpha_{\varp} < \alpha_{\bfm},$
            
            \item $a_{2 \varp} = (1, t_{\alpha_{\varp}}, s, \ell) \in A_{\alpha(\bfm)}^{1}$ and $\alpha_{\varp} = \alpha_{\bfm},$ 
            
            \item $a_{2 \varp + 1} = (\zeta, t_{\alpha_{\varp}}, s, \bfb, \varp) \in A_{\alpha_{\varp}}^{3}.$
        \end{itemize}
        
            
                
            
            
            
            
        
        \item[(e)] if $\varp < \zeta_{p}$ then one of the cases above occurs.
            
        \item[(f)] if $\varp \in [2 \zeta_{p}, \lambda)$ then $a_{\varp}$ is the triple $(2, \alpha_{\bfm}, \psi) \in A_{\alpha(\bfm)}^{2}$ where $\psi \in \bbP_{\bfm}[M_{\bfm}]$ is witnessed by $p.$
    \end{enumerate}
    
    2) We say $\bar{a}$ is a formal representative for $\bbP_{\bfm}$ \underline{when} for some $\bar{\alpha}$ the demands above holds (ignoring the existence of $p$).
    
    3) We say $\bar{a} \in {}^{\lambda }(\md(\bfm))$ is a formal representation of a member of $\bbP_{\bfn}[M_{\bfm}]$ similarly using $A_{\alpha}^{4}, A_{\alpha(\bfm)}^{4}. $
\end{definition}

\begin{claim}\label{h10}
    Here,
    
    \begin{enumerate}
        \item[(a)] every $p \in \bbP$ is represented by some $\bar{a} \in {}^{\lambda}\md(\bfm),$
        
        \item[(b)] every formal representative represent some member of $\bbP_{\bfm},$
        
        \item[(c)] there is a formula $\psi_{\rm{rep}}(\bar{x}_{[\lambda]})$ in the logic $\bbL_{\lambda_{1}^{+}, \lambda_{1}^{+}}$ in the vocabulary of $\md(\bfm)$ defining the set of formal representative,
        
        \item[(d)] similarly for  $p \in \bbP_{\bfm}[L_{\bfm}]$ more accurately $\psi \in \bbL_{\lambda^{+}}[Y_{\bfm}]$ not excluding contradictory ones. 
    \end{enumerate}
\end{claim}

\begin{PROOF}{\ref{h10}}
    Easy. 
\end{PROOF}

\begin{definition}\label{h13}
    We say $L$ is \emph{good} \underline{when}: 
    
    \begin{enumerate}
        \item[(a)] $L$ is a initial segment of $L_{\bfm},$
        
        \item[(b)] $L$  is $\bbL_{\lambda_{1}^{+}, \lambda_{1}^{+}}$-definable in $\md_{\bar{t}}(\bfm)$ (without parameters),
        
        \item[(c)] the following are definable in $\md_{\bar{t}}(\bfm)$ by a formula (without parameters) in $\bbL_{\lambda_{1}^{+}, \lambda_{1}^{+}}$:
        
        \begin{itemize}
            \item $\bar{a}$ represent some $p \in \bbP_{\bfm}[L],$
            
            \item $\bar{a}$ represent some $p \in \bbP_{\bfm}(L),$
            
            \item $\bar{a}_{1}, \bar{a}_{2}$ represent $p_{1}, p_{2} \in \bbP_{\bfm}(L)$ respectively, and $p_{1} \leq_{\bbP_{\bfm(L)}} p_{2},$
            
            \item $\bar{a}_{1}, \bar{a}_{2}$ represent $p_{1}, p_{2} \in \bbP_{\bfm}[L]$ respectively and $p_{1} \leq_{\bbP_{\bfm}[L]} p_{2}.$
            
        \end{itemize}
    \end{enumerate}
\end{definition}

\begin{claim}\label{h16}
    1) The set $L = \{ s \in L_{\bfm}: \ \text{for no} \ t \in M_{\bfm} \ \text{do we have} \ t \leq_{\bfm} s \}$ is good.
    
    2) If $L$ is good and $t_{*} \in M_{\bfm} \setminus L$ but $L_{\bfm(< t_{*})} \subseteq L,$ then $L \cup \{ t_{*} \}$ is good. 
    
    3) If $\langle L_{\alpha}: \alpha < \delta \rangle$ is an $\subseteq$-increasing sequence of good sets and $\delta < \lambda_{1}^{+}$ \underline{then} so is $\bigcup_{\alpha < \delta} L_{\alpha}.$
 
    4) If $L$ is good then $L^{+} = \{ s: \text{there is no} \ t \in M_{\bfm} \setminus L \ \text{such that} \ t \leq_{\bfm} s \}$ is good.
\end{claim}

\begin{PROOF}{\ref{h16}}
    Notice that 1), 2) and 3) are straightforward. Concerning part (2) the reader may wonder: how do you define $\bbP_{\bfm}(L)$ not using parameters if, say, multiple such $t_{\ast}$'s exists?. The answer is by \ref{h2}(3) clause (c); that is, each $t \in M_{\bfm}$ is definable without parameters. 
    
    4) The point is that we do not like \ to induct on $\rm{dp}(s, L_{\bfm})$ just on $\rm{dp}(t, M_{\bfm}).$ Note that the clauses on $\bbP_{\bfm}[L^{+}]$ follows by those on $\bbP_{\bfm}(L^{+}).$ What we do is noting:
    
    \begin{itemize}
        \item[$\oplus$] for $p, q \in \bbP_{\bfm}(L^{+}), p \leq q$ iff:
        
        \begin{enumerate}
            \item[(a)] $p \rest L \leq_{\bbP_{\bfm}(L)} q \rest L,$
            
            \item[(b)] if $s \in \dom(p) \setminus L$ \underline{then} necessarily $s \in L^{+} \setminus L$  and $s / E_{\bfm}''$ appears in $\langle t_{\alpha} / E_{\bfm}: \alpha < \alpha_{\bfm} \rangle$ and,
            
            \begin{itemize}
                \item[$\bullet$] $\bbP_{\bfm}[L \cup (s / E_{\bfm})] \models p \rest (L \cup (L^{+} \cap s / E_{\bfm})) \subseteq q \rest (L \cup (L^{+} \cap s / E_{\bfm})).$
            \end{itemize}
        \end{enumerate}
    \end{itemize}
    
    [Why? Just think.]
    
    Recalling \ref{h2}(3)(A)(d) it suffice to prove that: 
    
    \begin{enumerate}
        \item[$(\ast)$] Assume $s \in L^{+} \setminus L, p, q \in \bbP_{\bfm}$ and $\dom(p) \subseteq \dom(q) \subseteq L \cup (s / E_{\bfm}')$ \underline{then} $\bbP_{\bfm} \models$ ``$p \leq q$ iff (a) $+$ (b)'', where:
        
        \begin{enumerate}
            \item[(a)] $(p \rest L) \leq_{\bbP_{\bfm}(L)} (q \rest L),$
            
            \item[(b)] for some $\psi \in \bbL_{\lambda^{+}}[Y_{s / E_{\bfm}' \cap L}]$ we have: 
            
            \begin{enumerate}
                \item[$\bullet_{1}$] $\bbP_{\bfm}[L_{\bfm}] \models \psi \subseteq q,$
                
                \item[$\bullet_{2}$] $\psi \wedge q \wedge \neg p \notin \bbP_{\bfm}[L^{+} \cap (s / E_{\bfm}')].$
            \end{enumerate}
        \end{enumerate}
    \end{enumerate}
    
    [Why $(\ast)$ holds? As in \S3C.]
\end{PROOF}

\begin{claim}\label{h19}
    1) $\mathbb{P} _ \mathbf{m} $ is $ \mathbb{L} _{\lambda ^+_1, \lambda ^+_1}$-interpretable in $ \md(\mathbf{m}). $ 

    2) We have $ \mathbb{P} _\mathbf{m}  \lessdot   \mathbb{P} _ \mathbf{n} $ \underline{when}:  
    
    \begin{enumerate} 
        \item[(a)]  $ \mathbf{m} \le _ \mathbf{M}  \mathbf{n}, $
        
        \item[(b)] for every $ \zeta < \lambda ^+_1$ and $ t \in L_{\bfn} \setminus L_{\mathbf{m}} $ there are  at least $ \lambda ^+_1 $ elements $s \in L_{\bfm} \setminus M_{\bfm}$  such that, recalling \ref{h2}(3), the models  $\md( \mathbf{n} ) \rest (t / E_{\bfn}), \, \md( \mathbf{n} ) \rest (s / E_{n})$  are $ \mathbb{L} _{\beth^+ _ \zeta , \beth ^+ _ \zeta }$-equivalent.
    \end{enumerate}      
 
    3) If $\bfn \in \bfM_{\rm{ec}}$ is wide and full, $A \subseteq L_{\bfn}$ has cardinality $\leq \beth_{\lambda_{1}^{+}}$ \underline{then} there is $\bfm$ such that: 

    \begin{enumerate}
        \item[(a)] $\bfm \leq_{\bfM} \bfn,$
        
        \item[(b)] $L_{\bfm}$ has cardinality $\leq \beth_{\lambda_{1}^{+}},$
        
        \item[(c)] $\bfm \in \bfM_{\rm{ec}}.$
    \end{enumerate}
    
    4) Similarly to part (3) for $\bfM_{\rm{bec}}.$ 
\end{claim} 

\begin{PROOF}{\ref{h19}}
    1) Let $\langle s_{\zeta}: \zeta < \zeta_{\ast} \rangle$ lists the elements of $M_{\bfm}$ such that $s_{\varp} <_{M_{\bfm}} s_{\varp} \Rightarrow \varp < \zeta:$ exists as $L_{\bfm}$ is a (possibly partial) well order. Clearly $\zeta_{\ast} < \Vert M_{\bfm} \Vert^{+} \leq \lambda_{1}^{+}.$ We define $L_{\zeta}$ for $\zeta \leq 2 \zeta_{\ast} + 1$ as follows: 
    
    \begin{itemize}
        \item if $\varp \leq \zeta_{\ast},$ then $L_{2 \varp} = \{ t \in L_{\bfm}:$ for some $\zeta < \varp$ we have $t \leq_{L_{\bfm}} s_{\zeta} \},$
        
        \item for $\varp \leq \zeta_{\ast}$ we let $L_{2 \varp +1} = L_{2 \varp} \cup \{ t \in L_{\bfm}:$ if $\zeta \in [\varp, \zeta_{\ast})$ then $s_{\zeta} \nleq t$ \underline{or} for some $\xi < \varp$ we have $(\forall \zeta \in [\xi, \zeta_{\ast}))(s_{\zeta} \nleq t) \}.$  
    \end{itemize}

    Pedantically, the models   $ \langle  \md( \mathbf{m} ) \rest A_{\mathbf{m} , \alpha }: \alpha < \alpha _ \mathbf{m} \rangle $ are not pairwise disjoint but the common part consists  of $\lambda_{1}$ individual constants \underline{hence} this does not matter.
    
    Clearly $L_{2 \zeta_{\ast} +1} = L_{\bfm}, L_{\varp}$ is a definable initial segment of $L_{\bfm}.$ By Definition \ref{h13} it suffice to prove that $L_{\bfm}$ is good. 
    
    Now we prove by induction on $\varp$ that $L_{\varp}$ is good, so for $\varp = 2 \zeta_{\ast} +1$ we get the desired conclusion. 
    
    For $\zeta = 0,$ this holds by \ref{h16}(1).
    
    For $\zeta = 2 \varp +1$ we have $L_{\zeta} \setminus L_{2 \varp} = \{ s_{\varp} \}$ and $L_{\bfm(< s_{\alpha})} \subseteq L_{2 \varp};$ hence by \ref{h16}(2) we are done. 
    
    For $\zeta = 2 \varp + 2$ we apply \ref{h16}(4).
    
    Lastly, for $\zeta$ limit we apply \ref{h16}(3). Together we are done. 
    
    2) By part (1) and the addition theorem, (best formulated for the intermediate logic $ \mathbb{L} _{\infty , \lambda ^+_ 1, \zeta}$ for $ \zeta < \lambda ^+_1)$,  see \cite{Dic85}). 
        
    3), 4) As in the proof of \ref{h22} below. 
\end{PROOF}

\begin{claim}\label{h22}
    Assume $ \lambda _2 \ge \beth _{\lambda ^+_1}.$
    
    1) If $ \mathbf{n} \in \mathbf{M}$ satisfy $ L_\mathbf{n} = M_ \mathbf{n}$ (e.g. it is isomorphic to $ (\gamma, < ), \gamma < \lambda_{1}^{+})$  \underline{then}  there is $ \mathbf{m} \in \mathbf{M} _{\ec}$  of cardinality $  \lambda _2$ above $ \mathbf{n} $. 
    
    1A) Similarly for $\bfM_{\rm{bec}}.$ 
    
    1B) Moreover if $\bfn$ is strongly $(< \lambda^{+})$-directed (see \ref{b36}(2), if $L_{\bfn} = M_{\bfn} = (\gamma, <)$ for some $\gamma < \lambda_{1}^{+},$ this mean $\cf(\gamma) = \lambda$) \underline{then} (in part (1A)) $\bfm$ is strongly $(< \lambda^{+})$-directed, so $\{ \name{\eta}_{r}: r \in M_{\bfm}  \}$ is cofinal in $\Pi_{\varp < \lambda} \theta_{\varp}$ in $\bfV^{\bbP_{\bfm}},$ so $\bfm \in \bfM_{\rm{bec}}.$
    
    2) If $ \mathbf{m} _1 \in \mathbf{M} $ has cardinality $ \le \lambda _2$ \underline{then} we can demand in part (1) $\mathbf{m}_1 \le _ \mathbf{M} \mathbf{m}.$
    
    2A) If $\bfm_{1} \in \bfM_{\bd}$ has cardinality $\leq \lambda_{1},$ \underline{then} in part (1A) we can demand $\bfm_{1} \leq_{\bfM} \bfm.$ 
\end{claim} 

\begin{PROOF} {\ref{h22}}
    1) Let $ \mathbf{n}_{\bullet}$ be  very wide full  
    of cardinality  $ 2^{\lambda _2}$ such that $ \mathbf{n} \le _\mathbf{M} \mathbf{n}_{\bullet}$ and let $\bfn_{\ast} = \bfn_{\bullet}^{[\bd]}$.  We can find  $ \mathbf{m} \le _ \mathbf{M} \mathbf{n}_*$ of cardinality $ \lambda _2$  as in  \ref{h19}(2), because for every $ \zeta < \lambda ^+_ 1$  there are  $ < \beth _{\lambda ^+_1}$ theories in the relevant vocabulary and logic. So $ L_ \mathbf{m} $ has cardinality $\leq \lambda_{2}$ and $ \mathbf{n} \le _\mathbf{M} \mathbf{m} $ but why does it belong to
    $ \mathbf{M} _{\rm{bec}}$? Toward contradiction let   $ \mathbf{m} _1, \mathbf{m} _2 \in \bfM_{\rm{bec}}$ be such that $ \mathbf{m} \le _ \mathbf{M} \mathbf{m} _1 \le _ \mathbf{M} \mathbf{m} _2$ but $ \mathbb{P} _{\mathbf{m} _1 } \lessdot \mathbb{P} _{\mathbf{m} _2}$ fail. By the L.S.T. argument,  (see the proof of \ref{e37} third paragraph),  \wilog \,  $ \mathbf{m} _2 $ has cardinality $ \le 2^{\lambda _2}$,  Hence by the choice of $ \mathbf{m} , \mathbf{n} $ \wilog \, $\mathbf{m} _2 \le _ \mathbf{M} \mathbf{n}_{\ast}$.  Now for $ {\ell} = 1,2 $,  by \ref{h19}(2) applied to  $ (\mathbf{m} _ {\ell} , \mathbf{n}_{\ast})$ we have $ \mathbb{P} _ {\mathbf{m}_{\ell}}  \lessdot  \mathbb{P}_{\mathbf{n}_{\ast}}.$  But this implies $ \mathbb{P} _{\mathbf{m} _1 } \lessdot   \mathbb{P} _{\mathbf{m} _2}$  so we are done. 
    
    1A) Similarly.
    
    1B) In the proof of part (1) we restrict ourselves to strongly $(< \lambda^{+})$-directed $\bfm$-s (see \ref{c5}(10)) so we use the relevant criterion for being in $\bfM_{\rm{bec}},$ see \ref{b47}(7) i.e. consider bounded $\bfm$-s only: $\bfm \leq \bfm_{1} \leq \bfm_{2}, \bfm_{1}, \bfm_{2}$ strongly $\lambda^{+}$-directed $\Rightarrow \bbP_{\bfm_{1}} \lessdot \bbP_{\bfm_{2}}.$ The cofinality is by \ref{c37}(3).
    
    2), 2A) Similarly. 
\end{PROOF} 

\begin{conclusion}\label{h26}
    1) If $\mathbf m \in \mathbf M,M \subseteq M_{\mathbf m}$ and $\mathbf n = \mathbf m \rest M$ \then \, $\bbP^{\cor}_{\mathbf n} \lessdot \bbP^{\cor}_{\mathbf m}$.
    
    2) If $\mathbf m_\ell \in \mathbf M$ and $M_\ell \subseteq M_{\mathbf m_\ell}$ for $\ell=1,2$ and $h$ is an isomorphism from $\mathbf m_1 \rest M_1$ onto $\mathbf m_2 \rest M_2$ \then \, $h$ induces 
    an isomorphism from $\bbP^{\cor}_{\mathbf m_1}[M_1]$ onto  $\bbP^{\cor}_{\mathbf m_2}[M_2]$.

    3) If $\bfm \in \bfM_{\rm{bec}}$ is strongly $\lambda^{+}$-directed, $M \subseteq M_{\bfm}$ is cofinal in $M_{\bfm}$ \underline{then}  $\Vdash_{\bbP_{\bfm}}$``$\{ \name{\eta}_{s}: s \in M \}$ is cofinal in $(\Pi_{\varp < \lambda} \theta_{\varp}, <_{J_{\lambda}^{\rm{bd}}})$''.

\end{conclusion} 

\begin{PROOF}{\ref{h26}}
    For 1) and 2) it suffices to proceed exactly as the proof \ref{e50}, replacing quoting \ref{e37} by quoting \ref{h22}(2).  Also, 3) is easy by now. 
\end{PROOF}

\subsection{General $\bfm$'s}\label{3E}\

See Discussion \ref{g0} for our aim and \ref{g15} on the connection to \cite{Sh:945}.   

\begin{definition}\label{e51}
    Assume  $ \mathbf{m}$ is $\lambda_{0}$-wide.   Let $ \mathbb{P} ^\dagger _\mathbf{m} =  \mathbb{P} ^ \dagger_\mathbf{m} [M_ \mathbf{m} ]$     be the forcing notion $ \mathbb{P} _ \mathbf{m}  [M_ \mathbf{m}]$ restricted to the set of $ \psi \in \mathbb{P}_\mathbf{m} [M_\mathbf{m} ]$ such that there is $ p \in \mathbb{P} _\mathbf{m} (L_ \mathbf{m})$ witnessing it, which means that (it is the projection of $ p $ into $ \mathbb{P} _ \mathbf{m} [M_ \mathbf{m} ] $, that is): 
    
    \begin{enumerate}
        \item[$\bullet _1$] the condition $ \psi  $ is smaller or equal to $ p $  in the forcing notion $\mathbb{P} _ \mathbf{m} [L_\mathbf{m} ],$
        
        \item[$\bullet _2$]  if $ \mathbb{P} _ \mathbf{m} [M_\mathbf{m}] \models ``\psi \le \varphi" $ \then \, $ \varphi, p $ are compatible in the forcing notion $ \mathbb{P} _\mathbf{m} [L_ \mathbf{m}].$
    \end{enumerate}
\end{definition} 

\begin{claim}\label{e52} 
    Assume $ \mathbf{m} $ to be $\lambda_{0}$-wide.
      
    1) $\mathbb{P} ^ \dagger _ \mathbf{m} $ is a dense subset  of $ \mathbb{P} _ \mathbf{m} [M_\mathbf{m} ]$, hence  $ \mathbb{P} ^ \dagger _ \mathbf{m}  \lessdot \mathbb{P} _ \mathbf{m} [L_ \mathbf{m}]$.
    
    2) If $L$ is an initial segment of $L_{\bfm}$ and $\bfn = \bfm \rest L,$ \then \,  $ \mathbb{P} ^ \dagger_\mathbf{n} =  \mathbb{P} ^\dagger_ \mathbf{m}  \cap  \mathbb{P} _ \mathbf{n}[M_ \mathbf{n}]$. 
    
    3) If $L$ is a $\lambda_{0}$-wide initial segment of $L_ \mathbf{m}$, and $ \mathbf{n} = \mathbf{m} \upharpoonright L,$ \then:  
     
     \begin{enumerate} 
         \item[(a)]  $ \mathbb{P}_\mathbf{n} [M_\mathbf{n} ] \lessdot \mathbb{P} _ \mathbf{n} [L_ \mathbf{n} ]$ and $ \mathbb{P} _ \mathbf{n} [M_ \mathbf{n} ] = \mathbb{P} _ \mathbf{m} [M_\mathbf{n} ]  \lessdot \mathbb{P} _ \mathbf{m} [M_ \mathbf{m}]$,  
         
         \item[(b)] if $ p_1 \in \mathbb{P} _ \mathbf{n} (L_ \mathbf{n} )$ \underline{then}  there is $ \psi  \in \mathbb{P} _ \mathbf{n} [M_ \mathbf{n}]$ satisfying:  
         
         \begin{enumerate} 
             \item[($ \alpha $)]  $ \psi   \le p_1 \in \mathbb{P} _ \mathbf{n} [L_ \mathbf{n}]$, 
             \item[($ \beta  $)] if $ \psi  \le  \varphi  \in \mathbb{P} _ \mathbf{n} [M_ \mathbf{n} ] $   then $ p_1, \varphi $ are compatible  in $ \mathbb{P} _ \mathbf{n} [L_ \mathbf{n}]$,
             
             \item[($\gamma$)]  $ \psi $ being witnessed by $ p_1$, (see Definition \ref{e51} this follows).
        \end{enumerate} 
    \end{enumerate} 
\end{claim}

\begin{PROOF}{\ref{e52}}
    1) Let $\varphi \in \bbP_{\bfm}[M_{\bfm}]$ and we should find $\psi \in \bbP_{\bfm}^{\pigyon}$ above it,  this suffice. Clearly there is $p_{1} \in \bbP_{\bfm}$ such that $\varphi \leq p_{1},$ that is $p_{1} \Vdash$``$\varphi \in \name{\bfG}_{{\bbP_{\bfm}[L_{\bfm}]}}$''. Now let $\langle \psi_{i}: i < i_{*} \rangle$ be a maximal anti-chain of members of $\bbP_{\bfn}[M_{\bfn}]$ which are incompatible with $p_{1}$ in $\bbP_{\bfm}[L_{\bfm}].$ Clearly $i_{*} < \lambda^{+}$ hence without loss of generality $i_{*} \leq \lambda$ and let $\psi = \bigwedge_{i < i_{*}} \neg \psi_{i}.$  Clearly $p_{1}$ witnesses $\psi \in \bbP_{\bfm}^{\pigyon}$ hence $\varphi \leq \psi,$ see more details in the proof of \ref{e52}(3). 
    
    2) Trivially $\bbP_{\bfn}^{\dagger} \subseteq \bbP_{\bfn}[M_{\bfn}],$ so it suffice assume $\psi \in \bbP_{\bfn}[M_{\bfn}]$ and prove $\psi \in \bbP_{\bfn}^{\dagger} \Leftrightarrow \psi \in \bbP_{\bfm}^{\dagger}.$
    
    First assume $\psi \in \bbP_{\bfn}^{\dagger}$ is witnessed by $p \in \bbP_{\bfn}$ and we shall prove that $p$ witness $\psi \in \bbP_{\bfm}^{\dagger};$ we have to check the two conditions $\bullet_{1} + \bullet_{2}$ of Definition \ref{e51}. Now clearly $p \in \bbP_{\bfm}$ and $\psi \in \bbP_{\bfm}[M_{\bfm}]$ (the second because $\psi \in \bbP_{\bfm}[M_{\bfm}]$ and $\bbP_{\bfn} \lessdot \bbP_{\bfm}$ and $M_{\bfn} \subseteq M_{\bfm},$ hence $\bbP_{\bfm}[M_{\bfn}] \lessdot \bbP_{\bfm}[M_{\bfm}]$). Also $\bbP_{\bfn}[L_{\bfn}] \models \psi \leq p$ but $\bbP_{\bfn}[L_{\bfn}] \lessdot \bbP_{\bfm}[L_{\bfm}]$ hence $\bbP_{\bfm}[L_{\bfm}] \models \psi \leq p.$ So in Definition \ref{e51} condition $\bullet_{1}$ holds; for proving condition $\bullet_{2},$ assume that $\bbP_{\bfm}[M_{\bfm}] \models$``$\psi \leq \varphi$'' hence, by part (1), we can find $q \in \bbP_{\bfm}$ and $\vartheta \in \bbP_{\bfm}[M_{\bfm}]$ which is witnessed by $q$ such that $\vartheta$ is above $\varphi.$ Without loss of generality, $\dom(q) \cap  \dom(p) \subseteq M_{\bfm}$ and let $q_{1} = q \rest L,$ now $q_{1}, \psi$ are compatible in $\bbP_{\bfm}[L_{\bfm}],$ hence in $\bbP_{\bfn}[L_{\bfn}],$ also $\dom(q_{1}) \cap \dom(p)$ is included in $L = L_{\bfn}$ and is included in $\dom(q) \cap \dom(p)$ which is included in $M_{\bfm};$ together $\dom(q_{1}) \cap \dom(p) \subseteq L_{\bfn} \cap M_{\bfm} = M_{\bfn}.$  Therefore by \ref{e29}(1), $p, q_{1}$ are compatible in $\bbP_{\bfm}(L_{\bfm}),$ hence in $\bbP_{\bfn}(L_{\bfn})$, so let $r \in \bbP_{\bfn}$ be a common upper bound. As $q \rest L \leq_{\bbP_{\bfn}} r,$ clearly $r, q$ has a common upper bound $r_{1}$ (in $\bbP_{\bfm}$) and so $r_{1}$ is a common upper bound of $\varphi, p.$
    
    So we are done proving one implication (the ``first'' above)  and the second is easier: if $p \in \bbP_{\bfm}$ witness $\psi \in \bbP_{\bfm}^{\dagger},$ then $p \rest L$ witness $\psi \in \bbP_{\bfn}^{\dagger}.$ 
    
    3) Why?
    
    \underline{Clause (a)}:
    
    The first clause is obvious, the second recalling  $\bbP_{\bfn} \lessdot \bbP_{\bfm}$ it is clear.

    \underline{Clause (b)}:
    
    In $\bbP_{\bfn}[M_{\bfn}]$ let $\bar{\psi} = \langle \psi_{i}: i < i_{*} \rangle, \psi$ be as in the proof of part (1) for $p_{1}$ (and $\bfn$). 
    
    Now, 
    
    \begin{itemize}
        \item[$(*)$] $\bar{\psi}$ is maximal also  for $\bfm, p_{1}.$ 
    \end{itemize}
    
    Why $(*)$ holds? It means: if $\vartheta \in \bbP_{\bfm}[M_{\bfm}]$ is incompatible with $p_{1}$ in $\bbP_{\bfm}[L_{\bfm}]$ then $\vartheta$  is compatible with some $\psi_{i}$ ($i < i_{*}$) in $\bbP_{\bfm}[L_{\bfm}].$ But if $\vartheta$ is a counterexample then there is $p_{2} \in \bbP_{\bfm}(L_{\bfm})$ above $\vartheta,$ so $p_{2}$ is incompatible with $\psi_{i}$ for $i < i_{*}$ and with $p_{1}.$ Let $q_{2} \in \bbP_{\bfm}(L_{\bfn})$ be above $p_{2} \rest L_{\bfn}$ and decide ``does $\psi \in \name{\bfG}_{\bbP_{\bfm}[L_{\bfm}]}$''.  As $q_{2}, p_{2}$ are compatible necessarily $q_{2} \Vdash$``$\psi \in \name{\bfG}_{\bbP_{\bfm}[L_{\bfn}]}$'' hence $q_{2} \Vdash_{\bbP_{\bfn}[L_{\bfn}]}$``$\psi \in \name{\bfG}_{\bbP_{\bfn}[L_{\bfn}]}$''. Let $\vartheta_{*} \in \bbP_{\bfn}[M_{\bfn}] $ be witnessed by $ q_{2}$ (exists by part (1)) so $\bbP_{\bfn}[L_{\bfn}] \models$``$\psi \leq \vartheta_{*}$''. Also without loss of generality $\dom(q_{2}) \cap \dom(p_{1}) \subseteq M_{\bfn}$ (by \ref{e27}) and $p_{1}, q_{2}$ are incompatible in $\bbP_{\bfn}[L_{\bfn}],$ (otherwise $p_{1}, q_{1}$ would be compatible). 
    
    So by \ref{e27}, $\vartheta_{*}, p_{1}$ are incompatible in $\bbP_{\bfn}[L_{\bfn}]$ so $\vartheta_{*}$ contradicts the maximality of $\bar{\psi}.$ 
\end{PROOF} 

\begin{definition}\label{e53}
    1) Let $ \mathbf{R} $ be the class of objects 
    $ \mathbf{r} $ consisting of  (so $ N=N_\mathbf{r} ,\mathbf{m} = \mathbf{m} _ \mathbf{r} $, 
    but we may omit the subscript $ \mathbf{r} $ when its identity is clear from the context, also in other parts):
    
    \begin{enumerate}
        \item[(a)]  $\bfm \in \bfM$  which is $ \lambda^+ _2 $-wide    (actually $ \lambda_1^+$  suffices),  
        
        \item[(b)] a cardinal $ \chi $ such that 
        $\bfm \in \cH(\chi)$ and $2^{|L_{\bfm}|+\lambda_2} < \chi$, 
        
        \item[(c)]  $N \prec (\cH(\chi),\in)$ such that $ \mathbf{m} \in N $,  and $N \cap \Ord = N \cap \chi$ has order type $\chi_{\bfr}$ (a cardinal $< \lambda$),
        
        \item[(d)]  $N \cap \lambda$ is an inaccessible
        cardinal $< \lambda$ called $\lambda _\mathbf{r} = \lambda ( \mathbf{r} ) =\lambda_N = \lambda(N)$, 
        
        \item[(e)]  $\|N\| < \theta_{\lambda(\mathbf{r} )}$  and\footnote{we shall use just $ \mathbb{P} _\mathbf{m} [M_ \mathbf{m}]$ has cardinality $\le \lambda_1$ because $ \lambda _1 =  \lambda_1 ^{< \lambda_{0}} $ in the proof $ (*)_3$ in \ref{c41}(1).} $ [N]^{< \lambda (\mathbf{r} ) } \subseteq N$,
        
        \item[(f)] $M_{\bfm}$ is listed in non-decreasing order $\bar{s}_{\bfr} = \langle s_i = s(i): i < i(\bfm) = i_{\bfm}  \rangle$ and let $s_{i(\mathbf{m})}$ be $\infty(\in L^+_{\bfm})$,  (so $ s_{\mathbf{r} , i}  = s_{\bfm, i} = s_i $); let $\bfU_{\bfr} := \{ j < i_{\bfm}: s_{j} \in N \}$ and $\bfU_{\bfr, i} := \{ j < i: s_{j} \in N \},$ and $\bfU_{\bfr}^{+} = \bfU_{\bfr} \cup \{ i_{\bfm} \},$  
        
        \item[(g)]  for $ i \in \bfU_{\bfr}^{+} $ let $L_{\mathbf{r}, i } = \cup \{ L_{\mathbf{m} (\leq s_j)}: j < i \} \cap N,$ and $L_{\bfr} = L_{{\bfr, i(\bfm)}} \subseteq  N$ so if $s_{i}$ is $<_{\bfm}$-increasing, then $ i = j +1 \Rightarrow L_{\bfr, i} = L_{\bfm(\leq s_{j})} \cap N,$ 
         
        \item[(h)]$  \Xi _ \mathbf{r}^{+} \not= \emptyset $,  see  (2B) below. 
    \end{enumerate} 
    
    2) For $ \mathbf{r} \in \mathbf{R}$ and $ i \in \bfU_{\bfr}^{+}$ let $  \Xi _{i}^{\dagger}= \Xi _{\mathbf{r}, i }^{\dagger}$ be\footnote{Justified when $ \mathbf{r} $ is clear from the context.} the set of sequences $ \bar{ \nu }$ such that:  
 
    \begin{enumerate} 
        \item[(a)] $ \bar{ \nu } = \langle \nu _j: 
         j \in \bfU_{\bfr, i} \rangle $, 
         
        \item[(b)] $ \nu _j \in \Pi _{  \varepsilon < \lambda (\mathbf{r} )} \theta _ \varepsilon$, 
     
        \item[(c)] there is  $ \mathbf{G} $ 
         weakly witnessing $\bar{\nu}$  which means: 
        
        \begin{enumerate} 
            \item[$(\alpha)$] $\mathbf{G} \subseteq N \cap  \mathbb{P}^ \dagger_{\mathbf{m}  \rest L_{\bfr, i}}$  is generic over $ N $; 
    
            \item[$(\beta)$] if $ j \in \bfU_{\bfr, i} $ \then \,  $ \nu_j = \name{ \eta} _ {s(j)} [ \mathbf{G}]$, that is  for every $ \xi < \lambda _\mathbf{r},$ for some $ \psi  \in \mathbf{G}$ we have $\psi \Vdash _ { \mathbb{P} _{\mathbf{m} (< s(i) )  }[M_{\mathbf{m}(<  s(i))}]}``  (\name{\eta}_{s(j)} \rest \xi ) = (\nu _j \rest \xi) "$.  
        \end{enumerate} 
    \end{enumerate} 
    
    Note that, if $M_{\bfm} \cap N$ is not linearly ordered, then maybe $j < i$ and $s(j) \notin L_{\bfm(< s(i))}$ but $s(j) \in L_{\bfr, i}$ so these two may not coincide. 

    2A) We have:
    
    \begin{enumerate} 
        \item[(a)] For $\bfr \in \mathbf{R}, i \in  \bfU_{\bfr}^{+}$ and $\bar{\nu} \in \Xi_{\bfr, i}^{\dagger},$ let $\mathbf{G} ^\dagger_{ \bar{ \nu }} = \mathbf{G} ^ \dagger _{\mathbf{r}, \bar{ \nu }}$  weakly witness   $\bar{ \nu } $, see  \ref{e53}(2)(c) above, so (by \ref{e61} below)   uniquely determined (by $ \bar{ \nu } $ and $ \mathbf{r} $), unlike in (2B) below.

        \item[(b)] let $\Xi ^\pigyon = \Xi ^\pigyon_{\mathbf{r}}$  be $ \Xi ^ \pigyon_{ \mathbf{r}, i(\mathbf{m} )}$\footnote{Justified when $ \mathbf{r} $ is clear from the context.}.
    \end{enumerate} 
  
    2B) For $ \mathbf{r} \in \mathbf{R}$ and $ i \in \bfU_{\bfr}^{+}$ let $ \Xi^+  _i = \Xi^+ _{\mathbf{r}  ,i }$ be\footnote{Justified when $ \mathbf{r} $ is clear from the context.} the set of sequences $ \bar{ \nu }$ such that: 
 
    \begin{enumerate} 
        \item[(a)] $ \bar{ \nu } = \langle \nu _j: j < i \rangle $, 
 
        \item[(b)] $ \nu _j \in \Pi _{  \varepsilon < \lambda (\mathbf{r})} \theta _ \varepsilon$, 

        \item[(c)] there is  $ \mathbf{G} $
        strongly witnessing $\bar{\nu}$  which means:
 
        \begin{enumerate} 
            \item[$(\alpha)$] $ \mathbf{G} \subseteq N \cap  \mathbb{P} _{\mathbf{m} }( L_{\mathbf{r}, i })$ is generic over $ N $; (\underline{but} $ \mathbf{G}  \subseteq \mathbb{P} _ \mathbf{m} [M_\mathbf{m} ]$ is not sufficient), 
  
            \item[$(\beta)$]  if $ j < i $ \then \,  $ \nu _ j = \name{ \eta} _ {s(j)} [ \mathbf{G}], $ that is  for every $ \xi < \lambda _\mathbf{r},$ for some $ \psi  \in \mathbf{G}$ we have $\psi \Vdash_{ \mathbb{P} _{\mathbf{m}(L_{i})  }[M_{\mathbf{m}(L_{i})}]} ``  (\name{\eta}_{s(j)} \rest \xi ) = (\nu _j \rest \xi)",$
            
            \item[$(\gamma)$]  $ N[{\mathbf{G} }]$ is isomorphic to ${\mathscr H} (\chi ')$ for some $ \chi '$  in fact  $\chi' = \chi _ \mathbf{r} $;  see clause (c) of part (1).
        \end{enumerate} 
    \end{enumerate} 
 
  2C) We have:   
 
    \begin{enumerate} 
        \item[(a)] for $\bar{\nu} \in \Xi_{\bfr}^{+},$  let $  \mathbf{G}^+  _ { \bar{ \nu }} = \mathbf{G}^+ _{ \mathbf{r}, \bar{ \nu } }$  strongly  witness $ \bar{ \nu } $, see \ref{e53}(2B)(c) above, so not necessarily uniquely determined, 
    
        \item[(b)] let $ \Xi^+ = \Xi ^+ _ \mathbf{r} $  be $\Xi^+_{\mathbf{s}, i(\mathbf{r} ) }$\footnote{Justified when $ \mathbf{r} $ is clear from the context.}.
    \end{enumerate} 
\end{definition} 
 
\begin{remark}\label{e59}  
    Assume $ \mathbf{m} $ is $ \lambda $-wide.  The following Claim \ref{e61}  justifies \ref{e53}(2A)(a).
\end{remark} 
  
\begin{claim} \label{e61} 
    Let $ \mathbf{r} \in \mathbf{R} $ and $ i \in \bfU_{\bfr}^{+} $, and  $ N= N_ \mathbf{r}.$
    
    1) If $ \mathbf{G} \subseteq  \mathbb{P} _{\mathbf{m} }[ L_{\mathbf{r}, i }  \cap M_ \mathbf{m}] \cap N$ is generic over $N$ \underline{then} there is one and only one $\bar{\nu}  \in {}^{ i}(\Pi _{\varepsilon < \lambda _ \mathbf{r}} \theta _ \varepsilon  ) $ such that  for every $ j \in \bfU_{\bfr, i} $ we have  $\nu _j = \cup \{ \varrho:$  there is $\psi \in \mathbf{G} $ satisfying  $\psi$ forces $ \varrho \trianglelefteq \name{ \eta }_{s(j)}\} $.  
    
    1A) We can use $ \mathbb{P} ^ \dagger _ \mathbf{m} $ instead of $ \mathbb{P} _ \mathbf{m} [M_ \mathbf{m} ]$. 
  
    2) If $ \mathbf{G} _1, \mathbf{G} _2 \subseteq 
    \mathbb{P} _{\mathbf{m} }[ L_{\mathbf{r}, i }  \cap M_ \mathbf{m}] \cap N $ are generic over $ N$,  $i \in \bfU_{\bfr}^{+}$ and $ \bar{ \nu } = \langle \nu  _j: j < i \rangle $,  and the pair $(\mathbf{G} _ {\ell} , \bar{\nu })$ is as above  for $ {\ell} = 1,2$ \underline{then} $\mathbf{G} _1 = \mathbf{G} _2,$ (not essentially used). 
   
    3) In part (1), we have $\mathbf{G} \cap \bbP_{\bfm}^{\pigyon} = \mathbf{G}^\pigyon   _{\langle \nu _j:j \in \bfU_{\bfr, i} \rangle} = \mathbf{G}^\pigyon_{\mathbf{r} , \langle \nu _j : j \in \bfU_{\bfr, i} \rangle }$, see \ref{e53}(2A)(a).

    4) Similarly for $ \mathbb{P}_\mathbf{m}[L_ \mathbf{m} ], \langle \name{ \eta } _s : s \in L_ \mathbf{m} \cap N \rangle$  instead $ \mathbb{P}_\mathbf{m} [M_ \mathbf{m} ], \langle \name{ \eta } _ s : s \in M_ \mathbf{m} \cap N\rangle$.
\end{claim} 
   
\begin{PROOF}{\ref{e61}}  
    1) For $  \psi \in \mathbb{P} _{\mathbf{m} }[ L_{\mathbf{r}, i }\cap M_ \mathbf{m}] \cap N$
    and $  j < i $ let $ \varrho _{\psi , j } $ 
    be the $ \trianglelefteq  $-maximal $  \varrho $ 
    such that $\psi \Vdash _{\mathbb{P} _ \mathbf{ m} [M_{\mathbf{m} (< s_i)}]}``\varrho \trianglelefteq \name{ \eta }_{s(j)}".$  
    
    Clearly, 
   
    \begin{enumerate}
        \item[$(*)_1$] for $ \psi, j $  as above, $ \varrho _{\psi, j}$is well defined and belongs to  $ \cup \{ \Pi _{\varepsilon < \xi } \theta _\varepsilon : \xi < \lambda _ \mathbf{r} \}.$  
    \end{enumerate}
   
    [Easy, e.g. why $ \lg ( \varrho _{\psi, j})< \lambda _ \mathbf{r} $? because $ \Vdash  ``  \name{ \eta } _{s(j)} \notin  \mathbf{V}"$ and $ \psi \in N_\mathbf{r} $.]
   
    \begin{enumerate}
        \item[$(*)_2$] for $ j < i $ and $ \xi < \lambda _ \mathbf{r}$ for some $ \psi \in \mathbf{G}_{\bar{ \nu }}$ we have $ \lg ( \nu _{\psi, j}) \ge  \xi.$ 
    \end{enumerate}
  
    [Why? by genericity and the definition of $\mathbb{P} _{\mathbf{m} }[ L_{\mathbf{r}, i } \cap M_ \mathbf{m}]$]. 
      
    \begin{enumerate}
        \item[$(*)_3$] if $ j < i $ and $ \psi _1 \le \psi_2$ are from $\mathbb{P} _{\mathbf{m} }[L_{\mathbf{r}, i } \cap M_ \mathbf{m}]$  then  $ \varrho _{\psi _1 , j} \trianglelefteq  \varrho _{\psi _{2, j}}.$ 
    \end{enumerate}
       
    [Why? Obvious].
       
    \begin{enumerate}
        \item[$ (*)_4$] if $ j < i $ and  $ \psi _1, \psi _2 \in \mathbf{G} $ where  $ \mathbf{G} $ is a subset of $ \mathbb{P} _{\mathbf{m} }[L_{\mathbf{r}, i }\cap M_ \mathbf{m}]$ generic over $ N $  then $ \varrho _{\psi_1 ,j}, \varrho _{\psi _2 , j}$ are $\trianglelefteq $-comparable. 
    \end{enumerate}

    [Why? as $ \mathbf{G} $ is directed and $ ( * ) _3 $].
   
    Together we are done proving part (1).
   
    1A) Easy. 
   
    2) Toward contradiction $\bfG_{1} \neq \bfG_{2}$ so we can assume $\psi_{1} \in \bfG_{1} \setminus \bfG_{2},$ hence there is $\psi_{2} \in \bfG_{2}$ which is incompatible with $\psi_{1}.$ Without loss of generality $\psi_{1}, \psi_{2} \in \bbP^{\dagger}[M_{\bfm}].$ So there is $p_{l} \in \bbP_{\bfm}$ witnessing $\psi_{l}$ (for $l = 1, 2$), and without loss of generality $\dom(p_{1}) \cap \dom(p_{2}) \subseteq M_{\bfm}.$ 

    Now by induction on $n$ we choose \ $(\psi_{1, n}, p_{1, n}, \psi_{2, n}, p_{2, n})$ such that: 
    
    \begin{itemize}
        \item[$(*)_{n}$] for $l = 1, 2$:
        
        \begin{enumerate}
            \item[(a)] $p_{l, n} \in \bbP_{\bfm} \cap N_{\bfr}$ and $m < n \Rightarrow p_{l} \leq p_{l, m} \leq p_{l, n},$
            
            \item[(b)] $\psi_{l, n} \in \bfG_{l}$ is witnessed by $p_{l, n},$
            
            \item[(c)] $m < n \Rightarrow \bbP[M_{\bfm}] \models \psi_{l, m} \leq \psi_{l, n},$
            
            \item[(d)] $\dom(p_{1, n}) \cap \dom(p_{2, n}) \subseteq M_{\bfm},$
            
            \item[(e)] if $n = m +1, s \in \dom(p_{1, m}) \cap \dom(p_{2, m})$ then $ \max \{ \rm{lg}(\tr(p_{1, m}(s))), \rm{lg}(\tr(p_{2, m}(s))) \} <  \min \{ \rm{lg}(\tr(p_{1, n}(s))), \rm{lg}(\tr(p_{2, n}(s))) \},$
            
            \item[(f)] if $s \in \dom(p_{l, n}) \cap M_{\bfm}$ then $\eta^{p_{2, n}(s)} \lhd \nu_{s}.$
        \end{enumerate}
    \end{itemize}

    Why it is enough to carry the induction? Because for $l = 1, 2$ we can let $p_{l}$ be  the lub of the increasing sequence $\langle p_{l, n}: n < \omega \rangle,$ and now $p_{1}, p_{2}$ are  compatible (as $s \in \dom(p_{1}) \cap \dom(p_{2})$) implies $s \in \dom(p_{1, n}) \cap \dom(p_{2, n}) \cap M_{\bfm}$ for some  $n \in \omega$ which implies $\tr(p_{1}(s)) = \tr(p_{2}(s)).$
    
    Now if $q$ is a common upper bound \ $p_{1}, p_{2}$ in $\bbP_{\bfm},$ then it is a common upper bound of $\psi_{1}, \psi_{2}$ in $\bbP_{\bfm}[L_{\bfm}],$ contradicting the choice of $\psi_{2}.$
    
    Why can we carry the induction?

    In the induction step we use having enough automorphisms and  (reflecting to $N_{\bfr}$).
    
    \begin{itemize}
        \item[$(*)$] if $q_{1} \in \bbP_{\bfm} \cap N_{\bfr}$  witnesses $\vartheta \in G_{l}$ and $\zeta < \lambda_{\bfr}$ then there are $q_{2} \in \bbP_{\bfm}$ and $\vartheta_{2} \in \bbP_{\bfm}^{\dagger}$ such that $\vartheta_{1} \leq \vartheta_{2}, q_{1} \leq q_{2}, q_{2}$ witnesses $\vartheta_{2}$ and $s \in \dom(q_{1}) \cap M_{\bfm} \Rightarrow \lg(\tr(q_{2})) \geq \zeta.$
    \end{itemize}

    [Why? let $\cI = \{ \varphi \in \bbP_{\bfm}^{\dagger}:  \text{\underline{either}} \ \varphi, \vartheta_{1} \ \text{are incompatible in} \ \bbP_{\bfm}[M_{\bfm}] \ \text{\underline{or}} \ \bbP_{\bfm}[M_{\bfm}] \models  \text{``}\psi_{2} \leq \varphi \text{'' and there is} \ q_{2} \in \bbP_{\bfm} \ \text{above} \ q_{1}, s \in \dom(q_{1}) \cap M_{\bfm} \Rightarrow \lg(\tr(q_{2}(s))) \geq \zeta \ \text{and} \ q_{2} \text{witnessing} \ \varphi \}$. By \ref{e29}(1), $\cI$ is a dense subset of $\bbP_{\bfm}[M_{\bfm}]$ and it belongs to $N,$ so necessary $\cI \cap \bfG_{l} \neq \emptyset$ and we can finish.]
    
    3) Follows.   
   
    4) Similarly. 
\end{PROOF} 

We may note:    
    
\begin{definition}\label{e62}
    Assume that $ \bar{p} = \langle  p_i: i < i_*
    \rangle $ where $ p_ i \in \mathbb{P} _ \mathbf{m} $ for $ i < i_*$  and $ i_{*} < \theta _0 $  (or just $ i_* < \lambda $  and $ i_* < \theta _{\lg(\tr( p_i(s )))}$ for every $ i < i_*, s \in \dom(p_i)$).
    
    We define $q = \oplus (\bar{p})$ as the following function  $q$:

    \begin{enumerate} 
        \item[(A)] $ q $ is a  function with domain  $\cup \{\dom(p_i): i < i_* \}$,  
        
        \item[(B)] if $ s \in \dom (q)$ then $ q(s ) $
        is  defined as in Definition \ref{c6}, as follows  (see (c) on $ j_s  $): 
        
        \begin{enumerate} 
            \item[(a)] $ \tr(q( s ) ) = \cup \langle \tr(p_i(s ): i < j_ s$ satisfies $ s \in \dom(p_i)\rangle $,  on $ j_s $ see below, 
            \item[(b)] for $ \varepsilon \in [\lg( \tr(q( s  )),  \lambda )$ we let $f_{q}( \varepsilon ) = \sup \{p_i(s ) (\varepsilon ): i < j_s $ satisfies $ s \in \dom(p_i)\rangle\} $; pedantically we consider each of the \lqq components" of $ f_ q$; see Definition. \ref{c6}; \underline{where}:   
            
            \item[(c)]  $j_s = \sup \{j: j \leq i_{*} $ and  the set $ \{\tr(p_i(s)) : i < j $ and $s \in\dom (p_i) \} $ is a set of pairwise $ \trianglelefteq $-comparable sequences$\}$. 
        \end{enumerate} 
    \end{enumerate} 
\end{definition} 

\begin{claim}{\label{e63}}
    1) If (A) then (B) where:
    
    \begin{enumerate}  
        \item[(A)] 
        
        \begin{enumerate} 
            \item[(a)] $ p_i \in \mathbb{P} _ \mathbf{m} $ for $ i < i_*$, 
            
            \item[(b)]  $ i_* < \theta _0 $ at least $ i_* < \theta_{\lg (\tr(p_i ( s )))}$ whenever $ i < i_*, s \in \dom(p_i)$. 
        \end{enumerate} 
    
        \item[(B)]
    
        \begin{enumerate} 
            \item[(a)]  $ q = \oplus (\bar{ p })$ is a member of $ \mathbb{P} _ \mathbf{m}$,
            
            \item[(b)] if $ r \in \mathbb{P} _ \mathbf{m} $ is a common upper bound of $ \{ p_i:  i < i_* \} $ then $q \le r \in \mathbb{P}_\mathbf{m}$,
            
            \item[(c)]  $ q $ is  a common upper bound of  $ \{ p_i:  i < i_* \} $  \when \, in addition to (A) and (B)(a):
        
            \begin{enumerate} 
                \item[(*)] if  $ i_1, i_2 < i_*, s \in \dom(p_{i_1}) \cap \dom(p_{i_2})$ and $ \lg (\tr(p_{i_1})) \le \epsilon <  \lg (\tr(p_{i_2}))$ \then \, for some $ i_3 < i_*$ we have $ p_{i_{\ell} }  \le p_{i_3}$ for $ {\ell} = 1,2$.
            \end{enumerate} 
        \end{enumerate} 
    \end{enumerate} 
    
    2)  If $(A)^+$ then $(B)^+$ where:
    
    \begin{enumerate} 
        \item[$(A)^+$]  as in (A) above adding:  
        
        \begin{enumerate} 
            \item[(c)]  $ p_i $ witnesses $ \psi _i \in \mathbb{P}^\dagger _\mathbf{m} [M_ \mathbf{m} ]$,
            
            \item[(d)]  $ \psi = \bigwedge _{i < i_*} \psi _ i \in \mathbb{P}^\dagger_\mathbf{m} [M_ \mathbf{m}]$, 
        \end{enumerate} 
    
        \item[$(B)^+$]  as in (B) above adding:
        
        \begin{enumerate} 
            \item[(d)] if $ q $  is a common upper bound  of $ \{ p_i:  i < i_* \}$ \then \, $q$ witnesses $ \psi   \in \mathbb{P}^  \dagger _\mathbf{m} [M_ \mathbf{m} ]$. 
        \end{enumerate} 
    \end{enumerate} 
\end{claim} 

\begin{PROOF}{\ref{e63}}
    1) Clearly $ q \in \mathbb{P} _ \mathbf{m}.$   Also   if $ r \in \mathbb{P} _ \mathbf{m} $ is a common upper bound of $ \{ p_i:  i < i_* \}$ clearly $q \le r  \in \mathbb{P} _ \mathbf{m} $.
     
    2) Easy and will not be used. 
\end{PROOF}      

\subsection{Nicely existentially closed}\label{4C}\ 

\begin{discussion}\label{g0}
    1) In the main case we have $M \subseteq M_{\bfm}$ cofinal in $M_{\bfm}$ and $\bfm \rest M \cong \bfm \rest M_{\bfm}.$ In \S4A we proved that if $\bfm \in M_{\rm{ec}}$ then $\bbP_{\bfm}[M] \cong \bbP_{\bfm}[M_{\bfm}]$ even in the general case.  
    
    Our main aim is to prove more; e.g.
    
    \begin{itemize}
        \item[$(*)$]
        
        \begin{enumerate}
            \item[(a)] for $\bfm, M$ as above,  there is $\bfn \in \bfM$ such that $\bfn = \bfm \rest L_{\bfn}, M_{\bfn} = M$ and $\bbP_{\bfn} \lessdot \bbP_{\bfm},$ 
            
             \item[(b)] moreover there is $\bfn$ such that $M_{\bfn} = M, L_{\bfn} \subseteq L_{\bfm}, \bbP_{\bfn}[M] = \bbP_{\bfm}[M]$ and $\bbP_{\bfm}[L_{\bfm}]$ is isomorphic to $\bbP_{\bfn}[L_{\bfn}]$ over $\bbP_{\bfm}[M].$
        \end{enumerate}
    \end{itemize}
    
    2) In the main case ($\bfm \in \bfM_{\ec}$ is strongly by $\lambda^{+}$-directed even in the general case), 
    
    \begin{itemize}
        \item[$\bullet$] $\Vdash_{\bbP_{\bfm}}$``$\{ \name{\eta}_{s}: s \in M \}$ is cofinal in $(\Pi_{\varp < \lambda} \theta_{\varp}, <_{J_{\lambda}^{\rm{bd}}})$''. 
    \end{itemize}
    
    From Definition \ref{g2} we shall use $L' = \bigcup \rm{cmp}(M, \bfm)$ defined below. 
    
    Earlier in \S3D we doctored $\bfm \in \bfM_{\leq \lambda_{2}}$ to an equivalent $\bfn$ such that $M_{\bfn} = M_{\bfm}, \, E''_{\bfn}$ has one equivalent class ``glueing'' together all $t / E_{\bfn}, t \in L_{\bfm} \setminus M_{\bfm}.$ Here things are more elaborated. First, in \ref{g2} we define the set $\rm{cmp}(M, \bfm)$ of the $t / E'_{m}$  for which $M$ is enough and then  in \ref{g5} doctor $\bfm$ to an equivalent $\bfn$ with $M_{\bfn} = M,$ but glueing \ together all $t / E'_{\bfn}$ not in $\rm{cmp}(M, \bfm).$ Later we can treat $\bfn$ as earlier. 
\end{discussion}
 
\begin{definition}\label{g2}
    Assume $\bfm \in \bfM$ and $M \subseteq M_{\bfm}.$
    
    1) Let $\rm{cmp}(M, \bfm)$ be the set of $L$ such that for some $t \in L_{\bfm} \setminus M_{\bfm},$ we have:
    
    \begin{enumerate}
        \item[(a)] $L = t / E_{\bfm}' = \{ s \in L_{\bfm}: s E_{\bfm}' t \},$
        
        \item[(b)] $L \cap M_{\bfm} \subseteq M.$
    \end{enumerate}

    2) Let $\rm{cmp}^{+}(M, \bfm)  = \rm{cmp}(M, \bfm) \cup \{ M \}.$
    
    3) For $t_{*} \in M$ let $\rm{cmp}(t_{*}, M, \bfm)$ be the set of $L \in \rm{cmp}(M, \bfm)$ such that $L \subseteq L_{\bfm(\leq t_{*})}$ and similarly  $\rm{cmp}^{+}(t_{*}, M, \bfm) = \cmp(t_{\ast}, M, \bfm) \cup (M_{\bfm(\leq t_{\ast})} \cap M).$
    
    
        
\end{definition}

\begin{definition}\label{g5}
    Assume $\bfm \in \bfM, M \subseteq M_{\bfm}$ and $\mathscr{L} \subseteq \rm{cmp}(M, \bfm).$

    1) We define $\rm{rest}(\mathscr{L}, M, \bfm)$ as  the object $\bfn$ consisting of (intended to be in $\bfM$):
    
    \begin{enumerate}
        \item[(a)] the set of elements of $L_{\bfn}$ is $\bigcup \{ L: L \in \mathscr{L} \} \cup M,$
        
        \item[(b)] the order on $L_{\bfm}$ is $<_{\bfm} \rest L_{\bfn},$
        
        \item[(c)] $M_{\bfn} = M,$
        
        \item[(d)] $E_{\bfn}' = E_{\bfm}' \rest L_{\bfn},$
        
        \item[(e)] $u_{\bfn, s} = u_{\bfm, s} \cap L_{\bfn},$
        
        \item[(f)] $\cP_{\bfn, s} = \cP_{\bfm, s} \cap [u_{\bfn, s}]^{\leq \lambda}.$
    \end{enumerate}

    2) Further assume $\bfm \in \bfM_{\leq \lambda_{2}}.$ We define $\rm{Rest}(\mathscr{L}, M, \bfm)$ as the object $\bfn$ consisting of (intended to be in $\bfM$):

    \begin{enumerate}
        \item[(a)] $L_{\bfn} = L_{\bfm},$
        
        \item[(b)] $M_{\bfn} = M,$
        
        \item[(c)] $E_{\bfn}' = \{ (s_{1}, s_{2}): \text{for some} \ L \in \mathscr{L}, (s_{1}, s_{2}) \in E_{\bfm}' \rest L \ \text{\underline{or}} \ (\forall L \in \mathscr{L})({s_{1}, s_{2} \notin L \setminus M}), \\ \{s_{1}, s_{2} \} \subseteq L_{\bfn}, \{ s_{1}, s_{2} \} \nsubseteq M_{\bfn} \},$
        
        \item[(d)] $u_{\bfn, s} = u_{\bfm, s}$ for $s \in L_{\bfn},$
        
        \item[(e)] $\cP_{\bfn, s} = \cP_{\bfm, s}.$
    \end{enumerate}

    3) We may omit $\mathscr{L}$ when $\mathscr{L} = \rm{cmp}(M, \bfm).$
    
    4) For $\bfr \in \bfR,$ let $\Xi_{\bfr}^{\bullet}$ be the set of $\bar{\nu}$ such that some $\bfG^{\bullet}$ witness it, which means:
    
    \begin{enumerate}
        \item[(a)] $\bfG^{\bullet} \subseteq \bbP_{\bfm} \cap N_{\bfr}$ is generic over $N_{\bfr},$
        
        \item[(b)] $\bar{\nu} = \langle \nu_{s}: s \in M_{\bfm} \cap N_{\bfr} \rangle,$
        
        \item[(c)] $\nu_{s} = \name{\eta}_{s}[\bfG^{\bullet}]$ for $s \in M_{\bfm} \cap N_{\bfr}.$  
    \end{enumerate}
    
    (So compared to Definition \ref{e53}(2B) clause (c)($\gamma$) is not required here). 
    
    5) For $\bfr \in \mathbf{R}, M' \subseteq M_{\bfm}$ such that $ M' \in N_{\bfr}$ and $ M = M' \cap N_{\bfr},$ let $\Xi_{M}^{\bullet} = \Xi_{\bfr, M}^{\bullet}$ be the set of $\bar{\nu} = \langle \nu_{s}: s \in M \rangle$ such that some pair $(\bfn, \bfG)$ witnesses $\bar{\nu}$ which means: 
    
    \begin{enumerate}
        \item[$(\alpha)$]  $\bfn \leq_{\bfM} \rm{rest}(M', \bfm)$ and $\bfn \in N_{\bfr}$ and $\bbP_{\bfn} \lessdot \bbP_{\bfm},$
        
        \item[$(\beta)$] $\bfG$ is a subset of $\bbP_{\bfn} \cap N_{\bfr}$ generic over $N_{\bfr},$
        
        \item[$(\gamma)$] $\nu_{s} = \name{\eta_{s}}[\bfG]$ for $s \in M.$
    \end{enumerate}
\end{definition}

\begin{claim}\label{g8}
    Assume $\bfm \in \bfM, M \subseteq M_{\bfm},$  $\mathscr{L} \subseteq \rm{cmp}(M, \bfm),$ ${\bfn_{2}} := \rm{Rest}(\mathscr{L}, M, \bfm)$ (see \ref{g5}(2)) and $\bfn_{1} = \rm{rest}(\mathscr{L}, M, \bfm),$ (see \ref{g5}(1)):
    
    1) $\bfn_{1} \leq_{\bfM} \bfn_{2},$ 
    
    2) $L_{\bfn_{2}} = L_{\bfm}, \bbP_{\bfn_{2}} = \bbP_{\bfm}$ and $\rm{cmp}(M, \bfn_{1}) = \mathscr{L} \subseteq \cmp(M, \bfn_{2})$ and $\cmp(M, \bfn_{2}) \setminus \cL$ is empty or a singleton; note that this is the single $E_{\bfn_{2}}'$-class that is not in $\cL_{\bfn_{1}}.$
    
    3) If $M = M_{\bfm}$ and $\cL = \cmp(M, \bfm)$ then $\bfn_{1} = \bfn_{2} = \bfm.$
     
    4) if $\iota \in \{1, 2 \}$ and $\bfn_{\iota} \leq_{\bfM} \bfn_{*}$ and $L_{\bfn_{*}} \cap L_{\bfm} = L_{\bfn_{\iota}},$ \underline{then} we can find $\bfm_{*}$ such that: 
    
    \begin{enumerate}
        \item[(a)] $\bfm \leq_{\bfM} \bfm_{*}$ and $L_{\bfm_{*}} = L_{\bfn_{*}} \cup L_{\bfm},$
        
        \item[(b)] $\mathscr{L} \subseteq \rm{cmp}(M, \bfm_{*}),$
        
        \item[(c)] if $\iota = 2$ then  $ \bfn_{*} = \rm{Rest}(\mathscr{L}, M, \bfm_{*}),$
        
        \item[(d)] if $\iota = 1 $ letting $\mathscr{L}_{1} = \rm{cmp}(M, \bfn_{*})$ we have $\bfn_{*} = \rm{rest}(\mathscr{L}_{1}, M, \bfm_{*}),$
        
        \item[(e)] (choosing minimal $u$) if $t \in M_{\bfm} \setminus M$ and $\iota = 1,$ then $u_{\bfm_{*}, t} = u_{\bfm, t} \cup \{ s: s \in L_{\bfn_{*}} \setminus L_{\bfm}  \ s \in u_{\bfn, t} \}.$ 
    \end{enumerate}
\end{claim}

\begin{PROOF}{\ref{g8}}
    Straightforward, as in earlier proofs (in particular \ref{b44}) in particular for part (1), check the clause (e)$(\gamma)$ of \ref{e4}. 
\end{PROOF}

\begin{claim}\label{g11}
    Assume $\bfm \in \bfM.$
    
    1) If $t \in M_{\bfm}$ and $\bfn = \bfm(\leq t)$ \underline{then} $t = \max(M_{\bfn}) = \max(L_{\bfn}).$
    
    2) If  for every $t \in M_{\bfm}, \rm{rest}(M_{\bfm(\leq t)}, \bfm) \in \bfM_{\rm{{ec}}}$ \underline{then} $\bfm \in \bfM_{\rm{ec}}$ provided that  is strongly $\lambda^{+}$-directed.
    
    3) If $\bfm_{1} \leq_{\bfM} \bfm_{2}$ and $ t \in M_{\bfm_{1}}$ then $\rm{rest}(M_{\bfm(\leq t), \bfm}) \leq_{\bfM} \rm{rest}(M_{\bfm(\leq t)}, \bfm).$
\end{claim}

\begin{PROOF}{\ref{g11}}
    Easy. 
\end{PROOF}

\begin{observation}\label{g13}
    1) If $L_{\bfm} = M_{\bfm}$ \underline{then} $\bfm$ is  essentially $\lambda^{+}$-directed (see Definition \ref{b36}(1)) \underline{if} $\bigcup \{ u_{s}: s \in M_{\bfm} \} = M_{\bfm} = M_{\bfm}$ and $(\{ u_{s}: s \in M_{\bfm} \}, \subseteq)$ is $\lambda^{+}-$directed.

    2) Assume $\bfm$ is strongly $\mu$-directed and $\bfm_{1} \leq \bfm,$ \underline{then} $\bfm_{1}$  is strongly $\mu$-directed.

    3) if $\bfm_{1} \leq_{\bfM} \bfm$ and $\bfm$ is essentially $\mu$-directed \underline{then} so is $\bfm_{1}.$ 
\end{observation}

\begin{claim}\label{g14}
    1) Assume $\bfm \in \bfM_{\leq \lambda_{2}}$ is strongly $\lambda^{+}$-directed (hence bounded). There are $\bfn$ and $\bar{\cL}$ such that: 
    
    \begin{itemize}
        \item[$\boxplus$] 
        
        \begin{enumerate}
            \item[(a)] $\bfm \leq_{\bfM} \bfn \in \bfM_{\lambda_{2}},$
            
            \item[(b)] $\bfn \in \bfM_{\rm{bec}},$
            
            \item[(c)] $\bar{\cL} = \langle \cL_{M}: M \in \cP^{-}(M_{\bfm}) \rangle,$ recalling that for a set $X, \, \cP(X)^{-} := \{ Y \subseteq X: Y \neq \emptyset \},$ 
            
            \item[(d)] $\cL_{M} \subseteq \rm{cmp}(M, \bfn)$ for $M \in \cP^{-}(M_{\bfm})$ and $\cL_{M_{\bfm}} = \cmp(M_{\bfm}, \bfm),$ 
            
            \item[(e)] $\bfn_{M} = \bfn[M] := \rm{rest}(\cL_{M}, M, \bfm) \in \bfM_{\rm{bec}}$ for $M \in \cP^{-}(M_{\bfm}),$ so $\bfn_{M_{\bfm}} = \bfm,$ 
            
            \item[(f)] $\bbP_{\bfn_{M}} \lessdot \bbP_{\bfn}$ for $M \in \cP^{-}(M_{\bfm}),$
            
            \item[(g)] if $t_{2} <_{\bfm} t_{3}$ are from $M_{\bfm}$ and $t_{1} \in L_{{\bfn}} \setminus L_{\bfm}$ and $(\forall s \in t_{1} / E_{\bfn}'')$ ($s <_{\bfm} t_{2}$) then $t_{1} / E''_{\bfn} \subseteq u_{{\bfn}, t_{3}}$ ({yes! not $t_{1} / E_{\bfn}'$})
            
            \item[(h)] $\bfn$ is strongly $\lambda^{+}$-directed (see \ref{b36}(2)),
            
            \item[(i)] if $M_{1} \subseteq M_{2}$ are from $\cP^{-}(M_{\bfm})$ then $\mathscr{L}_{M_{1}} \subseteq \mathscr{L}_{M_{2}}.$
        \end{enumerate}
    \end{itemize}
    
    2) We can add: 
    
    \begin{itemize}
        \item[] 
        
        \begin{enumerate}
            \item[(j)] if $M_{1}, M_{2} \subseteq M_{\bfm}$ and $h$ is an isomorphism from $\bfm \rest M_{1}$ onto $\bfm \rest M_{2}$ \underline{then} $h$ can be extended to $\hat{h},$ an isomorphism from $\bfn_{M_{1}}$ onto $\bfn_{M_{2}},$ 
            
            \item[(k)] if $M \in \cP^{-}(M_{\bfm})$ and $L \in \mathscr{L}_{M}$ \underline{then} $L \cap L_{\bfm} \subseteq M_{\bfm}.$
            
            \item[(l)] If $M \in \cP^{-}(M_{\bfm})$ and $h$ is an isomorphism from $\bfm \rest M$ onto $\bfm \rest M_{\bfm},$ \underline{then}:
            
            \begin{itemize}
                \item[$\bullet$] there is an isomorphism $\hat{h}$ from $\bfn_{M}$ onto $\bfn_{M_{\bfm}} = \bfm$ embedding $h,$
                
                \item[$\bullet$] if $\bfG^{+}$ is a subset of $\bbP_{\bfm}$ generic over $\bfV,$ then there is $\bfG \in \bfV[\bfG^{+}],$ a generic subset of $\bbP_{\bfm}$ over $\bfV$ such that $s \in M \Rightarrow \name{\eta}_{s}[\bfG^{+}] = \eta_{h(s)}[\bfG].$
            \end{itemize}
        \end{enumerate}
    \end{itemize}
    
    3) Above $\bfn$ is reasonable (see Definition \ref{e36}).
\end{claim}

\begin{PROOF}{\ref{g14}}
    1) Let $\langle M_{\alpha}: \alpha < \alpha_{*} < \lambda_{1}^{+} \rangle$ list $\cP^{-}(M_{\bfm})$ so $\alpha_{*} <   \lambda_{2}$ such that $t_{\alpha} <_{\bfm} t_{\beta} \Rightarrow \alpha < \beta.$ 
    
    We choose by induction on $\alpha \leq \alpha_{*}, \bfm_{\alpha}$ and if $\alpha < \alpha_{*},$ also $\bfn_{\alpha}^{0}, \bfn_{\alpha}^{1}, \cL_{\alpha}$ such that:    
    
    \begin{itemize}
        \item[$(*)_{0}$] 
        
        \begin{enumerate}
            \item[(a)] $\bfm_{\alpha} \in \bfM_{\leq \lambda_{2}},$
            
            \item[(b)] $\bfm_{\alpha}$ is $\leq_{\bfM}$-increasing continuous,
            
            \item[(c)] $\bfm_{0} = \bfm,$
            
            \item[(d)] if $\alpha = \beta + 1,$ then: 
            \begin{enumerate}
                \item[$(\alpha)$] $\bfn^{0}_{\beta} = \bfm \rest M_{\beta},$
                
                \item[$(\beta)$] $\bfn^{0}_{\beta} \leq_{\bfM} \bfn^{1}_{\beta} \in \bfM_{\rm{bec}},$
                
                \item[$(\gamma)$] $\bfn_{\beta}^{1} \in \bfM_{\leq \lambda_{2}}$ and $L_{\bfn_{\alpha}^{1}} \cap L_{\bfm_{\beta}} = L_{\bfn_{\beta}^{0}},$
                
                \item[$(\delta)$] $\bfn^{1}_{\beta} = \rm{rest}(\cL_{\beta}, M_{\beta}, \bfm_{\alpha}),$
                
                \item[$(\varp)$] if $t \in M_{{\bfm}} \setminus \cup \{ M_{\gamma}: \gamma \leq \beta \},$ then $u_{\bfm_{\alpha}, t} = u_{\bfm_{\beta}, t} \cup \{ s \in L_{\bfm_{\alpha}}: s \notin L_{\bfm_{\beta}}, s\rest E_{\bfm_{\alpha}}' \subseteq L_{\bfm_{\alpha(\leq t)}} \}.$
            \end{enumerate}
        \end{enumerate}
    \end{itemize}
    
    Why can we carry the induction?
    
    First, arriving to $\alpha$ we choose $\bfm_{\alpha}$ as follows: 
    
    \begin{itemize}
        \item[$\bullet$] if $\alpha = 0$ let $\bfm_{\alpha} = \bfm$ (so $(*)(c)$ holds).
        
        \item if $\alpha$ is a limit ordinal then $\bfm_{\alpha} = \cup \{ \bfm_{\beta}: \beta < \alpha \},$ see \S1A, noting $\bfm_{\alpha} \in \bfM_{\leq \lambda_{2}}$ because $\alpha \leq \alpha_{*} < \lambda_{2}.$
        
        \item  if $\alpha = \beta + 1$ so $\bfn_{\beta}^{0}, \bfn_{\beta}^{1}$ have been chosen then choose $\bfm_{\alpha}$ by \ref{g8}(4) and $(*)(d)(\varp).$
    \end{itemize}
    
    Second assuming $\alpha < \alpha_{*},$ and $\bfm_{\alpha}$ has been defined we choose $\bfn_{\alpha}^{0}$ as $\rm{rest}(M_{\alpha}, \bfm_{\alpha})$ so $\bfn_{\alpha}^{0} \in \bfM_{\leq_{\lambda_{2}}}$ by \ref{g8}(1), so $(*)(d)(\alpha)$ holds.  Then choose $\bfn_{\alpha}^{1} \in M_{\leq \lambda_{2}}$ such that $\bfn_{\alpha}^{0} \leq_{\bfM} \bfn_{\alpha}^{1} \in \bfM_{\rm{ec}}$ (by claim \ref{e37}), without loss of generality $L_{\bfn_{\alpha}^{1}} \cap L_{\bfm_{\alpha}} = L_{\bfn_{\alpha}^{0}}$ so clause $(*)(d)(\beta), (\gamma)$ holds. 
    
    So we have carried the induction. 
    
    \begin{itemize}
        \item[$(*)_{1}$] Let $\bfn = \bfm_{\alpha_{*}}$ so $\bfn \in \bfM_{\leq\lambda_{2}}$ and $\bfm = \bfn_{0} \leq_{\bfM} \bfn.$ 
    \end{itemize}
    
    So clause $\boxplus(a)$ holds. Why $\bfn \in \bfM_{\rm{bec}}?$ (i.e.  clause $\boxplus(b)$). As clearly $M_{\bfn}$ is strongly $(< \lambda^{+})$-directed, by \ref{b47}(5) it  suffices to prove that $M_{\bfn(\leq t)} \in \bfM_{\rm{bec}}$ for every $t \in M_{\bfm}.$ But $M_{\bfn(\leq t)} \in \{ M_{\alpha}: \alpha < \alpha_{*} \}$ and if $M = M_{\bfm({\leq t_{\alpha}})}$ then $\rm{rest}(M, \bfn) \leq_{\bfM} \bfm(\leq t_{\alpha})$ so this follows if we prove clauses $\boxplus(d), (e).$
    
    Why clauses $\boxplus(d), (e)$ holds? So assume $M \in \cP^{-}(M_{\bfm})$ then for some $\alpha = \alpha(M) < \alpha_{*}$ we have $M = M_{\alpha}$ and let $\mathscr{L}_{M} = \rm{cmp}(M, \bfn_{\alpha}^{1})$ (so $\bigcup \{ L: L \in \mathscr{L}_{M} \} = L_{\bfn_{\alpha}^{1}}).$ Assume $\rm{rest}(\cL_{M}, M, \bfn) = \bfn_{0} \leq_{\bfM} \bfn_{1} \leq_{\bfM} \bfn_{2}$ and we shall prove that $\bbP_{\bfn_{1}} \lessdot \bbP_{\bfn_{2}},$ this suffices for $\boxplus(e).$ 
    
    Let $\alpha < \alpha_{*}$ be such that $M_{\alpha} = M,$ so clearly. 
    
    \begin{itemize}
        \item[$(*)_{2}$] 
        
        \begin{enumerate}
            \item[$\bullet_{1}$] $ \rm{rest}(\mathscr{L}_{M}, M, \bfn) = \bfn_{\alpha}^{1},$
            
            \item[$\bullet_{2}$] $\bfn_{\alpha}^{1} = \rm{rest}(M_{\alpha}, \bfm_{\alpha}),$
            
            \item[$\bullet_{3}$] $\rm{rest}(M_{\alpha}, \bfm_{\alpha}) \leq_{\bfM} \rm{rest}(M_{\alpha}, \bfn),$
            
            \item[$\bullet_{4}$] $\rm{rest}(M_{\alpha}, \bfn) \leq_{\bfM} \rm{Rest}(M_{\alpha}, \bfn) \leq_{\bfM} \rm{Rest}(M_{\alpha}, \bfn_{1}) \leq_{\bfM} \rm{Rest}(M_{\alpha}, \bfn_{2}).$
        \end{enumerate}
    \end{itemize}
    
    As $\bfn_{\alpha}^{1} \in \bfM_{\rm{bec}}$ it follows that $\rm{Rest}(M_{\alpha}, \bfn_{1}) \lessdot \rm{Rest}(M_{\alpha}, \bfn_{2}),$ but by \ref{g8}(1),
    
    \begin{itemize}
        \item[$(*)_{3}$]  $\bbP_{\bfn_{l}} = \bbP_{\rm{Rest}(M_{\alpha}, \bfn_{l})}$ for $l = 1, 2.$
    \end{itemize}
    
    So we are done proving $\boxplus(e).$
    
    Now $\boxplus(c)$ is just a choice, so we are left with $\boxplus(f)$ which says $\bbP_{\bfn_{M}} \lessdot \bbP_{\bfn}$ but it follows by $(*)_{2}$ and $(*)_{3}.$
    
    \noindent 
    2) We can find $\bfn^{+}$ such that: 
    
    \begin{itemize}
        \item[$(*)$]
        
        \begin{enumerate}
            \item[(a)] $\bfn \leq_{\bfM} \bfn^{+} \in \bfM_{\leq \lambda_{2}},$ ($\bfn$ is from part (1)),
            
            \item[(b)] if $M_{1}, M_{2} \in \cP^{-}(M_{\bfm})$ and $h$ is an isomorphism from $\bfm \rest M_{1}$ onto $\bfm \rest M_{2}$ and $L \in \rm{cmp}(M_{1}, \bfn)$ \underline{then} there is $\langle (L_{i}, {h_{i}}): i < \lambda_{2} \rangle$ such that: 
            \begin{enumerate}
                \item[$(\alpha)$] $L_{i} \in \rm{cmp}(M_{2}, \bfm),$
                
                \item[$(\beta)$] $h_{i}$ is an  isomorphism from $\bfm \rest L$ onto $\bfm \rest L_{i}$ which extends $h \rest (L \cap M_{i}),$
                
                \item[$(\gamma)$] $L_{i} \cap  L_{\bfn} \subseteq M_{\bfm}.$
            \end{enumerate}
            
            \item[(c)] if $t \in L_{\bfn}(t) \setminus L_{\bfn}$ \underline{then} for some $s \in L_{\bfn} \setminus M_{\bfm}$ we have $ \bfn \rest (t / E_{\bfn(t), \ast}), \bfn \rest (s / E_{\bfn, t})$ are isomorphic, 
            
            \item[(d)] if $s \in L_{\bfn} \setminus M_{\bfn}, M_{1} = s / E', M_{2} \in \cP^{-}(M_{\bfm})$ and $h$ is an isomorphism from $\bfm \rest M_{1}$ onto $\bfm \rest M_{2},$ \underline{then} there is a sequence $\langle t_{\varp}: \varp < \lambda_{2} \rangle$ of pairwise non-$E_{\bfn(t)}''$-equivalent  members of $L_{\bfm} \setminus M_{\bfm}$ such that $(t_{\zeta} / E_{\bfm}') \cap M_{\bfm} = M_{2}$ and there is an isomorphism from $\bfn \rest (s / E_{\bfn}')$ onto $\bfn^{\dagger} \rest (t_{\zeta} / E_{\bfn(+)}')$  
        \end{enumerate}
    \end{itemize}
    
    Now we can easily find the isomorphism promised in clause (j). Lastly, clause (k) holds because $\bfn_{\alpha}^{0} = \bfm \rest M_{\alpha}$ above. So $\bfn^{+}$ is as required recalling $t \in L_{\bfn^{+}} \setminus M_{\bfn^{+}} \Rightarrow t / E_{\bfn^{+}}' = \cup \{ L: L \in \cmp(M_{\bfm}, \bfm) \}$ as $\bfm$ is bounded.  
    
    3) Easy. 
\end{PROOF}

\begin{remark}\label{g15}
    How does this subsection help \cite{Sh:945}?

    1) Note in the family $\mathbf{R}$ of $\bfr'$s see Definition \ref{e53} we demand that there is $\bfG^{+} \subseteq N_{\bfr} \cap \bbP_{\bfm}$ generic over $N$ such that: (so $N \prec (\cH(\chi), \in), (\bfn, \lambda, \cdots) \in N$ and $\mathbf{j}_{N}$ is the Mostowski collapse)
    
    \begin{itemize}
        \item[$(*)$] $\mathbf{j}''_{N}(N)[\bfG^{+}] = \mathscr{H}(\chi_{\varp}), \chi_{\varp} = \otp(N \cap \chi).$ 
    \end{itemize}
    
    You can think of it as: in the preliminary forcing to get Laver diamond, in stage $\lambda_{N} = N \cap \lambda$ we force by $\mathbf{j}_{\varp} (\bbP_{\bfm} \cap N).$
    
    2) The present \ref{e70} tells us to use $\Xi_{\nu}^{\bullet}$ (defined in \ref{g5}(4)) that instead of using  $\Xi_{\nu}^{\dagger} = \{ \bfG: \bfG \subseteq N_{s} \cap \bbP_{\bfm}^{\dagger} \ \text{is generic over} \ N \ \text{such that} \ \nu_{\alpha} = \name{\eta}_{\alpha}[\bfG] \ \text{for} \ \alpha \in M \}$ which gives too may candidates \underline{or} $\Xi_{\bar{\nu}}^{+} = \{ \bfG: \bfG \subseteq N \cap \bbP_{\bfm} \ \text{is generic over} \ N \ \text{such that} \ \mathbf{j}''(N) [\mathbf{j}_{N}''(\bfG)] = \cH(\chi_{N}) \ \text{and} \ \nu_{\alpha} = \name{\eta}_{\alpha}[G]\}$ \color{black} which seems too restrictive. 
    
    Enough to use the middle ground $\Xi_{\bar{\nu}}^{\bullet} = \{ \bfG: \bfG \subseteq N \cap \bbP_{\bfm} \ \text{is generic over} \ N \ \text{and} \ \nu_{\alpha} = \name{\eta}_{\alpha}[\bfG] \ \text{for} \ \alpha \in M \}.$ 
    
    3) Now the original idea was that $\bfG \in \Xi_{\bar{\nu}}^{\dagger}$ is enough in \cite{Sh:945} but not so, however $\bfG \in \Xi_{\bar{\nu}}^{\bullet}$ is sufficient.
    
    4) Also we need that $\bfm$ is reasonable (see \ref{e36}(3)) so if $M \subseteq M_{\bfm}$  is cofinal then $\langle \name{\nu}_{\alpha}: \alpha \in M \rangle$ is cofinal for $\bfm.$  
    
    5) The point is that for $M \nsubseteq M_{\bfm}$ (or with $M \cap N, M_{\bfm} \cap N $ the reflection) we need stronger homogeneity of $\bbP_{\bfm},$ which is the aim of \ref{g0}-\ref{g17} relying on \ref{g14}. 
\end{remark}

\begin{conclusion}\label{g17}
    If $\bfm_{0} \in \bfM$ is strongly $\lambda^{+}$-directed, (so bounded) of cardinality $\leq \lambda_{2},$ \underline{then} there is $\bfm$ such that: 
    
    \begin{enumerate}
        \item[(a)] $\bfm \in \bfM$ of cardinality $\lambda_{2},$
        
        \item[(b)] $\bfm_{0} \leq_{\bfM} \bfm,$
        
        \item[(c)] $\bfm \in \bfM_{\rm{bec}},$
        
        \item[(d)] $(\alpha) \Rightarrow (\beta),$ where: 
        \begin{enumerate}
            \item[($\alpha$)]
            
            \begin{enumerate}
                \item[$\bullet_{1}$] $\bfr \in \mathbf{R}$ and $\bfm_{\bfr} = \bfm$ and $ \langle \eta_{s}: s \in M_{\bfm} \rangle \in \Xi_{\bfr}^{+}.$
                
                \item[$\bullet_{2}$] $M' \subseteq M_{\bfm},$ (in the main case is $M_{\bfm} \cong (\kappa, <), M_{*}$ a cofinal subset of $M_{\bfm}$), $M_{*} = M' \cap N_{\bfr}, M' \in N_{\bfr},$
                
                \item[$\bullet_{3}$] $h$ is an isomorphism from $\bfm \rest M_{*}$ onto $\bfm \rest M_{\bfm},$
                
                \item[${\bullet_{4}}$] so there is $\bfG^{+} \subseteq \bbP_{\bfm} \cap N_{\bfr}$ generic over $N_{\bfr}$ such that $s \in M_{\bfm} \Rightarrow \eta_{s} = \name{\eta}_{s}[\bfG^{+}]$ and $\cH(\chi_{\bfr}) = \mathbf{j}_{N}''(N_{\bfr})[\mathbf{j}_{N}''[\bfG^{+}]]$ (by the definition of $\mathbf{R}$).
            \end{enumerate}
            
            \item[($\beta$)] there is $\bfG \subseteq N \cap \bbP_{\bfm}$ generic over $N$ such that $s \in M_{*} \Rightarrow \name{\eta}_{h(s)}[\bfG] = \eta_{s}.$
        \end{enumerate}
    \end{enumerate}
    
    
        
            
 
\end{conclusion}

\begin{PROOF}{\ref{g17}}
    Let $(\bfn, \bar{\mathscr{L}})$ be as in \ref{g14} for $\bfm_{0}$ and we shall show that $\bfn$ can serve as $\bfm.$ 

    Clauses (a), (b) of \ref{g17} holds by clause (a) of \ref{g14}(1). 
    
    Clause (c) of \ref{g17} holds by clause (b) of \ref{g14}(1).  
    
    To prove clause (d) of \ref{g17} assume $(\alpha)$ there.
    Let $\bfn = \bfn_{M}$ from 
    \ref{g14}(1)(e)  and use \ref{g14}(1),(2). 
    
    2) Use \ref{g14}(2).
\end{PROOF}

\begin{claim}\label{e70}
    1) Assume $\bfm$ is as in \ref{g17}, $\bfr \in \mathbf{R}$ (and $\bfm = \bfm_{\bfr}$). If $M \subsetneqq M_{\bfm} \cap N_{\bfr},$ $M = M' \cap N_{\bfr}, M'  \in N_{\bfr}$ and $\bar{\nu} =  \bar{\eta}_{\bfr} \rest M$ and $h \in N$ is an isomorphism from $\bfm \rest M'$ onto $\bfm \rest M_{\bfm},$ \underline{then} there is $\bfG$ witnessing $\bar{\eta}_{h} = \langle \eta_{h^{-1}(s)}: s \in M \rangle$ belongs to $\Xi_{\bfr}^{\bullet}$ (see Definition \ref{g5}(4)). 
    
    2) Above, $\bfG$ has (in $\bbP_{\bfm}$) an upper bound $p^{+}$ which satisfies $s \in M_{\bfn} \Rightarrow \rm{tr}(p(s)) = \name{\eta}_{s}[\mathbf{G}]$.
    
    3) If $p^{+} \in \bbP_{\bfm}$ is as above \underline{then} $p^{+}$ is also an upper bound of $\bfG' = \mathbf{G} \cap \bbP_{\bfm}[M]$ in $\bbP_{\bfm}[L_{\bfm}]$
    
    4)  If  $\bfm$ is strongly $(< \lambda^{+})$-directed, then $\bar{\eta}_{h}$ is cofinal in $\left( \Pi_{\zeta < \lambda_{\bfr}} \theta_{\zeta}, <_{J_{\lambda_{\bfr}}^{\rm{bd}}} \right)^{\bfV[\bbP_{\bfn}]}.$ 
    
    5) If $\bfm$ is strongly $(< \lambda^{+})$-directed (or just essentially directed, see \ref{e35}) \underline{then} for every $p \in \bbP_{\bfm}$ and $s \in \dom(p) \cap M_{\bfm}$ for every large enough $t \in M_{\bfm}$ we have $p \Vdash_{\bbP_{\bfm}}$``$\name{f}_{p(s)} \leq \name{\eta}_{t} \mod J_{\lambda}^{\rm{bd}}$''
\end{claim}

\begin{PROOF}{\ref{e70}}
    1) Let $\bfn$ from \ref{g17}(1)(d)($\beta$) for our $M.$
    
    2) This is because $\Vert N_{2} \Vert < \theta_{\lambda(\bfr)}$ \ref{e53}(1)(e) and \ref{e63}.

    3) As $\bbP_{\bfn}[M_{\bfn}] = \bbP_{\bfm}[M]$ because $\bbP_{\bfn} \lessdot \bbP_{\bfm}.$
    
    4) Easy recalling  \ref{g14}(2)(j) and  \ref{g13}(4)). 
    
    5) Just check the definition (and) or see \ref{c37}.
\end{PROOF}

\bibliographystyle{amsalpha}
\bibliography{shlhetal}

\end{document}